%% file: survey.tex
\documentclass[11pt]{article}
\usepackage[T1]{fontenc}
\usepackage[utf8]{inputenc}
\usepackage{lmodern}
\usepackage{hyperref}
\usepackage{graphicx}
\usepackage{wrapfig}
\usepackage[english]{babel}
\usepackage{amsfonts,amssymb,mathrsfs,amscd,amsmath,amsthm,xcolor}
\usepackage[all]{xy}
\usepackage{makeidx}
\usepackage[margin=1.0in]{geometry}

\def\R{{\mathbb R}}
\def\RP{{\mathbb{R}P}}

\newtheorem{theorem}{Theorem}
\newtheorem{lemma}{Lemma}
\newtheorem{proposition}{Proposition}
\newtheorem{corollary}{Corollary}
\newtheorem{statement}{Statement}
\newtheorem{conjecture}{Conjecture}
\theoremstyle{definition}
\newtheorem{example}{Example}
\newtheorem{remark}{Remark}
\newtheorem{definition}{Definition}

\def\:{\colon}

\def\sign{\operatorname{sign}}

\def\relint{\operatorname{relint}}
\def\conv{\operatorname{conv}}

\newcommand{\eps}{\varepsilon}

\newcommand{\bx}{\mathbf{x}}

\def\R{{\mathbb R}}

\def\Z{{\mathbb Z}}

\def\0{{\mathbf 0}}
\def\1{{\mathbf 1}}

\newcommand{\wG}{\widetilde G}

\title{On Groups $G_{n}^{k}$ and $\Gamma_{n}^{k}$: A Study of Manifolds, Dynamics, and Invariants}
\author{\textsc{Vassily O. Manturov\footnote{Moscow Institute for Physics and Technology, Kazan Federal University, North-Eastern University (China)}, Denis A. Fedoseev\footnote{Moscow State University, Moscow Center for Fundamental and Applied Mathematics},} \\ \textsc{Seongjeong Kim\footnote{Moscow Institute for Physics and Technology}, Igor M. Nikonov\footnote{Moscow State University}}}

\date{}

\makeindex

\begin{document}
\maketitle

\begin{abstract}
\input{abstract.tex}
\end{abstract}

\vspace{5mm}


 \newpage

\tableofcontents

\part{Introduction}
\input{introduction.tex}

\part{Group diagrams and the Greendlinger theorem}
\label{part:group_diagrams}
\input{greendlinger.tex}

\part{The Groups $G_{n}^{k}$} \label{chap:gnk_definition}

Usually, invariants of mathematical objects are valued in numerical
or polynomial rings, rings of homology groups, etc. In the present
section we prove a general theorem about invariants of {\em dynamical
systems} which are valued in groups very close to pictures, the so-called {\em free
$k$-braids.}

Formally speaking, free  $k$-braids form a group presented by
generators and relations; this group has lots of picture-valued invariants. For
free $2$-braids, the following principle can be realised: \\

{\em If a braid diagram $D$ is complicated enough, then it realises
itself as a subdiagram of any diagram $D'$ equivalent to $D$}. \\

In topology, this principle was first demonstrated in terms of {\em
parity} for the case of {\em virtual knots}, see \cite{knot_parity}.

Our invariant of braids is constructed by using horizontal trisecant lines.
Herewith,  the set of {\em critical values} (corresponding to these
trisecants) leads to  a certain picture which appears in all
diagrams equivalent to the initial picture.

The main theorem of the present section has various applications in
knot theory, geometry, and topology. It is based on the following main principle: 

\begin{center}
{\em If dynamical systems describing the motion of $n$ particles possess a nice codimension one property, then these dynamical systems admit a topological invariant valued in $G_{n}^{k}$}. \\
\end{center}

\input{section_6.tex}

\section{Groups $G_n^k$ and their homomorphisms} \label{chap:gnk_homomorphisms}
\input{homomorphisms.tex}

\section{Brunnian braids and generalisations of the groups $G_n^k$} \label{sec:brunnian}
\input{generalisations.tex}

\section{Realisation of spaces with $G_{n}^{k}$ action}\label{section:realisation}
\input{realisation.tex}

\section{The groups $G_{n}^{k}$ and invariants of manifolds} \label{sec:manifold_invariants}
\input{manifolds_invariants.tex}

\section{The Kirillov-Fomin algebras}
\label{sec:kir-fom}

We culminate the present part of the survey with a short overview of an interesting algebraic structure which was introduced and studied by S.~Fomin and A.~Kirillov (see~\cite{KirFom}).

First, recall that an associative algebra defined by generators and relations is called {\em quadratic} if all the relations are of degree 2.

\begin{definition} \label{def:kirillov-fomin_algebra}
	The Kirillov-Fomin algebra\index{Kirillov-Fomin algebra} $\mathcal{E}_n$ is a quadratic algebra (say, over $\mathbb{Z}$) with generators $[ij]$, $1 \le i < j \le n$, and the following relations: \\ \\
	$[ij]^2=0 \;\; \text{for} \;\; i<j;$ \\
	$[ij][jk]=[jk][ik]+[ik][jk] \;\; \text{for} \;\; 1 \le i < j < k\le n$ \\
	$[jk][ij]=[ik][jk]+[jk][ik] \;\; \text{for} \;\; i<j<k; $ \\
	$[ij][kl]=[kl][ij] \;\; \text{for} \;\; \{i,j\}\cap\{k,l\}=\emptyset, i<j, k<l.$
\end{definition}

Alternatively, the $\mathcal{E}_n$ algebra may be defined as the algebra with generators $[ij]$, where $i$ and $j$ are any distinct elemants of the set $\bar{n}=\{1,\dots, n\}$ (not necessarily $i<j$), subject to relations: \\ \\
$[ij]+[ji]=0;$ \\
$[ij]^2=0;$ \\
$[ij][jk]+[jk][ki]+[ki][ij]=0;$ \\
$[ij][kl]=[kl][ij]$ \\ \\
for any distinct $i,j,k$ and $l$. \\

The Kirillov-Fomin algebra admits many different representations. One of them is interesting in the sense that it is related to Bruhat order.\index{Bruhat order}

Let $s_{ij}\in S_n$ denote the transposition of the elements $i$ and $j$. Consider the ``Bruhat operators'' $[ij]$ acting in the space of formal linear combinations of elements of $S_n$ defined by the following formula:

\begin{center}
	$[ij]w = \left\{\begin{array}{ll}ws_{ij}, \;\; \text{if} \;\; l(ws_{ij})=l(w)+1, \\
	0, \;\; \text{otherwise}. \end{array} \right.$
\end{center}

Here $l(w)$ denotes the length of a permutation $w$. In other words, the operator multiplies $w$ by $s_{ij}$ if it moves it up one in the Bruhat order, and annihilates it otherwise.

It is easy to check that the Bruhat operators satisfy the relations from Definition~\ref{def:kirillov-fomin_algebra}, thus we obtain a representation of the Kirillov-Fomin algebra. Note though, that the representation is not faithful: it produces additional realtions. For example, for any $i<j<k$ we have $$[ij][ik][ij]=[jk][ik][jk]=0.$$ \\

The Kirillov-Fomin algebra relations strongly resemble those of certain $G_n^k$ groups. To be precise, consider the generators $a_{ij}$ of the group $G_n^2$ and let us formally write $a_{ij}=1+A_{ij}$. Further, let us work over $\mathbb{Z}_2$. Then we shall write down the relations of the group $G_n^2$ disregarding all summands of the order greater than 2. That gives us exactly the relations of the Kirillov-Fomin algebra generated by $A_{ij}$ if we consider commutation as the multiplication in the algebra.

For that reason it seems an interesting and promising idea to study these similarities to find some deeper connections between those algebraic structures.

\part{The groups $\Gamma_{n}^{k}$}

Returning back to 1954, J.W.~Milnor in his first paper on link group~\cite{Milnor} formulated the
idea that manifolds which share lots of well known invariant (homotopy type, etc.)
may have different {\em link groups}.\index{Link groups}

We shall look at braids on a manifold $M^{n}$ as loops in some configuration space of a manifold $W$, more precisely, $W = M^{n} \backslash K^{n-2}$, where $K^{n-2} \subset M^{n}$ is a submanifold of $M^{n}$ of codimension $2$.
Following $G^k_n$-ideology, we mark singular configurations, while moving along a loop in the configuration space, to get a word --- an element of a ``$G^k_n$-like'' group $\Gamma_n^k$. Singular configurations will arise here from the moments of transformation of triangulations spanned by the configuration points; in this case  the good property of codimension one is

\begin{center}
{\em ``$k$ points of the configuration lie on a sphere of dimension $k-3$ and there are no points inside the sphere''}.
\end{center}

Below, we shall consider an arbitrary closed $n$-manifold $M^{n}$ and for $N$ large enough we shall construct the braid group $B_{N}(M^{n})=\pi_1(C_N(M^n))$, where $C_N(M^n)$ is the space of generic (in some sense) configurations of points in $M^n$ that is dense enough to form a triangulation of $M$. This braid group will have a natural map to the group $\Gamma_{N}^{n+2}$. Here $N$ will be estimated in terms of minimal number of triangles of a triangulation of $M$.

The very first example in this theory, the groups $\Gamma_{n}^{4}$ which will be defined in the next section, and which correspond to usual braids on 2-surfaces, turn out to be related to many areas of mathematics: Ptolemy equation, pentagon relation, and cluster algebra presentations.

\section{The group $\Gamma_{n}^{4}$}
\label{sec:gamma_n4}
\input{gamma_n4.tex}

\section{Representations of braids via triangulations} \label{sec:gamma_presentations}
\input{braid_inv_triangulation}

\section{The group $\Gamma_{n}^{5}$}
\label{sec:gamma_n_5}

In Section~\ref{sec:gamma_n4} we defined the groups $\Gamma_n^4$. Now we move on to define their ``elder siblings'', groups $\Gamma_n^5$. Just as $\Gamma_n^4$ were closely related to Delaunay triangulations of 2-dimensional surfaces, $\Gamma_n^5$ come from 3-dimensional Delaunay triangulations. Later on, we will give the general definition of $\Gamma_n^k$ for arbitrary $k$ (see Section~\ref{sec:higher_gamma}). \\

Consider a configuration of $n$ points in general position in $\R^3$. We can think of those points as lying in a fixed tetrahedron $ABCD$. The points induce a unique Delaunay triangulation\index{Delaunay triangulation} of the tetrahedron: four points form a simplex of the triangulation if and only if there are no other points inside the sphere circumscribed over these points. When the points move in the space, the triangulation transforms.

In order to avoid degenerate Delaunay triangulations we exclude configurations where four points lie on one circle (intersection of a plane and a sphere). Note that this exclusion resembles the case of excluding triples of points on a line, but is essentially different and yields different results.

\begin{figure}
\centering\includegraphics[width=300pt]{ 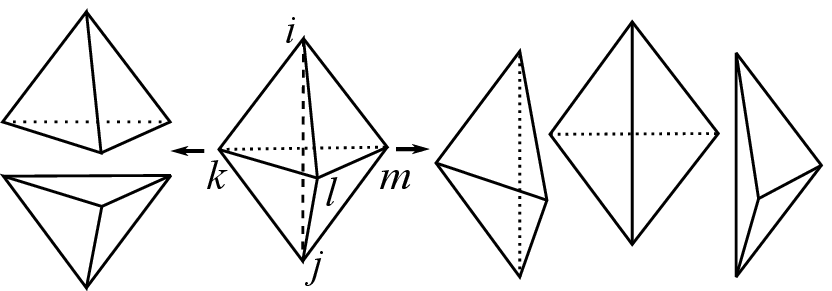}
\caption{A Pachner move}
\label{fig:pachner_move}
\end{figure}

Transformations of the combinatorial structure of the Delaunay triangulation correspond to configurations of codimension $1$ when five points lie on a sphere which does not contain any points inside. At this moment two simplices of the triangulation are replaced with three simplices as shown in Fig.~\ref{fig:pachner_move} (or vice versa). Such a transformation is called a {\em $(2-3)$ Pachner move} (or a $(3-2)$ Pachner move)\index{Pachner move}.

If we want to trace the evolution of triangulations that correspond to a dynamics of the points (and to describe it in algebraic language) we can mark each Pachner move with a letter. For the move that replaces the simplices $iklm, jklm$ with the simplices $ijkl, ijkm, ijlm$ in Fig.~\ref{fig:pachner_move} we use the generator $a_{ij,klm}$, which will be a generator of the group $\Gamma_n^5$ which we construct. Note that

\begin{enumerate}
	\item  we can split the indices into two subsets according to the combinatorics of the transformation;
	\item the generator $a_{ij,klm}$ is not expected to be involutive because it changes the number of simplices of the triangulation.
 \end{enumerate}

The relations on generators $a_{ij,klm}$ correspond to configurations of codimension $2$ which occurs when either

\begin{enumerate}
	\item six points lie on the same sphere with empty interior, or
	\item there are two spheres with five points on each of them, or
	\item five points on one sphere compose a codimension $1$ configuration with respect to configurations of codimension $1$.
\end{enumerate}

\begin{figure}[h!]
	\begin{minipage}{.45\textwidth}
		\begin{center}
 			\includegraphics[width =.63\textwidth]{ 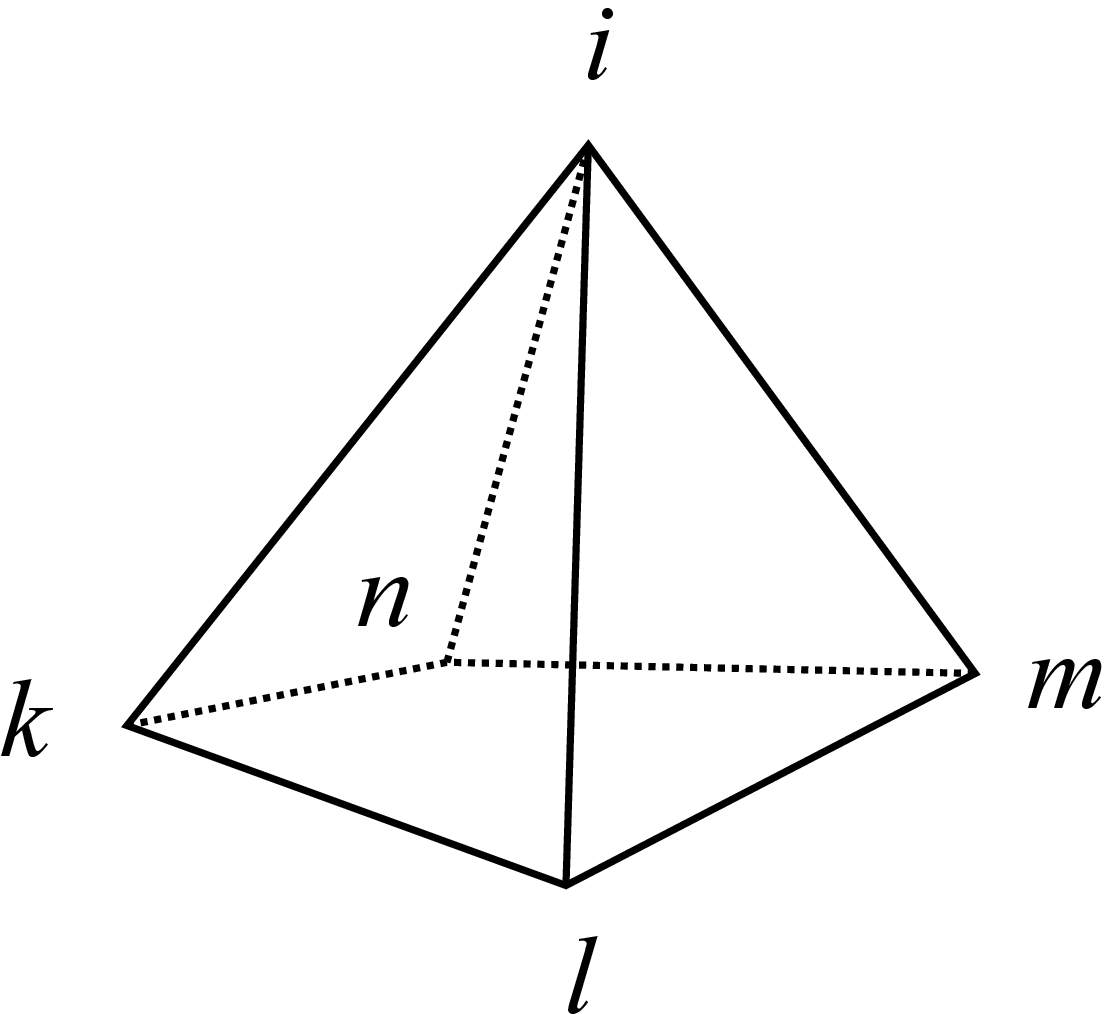}
		\end{center}
 		\caption{A quadrangular pyramid} \label{fig:delauney_pentahedron}
 	\end{minipage}
 	\hspace{5mm}
 	\begin{minipage}{.45\textwidth}
 		\centering\includegraphics[width = .45\textwidth]{ 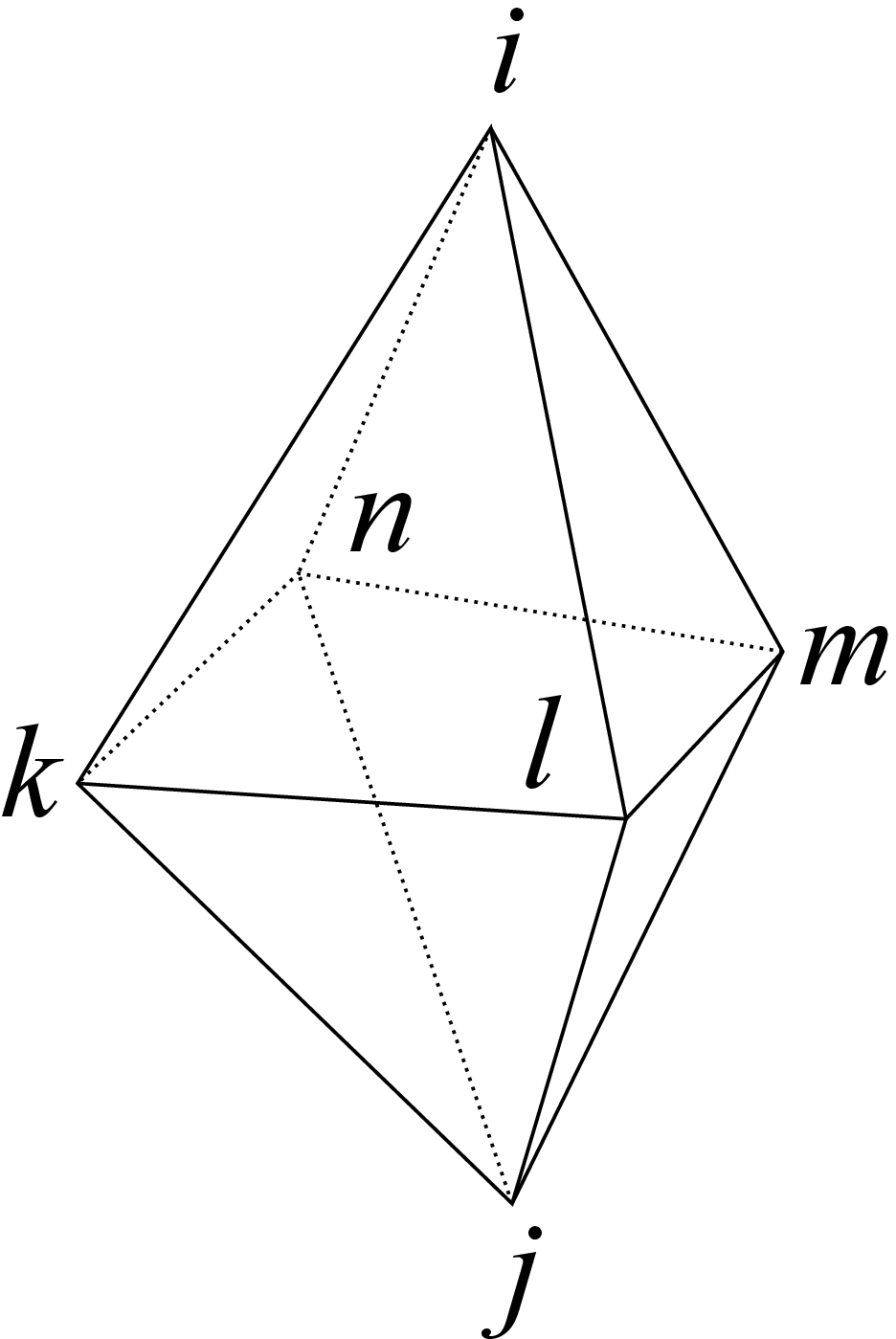}\quad\centering\includegraphics[width = .45\textwidth]{ 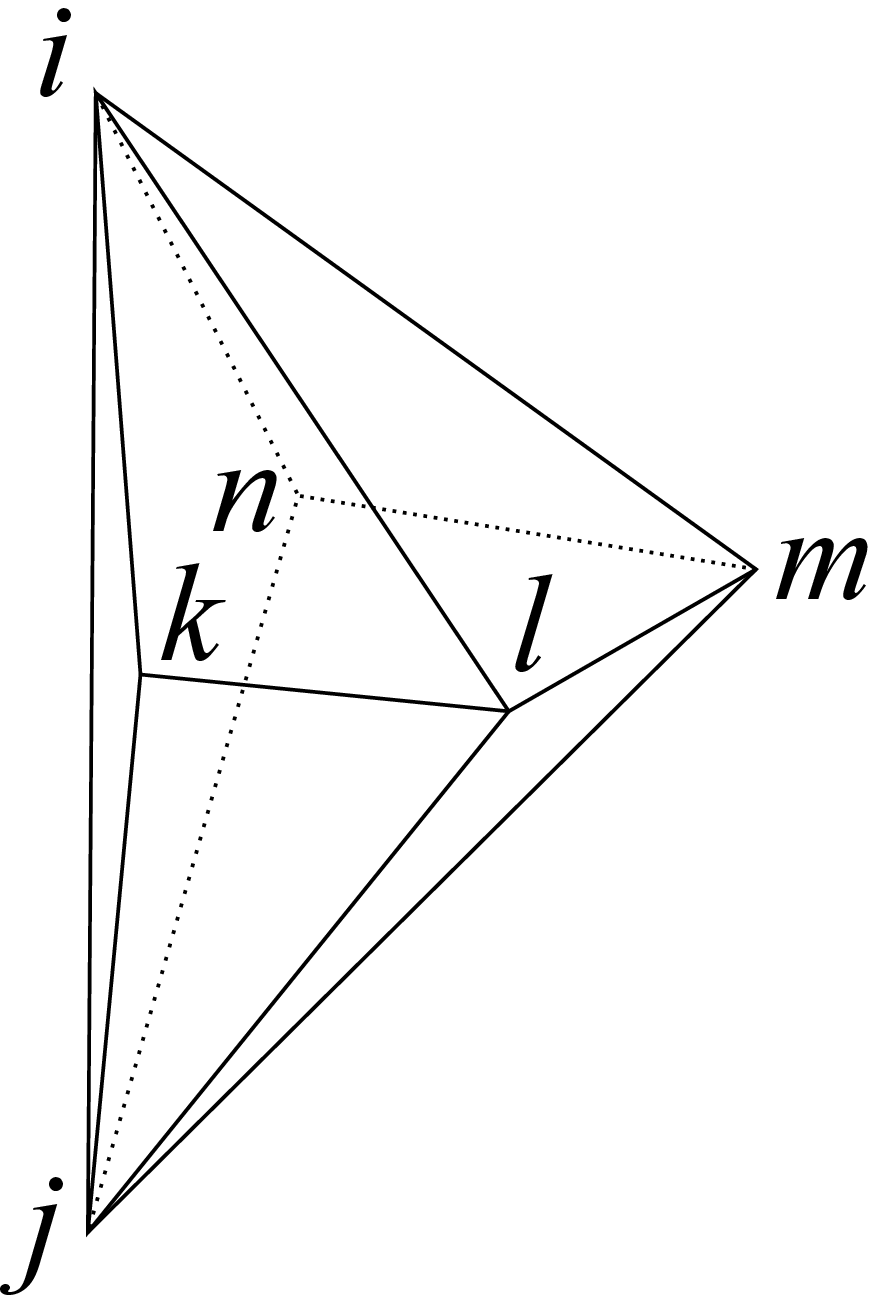}
		\caption{Convex polyhedra with triangular faces and $6$ vertices} \label{fig:delauney_octahedra}
 	\end{minipage}
 \end{figure}

The last case means the convex hull of the five points has a quadrangular face (Fig.~\ref{fig:delauney_pentahedron}). The vertices of this face lie on one circle so we exclude this configuration.

If there are two different spheres with five points on each of them, the inscribed simplices for the first sphere and for the second sphere are different. We can suppose that the inscribed convex hull has only triangular faces for each sphere; otherwise, the codimension should be greater than two. So the convex hulls can have one common face at most, and we can transform them independently. This gives us a commutation relation
$$a_{ij,klm} a_{i'j',k'l'm'}=a_{i'j',k'l'm'}a_{ij,klm},$$
where
$$|\{i,j,k,l,m\}\cap\{i',j',k',l',m'\}|<4,$$
$$|\{i,j\}\cap\{i',j',k',l',m'\}|<2,$$ and $$|\{i',j'\}\cap\{i,j,k,l,m\}|<2.$$

Consider now the case of six points on one sphere. The convex hull of these points is a convex polyhedron. The polyhedron must have only triangular faces; otherwise, there is an additional linear condition (four points lie on one plane) which rises the configuration of codimension $2$. There are two such polyhedra, see Fig.~\ref{fig:delauney_octahedra}.

\begin{figure}[h!]
	\begin{minipage}{.45\textwidth}
		\begin{center}
 			\includegraphics[width =.67\textwidth]{ 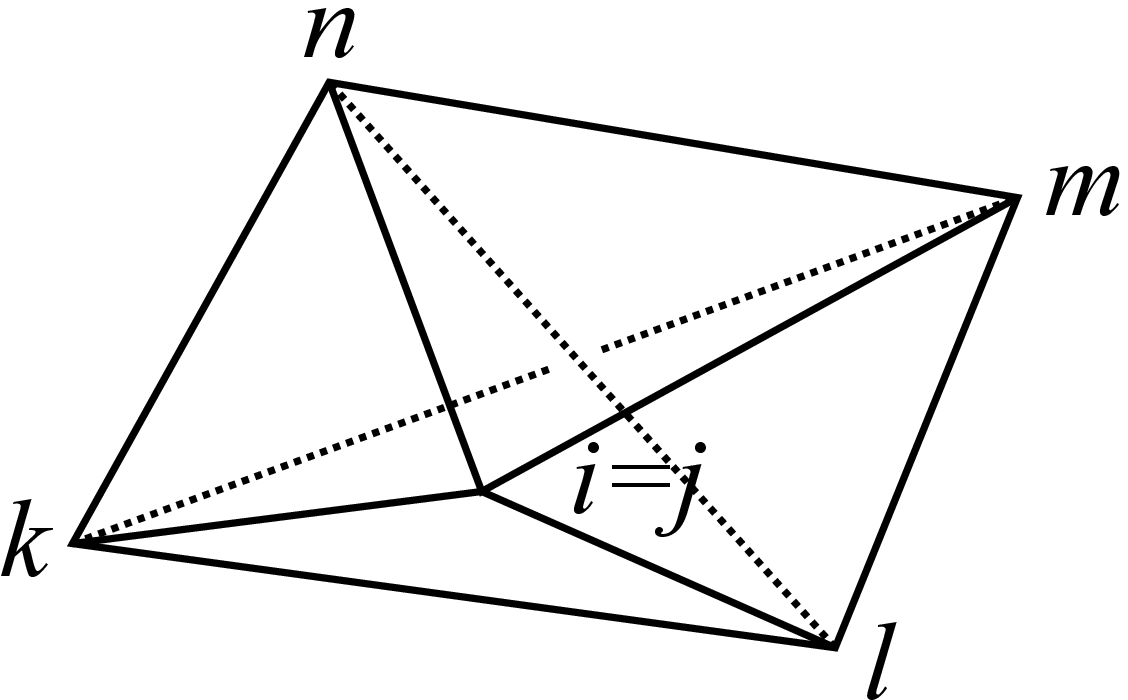}
		\end{center}
		\vspace{0.2cm}
 		\caption{An orthogonal projection of the octahedron} \label{fig:delauney_octahedron1_projection}
 	\end{minipage}
 	\hspace{5mm}
 	\begin{minipage}{.45\textwidth}
		\begin{center}
 			\includegraphics[width =.63\textwidth]{ 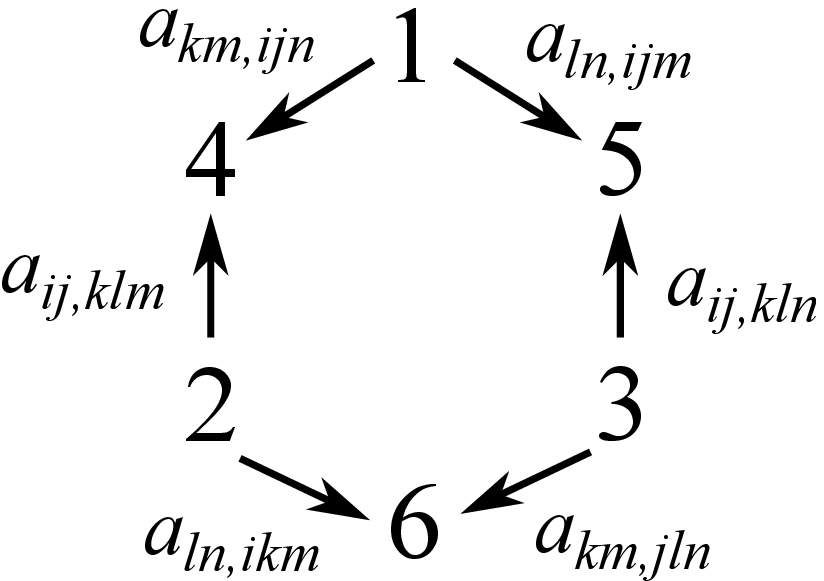}
		\end{center}
 		\caption{Triangulation graph for the octahedron} \label{fig:delauney_octahedron1_hasse}
 	\end{minipage}
 \end{figure}

For the octahedron on the left we specify the geometrical configuration assuming that the orthogonal projection along the direction $ij$ maps the points $i$ and $j$ near the projection of the edge $kl$ and for the line $ln$ to be higher than $km$ if you look from the vertex $i$ to the vertex $j$ (Fig.~\ref{fig:delauney_octahedron1_projection}). In this case we have six triangulations:
\begin{enumerate}
\item $ijkl, ijkn, ijlm, ijmn$;
\item $iklm, ikmn, jklm, jkmn$;
\item $ikln, ilmn, jkln, jlmn$;
\item $ijkl, ijkm, ijlm, ikmn, jkmn$;
\item $ijkl, ijkn, ijln, ilmn, jlmn$;
\item $ikln, ilmn, jklm, jkmn, klmn$.
\end{enumerate}

The first three have four simplices, the last three have five simplices. The Pachner moves between the triangulations are shown in Fig.~\ref{fig:delauney_octahedron1_hasse}. Thus, we have the relation
$$
a^{}_{km,ijn}a^{-1}_{ij,klm}a^{}_{ln,ikm}a^{-1}_{km,jln}a^{}_{ij,kln}a^{-1}_{ln,ijm}=1.
$$

\begin{figure}[h!]
	\begin{minipage}{.45\textwidth}
		\begin{center}
 			\includegraphics[width =.5\textwidth]{ 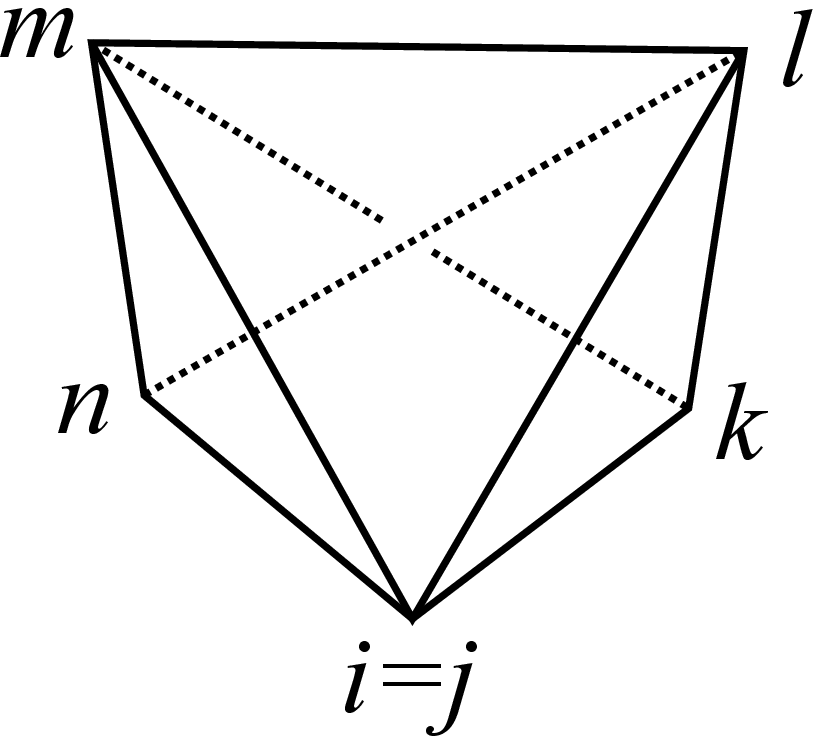}
		\end{center}
		\vspace{-0.1cm}
 		\caption{An orthogonal projection of the shifted octahedron} \label{fig:delauney_octahedron2_projection}
 	\end{minipage}
 	\hspace{5mm}
 	\begin{minipage}{.45\textwidth}
		\begin{center}
 			\includegraphics[width =.63\textwidth]{ 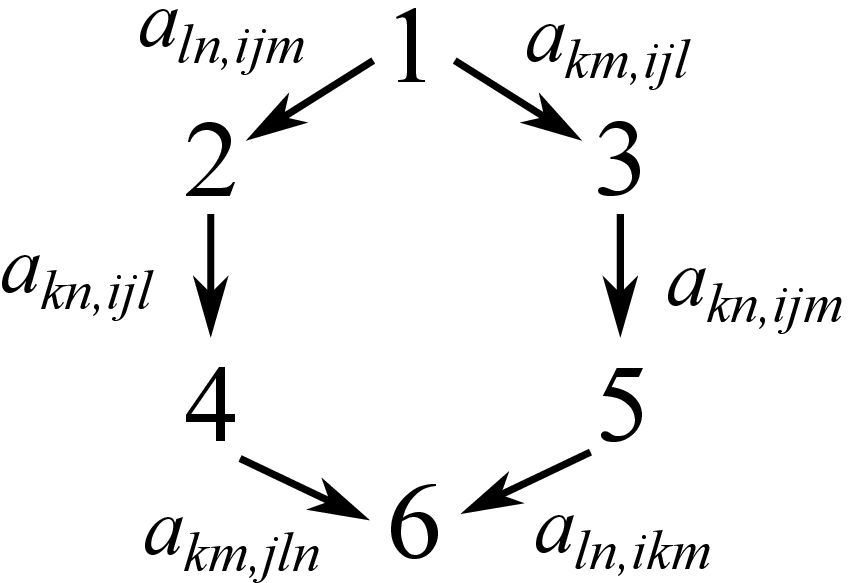}
		\end{center}
 		\caption{Triangulation graph for the shifted octahedron} \label{fig:delauney_octahedron2_hasse}
 	\end{minipage}
 \end{figure}

In the case of the shifted octahedron (Fig.~\ref{fig:delauney_octahedra} right) we assume that the line $ln$ lies higher than the line $km$ when one looks from the vertex $i$ to the vertex $j$ (Fig.~\ref{fig:delauney_octahedron2_projection}). Then we have six triangulations:
\begin{enumerate}
\item $ijkl, ijlm, ijmn$;
\item $ijkl, ijln, ilmn, jlmn$;
\item $ijkm, ijmn, iklm, jklm$;
\item $ijkn, ikln, ilmn, jkln, jlmn$;
\item $ijkn, iklm, ikmn, jklm, jkmn$;
\item $ijkn, ikln, ilmn, jklm, jkmn, klmn$.
\end{enumerate}

The Pachner moves between the triangulations are shown in Fig.~\ref{fig:delauney_octahedron2_hasse}. This gives us the relation
$$
a_{ln,ijm}a_{kn,ijl}a_{km,jln}=a_{km,ijl}a_{kn,ijm}a_{ln,ikm}.
$$


Now we can give the definition for $\Gamma_n^5$.

\begin{definition}\label{def:Gamma_n5}
The group $\Gamma_n^5$\index{Group!$\Gamma_n^5$} is the group with generators $$\{ a_{ij,klm}\,|\, \{i,j,k,l,m\}\subset \bar n, |\{i,j,k,l,m\}|=5\}$$ and relations
\begin{enumerate}
\item $a_{ij,klm}=a_{ji,klm}=a_{ij,lkm}=a_{ij,kml}$,
\item $a_{ij,klm} a_{i'j',k'l'm'}=a_{i'j',k'l'm'}a_{ij,klm}$, for $|\{i,j,k,l,m\}\cap\{i',j',k',l',m'\}|<4$, $|\{i,j\}\cap\{i',j',k',l',m'\}|<2$ and $|\{i',j'\}\cap\{i,j,k,l,m\}|<2$,
\item $a^{}_{km,ijn}a^{-1}_{ij,klm}a^{}_{ln,ikm}a^{-1}_{km,jln}a^{}_{ij,kln}a^{-1}_{ln,ijm}=1$ for distinct $i,j,k,l,m,n$,
\item $a_{ln,ijm}a_{kn,ijl}a_{km,jln}=a_{km,ijl}a_{kn,ijm}a_{ln,ikm}$ for distinct $i,j,k,l,m,n$.
\end{enumerate}
\end{definition}

The group $\Gamma_n^5$ can be used to construct invariants of braids and dynamical systems like the groups $G_n^k$.

Let us consider the configuration space
$\tilde C_n(\R^3)$ which consists of nonplanar
$n$-tuples $(x_1,x_2,\dots,x_n)$ of points in $\R^3$ such that for any distinct $i,j,k,l$ the points $x_i,x_j,x_k,x_l$ do not lie on the same circle.

Let us construct a homomorphism from $\pi_1(\tilde C_n(\R^3))$ to $\Gamma_{n}^5$. Let
$$\alpha=(x_1(t),\dots, x_n(t)),\quad t\in [0,1],$$
be a loop in $\tilde C_n(\R^3)$. For any $t$ the set
$\mathbf{x}(t)=(x_1(t), \dots, x_n(t))$
determines a Delaunay triangulation $T(t)$ of the polytope $\conv(\mathbf{x}(t))$. If $\alpha$ is in general position, then there will be a finite number of moments $0<t_1<t_2<\dots<t_N<1$ when the combinatorial structure of $T(t)$ changes, and for each $p$ the transformation of the triangulation at the moment $t_p$ will be the Pachner move on simplices with vertices $i_p,j_p,k_p,l_p,m_p$. We assign to this move the generator $a_{i_pj_p,k_pl_pm_p}$ at the moment $t_{p}$ if 2 simplicies $i_{p}k_{p}l_{p}m_{p}, j_{p}k_{p}l_{p}m_{p}$ are replaced by the 3 simplices $i_{p}j_{p}k_{p}l_{p}, i_{p}j_{p}k_{p}m_{p}, i_{p}j_{p}l_{p}m_{p}$, or $a^{-1}_{i_pj_p,k_pl_pm_p}$ at the moment $t_{p}$ in the converse situation.

Denote
$$\phi(\alpha)=\prod_{p=1}^N a^{\epsilon_{p}}_{i_pj_p,k_pl_pm_p}\in\Gamma_{n}^5,$$
where $\epsilon_{p} \in \{1,-1\}$.
\begin{theorem}\label{thm:Gamman5_invariant}
The homomorphism $\phi\colon \pi_1(\tilde C_n(\R^3))\to\Gamma^5_{n}$ is well defined.
\end{theorem}

\begin{proof}
We need to show that the element $\phi(\alpha)$ does not depend on the choice of the representative $\alpha$ in a given homotopy class. Given a generic homotopy $\alpha(\tau), \tau\in[0,1],$  of loops in $\tilde C_n(\R^3)$, the transformations of the words $\phi(\alpha(\tau))$ correspond to point configurations of codimension $2$, considered above the definition of $\Gamma_n^5$, and thus are counted by the relations in the group $\Gamma^5_{n}$. Therefore, the element $\phi(\alpha(\tau))$ of the group $\Gamma^5_{n}$ remains the same when $\tau$ changes.
\end{proof}

\section{The strategy of defining $\Gamma_{n}^{k}$ for arbitrary $k$}
\label{sec:higher_gamma}

The groups $\Gamma_n^4$ and $\Gamma_n^5$ defined above admit extrapolation to groups $\Gamma_n^k, k\ge 4$\index{Group!$\Gamma_n^k$}, that we define in this section.

Consider the configuration space $\tilde C_n(\R^{k-2})$, $4\le k\le n$, consisting of $n$-point configurations $$\mathbf{x}=(x_1, x_2,\dots, x_n)$$ in $\R^{k-2}$, such that $\dim (\conv\mathbf{x})=k-2$ and there are no $k-1$ points which lie on one $(k-4)$-dimensional sphere (the intersection of a $(k-3)$-sphere and a hyperplane in $\R^{k-2}$).

A configuration $\mathbf{x}=(x_1, x_2,\dots, x_n)\in\tilde C_n(\R^{k-2})$ determines a Delaunay triangulation of $\conv\mathbf{x}$ which is unique when $\mathbf x$ is in general position. If the vertices $x_{i_1},\dots,x_{i_{k-1}}$ span a simplex of the Delaunay triangulation, then the interior of the sphere circumscribed over these points does not contain other points of the configuration. The inverse statement is true when $\mathbf x$ is generic. The condition that no $k-1$ points lie on one $(k-4)$-dimensional sphere ensures that there are no degenerate simplices in the Delaunay triangulation.

Let $\alpha=(x_1(t),\dots, x_n(t)), t\in [0,1],$ be a loop in $\tilde C_n(\R^{k-2})$. For any $t\in [0,1]$ the configuration $\mathbf{x}(t)=(x_1(t), \dots, x_n(t))$ determines a Delaunay triangulation $T(t)$ of the polytope $\conv(\mathbf{x}(t))$. If $\alpha$ is generic, then there will be a finite number of moments $0<t_1<t_2<\dots<t_L<1$ when the combinatorial structure of $T(t)$ changes.

For each singular configuration $\mathbf{x}(t_i)$, either one simplex degenerates and disappears on the boundary $\partial\conv(\mathbf{x}(t))$ or the Delaunay triangulation is not unique, that is, there is a sphere in $\R^{k-2}$ which contains $k$ points of $\mathbf x$ on it and no points of $\mathbf x$ inside. Assuming $\mathbf x$ is generic, the span of these $k$ points is a simplicial polytope. Below we shall consider only the latter type of singular configurations.

We recall that the {\em join} of two sets $X,Y\subset\R^m$ is defined as
$$X*Y=\{ \lambda x+(1-\lambda)y\,|\, x\in X, y\in Y, \lambda\in [0,1]\},$$
and the {\em relative interior}\index{Relative interior} of a finite set $X\subset\R^m$ is defined as
$$\relint X=\left\{ \sum_{x\in X}\lambda_x x\,|\,\forall x\ \lambda_x>0, \sum_{x\in X}\lambda_x=1\right\}.$$
The polytope $\Delta_P*\Delta_Q$ has two triangulations:
$$T_P=\{\mathbf{x}_{P\cup Q}\setminus\{x_p\}\}_{p\in P} = \left\{x_{p_1}\cdots x_{p_{i-1}}x_{p_{i+1}}\cdots x_{p_l}x_{q_1}\cdots x_{q_{k-l}}\right\}_{i=1}^l$$
and
$$T_Q=\{\mathbf{x}_{P\cup Q}\setminus\{x_q\}\}_{q\in Q} = \left\{x_{p_1}\cdots x_{p_l} x_{q_1} \cdots x_{q_{i-1}} x_{q_{i+1}} \cdots x_{q_{k-l}}\right\}_{i=1}^{k-l}.$$
Here $P=\{p_1,\dots,p_l\}$, $Q=\{q_1,\dots,q_{k-l}\}$ and $\mathbf{x}_J=\{x_j\}_{j\in J}$ for any $J\subset\{1,\dots,n\}$.

The simplicial polytopes in $\R^{k-2}$ with $k$ vertices are described in~\cite{Gruenbaum}. Each of them is the join $\Delta_P*\Delta_Q$ of simplices $\Delta_P=\conv(x_{p_1},\dots,x_{p_l})$ of dimension $l-1\ge 1$ and $\Delta_Q=\conv(x_{q_1},\dots,x_{q_{k-l}})$ of dimension $k-l-1\ge 1$ such that the intersection $\relint(x_{p_1},\dots,x_{p_l})\cap\relint(x_{q_1},\dots,x_{q_{k-l}})$ consists of one point.

The condition $\relint(\mathbf{x}_P)\cap\relint(\mathbf{x}_Q)=\{z\}$ implies $P\cap Q=\emptyset$.
Thus, when the configuration $\bx(t)$ goes over a singular value $t_i, i=1,\dots, L$, in the Delaunay triangulation simplices $T_{P_i}$ are replaced by simplices $T_{Q_i}$ for some subsets $P_i,Q_i\subset\{1,\dots,n\}$, $P_i\cap Q_i=\emptyset$, $|P_i|,|Q_i|\ge 2$, $|P_i\cup Q_i|=k$. This transformation is called a {\em Pachner move}\index{Pachner move}. We assign the letter $a_{P_i,Q_i}$ to this transformation.

Hence, the loop $\alpha$ produces a word
\begin{equation}\label{eq:Gamma_nk_braid_invariant}
\Phi(\alpha)=\prod_{i=1}^L a_{P_i,Q_i}
\end{equation}
 in the alphabet
\begin{equation*}
\mathcal A_n^k=\left\{ a_{P,Q}\,|\, P,Q\subset\{1,\dots,n\}, P\cap Q=\emptyset, |P\cup Q|=k, |P|,|Q|\ge 2\right\}.
\end{equation*}

Now let us consider a generic homotopy $\alpha_s, s\in[0,1],$ between two generic loops $\alpha_0$ and $\alpha_1$. A loop $\alpha_s=\{\bx(s,t)\}_{t\in[0,1]}$ can contain a configuration of codimension $2$. This means that for some $t$ the configuration $\bx(s,t)$ has two different $k$-tuples of points, each of them lies on a sphere whose interior contains no points of $\bx(s,t)$. If these spheres do not coincide, then their intersection contains at most $k-2$ points (the intersection can not contain $k-1$ points because $\bx(s,t)\in\tilde C_n(\R^{k-2})$). Hence, the simplices involved in one Pachner move cannot be involved in the other one, so the Pachner moves can be performed in any order.

If  the $k$-tuples of points lie on one sphere, then there is a sphere with $k+1$ points of $\bx(s,t)$ on it and its interior contains no points of $\bx(s,t)$. This $k+1$ points span a simplicial polytope in $\R^{k-2}$. Such polytopes are described in~\cite{Gruenbaum}. The description uses the notion of {\em Gale transformation}.\index{Gale transformation}

Let $X=\{x_1,\dots,x_n\}$ be a set of $n$ points in $\R^d$ such that $\dim\conv X = d$. Then $n\ge d+1$. Let $x_i=(x_{1i},\dots,x_{di})\in\R^d$, $i=1,\dots,n$, be the coordinates of the points of $X$. The matrix
$$
M = \left(\begin{array}{cccc}
x_{11} & x_{12} & \cdots & x_{1n}\\
x_{21} & x_{22} & \cdots & x_{2n}\\
\cdots &\cdots &\cdots &\cdots \\
x_{d1} & x_{d2} & \cdots & x_{dn}\\
1 & 1 & \cdots & 1
\end{array}\right)
$$
has rank $d+1$. Then the dimension of the space $$\ker M = \{ b\in\R^n\,|\, Mb=0\}$$ of dependencies between columns of $M$ is equal to $n-(d+1)$. Take any basis $b_j = (b_{j1},b_{j2},\dots, b_{jn})$, $j=1,\dots, n-d-1,$ of $\ker M$ and write it in matrix form
\begin{equation}\label{eq:gale_b_matrix}
B = \left(\begin{array}{ccc}
b_{11} & \cdots & b_{1n}\\
\cdots &\cdots  &\cdots \\
b_{n-d-1,1} & \cdots & b_{n-d-1,n}
\end{array}\right).
\end{equation}

The columns of the matrix $B$ form a set $$Y=\{y_1,\dots, y_n\}, y_i=(b_{1i},\dots,b_{n-d-1,i}),$$ in $\R^{n-d-1}$. The set $Y$ is called a {\em Gale transform}\index{Gale transform} of the point set $X$. Gale transforms which correspond to different bases of $\ker M$ are linearly equivalent. The vectors of the Gale transform $Y$ may coincide.

\begin{example}
Let $X$ be a pentagon in $\R^2$ with vertices $x_1=(0,2)$, $x_2=(-2,1)$, $x_3=(-1,-1)$, $x_4=(1,-1)$, $x_5=(2,1)$. Then
$$
M = \left(\begin{array}{ccccc}
0 & -2 & -1 & 1 & 2\\
2 & 1 & -1 & -1 & 1\\
1 & 1 & 1 & 1 & 1
\end{array}\right)
$$
and
$$
B = \left(\begin{array}{ccccc}
-4 & 1 & 3 & -5 & 5\\
-4 & 6 & -7 & 5 & 0
\end{array}\right).
$$
The Gale transform $Y$ consists of vectors $y_1=(-4,-4)$, $y_2=(1,6)$, $y_3=(3,-7)$, $y_4=(-5,5)$, $y_5=(5,0)$, see Fig~\ref{fig:pentagon_gale_example}.

\begin{figure}
\centering\includegraphics[width=120pt]{ 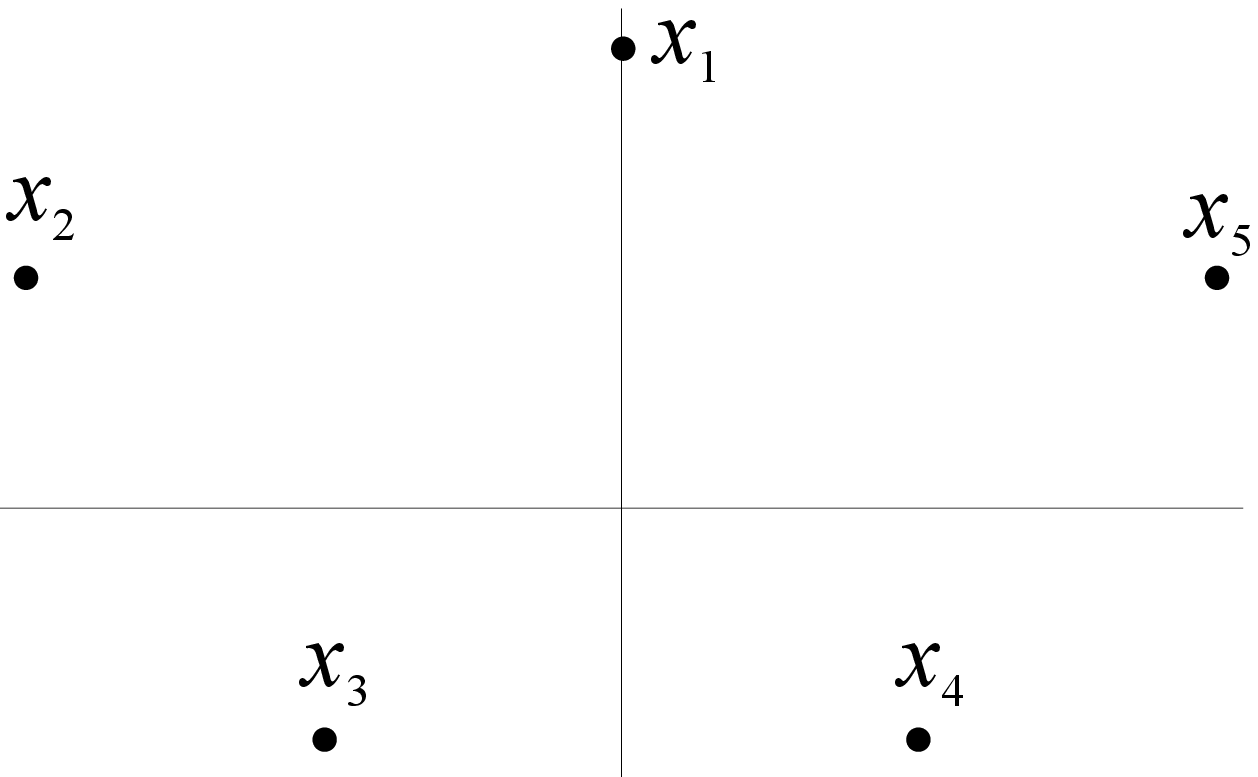}\hspace{15mm} \includegraphics[width=120pt]{ 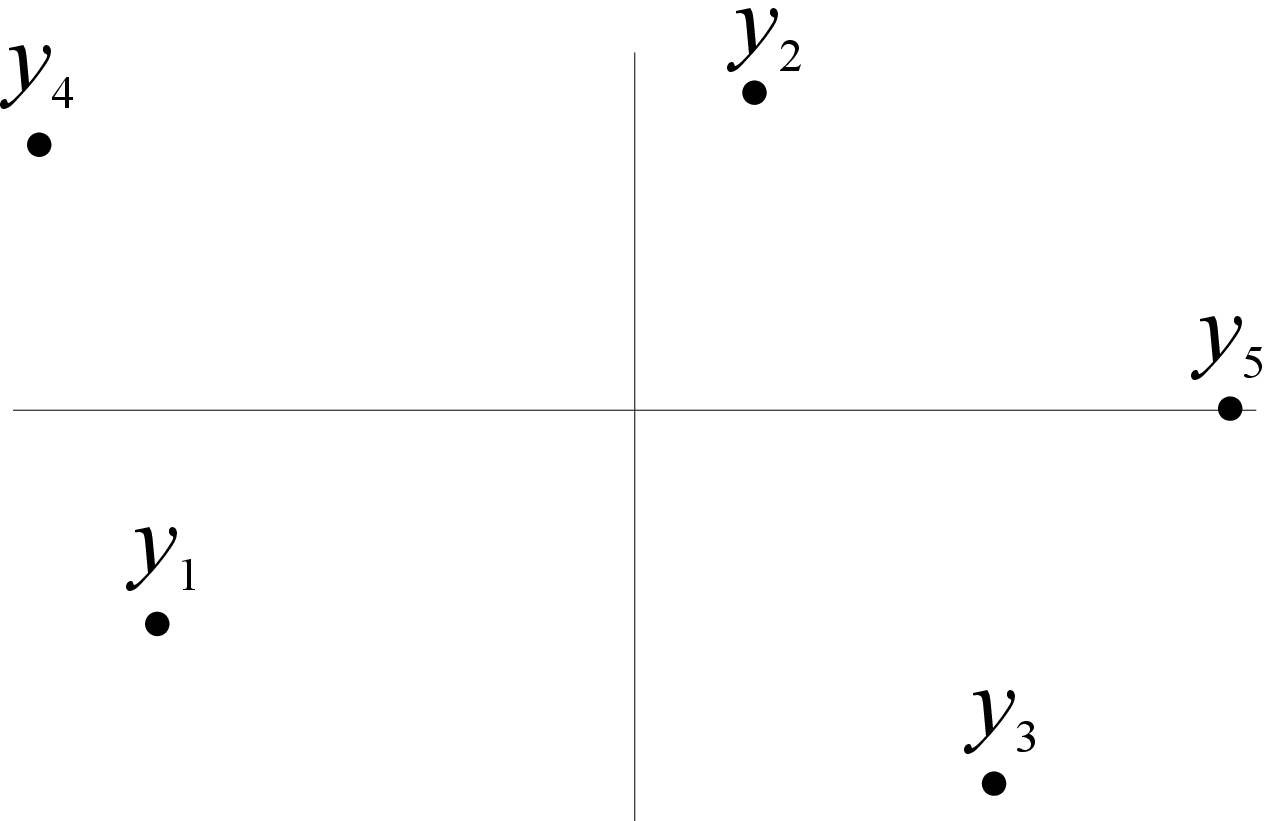}
\caption{A pentagon and its Gale transform}\label{fig:pentagon_gale_example}
\end{figure}
\end{example}

Let $Y=\{y_1,\dots, y_n\}$ be a Gale transform of $X$. The set
$$\bar Y = \{\bar y_1,\dots, \bar y_n\},\quad \bar y_i=\left\{\begin{array}{cl}\frac{y_i}{\|y_i\|}, & y_i\ne 0,\\ 0, & y_i=0 \end{array}\right.
$$
is called a {\em Gale diagram}\index{Gale diagram} of the point set $X$. It is a subset of $S^{n-d-2}\cup\{0\}$.

Two subsets $\bar Y = \{\bar y_1,\dots, \bar y_n\}$ and $\bar Y' = \{\bar y'_1,\dots, \bar y'_n\}$ in $S^{n-d-2}\cup\{0\}$ are called {\em equivalent} if there is a permutation $\sigma$ of $\bar{n}$ such that for any $J\subset \bar{n}$
$$
0\in\relint \bar Y_J \Longleftrightarrow 0\in\relint \bar Y'_{\sigma(J)}.
$$
Here we denote $\bar Y_J = \{\bar y_i\}_{i\in J}$ and $\bar Y'_J = \{\bar y'_i\}_{i\in J}$.
The properties of Gale diagrams can be summarized as follows~\cite{Gruenbaum}:

\begin{theorem}\label{thm:gale_diagram}
\begin{enumerate}
\item Let $X$ be a set of $n$ points, which are vertices of some polytope $P$ in $\R^d$, and $\bar Y$ be its Gale diagram. Then the set of indices $J\subset \bar{n}$ defines a face of $P$ if and only if $0\in\relint Y_{\bar{n} \setminus J}$.
\item Let $X$ and $X'$ be sets of vertices of polytopes $P$ and $P'$, $|X|=|X'|$, and $\bar Y$ and $\bar Y'$ be their Gale diagrams. Then $P$ and $P'$ are combinatorially equivalent if and only if $\bar Y$ and $\bar Y'$ are equivalent.
\item For any $n$-point set $\bar Y\in S^{n-d-2}\cup\{0\}$ such that $\bar Y$ spans $\R^{n-d-1}$ and $0$ lies in the interior of $\conv\bar Y$, there is an $n$-point set $X$ in $\R^d$ such that $\bar Y$ is a Gale diagram of $X$.
\end{enumerate}
\end{theorem}

The theorem implies (see~\cite{Gruenbaum})
that simplicial polytopes with $k+1$ vertices in $\R^{k-2}$ are in a bijection with standard Gale diagrams in $\R^2$ (defined uniquely up to isometries of the plane).

Here, a {\em standard Gale diagram}\index{Gale diagram!standard} of order $l=k+1$ is a subset $\bar Y$, $|\bar Y|=l$, of the vertices set $\{e^{\frac{\pi i p}{l}}\}_{p=0}^{2l-1}$ of the regular $2l$-gon inscribed in the unit circle $S^1$, such that:
\begin{enumerate}
\item any diameter of $S^1$ contains at most one point of $\bar Y$;
\item for any diameter of $S^1$, any open half-plane determined by it contains at least two points of $\bar Y$.
\end{enumerate}

The first property means that the corresponding polytope is simplicial, the second means that any of the $k+1$ vertices of the polytope form a face.

\begin{figure}
\centering
\begin{tabular}{cc}
\includegraphics[height=50pt]{ 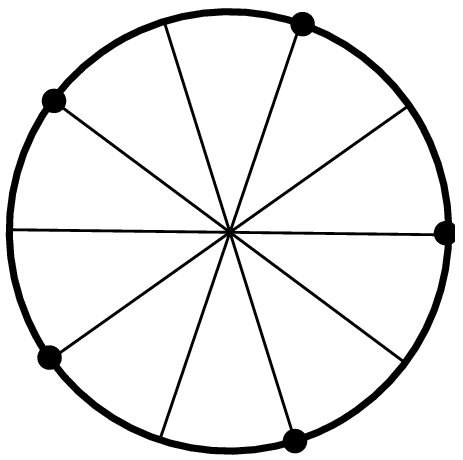} & \includegraphics[height=50pt]{ 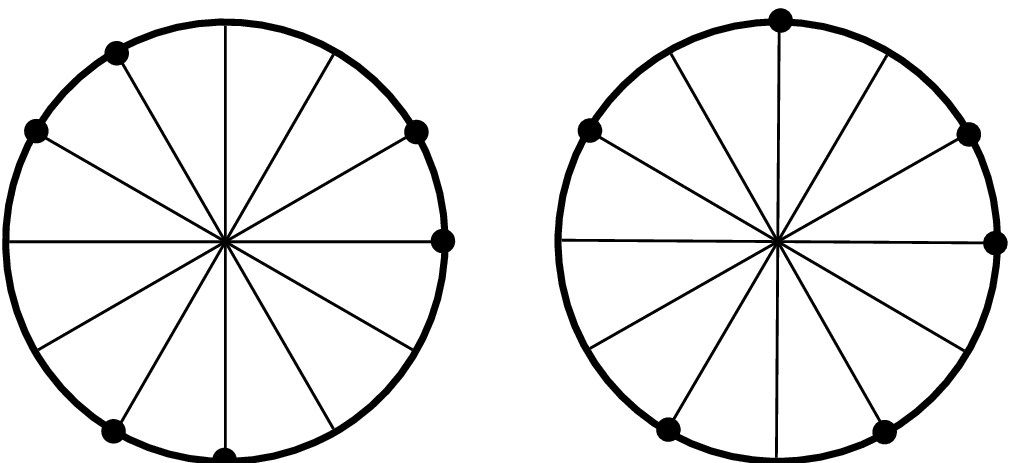}\\
 $l=5$ & $l=6$\\
 \multicolumn{2}{c}{\includegraphics[height=50pt]{ 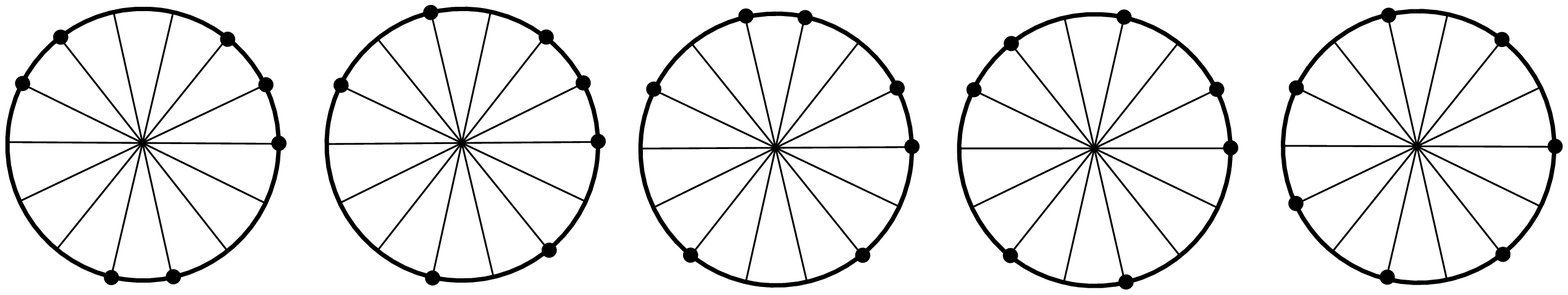}}\\
 \multicolumn{2}{c}{$l=7$}
\end{tabular}
\caption{Standard Gale diagrams of small order}\label{fig:gale_standard_diagrams}
\end{figure}

The number $c_l$ of standard Gale diagrams of order $l$ is equal to
$$
c_l= 2^{\left[\frac{l-3}{2}\right]}-\left[\frac{l+1}4\right]+\frac 1{4l}\sum_{h\colon 2\nmid h\mid l}\varphi(h)\cdot 2^{\frac lh},
$$
where $l=2^{a_0}\prod_{i=1}^{t}p_i^{a_i}$ is the prime decomposition of $l$ and $\varphi$ is Euler's function. For small $l$ the numbers are $c_5=1, c_6=2, c_7=5$, see Fig.~\ref{fig:gale_standard_diagrams}.

Let us describe the triangulations of the simplicial polytopes with $k+1$ vertices in $\R^{k-2}$.

Let $X=\{x_1,\dots,x_n\}$ be a subset of $\R^d$, so $x_i = (x_{i1},\dots, x_{id})$, $i=1,\dots, n$. Let $P=\conv X$ be the convex hull of $X$, assuming that $\dim P = d$. A triangulation $T$ of the polytope $P$ with vertices in $X$ is called {\em regular} if there is a height function $h\colon X\to\R$ such that $T$ is the projection of the lower convex hull of the lifting $X^h=\{x_1^h,\dots,x_n^h\}\subset\R^{d+1}$, where $x_i^h = (x_i, h(x_i))$. This means that a set of indices $J\subset\{1,\dots,n\}$ determines a face of $T$ if and only if there exists a linear functional $\phi$ on $\R^{d+1}$ such that $\phi(0,\dots,0,1)>0$ and $J=\{i\,|\, \phi(x^h_i) = \min_{x^h\in X^h}\phi(x^h)\}$. In case $T$ is regular we write $T=T(X,h)$. Any generic height function induces a regular triangulation.

The Delaunay triangulation of $X$ is regular with the height function $h\colon \R^d\to\R$, $h(z)=\|z\|^2= \sum_{i=1}^d z_i^2$ if $z=(z_1,\dots,z_d)\in\R^d$.

A height function $h\colon X\to\R$ can be regarded as a vector $h=(h_1,\dots, h_n)\in\R^n$ where $h_i=h(x_i)$. Denote $\beta(h)=Bh\in\R^{n-d-1}$ where $B$ is the matrix~\eqref{eq:gale_b_matrix} used to define a Gale transform of $X$. Let $\bar Y=\{\bar y_1,\dots,\bar y_n\}$ be a Gale diagram of $X$. Convex cones generated by the subsets of $\bar Y$ split the space $\R^{n-d-1}$ into a union of conic cells. A relation between the triangulation $T(X,h)$ and the Gale diagram $\bar Y$ can be described as follows.

\begin{theorem}\label{thm:gale_triangulation}
\begin{enumerate}
\item If $T(X,h)$ is a regular triangulation of $X$, then $\beta(h)$ belongs to a conic cell of maximal dimension in the splitting of $\R^{n-d-1}$ induced by $\bar Y$.
\item Let $J\subset \bar{n}=\{1,\dots,n\}$. The set $X_J$ spans cells of the triangulation $T(X,h)$ if and only if $\beta(h)\in\mathrm{concone}(\bar Y_{\bar{n}\setminus J})$, where $\mathrm{concone}(X_{\bar{n}\setminus J})$ is the convex cone spanned by the set $\bar Y_{\bar{n}\setminus J}$.
\end{enumerate}
\end{theorem}

Let $P$ be a simplicial polytope with $l=k+1$ vertices in $\R^{l-3}$ and $\bar Y=\{\bar y_1,\dots,\bar y_l\}$ be the corresponding standard Gale diagram. By Theorem~\ref{thm:gale_triangulation} there are $l$ different regular triangulations which correspond to open sectors between the rays spanned by the vectors of $\bar Y$. The graph, whose vertices are combinatorial classes of triangulations of $P$ and the edges are Pachner moves, is a cycle. Let us find which Pachner moves appear in this cycle.

We change the order of vertices of $P$ (and hence, the order of the points of $\bar Y$) so that the points $\bar y_1,\dots,\bar y_l$ appear in this sequence when one goes counterclockwise on the unit circle. For each $i$ denote by $R_{\bar Y}(i)$ ($L_{\bar Y}(i)$) the set of indices $j$ of vectors $\bar y_j$ that lie in the right (left) open half-plane incident to the oriented line spanned by the vector $\bar y_i$. Then the Pachner move which occurs when the vector $\beta(h)$ passes $\bar y_i$ from the right to the left, will be marked with the letter $a_{R_{\bar Y}(i),L_{\bar Y}(i)}\in \mathcal{A}_l^{l-1}$. The Pachner moves of the whole cycle of triangulation give the word $w_{\bar Y}=\prod_{i=1}^l a_{R(i),L(i)}$.

Thus,  in order to make the words like $\Phi(\alpha)$ independent of the resolutions of configurations of codimension $2$ we have to imply the that $w_{\bar Y}=1$, and hence we obtain a relation.

\begin{example}
Consider the standard Gale diagram of order $5$, see Fig.~\ref{fig:gale_d2_marked}. Then we have $R(1)=\{4,5\}$, $L(1)=\{2,3\}$, $R(2)=\{1,5\}$, $L(2)=\{3,4\}$ etc. The corresponding word is equal to
$$
w =a_{45,23}a_{15,34}a_{12,45}a_{23,15}a_{34,12}.
$$
\begin{figure}
\centering\includegraphics[width=80pt]{ 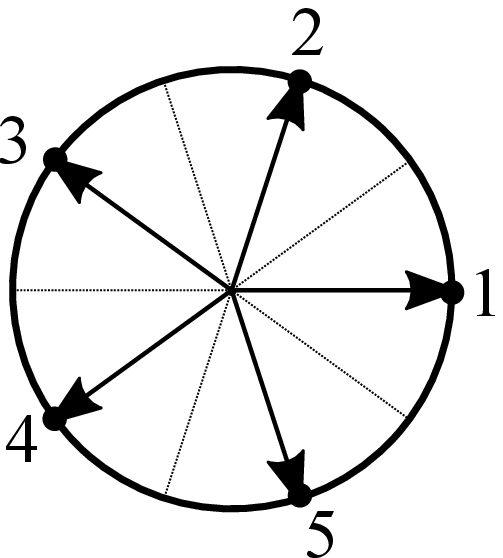}
\caption{Standard Gale diagram of order $5$}\label{fig:gale_d2_marked}
\end{figure}
\end{example}
We can give now the definition of $\Gamma_n^k$ groups.

\begin{definition}\label{def:Gamma_nk}
Let $4\le k\le n$. The group $\Gamma_n^k$ is the group with 
the generators
$$
\{ a_{P,Q}|P,Q\subset\{1,\dots,n\}, P\cap Q=\emptyset, |P\cup Q|=k, |P|,|Q|\ge 2\},
$$
and the relations:
\begin{enumerate}
\item $a_{Q,P}=a_{P,Q}^{-1}$;
\item (far commutativity) $a_{P,Q}a_{P',Q'}=a_{P',Q'}a_{P,Q}$ for each generators $a_{P,Q}$, $a_{P',Q'}$ such that
\begin{align*}
|P\cap(P'\cup Q')|&<|P|, & |Q\cap(P'\cup Q')|&<|Q|,\\
 |P'\cap(P\cup Q)|&<|P'|, & |Q'\cap(P\cup Q)|&<|Q'|;
\end{align*}
\item ($(k+1)$-gon relations) for any standard Gale diagram $\bar Y$ of order $k+1$ and any subset $M=\{m_1,\dots,m_{k+1}\}\subset \bar{n}$
$$ \prod_{i=1}^{k+1} a_{M_R(\bar Y,i), M_L(\bar Y,i)}=1, $$
where $M_R(\bar Y,i) = \{m_j\}_{j\in R_{\bar Y}(i)}$,  $M_L(\bar Y,i) = \{m_j\}_{j\in L_{\bar Y}(i)}$.
\end{enumerate}
\end{definition}

\begin{example}
Let us write the $(k+1)$-gon relations in $\Gamma_n^k$ for small $k$.

The group $\Gamma_n^4$ has one pentagon relation
$$
a_{m_4m_5,m_2m_3}a_{m_1m_5,m_3m_4}a_{m_1m_2,m_4m_5}a_{m_2m_3,m_1m_5}a_{m_3m_4,m_1m_2}=1.
$$

The group $\Gamma_n^5$ has two hexagon relations
\begin{multline*}
a_{m_5m_6,m_2m_3m_4}a_{m_1m_5m_6,m_3m_4}a_{m_1m_2,m_4m_5m_6}a_{m_1m_2m_3,m_5m_6}\cdot\\
\shoveright{a_{m_3m_4,m_1m_2m_6}a_{m_3m_4m_5,m_1m_2}=1,}\\
\shoveleft{a_{m_5m_6,m_2m_3m_4}a_{m_1m_5m_6,m_3m_4}a_{m_1m_2m_6,m_4m_5}a_{m_1m_2m_3,m_5m_6}\cdot}\\
a_{m_3m_4,m_1m_2m_6}a_{m_4m_5,m_1m_2m_3}=1.
\end{multline*}

The group $\Gamma_n^6$ has five heptagon relations
\begin{multline*}
a_{m_6m_7,m_2m_3m_4m_5}a_{m_1m_6m_7,m_3m_4m_5}a_{m_1m_2m_6m_7,m_4m_5}a_{m_1m_2m_3,m_5m_6m_7}\cdot\\
\shoveright{a_{m_1m_2m_3m_4,m_6m_7}a_{m_4m_5,m_1m_2m_3m_7}a_{m_4m_5m_6,m_1m_2m_3}=1,}\\
\shoveleft{a_{m_6m_7,m_2m_3m_4m_5}a_{m_1m_6m_7,m_3m_4m_5}a_{m_1m_2m_6m_7,m_4m_5}a_{m_1m_2m_3m_7,m_5m_6}\cdot}\\
\shoveright{a_{m_1m_2m_3m_4,m_6m_7}a_{m_4m_5,m_1m_2m_3m_7}a_{m_5m_6,m_1m_2m_3m_4}=1,}\\
\shoveleft{a_{m_6m_7,m_2m_3m_4m_5}a_{m_1m_6m_7,m_3m_4m_5}a_{m_1m_2m_7,m_4m_5m_6}a_{m_1m_2m_3m_7,m_5m_6}\cdot}\\
\shoveright{a_{m_1m_2m_3m_4,m_6m_7}a_{m_3m_4m_5,m_1m_2m_7}a_{m_5m_6,m_1m_2m_3m_4}=1,}\\
\shoveleft{a_{m_6m_7,m_2m_3m_4m_5}a_{m_1m_6m_7,m_3m_4m_5}a_{m_1m_2m_7,m_4m_5m_6}a_{m_1m_2m_3,m_5m_6m_7}\cdot}\\
\shoveright{a_{m_1m_2m_3m_4,m_6m_7}a_{m_3m_4m_5,m_1m_2m_7}a_{m_4m_5m_6,m_1m_2m_3}=1,}\\
\shoveleft{a_{m_5m_6m_7,m_2m_3m_4}a_{m_1m_6m_7,m_3m_4m_5}a_{m_1m_2m_7,m_4m_5m_6}a_{m_1m_2m_3,m_5m_6m_7}\cdot}\\
a_{m_2m_3m_4,m_1m_6m_7}a_{m_1m_4m_5,m_1m_2m_3m_7}a_{m_4m_5m_6,m_1m_2m_3}=1.
\end{multline*}
\end{example}

For groups $\Gamma_n^k$ we have a generalisation of Theorem~\ref{thm:Gamman5_invariant}.
\begin{theorem}\label{thm:Gamma_nk_invariant}
The map defined by the formula~\eqref{eq:Gamma_nk_braid_invariant} is a homomorphism $$\Phi\colon\pi_1(\tilde C_n(\R^{k-2}))\to\Gamma_n^k.$$
\end{theorem}
The proof repeats the arguments of the proof of Theorem~\ref{thm:Gamman5_invariant}. 

\section{The groups $\tilde{\Gamma}_n^k$} \label{sec:tilde_gamma}

We finish this part with a slight variation on the groups $\Gamma_n^k$. Geometrically speaking, here we consider oriented triangulations. Therefore, the indices of the generators of the groups are not independent and do not freely commute as was seen, for example, in the groups $\Gamma_n^5$ (see Definition \ref{def:Gamma_n5}, first relation).

To be precise, we introduce the following:

\begin{definition}\label{def:tildeGamma_nk}
Let $5\le k\le n$. The group $\tilde{\Gamma}_n^k$ is the group with the generators
$$
\{ a_{P,Q}|P,Q \; \hbox{--- cyclically oriented subsets of} \; \{1,\dots,n\}, P\cap Q=\emptyset, |P\cup Q|=k, |P|,|Q|\ge 2\},
$$
and the relations:
\begin{enumerate}
\item $a_{Q,P}=a_{Q',P'}$, where $Q=Q', P=P'$ as unordered sets, and as cyclically ordered sets $Q$ differs from $Q'$ by one transposition and $P$ differs from $P'$ by one transposition;
\item $a_{Q,P}=a_{P,Q}^{-1}$;
\item (far commutativity) $a_{P,Q}a_{P',Q'}=a_{P',Q'}a_{P,Q}$ for each generators $a_{P,Q}$, $a_{P',Q'}$ such that
\begin{align*}
|P\cap(P'\cup Q')|&<|P|, & |Q\cap(P'\cup Q')|&<|Q|,\\
 |P'\cap(P\cup Q)|&<|P'|, & |Q'\cap(P\cup Q)|&<|Q'|;
\end{align*}
\item ($(k+1)$-gon relations) for any standard Gale diagram $\bar Y$ of order $k+1$ and any subset $M=\{m_1,\dots,m_{k+1}\}\subset \bar{n}$
$$ \prod_{i=1}^{k+1} a_{M_R(\bar Y,i), M_L(\bar Y,i)}=1, $$
where $M_R(\bar Y,i) = \{m_j\}_{j\in R_{\bar Y}(i)}$,  $M_L(\bar Y,i) = \{m_j\}_{j\in L_{\bar Y}(i)}$.
\end{enumerate}
\end{definition}
\begin{wrapfigure}{R}{0.45\textwidth}
\begin{center}
 \includegraphics[width = 0.4\textwidth]{ 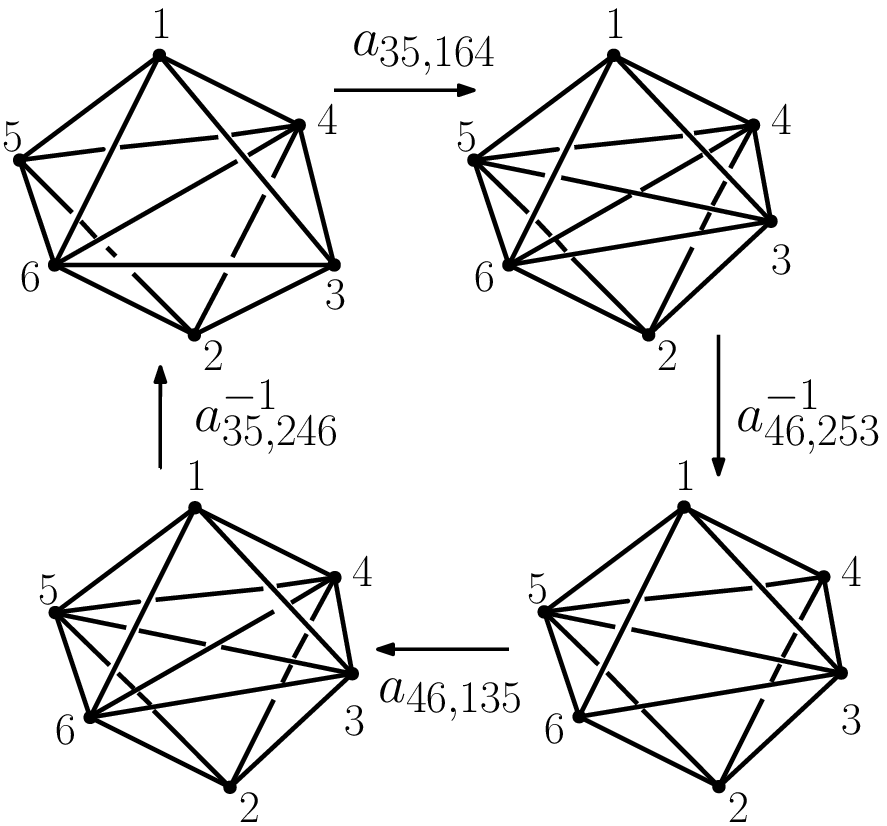}
\end{center}
\caption{A movement of a point around the configuration of four points on one circle}\label{exa_4word}
\end{wrapfigure}

These groups have the same connection to geometry and dynamics as the groups $\Gamma$ defined above. To illustrate that, consider the following dynamical system.
\begin{example}
	Let us have a dynamical system describing a movement of a point around the configuration of four points on one circle, see Fig.~\ref{exa_4word}. Such system may be presented as the word in the group $\tilde{\Gamma}_6^5$ $$w=a_{35, 164} a^{-1}_{46, 253} a_{46, 135} a^{-1}_{35, 246}.$$

	We can show that $\phi(\alpha)$ is nontrivial in the abelianisation $(\tilde{\Gamma}_6^5)_{ab}=\tilde{\Gamma}_6^5/[\tilde{\Gamma}_6^5,\tilde{\Gamma}_6^5]$ of the group $\tilde{\Gamma}_6^5$. Computer calculations show the group $(\tilde{\Gamma})_{ab}$ can be presented as the factor of a free commutative group with $120$ generators modulo $2\cdot 6!=1440$ relations that span a space of rank $90$ if we work over $\mathbb{Z}_2$. Adding the element $w$ to the relations increases the rank to $91$, so the element $w$ is nontrivial in $(\tilde{\Gamma}_6^5)_{ab}$, therefore it is nontrivial in $\tilde{\Gamma}_6^5$.

	Hence, we have encountered a peculiar new effect in the behaviour of $\tilde{\Gamma}_{n}^{5}$ which is not the case of $G_{n}^{k}$. Certainly, the abelianisation of $G_{n}^{k}$ is non-trivial and very easy to calculate since any generator enters each relation evenly many times. However, this happens not only for relations but for any words which come from braids. Thus, any abelianisations are trivial in interesting cases.
	
	For $\tilde{\Gamma}_{n}^{k}$, it is an interesting new phenomenon, and the invariants we have demonstrated so far are just the tip of the iceberg to be investigated further.
	
	Note that both for the groups $\Gamma_n^k$ and $\tilde{\Gamma}_n^k$ we have many invariants since the corank of those groups is big.
\end{example}


\input{unsolved.tex}

\input{appendum.tex}

\printindex

\bibliographystyle{hplain}
\bibliography{bibliography_new}

\end{document}

%% file: abstract.tex
Recently the first named author defined a 2-parametric family of groups $G_n^k$ \cite{gnk}. Those groups may be regarded as analogues of braid groups. 

Study of the connection between the groups $G_n^k$ and dynamical systems led to the discovery of the following fundamental principle: ``If dynamical systems describing the motion of $n$ particles possess a nice codimension one property governed by exactly $k$ particles, then these dynamical systems admit a topological invariant valued in $G_{n}^{k}$''.

The $G_n^k$ groups have connections to different algebraic structures, Coxeter groups, Kirillov-Fomin algebras, and cluster algebras, to name three. Study of the $G_n^k$ groups led to, in particular, the construction of invariants, valued in free products of cyclic groups. All generators of the $G_{n}^{k}$ groups are reflections which makes them similar to Coxeter groups and not to braid groups. Nevertheless, there are many ways to enhance $G_{n}^{k}$ groups to get rid of this $2$-torsion.

Later the first and the fourth named authors introduced and studied the second family of groups, denoted by $\Gamma_n^k$, which are closely related to triangulations of manifolds.

The spaces of triangulations of a given manifolds have been widely studied. The celebrated theorem of Pachner~\cite{pach} says that any two triangulations of a given manifold can be connected by a sequence of bistellar moves or Pachner moves. See also~\cite{GelKapZel,nab}; the $\Gamma_n^k$ naturally appear when considering the set of triangulations with the fixed number of points.

There are two ways of introducing the groups $\Gamma_n^k$: the geometrical one, which depends on the metric, and the topological one. The second one can be thought of as a ``braid group'' of the manifold and, by definition, is an invariant of the topological type of manifold; in a similar way, one can construct the smooth version.

In the present paper we give a survey of the ideas lying in the foundation of the $G_n^k$ and $\Gamma_n^k$ theories and give an overview of recent results in the study of those groups, manifolds, dynamical systems, knot and braid theories. 

%% file: introduction.tex
This survey paper  grew out of the lecture series in the Moscow State University in 2015--2016 in which Fedoseev, Kim, and Nikonov were attendees. As an introduction to the conceptual basis of the present paper and its background, Manturov reports the following.
\vspace{3mm}

A long time ago, when I first encountered knot tables and started unknotting knots ``by hand'', I was quite excited with the fact that some knots may have more than one minimal representative. 

\begin{figure}[h!]
\begin{center}
 \includegraphics[width = 4cm]{ 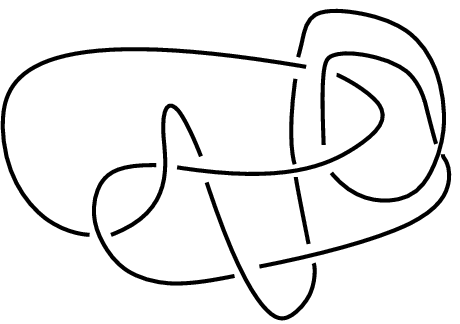}

\end{center}

\caption{Culprit knot}\label{culprit_knot}
\end{figure}


As I grew from an undergraduate student through my career, I learned, at least, the
following two things. 

1. In order to simplify diagrams of the unknot, it may be
necessary to make the diagrams more complicated, or {\em introduce} crossings before
{\em reducing the net number of crossings}, see Fig.~\ref{culprit_knot}.

\begin{figure}[t!]
\begin{center}
 \includegraphics[width = 40mm]{ 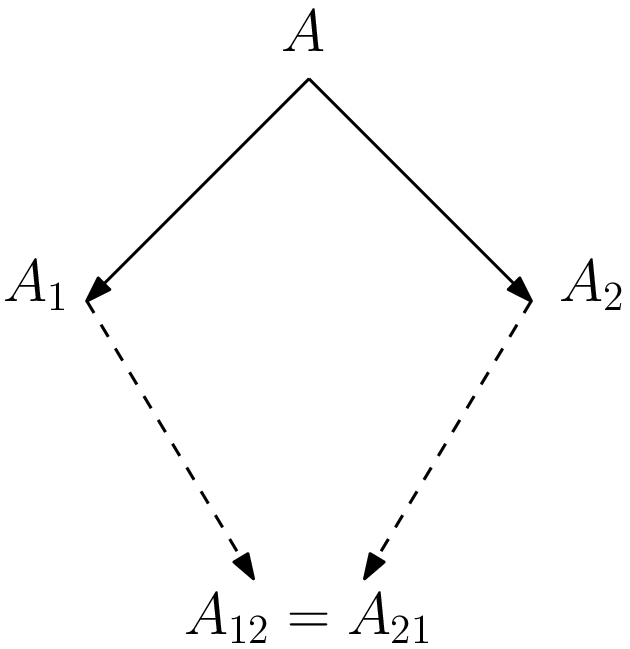}
\end{center}
\caption{The Diamond lemma}\label{diamond}
\end{figure}

2. In contrast, the word and
conjugation problems in free groups can be solved via a gradient descent algorithm.
Furthermore, the Diamond lemma (Fig.~\ref{diamond}) is a simple condition that guarantees the
uniqueness of minimal objects. As I learned and taught knot theory, I puzzled about
these contrasts between knot theory and (free) group theory.


Being a last year undergraduate and teaching a knot theory course for the first time,  I thought: ``Why do not we have a gradient descent (at least partially) in knot theory?''


By that time I knew about the Diamond lemma and solvability of many problems like word problem in groups by gradient descent algorithm. The van Kampen lemma and the Greendlinger theorem (see Sections~\ref{sec:van_kampen} and~\ref{sec:small_cancellations}) came to my knowledge much later.

Virtual knot theory \cite{Kauffman_virtual_knots, Kam, abstract, GPV}, which can be formally defined via Gau{\ss} diagrams (not necessarily planar), is a theory about knots in thickened surfaces $S_{g}\times I$ considered up to addition/removal of nugatory handles. It contains classical knot theory as a proper part: classical knots can be thought of as knots in the thickened sphere. Hence, virtual knots have a lot of additional information coming from the topology of the ambient space ($S_{g}$).
The topology of the ambient space gives more ``degrees of freedom'', more ``fundamental group structure'' which makes the topological theory closer to the theory of free groups.

In the study of Gau{\ss} diagrams, I rediscovered \cite{ManFreeKnot,knot_parity} the notion of free knots\footnote{A {\em free knot} is an equivalence class of Gau{\ss} diagrams without additional framing (arrows or signs) modulo formal Reidemeister moves} (introduced by Turaev \cite{tur} under the name of ``Homotopy classes of Gau{\ss} words''). Turaev conjectured that the theory was trivial. But the parity bracket invariant\footnote{The parity bracket invariant takes a diagram $D$ of a knot into the sum of all diagrams obtained from $D$ by all smoothings of its {\em even} crossings, taken modulo 2 and modulo the second Reidemeister move. For a rigorous definition, see for example~\cite{knot_parity}.} 
detects some non-trivial free knots. This invariant was constructed in such a way that all odd chord persisted and I got the formula~$$[K]=K$$ whenever $K$ is a diagram of a free knot with all chords odd where no two chords can be cancelled by a second Reidemeister move. The flat virtual that corresponds to the Kishino virtual is a knot of the type described above. So the parity bracket invariant is among the simplest ways to detect that the Kishino knot is, in fact, knotted.

Free knots are much simpler objects than virtual knots, nevertheless they admit powerful invariants of a new type; this motivates the whole construction. 

The deep sense of the formula $[K]=K$ can be expressed as follows: 

\begin{center}
{\em If a virtual diagram is complicated enough, then it realises itself.} \\
\end{center}

\begin{figure}
\centering\includegraphics[width=200pt]{ 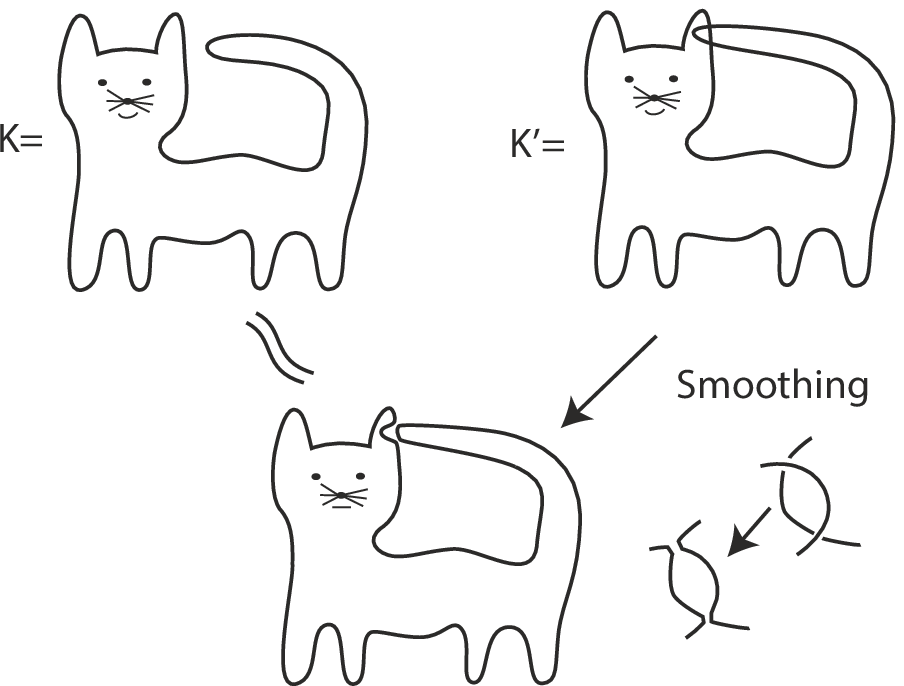}
\caption{A picture which is its own invariant}
\label{cat}
\end{figure}

Namely, $K$ on the left-hand side is {\em some} free knot diagram $K$, i.e. the knot itself, in other words, an equivalence class; and $K$ on the right-hand side is a {\em concrete graph}. Hence, if $K'$ is another diagram of the same knot, we shall have $[K']=K$ meaning that $K$ is obtained as a result of smoothing of $K'$. This principle is depicted on Fig. \ref{cat}.

In classical knot theory we lose the fact that {\em local minimality yields global minimality} or, in other words, {\em if a diagram is odd and irreducible} then it is {\em minimal} in a very strong sense. Not only one can say that $K'$ has larger crossing number than $K$, we can say that $K$ ``lives'' inside $K'$. Having these ``graphical'' invariants, we get immediate consequences about many characteristics of $K$. Actually, there are many other parities and many other brackets; for more details, see \cite{akimova,cob,IMN1,upcoming_book}.

Something similar to the above-mentioned principle ``local minimality yields global minimality'' appears in other situations as well: free groups or free products of cyclic groups, cobordism theory for free knots, or while considering other geometrical problems. For example, if we want to understand the genus of a surface where the knot $K$ can be realised, it suffices for us to look at the {\em minimal} genus where the concrete diagram $K$ can be realised, for the genus of any other $K'$ is a priori larger than that of $K$ \cite{Manturov_virtual_number}, see also \cite{IMN1}.

The diagrammatic approach takes me back to the time of my habilitation thesis. Once writing a knot theory paper and discussing it with Oleg Yanovich Viro, I wrote ``a classical knot is an equivalence class of classical knot diagrams modulo Reidemeister moves''. ``Well, --- said Viro, --- you are restricting yourself very much. Staying at this position, how can you prove that a non-trivial knot has a quadrisecant?''



The lack of homotopical information for {\em classical knots} can be easily recompensated at the level of {\em classical braids}: braid strands naturally capture a lot of homotopy information, and homotopy class of the braid naturally lives in the free group, the fundamental group of complement to the remaining strands. To make it fit into the above framework and capture free group information, one should change the notion of {\em braid crossing:} instead of looking at pairs of points having the same ordinate, one should look at triples of collinear points. Then we get closer to the free group theory.

\begin{figure}
\centering\includegraphics[width=400pt]{ 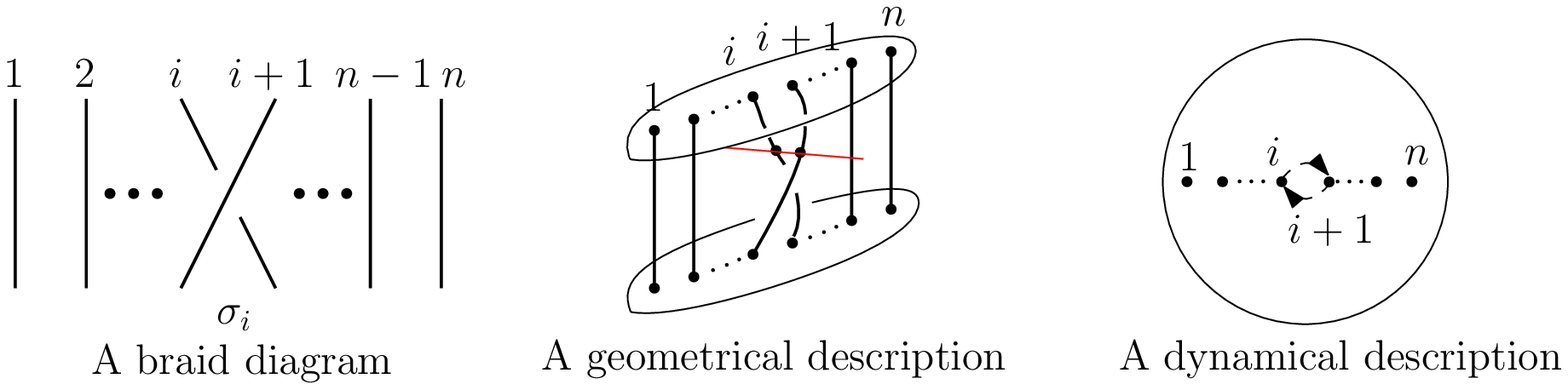}
\caption{Different interpretations of a braid}
\label{fig:braid_descriptions}
\end{figure}
\begin{figure}
\centering\includegraphics[width=300pt]{ 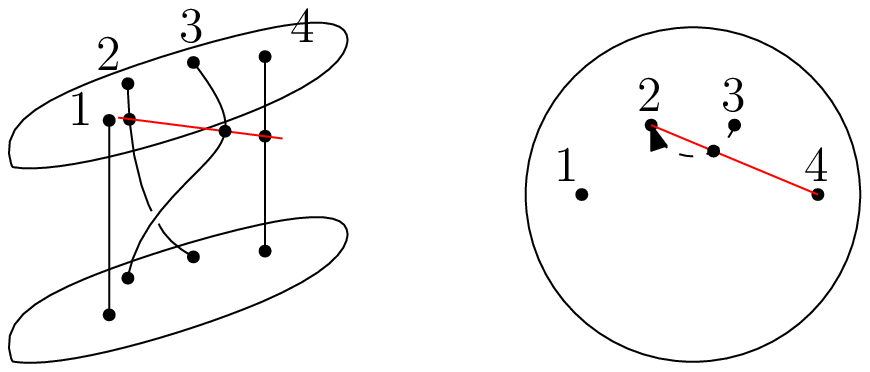}
\caption{A moment when three points appear on the same line}
\label{fig:exa_triple}
\end{figure}

Braids can be defined in many ways, see Fig.~\ref{fig:braid_descriptions}. Maybe the most familiar to the reader is the geometrical description of a braid as a set of strands in the three-dimensional space. But there is another, {\em dynamical} description. Given a (geometrical) braid, that is a set of $n$ strands connecting two parallel planes, take the bottom plane and start to lift it in thought in parallel. At each moment the plane will intersect the strands in $n$ points and the points will move as the plane ascends. Thus, the braid can be described as a movement of point in the plane, a dynamical system with $n$ points-particles. At different moments the points can form special configurations we call critical: some three points may become collinear, or four points lay on one circle, etc; see Fig.~\ref{fig:exa_triple}. The subject of study of this paper are chronicles of such critical events of dynamical systems. These chronicles are recorded in letters of some special groups known as $G_n^k$. 

These considerations led me to my initial preprint~\cite{gnk} 
and to an extensive study of braids and dynamical systems. This happened around New Year 2015.

Consider Artin's braids as dynamical systems of points in $\R^{2}$, but, instead of creating Reidemeister's diagram by projecting braids to a screen (say, the plane $Oxz$), I decided to look at those moments when some three points are {\em collinear}. This is quite a good property of a ``node'' which behaves nicely under generic isotopy. Let us denote such situations by letters $a_{ijk}$ where $i,j,k$ are numbers of points (this triple of numbers is unordered).

The situation when three points are collinear is a singular moment in the dynamics described by the braid. Such singularities form a codimension 1 set. In that sense, they are ``pretty common'' and give the generators of the group which is being constructed. To get the relations, one needs to consider the codimension 2 singularities. In particular, in case of points moving in $\R^2$, we shall look at the moments when {\em four points are collinear}. Considering them, we see that {\em the tetrahedron~(Za\-mo\-lod\-chi\-kov, {\rm see, for example, \cite{etin}}) equation} emerges. Namely, having a dynamics with a quadruple point and slightly perturbing it, we get a dynamics, where this quadruple point splits into four triple points.

Writing it algebraically, we get the Zamolodchikov equation:

$$S_{123}S_{124}S_{134}S_{234}=S_{234}S_{134}S_{124}S_{123}.$$

Taking arbitrary indices and passing to our standard generators, we get the following relation:

$$a_{ijk}a_{ijl}a_{ikl}a_{jkl}=a_{jkl}a_{ikl}a_{ijl}a_{ijk}.$$

\begin{definition}
The groups $G_{n}^{k}$ are defined as follows.
$$G_{n}^{k} = \langle a_{m} | (1),(2),(3) \rangle,$$
where the generators $a_{m}$ are indexed by all $k$-element subsets of $\{1,\dots, n\}$, the relation (1) means
\begin{center}
$(a_{m})^{2}=1$ for any unordered sets $m \subset \{1,\dots, n\}, Card(m)=k$;
\end{center}
(2) means
\begin{center}
$a_{m}a_{m'} = a_{m'}a_{m},$ if $Card(m\cap m') < k-1$;
\end{center}
and, finally, the relation (3) looks as follows. For every set $U \subset \{1,\dots, n\}$ of cardinality $(k+1)$, let us order all its $k$-element subsets arbitrarily and denote them by $m^{1}, \dots, m^{k+1}$. Then (3) is:

$$a_{m^{1}}\dots  a_{m^{k+1}} = a_{m^{k+1}}\dots  a_{m^{1}}.$$
\end{definition}

This situation with the Zamolodchikov equation happens almost everywhere, hence, I formulated the following principle: \\

\begin{center}
{\em If dynamical systems describing the motion of $n$ particles possess a nice codimension one property governed by exactly $k$ particles, then these dynamical systems admit a topological invariant valued in $G_{n}^{k}$}. \\
\end{center}

In topological language, it means that we get a certain homomorphism from some fundamental group of a configuration space to the groups $G_{n}^{k}$.

Note that somewhat later, we shall see another family of groups, $\Gamma_{n}^{k},k\ge 4$, satisfying a similar principle. Unfortunately, the relations for $k>4$ can not be written that easily (see Section~\ref{sec:higher_gamma} for full detail). 

Collecting all results about the $G_{n}^{k}$ groups, I taught a half-year course of lectures in the Moscow State University entitled ``Invariants and Pictures'' and a 2-week course in Guangzhou. 

Since that time, my seminar in Moscow, my students and colleagues in Moscow, Novosibirsk, Beijing, Guangzhou, and Singapore started to study the groups $G_{n}^{k}$, mostly from the following two points of view. \\

{\bf From the topological point of view,} which spaces can we study?

Besides the homomorphisms from the pure braid group $PB_{n}$ to $G_{n}^{3}$ and $G_{n}^{4}$, we can work with braids for different manifolds, real and projective spaces earlier studied for example by Birman~\cite{Bir}, van Buskirk~\cite{JvB}, and Fadell and Neuwirth~\cite{FN}. In our setting, they arise as follows. 

Of course, the configuration space $C(\R^{k-1},n)$ of ordered $n$-element subsets in $\R^{k-1}$ is simply connected for $k>3$ but if we take some {\em restricted} configuration space $C'(\R^{k-1},n)$ (see Section~\ref{section:realisation} for details), it will not be simply connected any more and will lead to a meaningful notion of higher dimensional braid.

What sort of the restriction do we impose? On the plane, we consider just braids, so we say that no two points coincide. In $\R^{3}$ we forbid collinear triples, in $\R^{4}$ we forbid coplanar quadruples of points.


As was mentioned by Jie Wu, consideration of restricted configuration spaces describing $k$-regular embeddings goes back to Carol Borsuk \cite{Borsuk} and P.\,L.~Chebysheff \cite{Chevyshev_poly}.


An interested reader may ask whether such braids exist not only for $\R^{k}$ (or $\R{}P^{k}$), but also for other manifolds, the question we shall touch on later, see Section~\ref{sec:manifold_invariants}. \\

{\bf From the algebraic point of view,} why are these groups good, how are they related to other groups, how to solve the word and the conjugacy problems, etc.?

It is impossible to describe all directions of the $G_{n}^{k}$ group theory in the survey, I just mention some of them.

For properties of $G_{n}^{k}$, we can think of them as $n$-strand braids with $k$-fold strand intersection. There are nice ``strand forgetting'' and ``strand deleting'' maps to $G_{n-1}^{k}$ and $G_{n-1}^{k-1}$, which resemble the formula ${n\choose k}={{n-1}\choose k}+{{n-1}\choose {k-1}}$. For details, see Fig.~\ref{(n,k)_(n-1,k)_(n-1,k-1)} in Section~\ref{sec:gnk_pictures}.

The groups $G_{n}^{k}$ have lots of epimorphisms onto free products of cyclic groups; hence, various invariants constructed from them are powerful enough and easy to compare. For example, the groups $G_{n}^{2}$ are commensurable with some Coxeter groups of special type (see Fig.~\ref{Gn2_Cox_Cay} in Section~\ref{sec:diagrams_prelim} for an example), which immediately solves the word problem for them. Here the groups are different, but there is a certain bijection, which allows to study the properties of one group looking at the other. This bijection is the result of the ``rewriting algorithm'', see Section~\ref{chap:gnk_homomorphisms} for details.

As Diamond lemma works for Coxeter groups, it works for $G_{n}^{2}$, and in many other places. \\

In general, we may say that the words from $G^k_n$ for $n=k+1$ can be ``realised'' as braids on $\RP^{k-1}$ describing dynamics of $k+1$ points, and for the case of $n>k+1$ we should investigate partial flag manifolds. That leads to Grassmanians and resolution of singularities.

After a couple of years of study of $G_{n}^{k}$, I understood that this approach was still somewhat restrictive. We can study braids, we can invent braids in $\R^{n}$ and $\R{}P^{n}$, but what if we consider just braids on a $2$-surface? What can we study then? The property ``three points belong to the same line'' is not quite good even in the metrical case because even if we have a Riemannian metric on a $2$-surface of genus $g$, there may be infinitely many geodesics passing through two points. Irrational cables may destroy the whole construction.

It turns out that the $G_{n}^{k}$-main principle is not very much universal in non-simply connected  when we have geodesics having multiple intersections and so on. In those situations it is natural to replace this principle by its ``local'' counterpart: we look at ``closest'' points satisfying some codimension one condition. When we start studying the corresponding groups, $\Gamma_{n}^{k}$, we see that they are naturally related to spaces of triangulations, polytopes, cluster algebras, and many other areas of mathematics.


Assume we have a collection of points in a $2$-surface and seek the $G^4_n$-property: four points belong to the same circle. We add just one condition: there must be no other points inside the circle. This minor modification reveals an unexpected connection to combinatorial topology through Vorono\"{\i} tilings and Delaunay triangulations. The potential of this connection is still to be comprehended.

Consider a $2$-surface of genus $g$ with $N$ points on it. We choose $N$ to be sufficiently large and put points in a position to form the centers of Vorono\"{\i} cells. It is always possible for the sphere $g=0$, and for the plane we may think that all our points live inside a triangle forming a Vorono\"{\i} tiling of the latter.

\begin{figure}[h]
\begin{center}
 \includegraphics[width = 80mm]{ 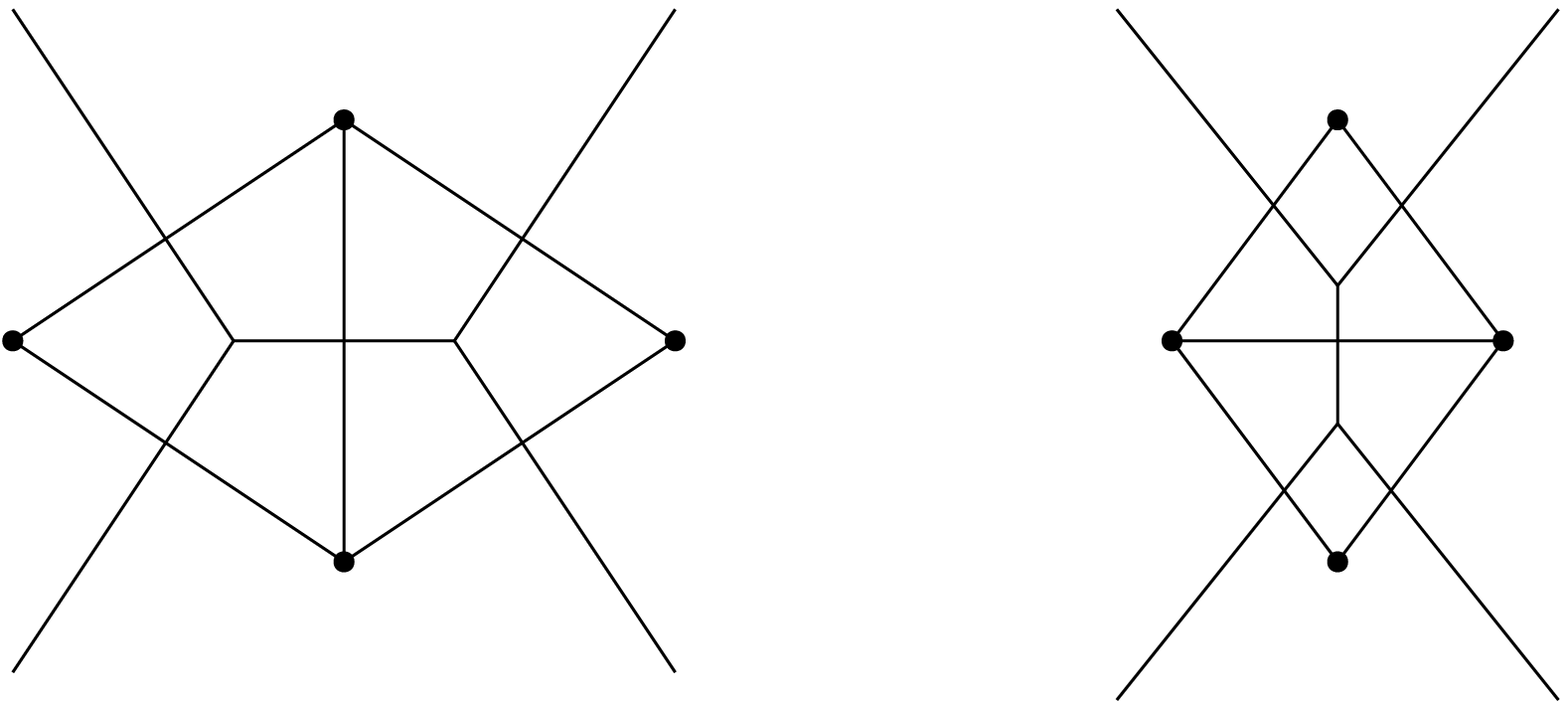}
\end{center}
\caption{Vorono\"{\i} tiling change}\label{tiling_change}
\end{figure}

We are interested in those moments when the combinatorics of the Vorono\"{\i} tiling changes, see~Fig.~\ref{tiling_change} (here the bold points denote the vertices of Delaunay triangulation, dual to the Vorono\"{\i} tiling shown in thin lines; the change of the tiling corresponds to the moment when the bold vertices appear on the same circle).

This corresponds to a flip, the situation when four {\em nearest points belong to the same circle}. This means that no other point lies inside the circle passing through these four ones. The most interesting situation of codimension $2$ corresponds to five points belonging to the same circle.

This leads to the relation:~$$d_{1234}d_{1245}d_{2345}d_{1235}d_{1345}=1.$$

Note that unlike the case of $G_{n}^{4}$, here we have five terms, not ten (as in the group $G_n^4$). What is the crucial difference? The point is that if we have five points in the neighbourhood of a circle, then every quadruple of them appears to be on the same circle {\em twice}, but one time the fifth point is outside the circle, and one time it is inside the circle. We denote the set $\{1,\dots, n\}$ by $\bar{n}$ and introduce the following

\begin{definition}
The group $\Gamma_{n}^{4}$ is the group given by group presentation generated by $\{ d_{(ijkl)}~|~ \{i,j,k,l\} \subset \bar{n}, |\{i,j,k,l\}| = 4\}$ subject to the following relations:

\begin{enumerate}
\item $d_{(ijkl)}^{2} = 1$ for $\{ i,j,k,l\} \subset \bar{n}$,
\item $d_{(ijkl)}d_{(stuv)} = d_{(stuv)}d_{(ijkl)}$, for $| \{i,j,k,l\} \cap \{s,t,u,v\} | < 3$,
\item $d_{(ijkl)}d_{(ijlm)}d_{(jklm)}d_{(ijkm)}d_{(iklm)} = 1$ for distinct $i,j,k,l,m$,
\item $d_{(ijkl)}=d_{(kjil)}=d_{(ilkj)}=d_{(klij)}=d_{(jkli)}=d_{(jilk)}=d_{(lkji)}=d_{(lijk)}$, for distinct $i,j,k,l$.
\end{enumerate}
\end{definition}

Just like we formulated the $G_{n}^{k}$ principle, here we could formulate the $\Gamma_{n}^{k}$
principle in whole generality, but we restrict ourselves with several examples.



\begin{figure}[h]
\begin{center}
 \includegraphics[width = 90mm]{ 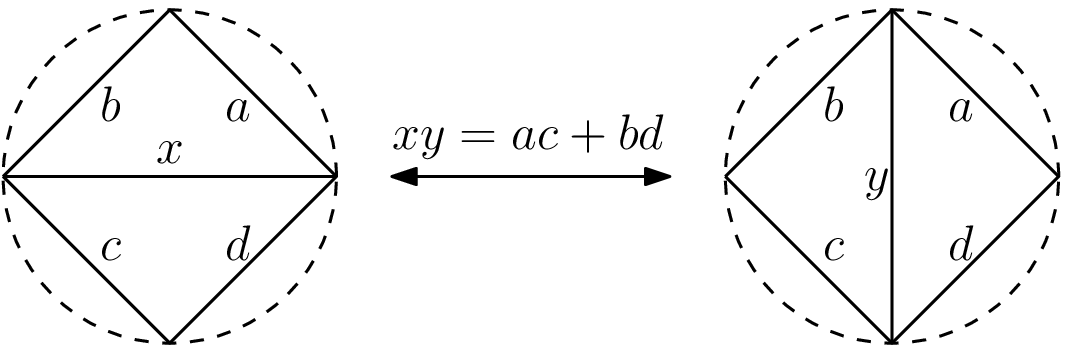}
\end{center}
\caption{The Ptolemy relation}\label{flip_ind}
\end{figure}

It turns out that groups $\Gamma_{n}^{4}$ have nice presentation coming from
{\em the Ptolemy relation} and cluster algebras. The Ptolemy relation
$$xy=ac+bd$$ says that the product of diagonals of an inscribed quadrilateral
equals the sum of products of its opposite faces, see Fig. \ref{flip_ind}.

We can use it when considering
triangulations of a given surface: when performing a flip, we replace one diagonal $(x)$ with the other diagonal $(y)$ by using this relation.
It is known that if we consider all five triangulations of the pentagon and
perform five flips all the way around, we return to the initial triangulation
with the same label.

This well known fact gives rise to presentations of $\Gamma_{n}^{4}$.

Thus, by analysing the groups $\Gamma_{n}^{4}$, we can get

\begin{enumerate}
\item Invariants of braids on $2$-surfaces valued in polytopes;
\item Relations to groups $G_{n}^{4}$;
\item Braids on 3-manifolds.
\end{enumerate}

Going slightly beyond, we can investigate braids which correspond to dynamics of $n$ points in $\R^{3}$, and the configuration space of polytopes in $\R^{3}$ (note that the usual Artin braids describe the dynamics of $n$ points on the plane; here we understand the word ``braid'' in a broader sense). Indeed, when studying bifurcations of simplicial polytopes in tree-space is similar to studying braids: crucial moments when four points belong to the same facet and the combinatorial structure changes are similar to those moments when four points belong to the same circle with no points inside.

We will not say much in the introduction about the groups $\Gamma_{n}^{k}$ for $k>4$ (they are discussed and defined in detail in Section~\ref{sec:higher_gamma}), but the main idea behind them is the following.
\begin{enumerate}
\item Generators (codimension 1) correspond to simplicial $(k-2)$-polytopes with $k$ vertices;
\item The most interesting relations (codimension 2) correspond to $(k-2)$-polytopes with $k+1$ vertices.
\item Roughly speaking, the most interesting relation in these groups deals with ``higher associativity''.
\end{enumerate}

It would be extremely interesting to establish the connection between $G_{n}^{k}$ with the Manin--Schechtmann ``{\em higher braid groups}'' \cite{ManinSchechtmann} where the authors study the fundamental group of complements to some configurations of {\em complex hyperplanes}.

It is also worth mentioning, that the relations in the group $G_n^2$ resemble the relations in {\em Kirillov--Fomin algebras} \cite{KirFom}. For that reason it seems interesting to study the interconnections between those objects.


The survey culminates with the {\em manifold of triangulations} and invariants of arbitrary manifolds.

The space of triangulations of a manifold with $N$ points has been considered since a long ago. It is known that triangulations are related to each other by Pachner moves.

The structure on the set of triangulations, its relations to associahedra, permutohedra, etc. has been widely studied, see \cite{nab,GelKapZel} and references therein. This structure is related to lots of different areas of mathematics. Say, Stasheff's polytopes were first introduced by Stasheff (following Milnor and Adams) in order to construct H-associativity on topological spaces.

Interestingly, the space of all triangulations of a sphere, the space of all simplicial polytopes of a given dimension, and similar objects were considered by many authors in many directions, but, to the best of my knowledge, the ``canonical'' group action on such objects has not been introduced.

In the Appendix of the present survey we introduce {\em the braid group for manifolds of arbitrary dimension.} If the manifold were just $\R^{k-1}$ ($\RP^{k-1}$) we would consider the restricted configuration space by deleting local codimension $2$ set from the set of whole triangulations. We get an (open) manifold, whose fundamental groups can be studied by means of groups $\Gamma_{n}^{k}$.

Depending on the point of view, {\em geometrical} or {\em topological}, we can define different manifolds of triangulations with different braid groups. These groups are naturally mapped to the corresponding $\Gamma_{n}^{k}$. Such spaces, triangulations, flips, cluster algebras were studied in \cite{KirFom} (see also references therein) but the group structure was missing.

Formulating this principle concisely, we seek a nice codimension $2$ property which would create a non-trivial group. A nice codimension $1$ property gives rise to the local condition, which guarantees that these groups have a nice map to $\Gamma_{n}^{k}$; the algebraic nature of $\Gamma_{n}^{k}$ allows one to ``take a picture'' of the discriminant set to get various presentations of the groups into algebraic objects (keep in mind the Pachner move).
\vspace{3mm}

In summary, we have constructed two families of groups $G_{n}^{k}$ and $\Gamma_{n}^{k}$ which are friendly to various well-known objects (braid groups, Coxeter groups, cluster algebras, configuration spaces, spaces of polytopes), touch different theories and areas of study (as diverse as picture-valued invariants and combinatorics of the Diamond lemma), lead to the construction of new invariants of manifolds, to the invention of new objects (braids for higher dimensional spaces), and motivate new areas of research in topology, algebra and geometry. \newline

\begin{figure}[t!]
\centering\includegraphics[width=440pt]{ 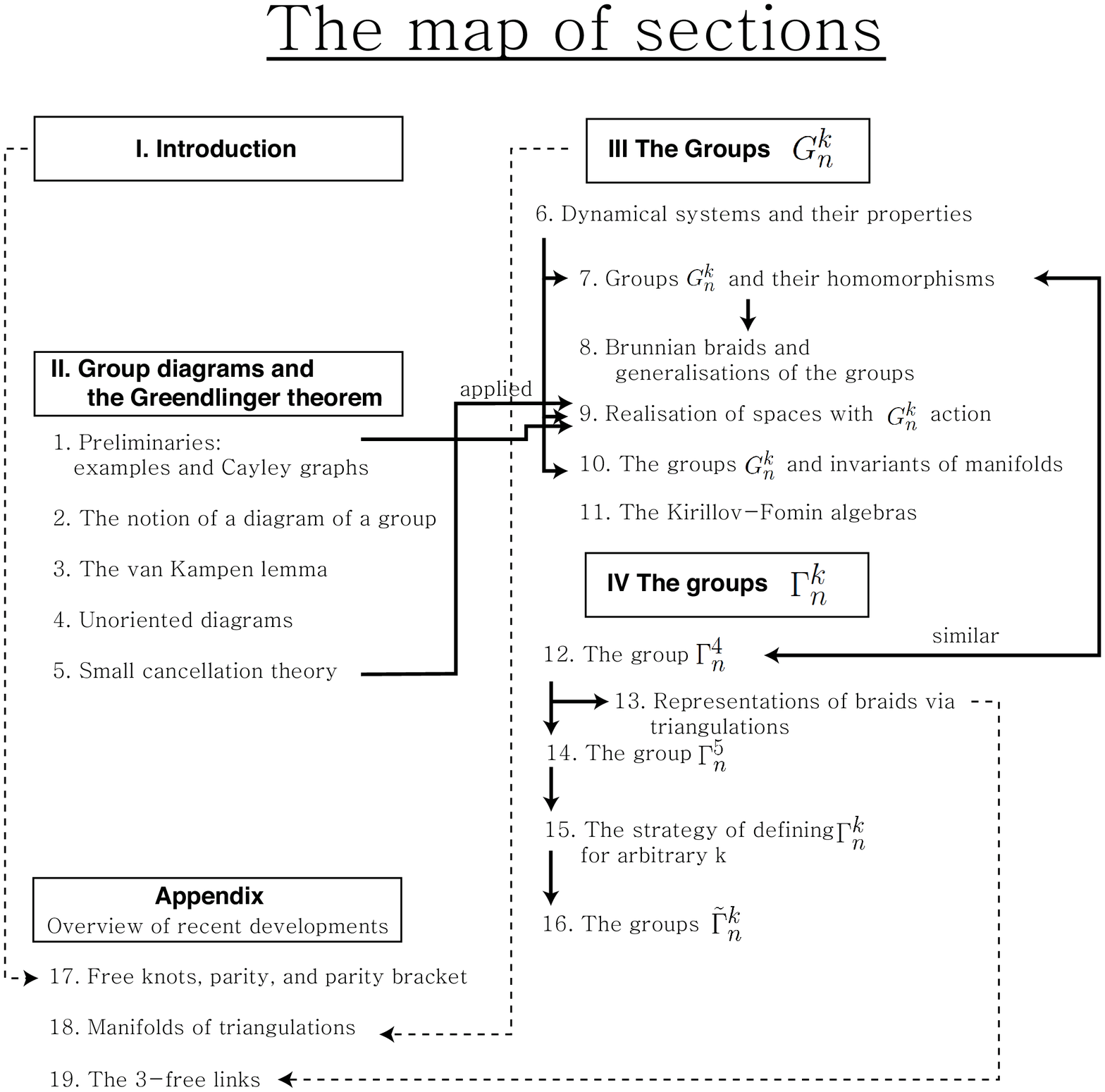}
\caption{The graph of the sections of the paper}
\label{fig:section_graph}
\end{figure}

The structure of the present text is as follows (it is summarised and visualised in Fig.~\ref{fig:section_graph}). The text consists of four parts and an appendix. 

The first part is actually the present Introduction, which gives an overview of what the present paper is about.

In the second part of the text, we give an overview of the diagrammatic language of group presentation, and provide certain approaches to working with groups using the geometric language, most prominently --- the Greendlinger theorem (see Theorem~\ref{th:greendlinger}), van Kampen lemma (see Lemma~\ref{lem:van_kampen}), and the small cancellation theory. This part may be omitted by a reader fluent in the group-diagrammatic language and wishing to delve into the $G_n^k$ and $\Gamma_n^k$ theory. This part breaks down into the following sections.

In Section~\ref{sec:diagrams_prelim} we give several examples, which illustrate the principles defined later in Part~\ref{part:group_diagrams}. 
In Section~\ref{sec:oriented} one can find explicit definitions of group diagrams and necessary results. 
In Section~\ref{sec:van_kampen}, certain theorems are stated, which give criteria to distinguish elements of a group from the identity by using diagrammatic descriptions for groups. In particular, here the celebrated van Kampen lemma~\ref{th:van_kampen_strong} is given.
In Section~\ref{sec:unoriented}, a modification of the van Kampen lemma for unoriented diagrams is introduced.
In Section~\ref{sec:small_cancellations} the notion of small cancellation conditions are formulated, which can be applied to solve ``the word problem'' for groups, which will be applied to the $G_n^k$-theory later in Section~\ref{section:realisation}. 

The third part of the text is dedicated to the $G_{n}^{k}$ theory, which was defined by the first named author in 2015. 

In Section~\ref{sect:dyn_sys} we discuss dynamical systems and define the $G_{n}^{k}$ groups, which are closely related to dynamical systems. In addition, the Main principle is formulated, on which the $G_{n}^{k}$-theory is based. 
In Section~\ref{chap:gnk_homomorphisms}, basing on the Main principle from Section~\ref{sect:dyn_sys}, group homomorphisms from pure braid groups to $G_{n}^{3}$ and $G_{n}^{4}$ are constructed. Moreover, mappings from $G_{n}^{k}$ to a free product of copies of $\mathbb{Z}_{2}$ are presented, which give an invariant for pure braids. In the end of Section~\ref{chap:gnk_homomorphisms}, crossing numbers are studied by using the invariant referred above. 
In Section~\ref{sec:brunnian}, we formulate a shortcoming of the invariant defined in Section~\ref{chap:gnk_homomorphisms} appearing for Brunnian braids. To overcome this difficulty, some modifications of groups $G_{n}^{2}$ and $G_{n}^{3}$ are constructed. By using these modifications one can define invariants for pure braids valued in a free product of several copies of $\mathbb{Z}_{2}$ and it is shown that the new invariant more or less abolishes the weakness of the earlier invariant for Brunnian braids. 
In Section~\ref{section:realisation}, the realisability of word from $G_{n}^{k}$ by dynamical systems is studied by applying topological and algebraic principles, which were introduced in Part~\ref{part:group_diagrams} of the paper. 
In Section~\ref{sec:manifold_invariants}, the $G_{n}^{k}$ groups are applied to study higher dimensional manifolds. 
In Section~\ref{sec:kir-fom} one can find short overview of an algebraic structure studied by S.~Fomin and A.~Kirillov and its similarities with $G_n^k$ groups. 

In the fourth part of the text, another family of groups $\Gamma^k_n$ is introduced, which is motivated by studying triangulations of surfaces and attempting to look at their transformations locally. In contrast, the generators of the groups $G^k_n$ measure interactions among all particles at given instances of their dynamics, and hence these groups are constructed from more global information.

In Section~\ref{sec:gamma_n4}, the group $\Gamma_{n}^{4}$ is defined by group presentation and several homomorphisms from pure braid groups to $\Gamma_{n}^{4}$ and products thereof are given. 
In Section~\ref{sec:gamma_presentations} some representations of pure braids by using certain ``singular configurations'' are introduced, which are related to the $\Gamma_{n}^{4}$ structure. 
In Section~\ref{sec:gamma_n_5}, the group $\Gamma_{n}^{5}$ is defined by group presentation and a homomorphism from $\pi_1(\tilde{C}_n(\mathbb{R}^3))$ to $\Gamma_{n}^{5}$ is defined. 
In Section~\ref{sec:higher_gamma}, the notion of $\Gamma_{n}^{4}$ and $\Gamma_{n}^{5}$ are generalised to $\Gamma_{n}^{k}$ for any pair $n$ and $k$ with $n \ge k$ and homomorphisms from $\pi_1(\tilde{C}_n(\mathbb{R}^{k-2}))$ to $\Gamma_{n}^{k}$ are defined. 
In Section~\ref{sec:tilde_gamma}, the groups $\Gamma_{n}^{k}$ are slightly modified in context of oriented triangulations. 

In the end of the main text of the paper one can find some unsolved problems and directions of the further research.

Since this paper was originally submitted, a number of new directions and results in the fields within the scope of the present survey were obtained and explored. We touch some of them in the Appendix of the present survey.  

In Section~\ref{sec:free_knots}, we give a brief overview of free knot theory and parity theory which are referred to throughout the text, and which demonstrate a number of effects appearing in the theory of $G_n^k$ groups and their geometry and topology.
In Section~\ref{sec:manifold_of_triang}, the space of triangulations of a given manifold is considered. A homomorphism from the fundamental group of the space of triangulations of a manifold of dimension $d$ to $\Gamma_{n}^{d+2}$. 
In Section~\ref{sec:3-free} we study the so-called {\em 3-free links}, which are the geometric counterpart of a modification of the group $G_{n}^{3}$. A picture-valued invariant of conjugacy classes of closed braids is constructed using the 3-free links.
In Section~\ref{sec:domino_tiling} we describe {\em domino tilings}, a type of surface partitioning resembling triangulations. Then we introduce two groups, arising from domino tiling, and discuss their similarities with the groups $G_n^k$ and $\Gamma_n^k$. \newline

The authors would like to express their heartfelt gratitude to L.\,A~Bokut', H.~Boden, J.\,S.~Carter, A.\,T.~Fomenko, S.\,G.~Gukov, Y.~Han, D.\,P.~Ilyutko, \fbox{V.\,F.\,R.~Jones}, A.\,B.~Karpov, R.\,M.~Kashaev, L.\,H.~Kauffman, M.\,G.~Khovanov, A.\,A.~Klyachko, P.\,S.~Kolesnikov, I.\,G.~Korepanov, Fenglin Ling, S.\,V.~Matveev, A.\,Yu.~Olshanskii, W.~Rushworth, G.\,I.~Sharygin, V.\,A.~Vassiliev, Jun Wang, J.~Wu and Zerui~Zhang for their interest and various useful discussions on the present work. We are grateful to Efim~I.~Zelmanov for pointing out the resemblance between the groups $G_{n}^{2}$ and Kirillov-Fomin algebras. \\

The first named author was supported by the development program of the Regional Scientific and Educational Mathematical Center of the Volga Federal District, agreement No. 075-02-2020 and by the Russian Foundation for Basic Research (grants No. 20-51-53022, 19-51-51004).  The second named author was supported by the Russian Foundation for Basic Research (grants No. 20-51-53022, 19-51-51004 and 19-01-00775-a). The third named author was supported by the Russian Foundation for Basic Research (grants No. 20-51-53022, 19-51-51004). The fourth named author was supported by the Russian Foundation for Basic Research (grants No. 20-51-53022 and 19-51-51004). 

%% file: greendlinger.tex
In the present section we discuss diagrammatic language of group description. This approach was first discovered and used by van Kampen \cite{vK}. The essence of his discovery was an interconnection between combinatorially-topological and combinatorially-group-theoretical notions. The gist of this approach is a presentation of groups by flat diagrams (that is, geometrical objects, --- flat complexes, --- on a plane or other surface, such as a sphere or a torus). We review this theory following \cite{Olshanskii}.

\section{Preliminaries: examples and Cayley graphs}
\label{sec:diagrams_prelim}

\subsection*{Preliminary examples}
First, let us present several examples of the principle, which will be rigorously defined later in this section.

\begin{example}
	Consider a group $G$ with relations $a^3=1$ and $bab^{-1}=c$. Clearly in such group we have $c^3=1$. This fact can be seen on the following diagram, see Fig. \ref{fig:triangle_group_diagram}.
	
\begin{figure}[h!]
\begin{minipage}{.45\textwidth}
	\begin{center}
			\includegraphics[width =.7\textwidth]{ 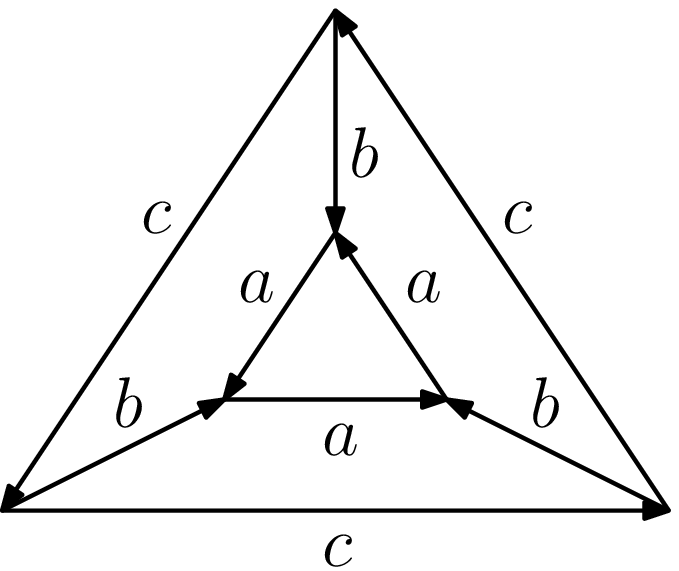}
	\end{center}
	\vspace{-0.1cm}
		\caption{Diagrammatic view of the $c^3=1$ relation} \label{fig:triangle_group_diagram}
	\end{minipage}
	\hspace{5mm}
	\begin{minipage}{.45\textwidth}
	\begin{center}
			\includegraphics[width =1\textwidth]{ 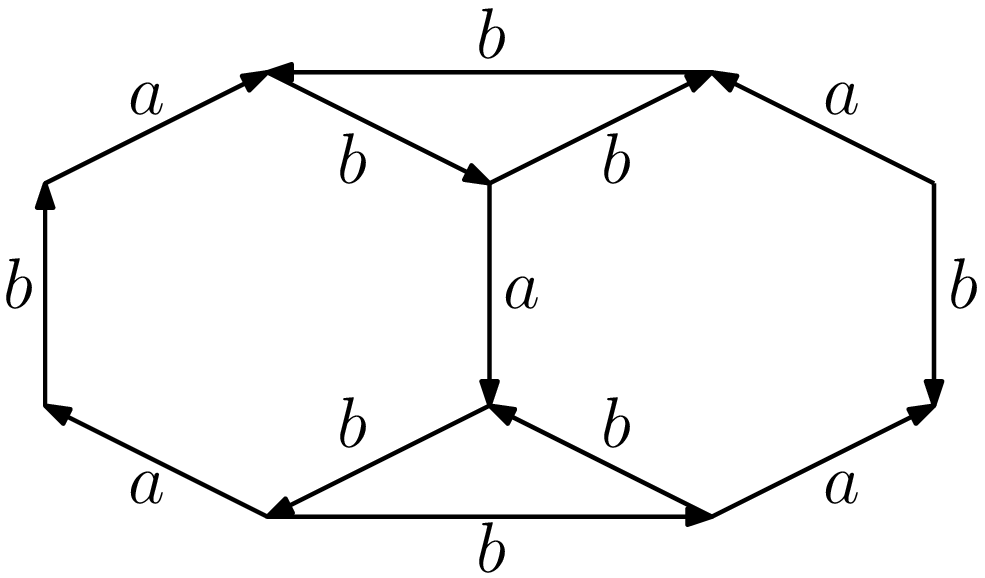}
	\end{center}
		\caption{Proof of the claim given in Example \ref{ex:second_group}} \label{fig:second_group_diagram}
	\end{minipage}
\end{figure}

	
	 Indeed, if we go around the inner triangle of the diagram, we get the relation $a^3=1$ (to be precise, the boundary of this cell gives the left-hand side of the relation; if we encounter an edge whose orientation is compatible with the direction of movement, we read the letter the edge is decorated with, otherwise we read the inverse letter; in this example we fix the counterclockwise direction of movement). Similarly, the quadrilaterals glued to the triangle all read the second given relation $c^{-1}bab^{-1}=1$. Now if look at the outer boundary of the diagram, we read $c^3=1$, and that is what we need to prove.
\end{example}

This simple example gives us a glimpse  of the general strategy: we produce a diagram, composed of cells, along the boundary of which given relations can be read. Then the outer boundary of the diagram gives us a new relation which is a consequence of the given ones.

Let us consider a bit more complex example of the same principle.

\begin{example} \label{ex:second_group}
	Consider a group $G$ where the relation $x^{3}=1$ holds for every $x\in G$. It is a well-known theorem that in such a group every element $a$ lies in some commutative normal subgroup $N\subset G$. Such situation arises, for example, in link-homotopy.
	
	To prove this fact it is sufficient to prove that any element $y=bab^{-1}$ conjugate to $a$ commutes with $a$. If that were the case, the subgroup $N$ could be constructed as the one generated by all the conjugates of $a$.
	
	So we need to prove that for every $b\in G$ the following holds: $$a(bab^{-1})=(bab^{-1})a$$ or, equivalently $$abab^{-1}a^{-1}ba^{-1}b^{-1}=1.$$
	This equality can be read walking clockwise around the outer boundary of the diagram in Fig. \ref{fig:second_group_diagram} composed of the relations $b^3=1, (ab)^3=1$, and $(a^{-1}b)^3=1$. 
\end{example} 
	

\subsection*{Cayley graphs}

Before delving into the notion of group diagrams, let us briefly discuss another geometrical description of groups: the Cayley graphs. Consider a group $G$ with a presentation $$G=\langle\mathcal{A}|R=1, R\in\mathcal{R}\rangle.$$ We can represent this group by an oriented graph $\Gamma$, where each vertex of the graph corresponds to an element of the group, and each oriented edge corresponds to a generator of the group. One vertex is distinguished (we denote it by $O$) and corresponds to the neutral element of the group. By definition, two vertices $v, v'$ of a Cayley graph of the group $G$ are connected by an edge, corresponding to a generator $a$, and oriented from the vertex $v$ to the vertex $v'$, if and only if for the elements of the group $g_v, g_{v'}\in G$, corresponding to the vertices $v$ and $v'$, respectively, the equality $g_{v'}=g_v a$ holds.

Now consider a path from the distinguished vertex $O$ to some vertex $v$ corresponding to an element $g_v\in G$. Walking along this path we ``read'' a word $w$ from $G$ (each edge corresponds to a generator; if the orientation of the edge is compatible with the direction of walking, we record the generator, otherwise, we record its inverse). This word $w$ represents the element $g_v$. Naturally, many different paths may connect the vertices $O$ and $v$. The words given by them are {\em equal} in the group $G$. In particular, trivial words (including the relations) correspond to loops connecting the vertex $O$ with itself.

\begin{figure}[h!]
\centering\includegraphics[width=260pt]{ 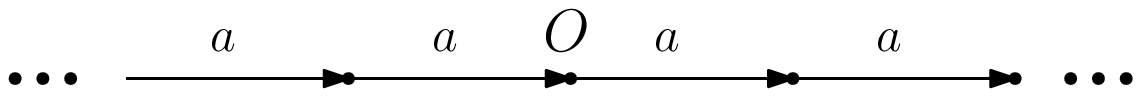}
\caption{An example of a Cayley graph}
\label{fig:cayley_example}
\end{figure}

\begin{example} 
Consider the graph depicted in Fig.~\ref{fig:cayley_example}. That is the graph of the group $\mathbb{Z} = \langle a \rangle$: each edge corresponds to the only generator of the group, $a$. Each vertex of the graph may be connected with the distinguished vertex $O$ by a path of the form $a^{\varepsilon_1}\dots a^{\varepsilon_n}$, where each $\varepsilon_i\in\{1,-1\}$. This path naturally corresponds to an element $\sum_{i=1}^n \varepsilon_i\in\mathbb{Z}$. 	
\end{example}

\begin{figure}[h]
\begin{center}
\includegraphics[width = 8cm]{ 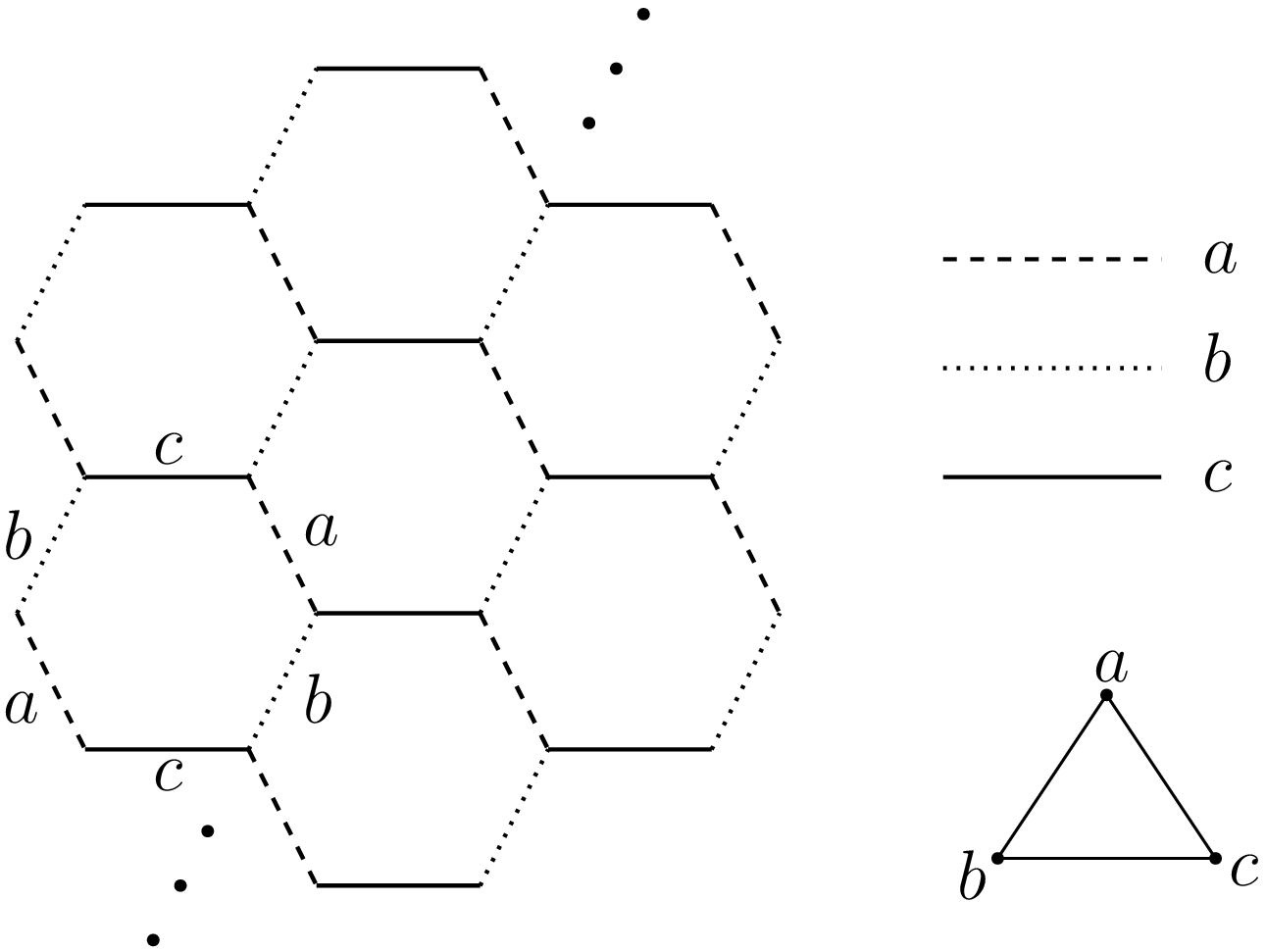}
\end{center}
\caption{The Cayley graph of the group $G_{3}^{2}$ and the Coxeter group $C(3,2)$}\label{Gn2_Cox_Cay}
\end{figure}

Note that one can ``see'' different groups in one and the same Cayley graph. For example, Fig.~\ref{Gn2_Cox_Cay} depicts the Cayley graph of the group $G_3^2$ (to be defined in fullness in Section~\ref{sec:gnk} and in example~\ref{ex:g32} in particular; for now we just need to know that this group has a presentation $G_3^2=\langle a,b,c\,|\,a^2=b^2=c^2=1, \, abcabc=1 \rangle$) and one of the Coxeter group $C(3,2)$. Note that here we use {\em unoriented} graphs; that is useful when dealing with groups with the relations $a^2=1$ for all generators: since $a=a^{-1}$ in such groups, there is no difference, whether we walk along the orientation of an edge, or in the other direction.

Recall, that a {\em Coxeter group $C(n,2)$} is defined as follows. Consider ${n \choose 2}$ generators $b_{ij}$, where $i$ and $j$ go over all unordered pairs of distinct numbers from $\{1,\dots,n\}$, and the following set of relation: $$(b_{ij}b_{ik})^3=1, \;\;\; |\{i,j,k\}|=3,$$ $$b_{ij}b_{kl}=b_{kl}b_{ij}, \;\;\; |\{i,j,k,l\}|=4,$$ $$b_{ij}^2=1, \;\;\; i\neq j.$$ Hence, the presentation of the simplest Coxeter group $C(3,2)$ may be written in the form $C(3,2)=\langle a,b,c \, |  \, a^2=b^2=c^2=1, (ab)^3=(ac)^3=(bc)^3=1 \rangle$.

These graphs are different, but the graph of the Coxeter group is ``inside'' the graph of the group $G_3^2$. \\

In a sense, Cayley graphs present a {\em 1-dimensional} diagrammatic method for describing a group. The group diagrams we discuss in the next section give a {\em 2-dimensional} formalism of group description.

\section{The notion of a diagram of a group}
\label{sec:oriented}

Now we can move on to the explicit definitions of group diagrams and the overview of necessary results in that theory.

In this section we follow closely the presentation \cite{Olshanskii}. A partition, $\Delta$, of a surface $S$ into cells (each homeomorphic to a disc) will be called a map on $S$ for short. For some particular surfaces we will also use special names; for example, a map on a disc will be called a {\em disc map}, on an annulus an {\em annular map}\index{Map!annular}, on a sphere or a torus --- {\em spherical}\index{Map!spherical} or {\em toric}\index{Map!toric}, respectively. Oriented sides of the partitioning are called {\em edges} of the map\index{Edge of the map}. Note, that if $e$ is an edge of a map $\Delta$, then $e^{-1}$ is also its edge with the opposite orientation (consisting of the same points of the surface $S$ as a side of the partitioning $\Delta$), see Definition~\ref{def:diagram} below for details.

Now consider an oriented surface $S$ and a given map $\Delta$ and let us fix an orientation on its cells --- e.g. let us walk around the boundary of each cell counterclockwise. In particular, the boundary of a disc map will be read clockwise and for an annular map, one boundary component (``exterior'') will be read clockwise, and the other (``interior'') --- counterclockwise.

Let a boundary component $Y$ of a map or a cell consist of $n$ sides. Walking around this component in accordance with the chosen orientation, we obtain a sequence of edges $e_1\dots e_n$ forming a loop. This loop will be called a {\em contour} of the map\index{Contour of the map} or the cell. In particular, a disc map has one contour, and an annular map has two contours (exterior and interior). Contours are considered up to a cyclic permutation, that is every loop $e_i\dots e_n e_1 \dots e_{i-1}$ gives the same contour. A contour of a cell $\Pi$ will be denoted by $\partial \Pi$ and we will write $e\in\partial\Pi$ if an edge $e$ is a part of the contour $\partial\Pi$ and we will call this situation ``the edge $e$ lies in the contour $\partial\Pi$''. Note, that if an edge $e$ lies in a contour $\partial\Pi$ its inverse $e^{-1}$ does not necessarily lie in that contour. For example, in the situation depicted in Fig. \ref{fig:triangle_group_diagram} an edge $a$ lies in the innermost triangular contour, but $a^{-1}$ does not lie there.

A {\em path} is a sequence of contiguous edges of the cells in the map $\Delta$ on the surface $S$. Given a path $p$ we may define a {\em subpath} in a natural way: a path $q$ is a subpath of the path $p$ if there exist two paths $p_1,p_2$ such that $p=p_1 q p_2$.\\

Given an alphabet $\mathcal{A}$ we denote by $\mathcal{A}^1=\mathcal{A}\cup\mathcal{A}^{-1}\cup\{1\}$, that is, the alphabet $\mathcal{A}^1$ consists of the letters from the alphabet $\mathcal{A}$, their inverses and the symbol ``$1$''. Let $\Delta$ be a map and for each edge $e$ of the map $\Delta$ a letter $\varphi(e)\in\mathcal{A}^1$ is chosen (edges with $\varphi(e)\equiv 1$ are called {\em 0-edges} of the map; other edges are called {\em $\mathcal{A}-$edges}). Note that ${\mathcal A}$-edges that are not $0$-edges may receive labels from ${\mathcal A}^1$. The mapping $\varphi\colon\{\text{edges of a map}\}\to \mathcal{A}^1$ is called {\em labelling}; for each edge $e$ we call the letter $\varphi(e)$ the {\em label} of the edge $e$.

\begin{definition}
	\label{def:diagram}
	Let $\Delta$ be a map over ${\mathcal A}$. If for each edge $e$ of a map $\Delta$ the following relation holds: $$\varphi(e^{-1})\equiv(\varphi(e))^{-1},$$ the map pair $(\Delta,\varphi)$ is called a {\em diagram over $\mathcal{A}$}\index{Diagram!over $\mathcal{A}$}. For short we will denote diagrams with the same letter as the underlying map.
\end{definition}

Here the symbol ``$\equiv$'' denotes the {\em graphical equality} of the words in the alphabet $\mathcal{A}$. In other words the notation $V\equiv W$ means that the words $V$ and $W$ are the same as a sequences of letters of the alphabet. By definition we set $1^{-1}\equiv 1$. Just like in case of paths and subpaths, a {\em subword} $w'$ of a word $w$ is such word that there exist words $v_1, v_2$ such that $w=v_1w'v_2$.

If $p=e_1\dots e_n$ is a path in a diagram $\Delta$ over $\mathcal{A}$ let us define its {\em label} the word $\varphi(p)=\varphi(e_1)\dots\varphi(e_n)$. If the path is empty, that is $|p|=0$, we set $\varphi(p)\equiv 1$ by definition. As before, the label of a contour is defined up to a cyclic permutation (and thus forms a cyclic word). \\

Consider a group $G$ with presentation

\begin{equation}
	\label{eq:presentation}
G=\langle\mathcal{A}|R=1, R\in\mathcal{R}\rangle.
\end{equation}

That means that $\mathcal{A}$ is a basis of a free group $F=F(\mathcal{A})$, $\mathcal{R}$ is a set of words in the alphabet $\mathcal{A} \cup \mathcal{A}^{-1}$ and there exists an epimorphism $\pi: F(\mathcal{A})\to G$ such that its kernel is the normal closure of the subset $\{[r]|r\in \mathcal{R}\}$ of the set of words $F(\mathcal{A})$. Elements of $\mathcal{R}$ are called the {\em relations}\index{Relation} of the presentation $\langle \mathcal{A}|\mathcal{R}\rangle$. We will always suppose that every element $r\in \mathcal{R}$ is a non-empty {\em cyclically-irreducible} word, that is $r$ or any of its cyclic permutations cannot include subwords of the form $ss^{-1}$ or $s^{-1}s$ for a $s\in \mathcal{A}$.

Note that here we presume that if a presentation of a group has a relation $R$, it has all its cyclic permutations as relations as well. The normal closure of the subset $\{[r]\colon r \in {\mathcal R}\}$ of the set of words or relators in $F({\mathcal A})$. Since the kernel is the normal closure, it contains all conjugates of elements in ${\mathcal R}$, and in particular, cyclic permutations of the elements $R$.


\begin{definition}
	\label{def:r-cell}
	A cell of the diagram $\Delta$ over $\mathcal{A}$ is called an {\em $\mathcal{R}$-cell}\index{Cell!$\mathcal{R}-$} if the label of its contour is graphically equal (up to a cyclic permutation) either to a word $R\in \mathcal{R}$, or its inverse $R^{-1}$, or to a word, obtained from $R$ or from $R^{-1}$ by inserting several symbols ``$1$'' between its letters.
\end{definition}

This definition effectively means that choosing direction and the starting point of reading the label of the boundary of any cell of the map and ignoring all trivial labels (the ones with $\varphi(e)\equiv 1$) we can read exactly the words from the set of relations of the group $G$ and nothing else.

Sometimes it proves useful to consider a cell with effectively trivial labels. To be precise, we give the following definition.

\begin{definition}
	\label{def:0-cell}
	A cell $\Pi$ of a map $\Delta$ is called a {\em 0-cell}\index{Cell!$0-$} if the label $W$ of its contour $e_1\dots e_n$ graphically equals $\varphi(e_1)\dots\varphi(e_n)$, where either $\varphi(e_i)\equiv 1$ for each $i=1,\dots ,n$, or for some two indices $i\neq j$ the following holds: $$\varphi(e_i)\equiv a\in\mathcal{A},$$ $$\varphi(e_j)\equiv a^{-1},$$ and $$\varphi(e_k)\equiv 1 \;\; \forall k\neq i,j.$$
\end{definition}

Finally, we can define a diagram of a group.

\begin{definition}
	\label{def:group_diagram}
	Let $G$ be the group given by a presentation (\ref{eq:presentation}). A diagram $\Delta$ on a surface $S$ over the alphabet $\mathcal{A}$ is called a {\em diagram on a surface $S$ over the presentation (\ref{eq:presentation})} (or a {\em diagram over the group $G$} for short)\index{Diagram!van Kampen} if every cell of this map is either an $\mathcal{R}$-cell or a $0$-cell.
\end{definition}

Note that we say ``a diagram over alphabet $\mathcal{A}$'' but the labelling is valued in the set $\mathcal{A}^1$. Note further, that while {\em drawing the diagram} usually the following convention is used. Each edge is oriented with an arrow and endowed with an element of the set $\mathcal{A}$. That means that the edge oriented in accordance with the arrow is labelled with this element, and the inverse edge (with opposite orientation) is labelled with the inverse elements. Hence, the labellings on the diagram are taken from the alphabet $\mathcal{A}$ only.

\section{The van Kampen lemma}
\label{sec:van_kampen}

Earlier we gave two examples of diagrams used to show that a certain equality of the type $W=1$ holds in a group given by its presentation. In fact this process is made possible by the following lemma due to van Kampen:

\begin{lemma}[van Kampen \cite{vK}]
\label{lem:van_kampen}	
Let $W$ be an arbitrary non-empty word in the alphabet $\mathcal{A}^1$. Then $W=1$ in a group $G$ given by its presentation (\ref{eq:presentation}) if and only if there exists a disc diagram over the presentation (\ref{eq:presentation}) such that the label of its contour graphically equals $W$.
\end{lemma}

\begin{proof}
1) First, let us prove that $\Delta$ is a disc diagram over the presentation (\ref{eq:presentation}) with contour $p$, its label $\varphi(p)=1$ in the group $G$.

If the diagram $\Delta$ contains exactly one cell $\Pi$, then in the free group $F$ we have either $\varphi(p)=1$ (if $\Pi$ is a 0-cell) or $\varphi(p)=R^{\pm 1}$ for some $R\in\mathcal{R}$ (if $\Pi$ is an $\mathcal{R}-$cell). In any case, $\varphi(p)=1$ in the group $G$.

If $\Delta$ has more than one cell, the diagram can be cut by a path $t$ into two disc diagram $\Delta_1,\Delta_2$ with fewer cells. We can assume that their contours are $p_1 t$ and $p_2 t^{-1}$ where $p_1 p_2=p$. By induction $\varphi(p_1 t)=1$ and $\varphi(p_2 t^{-1})=1$ in the group $G$. Therefore $$\varphi(p)=\varphi(p_1 p_2)=\varphi(p_1 tt^{-1} p_2)=\varphi(p_1 t)\varphi(t^{-1} p_2)=1$$ in the group $G$. \\

2) Now let us prove the inverse implication. To achieve that we need for a given word $W$ such that $W=1$ in the group $G$ construct a diagram $\Delta$ with contour $p$ such that $\varphi(p)=W$.

Because $W$ is in the normal closure of the set of relators in the free group $F= F({\mathcal A})$, the word $W$ equals a word $V=\prod^{n}_{i=1} X_i R^{\pm 1}_i X_i^{-1}$ for some $R_i\in \mathcal{R}$ and for some words $X_i$.

Construct a polygonal line $t_1$ on the plane and mark its segments with letters so that the line reads the word $X_1$. Connect a circle $s_1$ to the end of this line and mark it so it reads $R_1^{\pm 1}$ if we walk around it clockwise. Now we glue 0-cells to $t_1, s_1$ and $t_1^{-1}$ to obtain a set homeorphic to a disc. We obtain a diagram with the contour of the form $e_1\dots e_k$ with $\varphi(e_1)\equiv 1 \equiv \varphi(e_k)$ and $\varphi(e_2\dots e_{k-1})\equiv X_1 R_1^{\pm 1} X_1^{-1}$.

Construct the second diagram analogously for the word $X_2 R_2^{\pm 1} X_2^{-1}$ and glue it to the first diagram by the edge $e_k$.

Continue the process until we obtain a diagram $\Delta'$ such that $\varphi(\partial\Delta')\equiv V$, see Fig. \ref{fig:van_kampen}.

\begin{figure}
\centering\includegraphics[width=250pt]{ 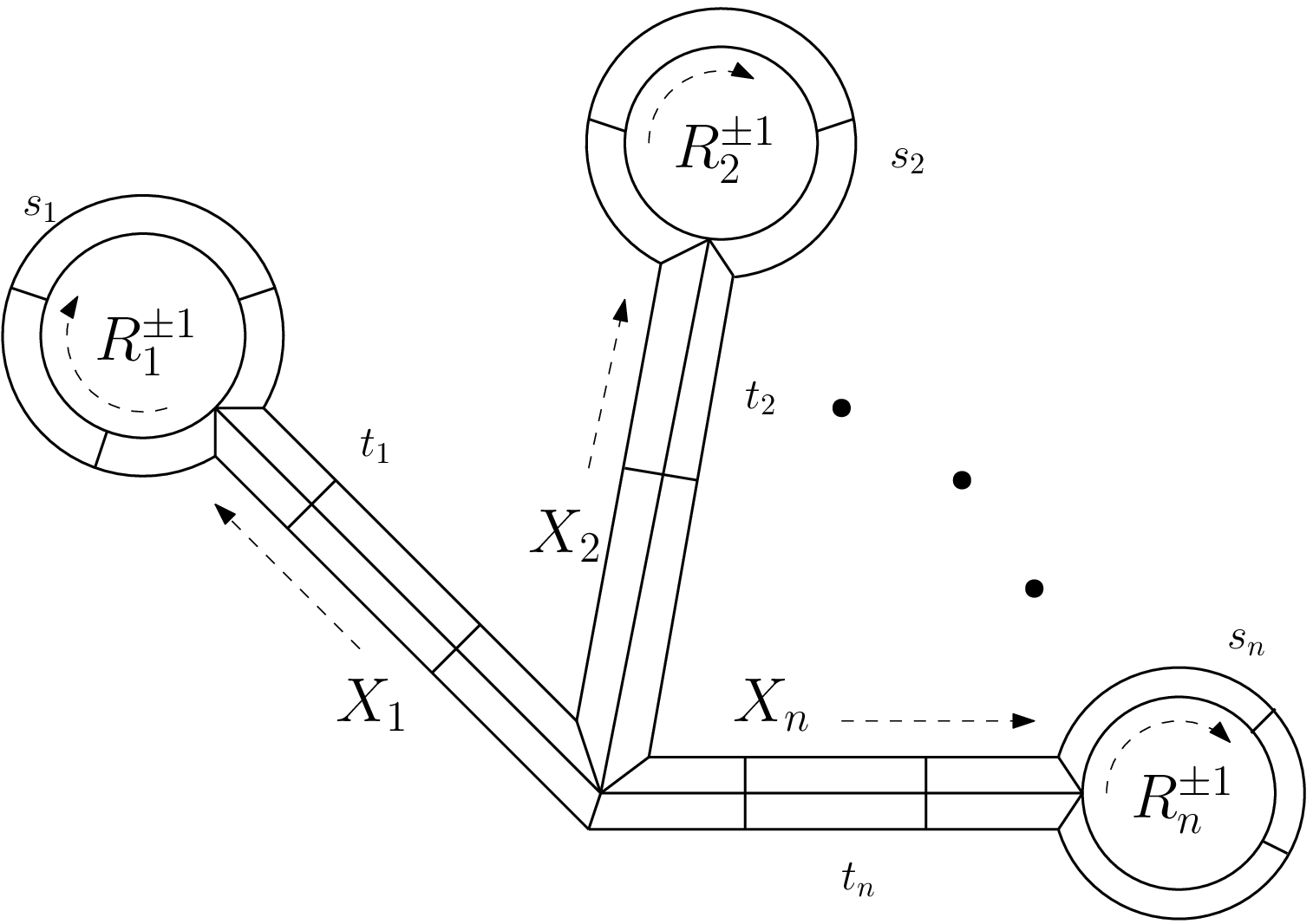}
\caption{The diagram $\Delta'$ with boundary label $\varphi(\partial\Delta')\equiv V$}
\label{fig:van_kampen}
\end{figure}

Finally, glueing some 0-cells to the diagram $\Delta'$ we can transform the word $V$ into the word $W$ getting a diagram $\Delta$ such that $\varphi(\partial\Delta)\equiv W$. That completes the proof.
\end{proof}

This lemma means that disc diagrams can be used to describe the words in a group which are equal to the neutral element of the group. It turns out that annular diagrams can be used to relate conjugate elements.

\begin{lemma}[Schupp \cite{Schupp}]
	\label{lem:schupp}
	Let $V, W$ be two arbitrary non-empty words in the alphabet $\mathcal{A}^1$. Then they are conjugate in a group $G$ given by its presentation (\ref{eq:presentation}) if and only if there exists an annular diagram over the presentation (\ref{eq:presentation}) such that it has two contours $p$ and $q$ with the labels $\varphi(p)\equiv V$ and $\varphi(q)\equiv W^{-1}$.
\end{lemma}


Let $p$ be a loop on a surface $S$ such that its edges form a boundary of some subspace $\tilde{S}\subset S$ homeomorphic to a disc. Then the restriction of a cell partitioning $\Delta$ to the subspace $\tilde{S}$ is a cell partitioning on the space $\tilde{S}$ which is called a {\em submap} $\Gamma$ of the map $\Delta$. Note that by definition a submap is always a {\em disc} submap.

A {\em subdiagram} of a given diagram $\Delta$ is a submap $\Gamma$ of the map $\Delta$ with edges endowed with the same labels as in the map $\Delta$. Informally speaking, a subdiagram is a disc diagram cut out from a diagram $\Delta$.

Let us state an additional important result about the group diagrams.

\begin{lemma}
\label{lem:homotopic_paths_diagram}
Let $p$ and $q$ be two (combinatorially) homotopic paths in a given diagram $\Delta$ over a presentation (\ref{eq:presentation}) of a group $G$. Then $\varphi(p)=\varphi(q)$ in the group $G$.	
\end{lemma}

In particular, the diagrammatic approach is used to deal with {\em groups satisfying the small cancellation conditions}. In that theory a process of cancelling out pairs of cells of a diagram is useful (in addition to the usual process of cancelling out pairs of letters $a$ and $a^{-1}$ in a word). The problem is that two cells which are subject to cancellation do not always form a disc submap, so to define the cancellation process correctly we need to prepare the map prior to cancelling a suitable pair of cells. Let us define those notions in detail.

First, for a given cell partitioning $\Delta$ we define its {\em elementary transformations} (note that elementary transformations are defined for any cell partitionings, not necessarily diagrams, which are defined as cell partitionings satisfying specific conditions, see Definition~\ref{def:diagram}).

\begin{figure}
\centering\includegraphics[width=300pt]{ 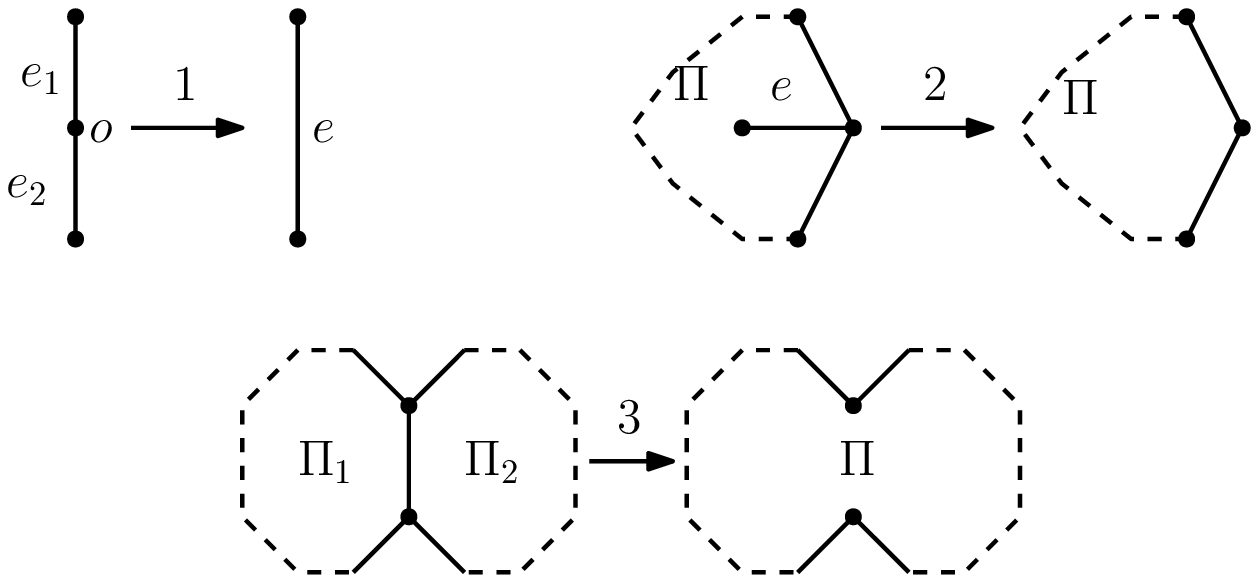}
\caption{Elementary transformations of a cell partitioning}
\label{fig:elementary_tranformations}
\end{figure}

\begin{definition}
	\label{def:elem_transform}
	The following three procedures are called the {\em elementary transformations} of a cell partitioning $\Delta$ (see Fig.~\ref{fig:elementary_tranformations}):
	\begin{enumerate}
		\item If the degree of a vertex $o$ of $\Delta$ equals 2 and this vertex is boundary for two different edges $e_1, e_2$, delete the vertex $o$ and replace the edges $e_1, e_2$ by a single side $e=e_1\cup e_2$;
		\item If the degree of a vertex $o$ of a cell $\Pi$ with $n$ edges ($n\ge 3$) equals 1 and this vertex is boundary for a side $e$, delete the side $e$ and the vertex $o$ (the second boundary vertex of the side $e$ persists);
		\item If two different cells $\Pi_1$ and $\Pi_2$ have a common edge $e$, delete the edge $e$ (leaving its boundary vertices), naturally replacing the cells $\Pi_1$ and $\Pi_2$ by a new cell $\Pi=\Pi_1\cup\Pi_2$.
	\end{enumerate}
\end{definition}

Now we can define a {\em 0-fragmentation} of a diagram $\Delta$. First, consider a diagram $\Delta'$ obtained form the diagram $\Delta$ via a single elementary transformation. This transformation is called an {\em elementary 0-fragmentation} if one of the following holds:

\begin{enumerate}
	\item The elementary transformation is of type 1 and either $\varphi(e_1)\equiv\varphi(e), \, \varphi(e_2)\equiv 1$ or $\varphi(e_2)\equiv\varphi(e), \, \varphi(e_1)\equiv 1$ and all other labels are left unchanged;
	\item The elementary transformation is of type 2 and $\varphi(e)\equiv 1$;
	\item The elementary transformation is of type 3 and the cell $\Pi$ became a 0-cell.
\end{enumerate}

\begin{definition}
	\label{def:0-fragmentation}
	A diagram $\Delta'$ is a {\em 0-fragmentation} of a diagram $\Delta$ if it is obtained from the diagram $\Delta$ by a sequence of elementary 0-fragmentations.
\end{definition}

Note that 0-fragmentation does not change the number of $\mathcal{R}-$cells of a diagram (see Definition~\ref{def:r-cell}).

Now consider an oriented diagram over a presentation (\ref{eq:presentation}). Let there be two $\mathcal{R}-$cells $\Pi_1, \Pi_2$ such that for some 0-fragmentation $\Delta'$ of the diagram $\Delta$ the copies $\Pi_1', \Pi_2'$ of the cells $\Pi_1, \Pi_2$ have vertices $O_1,O_2$ with the following property: those vertices can be connected by a path $\xi$ without selfcrossings such that $\varphi(\xi)=1$ in the free group $F$ and the labels of the contours of the cells $\Pi_1', \Pi_2'$ beginning in $O_1$ and $O_2$ respectively are mutually inverse in the group $F$. In that case the pair $\{\Pi_1, \Pi_2\}$ is called {\em cancelable} in the diagram $\Delta$.

Such pairs of cells are called cancelable because if a diagram $\Delta$ over a group $G$ on a surface $S$ has a pair of cancelable $\mathcal{R}-$cells, there exists a diagram $\Delta'$ over the group $G$ with two fewer $\mathcal{R}-$cells on the same surface $S$. Moreover, if the surface $S$ has boundary, then the cancellation of cells of the diagram $\Delta$ leaves the labels of its contours unchanged.

Given a diagram $\Delta$ and performing cell cancellation we get a diagram $\Delta'$ with no cancelable pairs of cells. Such diagrams are called {\em reduced}. Since this process of reduction does not change the boundary label of a diagram, we obtain the following enhancements of Lemma \ref{lem:van_kampen} and Lemma \ref{lem:schupp}:

\begin{theorem}
	\label{th:van_kampen_strong}	
Let $W$ be an arbitrary non-empty word in the alphabet $\mathcal{A}^1$. Then $W=1$ in a group $G$ given by its presentation (\ref{eq:presentation}) if and only if there exists a reduced disc diagram over the presentation (\ref{eq:presentation}) such that the label of its contour graphically equals $W$.
\end{theorem}

\begin{theorem}
	\label{th:schupp_strong}
	Let $V, W$ be two arbitrary non-empty words in the alphabet $\mathcal{A}^1$. Then they are conjugate in a group $G$ given by its presentation (\ref{eq:presentation}) if and only if there exists a reduced annular diagram over the presentation (\ref{eq:presentation}) such that it has two contours $p$ and $q$ with the labels $\varphi(p)\equiv V$ and $\varphi(q)\equiv W^{-1}$.
\end{theorem}

\section{Unoriented diagrams}
\label{sec:unoriented}

In the present section we introduce a notion of {\em unoriented diagrams}\index{Diagram!van Kampen!unoriented} --- a slight modification of van Kampen diagrams which is useful in the study of a certain class of groups.

Consider a diagram $\Delta$ over a group $G$ with a presentation (\ref{eq:presentation}). With the alphabet $\mathcal{A}$ we associate an alphabet $\bar{\mathcal{A}}$ which is in a bijection with the alphabet $\mathcal{A}$ and is an image of the natural projection $\pi: \mathcal{A}^{1} \to \bar{\mathcal{A}}$ defined as $$\pi(a)=\pi(a^{-1})=\bar{a}$$
for all $a\in \mathcal{A}$ with $\bar{a}$ being the corresponding element of $\bar{\mathcal{A}}$.

Now we take the diagram $\Delta$ and ``forget'' the orientation of its edges. The resulting 1-complex will be called an {\em unoriented diagram over the group $G$} and denoted by $\bar{\Delta}$. All definitions for diagrams (such as disc and annular diagrams, cells, contours, labels, etc.) are repeated verbatim for unoriented diagrams.

Each edge of the diagram $\bar{\Delta}$ is decorated with an element of the alphabet $\bar{\mathcal{A}}$. Therefore, walking around the boundary of a cell in a chosen direction we obtain a sequence of letters but, unlike the oriented case, we do not have an orientation of the edges to determine the sign of each appearing letter. Therefore we shall say that the label $\varphi(\partial \Pi)$ of the contour of a cell $\Pi$ of the diagram $\bar{\Delta}$ is a cyclic word in the alphabet $\bar{\mathcal{A}}$.

Given a word $\bar{w}=\bar{a_1}\dots\bar{a_n}$ in the alphabet $\bar{\mathcal{A}}$ we can produce $2^n$ words in the alphabet $\mathcal{A}^1$ of the form $a_1^{\varepsilon_1}\dots a_n^{\varepsilon_n}$ with $\varepsilon_i\in\{+1,-1\}$. We shall call each of those words a {\em resolution} of the word $\bar{w}$. \\

Unoriented diagrams are very useful when describing such group presentations that the relation $a^2=1$ holds for all generators of the group, in other words groups with a presentation

\begin{equation}
\label{eq:presentation_unoriented}
	G=\langle\mathcal{A}|R=1, R\in\mathcal{R}; a^2=1, a\in\mathcal{A}\rangle.
\end{equation}

In fact, the following analog of the van Kampen lemma holds:

\begin{lemma}
	\label{lem:van_kampen_unoriented}
	Let $W$ be an arbitrary non-empty word in the alphabet $\mathcal{A}^1$. Then $W=1$ in a group $G$ given by its presentation (\ref{eq:presentation_unoriented}) if and only if there exists a reduced unoriented disc diagram over the presentation (\ref{eq:presentation_unoriented}) such that there is a resolution of the label of its contour which graphically equals~$W$.\index{van Kampen lemma!unoriented}
\end{lemma}

\begin{proof}
	First let $W$ be a non-empty word in the alphabet $\mathcal{A}^1$ such that $W=1$ in the group $G$. Let us show that there exists an unoriented diagram with the corresponding label of its contour.
	
	Since $W=1$, due to the strong van Kampen lemma (Theorem \ref{th:van_kampen_strong}) there exists a reduced disc diagram over the presentation (\ref{eq:presentation_unoriented}) such that the label of its contour graphically equals $W$. Denote this diagram by $\Delta$. Now transform every edge $e_i$ of this diagram into a bigon with the label $\varphi(e_i)^2$. Note that the result of this transformation is still a disc diagram. Indeed, we replaced every edge with an $\mathcal{R}$-cell (since $\varphi(e_i)^2=1$ in the group $G$) and reversed the orientation of some of the edges in the boundary contours of the cells of the diagram $\Delta$ but they remain $\mathcal{R}-$cells due to the relations $a^2=1$.
	
	Now we ``collapse'' those 0-cells: replace every bigonal cell with label of the form $\varphi(e)^2$ with an unoriented edge. Thus we obtain an unoriented diagram and the label of its contour by construction has a resolution graphically equal to $W$. \\
	
	To prove the inverse implication, consider an unoriented diagram $\bar{\Delta}$ with the label of its contour $\bar{W}\equiv\varphi(\partial\bar{\Delta})$ in the group $G$ for which $a^2=1$ for $a \in {\mathcal A}$.. We need to show that for each resolution $W$ of the word $\bar{W}$ the relation $W=1$ holds in the group $G$.
	
	First, note that if this relation holds for one resolution of the word $\bar{W}$, it holds for every other resolution of this word. Indeed, due to the relations $a^2=1$ we may freely replace letters with their inverses: the relation $uav=ua^{-1}v$ holds in the group $G$ for any subwords $u, v$ and any letter $a$.
	
	But by definition the diagram $\bar{\Delta}$ is obtained from some diagram $\Delta$ over the group $G$ by forgetting the orientation of its edges. Therefore there exists a van Kampen diagram over the group $G$ with the label of its contour graphically equal to some resolution $R$ of the word $\bar{W}$. Therefore due to the van Kampen lemma $R=1$ in the group $G$. And thus for every other resolution $W$ of the word $\bar{W}$ the relation $W=1$ holds in the group $G$.
\end{proof}

It is worth noting that the approach described in this section may be applied to the $G_n^k$ groups.

\section{Small cancellation theory}
\label{sec:small_cancellations}

\subsection{Small cancellation conditions}

We will introduce the notion of {\em small cancellation conditions $C'(\lambda)$, $C(p)$ and $T(q)$}.\index{Small cancellation conditions} Roughly speaking, the conditions $C'(\lambda)$ and $C(p)$ mean that if one takes a free product of two relations, one gets ``not too many'' cancellations. To give exact definition of those objects, we need to define {\em symmetrisation} and {\em a piece}.

As before, let $G$ be a group with a presentation (\ref{eq:presentation}). Given a set of relations $\mathcal{R}$ its {\em symmetrisation} $\mathcal{R}_{*}$ is a set of all cyclic permutations of the relations $r\in \mathcal{R}$ and their inverses. 

\begin{definition}
\label{def:piece}
	A word $u$ is called {\em a piece with respect to $\mathcal{R}$} if there are two distinct elements $w_1, w_2\in \mathcal{R}_{*}$ with the common beginning $u$, that is $w_1=uv', w_2=uv''$. 
\end{definition}
The length of a word $w$ (the number of letters in it) will be denoted by $|w|$.

\begin{definition} \label{C'(lambda)}
Let $\lambda$ be a positive real number. A set of relations $\mathcal{R}$ is said to {\em satisfy a small cancellation condition $C'(\lambda)$} if
$$|u|<\lambda |r|$$
for every $r\in \mathcal{R}_{*}$ and for all its beginnings $u$ such that $u$ is a piece with respect to $\mathcal{R}$.
\end{definition}

\begin{definition} \label{C(p)}
Let $p$ be a natural number. A set of relations $\mathcal{R}$ is said to {\em satisfy a small cancellation condition $C(p)$} if every element of $\mathcal{R}_{*}$ is a product of at least $p$ pieces.
\end{definition}

The small cancellation conditions are given as conditions on the set of relations $\mathcal{R}$. If a group $\Gamma$ admits a presentation $\langle S|\mathcal{R}\rangle$ with the set of relations satisfying a small cancellation condition, then the group $\Gamma$ is said to satisfy this condition as well.

The condition $C'(\lambda)$ is sometimes called {\em metric} and the condition $C(p)$ --- non-metric. Note, that $C'(\frac{1}{n})$ always yields $C(n+1)$.

There exists one more small cancellation condition. Usually it is used together with either of the conditions $C'(\lambda)$ or $C(p)$.

\begin{definition} \label{T(q)}
Let $q$ be a natural number, $q>2$. A set of relations $\mathcal{R}$ is said to {\em satisfy a small cancellation condition $T(q)$} if for every $l\in \{3, 4, \dots, q-1\}$ and every sequence $\{r_1, r_2, \dots, r_l\}$ of the elements of $\mathcal{R}_{*}$ the following holds: {\em if $r_1 \neq r_2^{-1}, \dots, r_{l-1} \neq r_l^{-1}, r_l \neq r_1^{-1},$ then at least one of the products $r_1 r_2, \dots, r_{l-1} r_l, r_l r_1$ is freely reduced.}
\end{definition}

\begin{remark}
Those conditions have a very natural geometric interpretation in terms of van Kampen diagrams. Namely, the $C(p)$ condition means that every interior cell of the corresponding disc partitioning has at least $p$ sides; the $T(q)$ condition means that every interior vertex of the partitioning has the degree of at least $q$.
\end{remark}

Note that every set $\mathcal{R}$ satisfies the condition $T(3)$. Indeed, no interior vertex of a van Kampen diagram has degree 1 or 2.

\subsection{The Greendlinger theorem}

An important problem of combinatorial group theory is the {\em word problem}: the question whether for a given word $W$ in a group $G$ holds the equality $W=1$ (or, more generally, whether two given words are equal in a given group). Usually the difficult question is to construct a group where a certain set of relations holds but a given word is nontrivial. In other words, to prove that a set of relations $\mathcal{R}$ does not yield $W=1$ (note that there may exist groups where both the relations $\mathcal{R}$ and $W=1$ hold due to the presence of {\em additional} relations). Small cancellation theory proves to be a very powerful and useful instrument in that situation. In particular, an important place in solution of that kind of problems plays the Greendlinger theorem which we will formulate in this section.
We will always assume that the presentation (\ref{eq:presentation}) is symmetrised. \\

The Greendlinger theorem deals with the length of the common part of cells' boundary. Sometimes two cells are separated by 0-cells. Naturally, we should ignore those 0-cells. To formulate that accurately we need some preliminary definitions.

Let $\Delta$ be a reduced diagram over a symmetrised presentation (\ref{eq:presentation}). Two $\mathcal{A}-$edges $e_1,e_2$ are called {\em immediately close} in $\Delta$ if either $e_1=e_2$ or $e_1$ and $e_2^{-1}$ (or $e_1^{-1}$ and $e_2$) belong to the contour of some 0-cell of the diagram $\Delta$. Furthermore, two edges $e$ and $f$ are called {\em close} if there exists a sequence $e=e_1, e_2, \dots , e_l=f$ such that for every $i=1, \dots , l-1$ the edges $e_i, e_{i+1}$ are immediately close.

Now for two cells $\Pi_1, \Pi_2$ a subpath $p_1$ of the contour of the cell $\Pi_1$ is a {\em boundary arc between $\Pi_1$ and $\Pi_2$} if there exists a subpath $p_2$ of the contour of the cell $\Pi_2$ such that $p_1=e_1u_1e_2\dots u_{n-1}e_n, \, p_2^{-1}=f_1v_1f_2\dots v_{n-1}f_n$, the paths $u_i, v_i$ consist of 0-edges, $e_i, f_i$ are $\mathcal{A}-$edges such that for each $i=1, \dots , n$ the edge $f_i$ is close to the edge $e_i$. In the same way a boundary arc between a cell and the contour of the diagram $\Delta$ is defined.

Informally we can explain this notion in the following way. Intuitively, a boundary arc between two cells is the common part of the boundaries of those cells. The boundary arc defined above becomes exactly that if we collapse all 0-cells between the cells $\Pi_1$ and $\Pi_2$.

Now consider a {\em maximal} boundary arc, that is a boundary arc which does not lie in a longer boundary arc. It is called {\em interior} if it is a boundary arc between two cells, and {\em exterior} if it is a boundary arc between a cell and the contour $\partial\Delta$.

\begin{remark}
	\label{rem:geometrc_interpretation}
	It is easy to see that the small cancellation conditions $C'(\lambda)$ and $C(p)$ have a natural geometric interpretation in terms of boundary arcs. Namely, an interior arc of a cell $\Pi$ of a reduced diagram of a group satisfying the $C'(\lambda)$ condition has length smaller than $\lambda|\partial\Pi|$. Likewise, the $C(p)$ condition means that the boundary of every cell of the corresponding diagram consists of at least $p$ arcs.
\end{remark}

Now we can formulate the Greendlinger theorem.\index{Greendlinger theorem}

\begin{theorem}
	\label{th:greendlinger}
	Let $\Delta$ be a reduced disc diagram over a presentation of a group $G$ satisfying a small cancellation condition $C'(\lambda)$ for some $\lambda\le\frac{1}{6}$ and let $\Delta$ have at least one $\mathcal{R}-$cell. Further, let the label $\varphi(q)$ of the contour $q=\partial \Delta$ be cyclically irreducible and such that the cyclic word $\varphi(q)$ does not contain any proper subwords equal to 1 in the group $G$.
	
	Then there exists an exterior arc $p$ of some $\mathcal{R}-$cell $\Pi$ such that $$|p|>\frac{1}{2}|\partial\Pi|.$$
\end{theorem}

\begin{remark}
	If we formulate the Greendlinger theorem for unoriented diagrams, the theorem still holds.
\end{remark}

Before proving this theorem, let us interpret it in terms of group presentation and the word problem. Consider a group $G$ with a presentation (\ref{eq:presentation}) satisfying a small cancellation condition $C'(\lambda), \, \lambda\le\frac{1}{6}$. Due to van Kampen lemma \ref{lem:van_kampen} (and its strengthening, Theorem \ref{th:van_kampen_strong}) for a word $W=1$ in the group $G$ there exists a diagram with boundary label $W$. The boundary label of every $\mathcal{R}-$cell of the diagram by definition is some relation from the set $\mathcal{R}$ (or its cyclic permutation). Then, due to Greendlinger theorem \ref{th:greendlinger} there is an $\mathcal{R}-$cell such that at least half of its boundary ``can be found'' in the boundary of the diagram. Thus we obtain the following corollary (sometimes it is also called the Greendlinger theorem):

\begin{theorem}[Greendlinger \cite{Greendlinger}]
	\label{cor:greendlinger_group}
	Let $G$ be a group with a presentation (\ref{eq:presentation}) satisfying a small cancellation condition $C'(\lambda), \, \lambda\le\frac{1}{6}$. Let $W\in F$ be a nontrivial freely reduced word such that $W = 1$ in the group $G$. Then there exists a subword $V$ of $W$ and a relation $R\in\mathcal{R}$ such that $V$ is also a subword of $R$ and such that $$|V|>\frac{1}{2}|R|.$$
\end{theorem}

This theorem is very useful in solving the word problem.

\begin{example}[A.A. Klyachko]
Consider a relation $R=[x,y]^2$ and a word $$W=[x^{1000},y^{1000}]^{1000}.$$ The question is, whether in every group with the relation $R$ the equality $W=1$ holds.

First, it is easy to see that the set of relations $\mathcal{R}_*$ obtained from $R$ by symmetrisation satisfies the condition $C'(\frac{1}{6})$. Therefore, due to the Greendlinger theorem every word $V$ such that $V=1$ in the group $G$ has a cyclic permutation such that both its irreducible form and some relation $\tilde{R}\in\mathcal{R}_{*}$ have a common subword $p$ of length $|p|>\frac{1}{2}|\tilde{R}|=\frac{8}{2}=4$.

On the other hand, the longest common subwords of the word $W=[x^{1000},y^{1000}]^{1000}$ and any of the relations have length 2: those are $$x^{-1}y^{-1}, xy, y^{-1}x, yx^{-1}.$$ Therefore we can state that there exists a group $G$ where for some two elements $a,b\in G$ $[a,b]^2=1$ but $[a^{1000},b^{1000}]^{1000}\neq 1$.
\end{example}

Now let us prove Theorem \ref{th:greendlinger}.

\begin{proof}
	Let $G$ be a group with a presentation (\ref{eq:presentation}) and $\Delta$ be a diagram as in the statement of the theorem. Since $\Delta$ is a disc diagram, we can place it on the sphere. We construct the dual graph $\Phi$ to the 1-skeleton of the diagram $\Delta$ in the following way.
	
	Place a vertex of the new graph $\Phi$ into each $\mathcal{R}-$cell of the diagram $\Delta$ and one vertex $O$ into the outer region on the sphere. To construct the edges of the graph $\Phi$ (note that they are not edges of a diagram) we perform the following procedure. For every interior boundary arc between cells $\Pi_1$ and $\Pi_2$ we chose an arbitrary $\mathcal{R}-$edge $e_1$ and an edge $e_2$ close to it such that $e_2^{-1}$ lies in $\partial\Pi_2$. Now we connect the vertices of the graph $\Phi$, lying inside the cells $\Pi_1,\Pi_2$ by a smooth curve $\gamma$ transversally intersecting the interior of the edges $e_1$ and $e_2$ once and crossing all 0-cells lying between them. This curve $\gamma$ is now considered as an edge of the graph $\Phi$. This procedure is performed for every vertex of the graph $\Phi$ and every boundary arc (if the arc is exterior, we connect the vertex with the ``exterior'' vertex $O$).
	
	The main idea of the proof is to study the dual graph $\Phi$ and to prove the necessary inequality by the reasoning of Euler characteristic of sphere.
	
	First, note that among the faces of the graph $\Phi$ there are no 1-gons (because the boundary labels of the cells of the diagram $\Delta$ are cyclically irreducible) or 2-gons (because in the previous paragraphs we have constructed exactly one edge crossing every maximal boundary arc of $\Delta$). Therefore every face of $\Phi$ is at least a triangle. Thus, denoting the number of vertices, edges and faces of the graph $\Phi$ by $V,E$ and $F$ respectively, and considering the usual Euler formula $V-E+F=2$ we get the following inequality:
	
	\begin{equation}
	\label{eq:euler_phi}	
		V>\frac{1}{3}E+1.
	\end{equation}
	
	Now, suppose the statement of the theorem does not hold. For each edge of the graph $\Phi$ connecting the vertices $o_i, o_j$ lying in the cells $\Pi_i,\Pi_j$ of the diagram $\Delta$ we attribute this edge to each of those vertices with a coefficient $\frac{1}{2}$; if an edge connects a vertex $o_k$ with the vertex $O$, we attribute it to the vertex $o_k$ with a coefficient 1.
	
	Fix an arbitrary vertex $o\neq O$. Denote the cell this vertex lies in by $\Pi$ and consider the following possibilities.
	
	\begin{enumerate}
		\item Let the contour $\partial\Pi$ consist only of interior arcs. Due to the $C'(\frac{1}{6})$ condition the length of each of those arcs is smaller than $\frac{1}{6}|\partial\Pi|$ (see Remark \ref{rem:geometrc_interpretation}). Therefore, their number is not smaller than seven and thus at least $\frac{7}{2}$ edges of the graph $\Phi$ is attributed to the vertex $o$.
		\item Let there be exactly one exterior arc $p$. Since we suppose that $|p|\le\frac{1}{2}|\partial\Pi|$, the number of interior arcs of the contour of the cell $\Pi$ is at least four. Therefore we attribute at least $1+4\cdot\frac{1}{2}=3$ edges to the vertex $o$.
		\item Let the cell $\Pi$ have two exterior edges $p_1,p_2$. Note that the end of the arc $p_1$ can't be the beginning of the arc $p_2$ (and vice versa) and they can not be separated by 0-edges only because otherwise we could cut the diagram $\Delta$ with edges $f_1,\dots,f_l$ such that $\varphi(f_i)=$ for all $i=1,\dots , l$ and due to van Kampen lemma obtain two proper subwords $\varphi(q_1), \varphi(q_2)$ of the word $\varphi(q)$ (where $q$ is the contour of the diagram $\Delta$) equal to 1 in the group $G$. That contradicts the condition of $\varphi(q)$ being cyclically irreducible.
		
			Therefore, the cell $\Pi$ has at least two distinct interior boundary arcs and there are at least $2\cdot 1+2\cdot \frac{1}{2}=3$ edges attributed to the vertex $o$.
		\item If the cell $\Pi$ has at least three exterior arcs, there are at least three edges of the graph $\Phi$ attributed to the vertex $o$.
	\end{enumerate}
	
	Considering those possibilities for all vertices of the graph $\Phi$ except the ``exterior'' vertex $O$ we see that $$V-1\le \frac{1}{3}E$$ and that contradicts the inequality (\ref{eq:euler_phi}). This contradiction completes the proof.
\end{proof}

%% file: section_6.tex
\begin{figure}[h]
\begin{center}
 \includegraphics[width =0.3\textwidth]{ 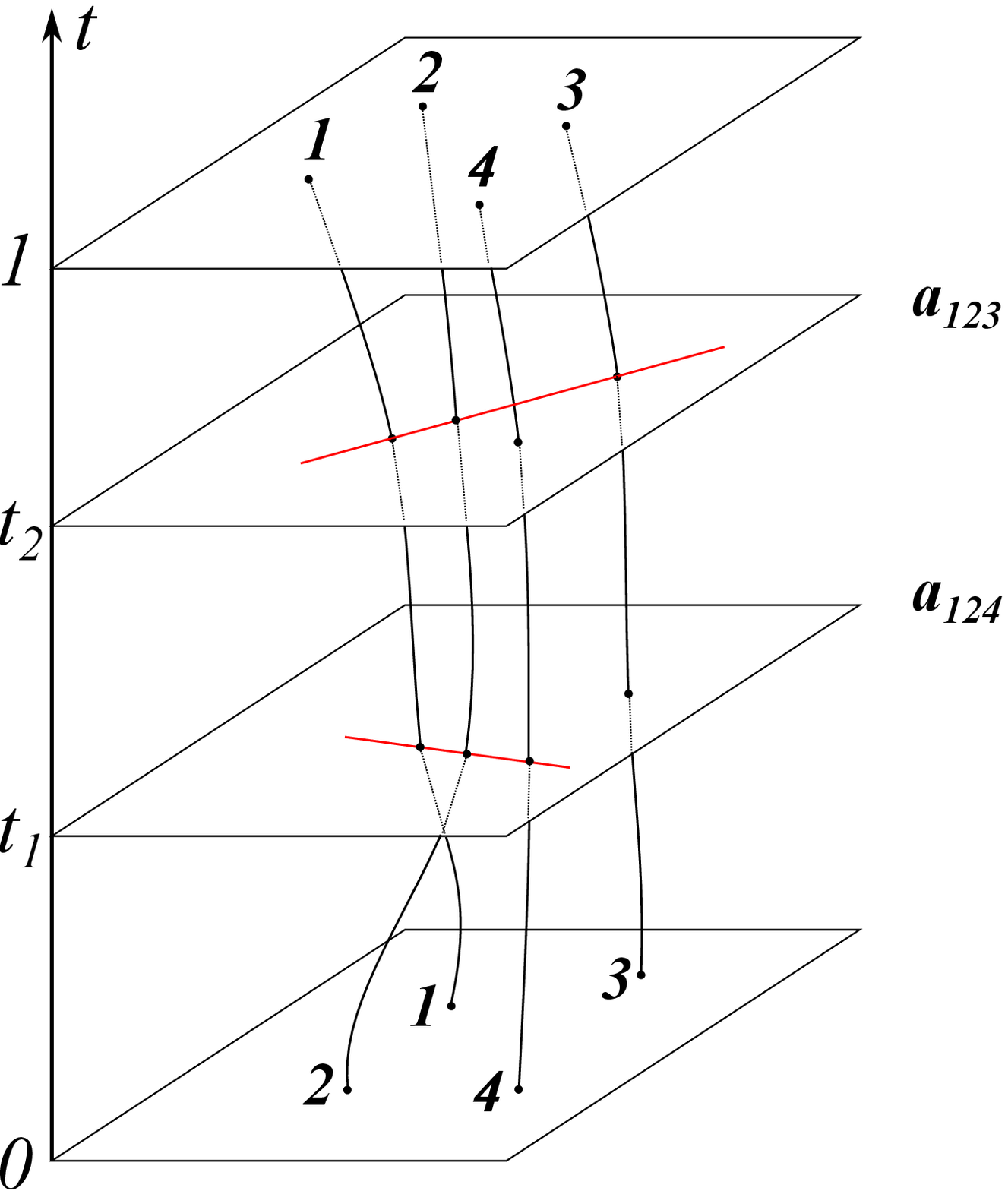}
\end{center}
\caption{Horizontal trisecants of a braid $\beta$}\label{fig:trisecants}
\end{figure}

The $G_n^k$-theory we consider below grows on the following basic example. Consider a geometrical braid $\beta$, i.e. a set of $n$ numbered strands connecting point on two parallel planes in $\R^3$, see~Fig.~\ref{fig:trisecants}. The braid determines a dynamics of $n$ points in the plane: take the bottom plane and move it in parallel to the top. The strands will intersect the plane in $n$ points which will move while the plane go up. We mark the moments when the point configuration becomes critical (``possesses a good property''). Here it means that some three of the intersection points are collinear. The line which contains these point is called a {\em horizontal trisecant}\label{def:horizontal_trisecant}. We assign the letter $a_{ijk}$ to the critical moment, $i,j,k$ are the numbers of strands intersecting the trisecant. The sequence of the letters assigned to the critical moments during ascension of the plane form a word $w(\beta)$. In our example $w(\beta)=a_{124}a_{123}$. Isotopies of the braid lead to transformations of  $w(\beta)$, and these transformations can be incorporated into relations of the group $G_n^3$ we define below. Thus, the word $w(\beta)$ considered as an element of $G_n^3$ becomes an invariant of the braid $\beta$.

\section{Dynamical systems and their properties}\label{sect:dyn_sys}

Given a topological space $\Sigma \simeq \Sigma^{N}$ of some high dimension $N$ (the model example we have in mind is $\Sigma=\R^N$), we shall fix some number $n$ and construct some {\em configuration space} to be a certain subset of the Cartesian power $\Sigma^{\times n}\simeq (\Sigma^{N})^{\times n}$ having dimension $Nn$.
The components $x_{i}$ of $\Sigma \; (i=1,\dots, n)$ of an element ${\mathbf x}=(x_1,\dots,x_n)\in\Sigma^{\times n}$ will be referred to as {\em particles}. Later on $N$ should often be omitted.

The topology on $\Sigma^{\times n}$ defines a natural topology
on the space of all continuous mappings $[0,1]\to \Sigma^{\times n}$; we
shall also study mappings $S^{1}\to \Sigma^{\times n}$, where
$S^{1}=[0,1]/\{0=1\}$ is the circle; these mappings will naturally
lead to {\em closures} of $k$-braids.

Let us fix positive integers $n$ and $k$.


The   {\em space of admissible dynamical systems} ${\cal D}$ is a
closed subset in the space of all maps $[0,1]\to \Sigma^{n}$ (later on we shall justify which subsets will be taken); here $n$ will denote the number of particles in $\Sigma$,
$k$ will be some parameter relating the way how to describe some good subset of $\Sigma^{n}$ to deal with.

A {\em dynamical system}\index{Dynamical system} is an element $D$ of $\mathcal D$. By a {\em state} of $\mathcal D$ we mean the ordered
set of particles $D(t)\in
\Sigma^{n}$. Herewith, $D(0)$ and $D(1)$ are called {\em the initial
state} and {\em the terminal state}.

Later on, we shall justify that we take $D(t)$ not barely in $\Sigma^{\times n}$ but rather in a certain subspace
$\Sigma^{*} = \Sigma^{\times n} \backslash \Delta$; in some cases (when $\Sigma$ is one-dimensional) $\Delta$ may be empty,
for the classical braids in $2$-surfaces, $\Delta$ is the {\em big diagonal}: $x_{i}=x_{j}$ for at least one
pair of indices $i\neq j$.

The parameter $k$ will be used to define the structure of $\Delta$; we shall discuss every specific case
separately.


As usual, we shall fix the initial state and the terminal state and
consider the set of admissible dynamics with such initial and
terminal states.

We shall also deal with {\em cyclic dynamical system}, where $D(0)=D(1)$.

Given $\Sigma^{*}\subset \Sigma^{\times n}$; a point from $\Sigma^{*}$ is
a set of particles from $\Sigma$ satisfying a certain condition.
We say that a property ${\cal P}$ defined for elements $(x_{1},\dots, x_{n})\in \Sigma^{*}$
is {\em $k$-good}\index{$k$-good property}, if the following
conditions hold:

\begin{enumerate}

\item if this property holds for some set of particles then it holds for every subset of this set;

\item this property holds for every set consisting of $k-1$ particles among $n$ ones (hence, for every smaller set);

\item fix $k+1$ pairwise distinct numbers $i_{1},\dots,i_{k+1},$ where each $i_{j}\in\{1,\dots,n\}$;
if the property ${\cal P}$ holds for particles with numbers
 $i_{1},\dots, i_{k}$ and for the set of particles with numbers $i_{2},\dots, i_{k+1}$, then it holds
 for the set of all
$k+1$ particles $i_{1},\dots, i_{k+1}$.
\end{enumerate}

The basic example of a $k$-good property of points in $\Sigma=\R^{k-1}$ is {\em points belong to one $(k-1)$-dimensional (affine) plane}. For instance, collinearity property in $\R^2$ should be $3$-good, coplanarity property in $\R^3$ should be $4$-good etc. But here we have a problem. For example, consider the coplanarity property and take three collinear points $1,2,3$ in $\R^3$. Then adding any point $4$ (or $5$) of $\R^3$ to them, we obtain a coplanar quadruple. On the other hand, we can add two points $4$ and $5$ so that the set $1,2,3,4,5$ is not coplanar. Thus, the third condition of the good property definition is violated.

One solution to keep coplanarity as a good condition is to exclude degenerated configurations as the above. That is why we restrict the consideration to a subset $\Sigma^{*}$ of the configuration space.

\begin{remark}
Besides {\em statically good properties} ${\cal P}$,
which are defined for subsets of the set of particles regardless the
state of the dynamics, one can talk about {\em dynamically good}
properties, which can be defined for states considered in time.
\end{remark}


\begin{remark}
Our main example deals with the case $k=3$, where for particles we
take different points on the plane $\Sigma = \R^2$, and ${\cal P}$ is the property
of points to belong to the same line. In general, particles may be
more complicated objects than just points. More examples are given on page~\pageref{exa:good_properties}.
\end{remark}

Let ${\cal P}$ be a $k$-good property defined on a set of $n$
particles, $n>k$. For each  $t\in [0,1]$ we shall fix the
corresponding state of particles, and pay attention to those $t$ for
which there is a set $k$ particles possessing ${\cal P}$; we shall
refer to these moments as  ${\cal P}$-critical (or just {\em
critical}).

\begin{definition}
We say that a dynamical system $D$ is {\em pleasant}\index{Dynamical system!pleasant}, if the set of its critical moments is finite whereas for each critical moment there exists exactly one $k$-index set for which the condition ${\cal P}$ holds (thus, for larger sets the property ${\cal P}$ does not hold). Such an unordered $k$-tuple of indices will be called a {\em multiindex} of critical moments.
\end{definition}

\begin{figure}[h]
\begin{center}
 \includegraphics[width =0.3\textwidth]{ 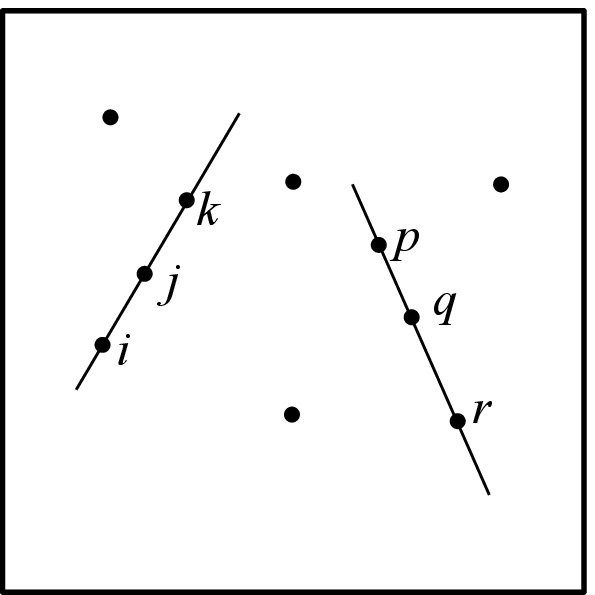}
\end{center}
\caption{An unpleasant critical moment for the collinearity property: there are two triples of collinear points.}\label{fig:unpleasant}
\end{figure}


For each dynamics $D$ and each multiindex $m=(m_{1},\dots, m_{k})$,
let us define the {\em $m$-type} of $D$ as the ordered set $t_{1}<
t_{2}< \dots$ of values of $t$ for which the set of particles
$m_{1},\dots, m_{k}$ possesses the property ${\cal P}$.

Notation: $\tau_{m}(D)$. If the number $l$ of values $t_{1}<\dots<
t_{l}$ is fixed, then the type can be thought of as a point in
$\R^{l}$ with coordinates $t_{1},\dots, t_{l}$.

\begin{definition}
By the {\em type} of a dynamical system, we mean the set of all its
types $\tau_{m}(D)$. Notation: $\tau(D)$.
\end{definition}

For the example in Fig.~\ref{fig:trisecants} we have $\tau_{123}=\{t_2\}$, $\tau_{124}=\{t_1\}$ and $\tau_m=\emptyset$ for the other index triples.

If $D$ is pleasant, then these set are pairwise disjoint for different $m$.

\begin{definition}
Fix a number $k$ and a $k$-good property  ${\cal P}$. We say that  $D$
is  {\em ${\cal P}$-stable}\index{Dynamical system!stable}, if there is a neighbourhood $U(D)$,
where each dynamical system $D'\in U\setminus\{D\}$ is pleasant, whereas for each
multiindex $m$ (consisting of $k$ indices), the number $l$ of
critical values corresponding to this multiindex $U(D)$ is the same, the type
$\tau_{m}$ is a continuous mapping $U(D)\to \R^{l}$.
\end{definition}

\begin{figure}[h]
\begin{center}
 \includegraphics[width =0.3\textwidth]{ 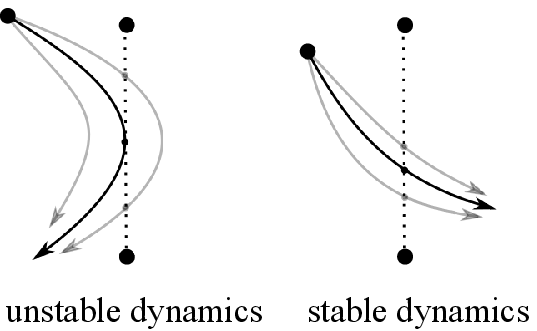}
\end{center}
\caption{Stable and unstable dynamical systems for the collinearity property.}\label{fig:stable_dynamics}
\end{figure}

For a non-stable dynamical system we can consider the following, Let $n-1$ points $z_{k}(t)$ be uniformly distributed along the unit circle, $z_{k}(t) = e^{\frac{2 \pi i k}{n-1}},k=1,\dots n-1$
and let the $n$-th point $z_{n}(t)$ have a trajectory which is tangent to the line passing through $z_{1}$ and $z_{2}$ for $t=t_{0}$. Then this dynamical system is not stable with respect to the property {\em three points are collinear}.

We shall often say {\em pleasant dynamical systems} or {\em stable
dynamical systems} without referring to ${\cal P}$ if it is clear
from the context which ${\cal P}$ we mean.

\begin{definition}
A {\em deformation} is a continuous path $s:[0,1]\to {\cal D}$ in
the space of admissible dynamical systems, from a stable pleasant
dynamics $s(0)$ to another pleasant stable dynamical system $s(1)$.
\end{definition}

\begin{definition}
We say that a deformation $s$ is {\em admissible} if:
\begin{enumerate}

\item the set of values $u$, where $s(u)$ is not pleasant or is not stable, is finite,
and for those $u$ where $s(u)$ is not pleasant, $s(u)$ is stable.

\item inside the stability intervals, the $m$-types are
continuous for each multiindex $m$.

\item for each value $u=u_{0}$, where $s(u_{0})$ is not pleasant, exactly one of the two following
cases occurs:

\begin{enumerate}
\item There exists exactly one $t=t_{0}$ and exactly one $(k+1)$-tuple
$m=(m_{1},\dots, m_{k+1})$ satisfying ${\cal P}$ for this $m$
(hence, ${\cal P}$ does not hold for larger sets).

Let ${\tilde m}_{j}=m\setminus \{m_{j}\},j=1,\dots, k+1$.  For types
$\tau_{{\tilde m}_{j}}$, choose those coordinates $\zeta_{j}$, which
correspond to the value $t=t_{0}$, i.e. $\zeta_j(u_0)=t_0$. It is required that for all these
values  $u_{0}$ all functions $\zeta_{j}(u)$ are smooth, and all
derivatives $\frac{\partial \zeta_{j}}{\partial{u}}$ are pairwise
distinct;

\item there exists exactly one value $t=t_{0}$ and exactly two
multiindices  $m=\{m_{1},\dots, m_{k}\}$ and $m'=\{m'_{1},\dots,
m'_{k}\}$ for which  ${\cal P}$  holds; we require that $Card(m\cap
m')<k-1$.
\end{enumerate}

\item For each value $u$, where the dynamical system $s(u)$ is not
stable, there exists a value $t=t_{0}$, which is not critical for
$D_{u}$, and a multiindex $\mu=(\mu_{1},\dots,\mu_{k})$, for which
the following holds. For some small  $\eps$ all dynamical systems
$D_{u_0}$ for $u_{0}\in (u-\eps,u)$ and $u_{0}\in (u,u+\eps)$ are
stable (for $\delta<\eps$), and the type $\tau_{\mu}(D_{u+\delta})$
differs from the type $\tau_{\mu}(D_{u-\delta})$ by an
addition/removal of two identical multiindices  $\mu$ in positions close to $t_{0}$.
\end{enumerate}
\end{definition}

\begin{figure}[h]
\begin{center}
 \includegraphics[width =0.5\textwidth]{ 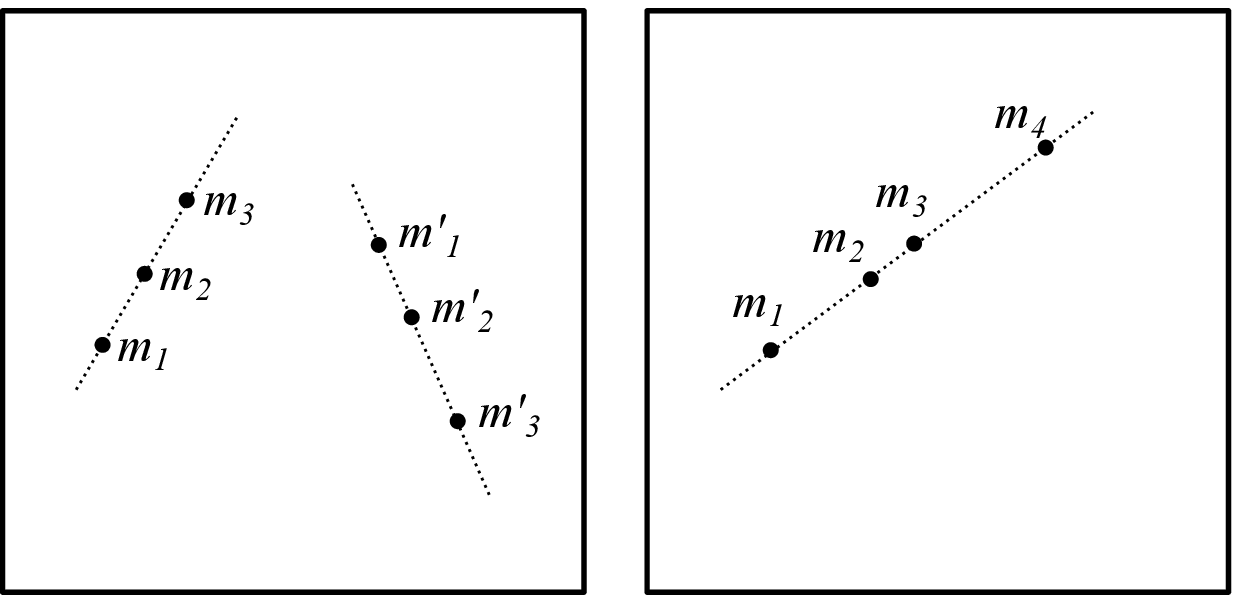}
\end{center}
\caption{Stable unpleasant configurations for the collinearity property.}\label{fig:unpleasant_stable}
\end{figure}

For the collinearity property a stable unpleasant configurations are four collinear points (Fig.~\ref{fig:unpleasant_stable} right) or two collinear triples (Fig.~\ref{fig:unpleasant_stable} left). An unstable configuration appears when the trajectory of one point is tangent to the line containing other two points (Fig.~\ref{fig:stable_dynamics} left).

Admissible deformation can be thought of as deformation in general position. Confining ourselves to admissible deformations, we evade complex critical configurations like five collinear points in the plane etc.

For the space of deformation $\mathfrak{D}$, one defines an induced topology.

\begin{definition} We say that a $k$-good property ${\cal P}$ is $k$-{\em correct} for the space of admissible dynamical
systems, if the following conditions hold:

\begin{enumerate}

\item In each neighbourhood of any dynamical system $D$ there exists a pleasant dynamical system $D'$.

\item For each deformation $s$ there exists an {\em admissible} deformation with the same ends $s'(0)=s(0),s'(1)=s(1)$.

\end{enumerate}

Correctedness of a good property depends essentially on the configuration space $\Sigma^*$ we work with. For example, planarity is a correct $3$-good property for sets of distinct points in $\R^2$, but it is not correct if the points-particles move in the grid
$$
\Gamma =\{(x,y)\in\R^2\,|\, x\in\Z\mbox{ or }y\in\Z\}.
$$
Indeed, in this case the pleasant configurations are not dense, see Fig.~\ref{fig:incorrect_planarity}.

\begin{figure}[h]
\begin{center}
 \includegraphics[width =0.25\textwidth]{ 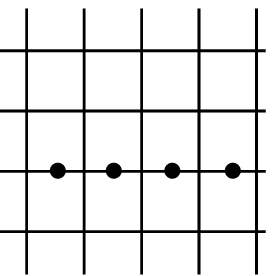}
\end{center}
\caption{Unpleasant configuration whose small neighbourhood contains no pleasant configurations.}\label{fig:incorrect_planarity}
\end{figure}

\end{definition}

\begin{definition}
We say that two dynamical systems $D_{0},D_{1}$ are {\em equivalent}, if there exists a deformation $s, s(0)=D_{0},s(1)=D_{1}$.
\end{definition}

Thus, if we talk about a correct  ${\cal P}$-property, we can talk
about an admissible deformation when defining the equivalence.

\subsection{The group $G_{n}^{k}$}
\label{sec:gnk}

Let us now pass to the definition of the {\em $n$-strand free
$k$-braid group $G_n^k$}\index{$G_n^k$}\index{Braid!free $k-$}.

Consider the following $n \choose k$ generators
$a_{m}$, where $m$ runs the set of all unordered $k$-tuples
$m_{1},\dots, m_{k},$ whereas each $m_{i}$ are pairwise distinct
numbers from $\left\{1,\dots, n\right\}$.

For each unordered $(k+1)$-tuple $U$ of distinct indices $u_{1},\dots, u_{k+1} \in
\{1,\dots, n\}$, consider the $k+1$ sets $m^{j}=U\setminus
\{u_{j}\}, j=1,\dots, k+1$. With $U$, we associate the relation

\begin{equation}\label{eq:gnk_tetrahedron_relation}
a_{m^1}\cdot a_{m^2}\cdots a_{m^{k+1}}= a_{m^{k+1}}\cdots a_{m^2}\cdot a_{m^1}
\end{equation}
for two tuples $U$ and ${\bar U}$, which differ by order reversal,
we get the same relation.

Thus, we totally have
$\frac{(k+1)! {n \choose k+1}}{2}$ relations.

We shall call them the {\em tetrahedron relations}\index{Tetrahedron relation}. Note that those relations appear in physics, see \cite{etin}.

For $k$-tuples  $m,m'$ with $Card(m\cap m')<k-1$, consider the {\em
far commutativity relation:}\index{Far commutativity}

\begin{equation}\label{eq:gnk_far_commutativity_relation}
a_{m}a_{m'}=a_{m'}a_{m}.
\end{equation}

Note that the far commutativity relation can occur only if $n>k+1$.

In addition, for all multiindices $m$, we write down the following
relation:

\begin{equation}\label{eq:gnk_order2_relation}
a_{m}^{2}=1.
\end{equation}

\begin{definition}\label{def:gnk}
The {\em $k$-free braid group} $G_{n}^{k}$ is defined as the quotient group of the free group generated
by all $a_{m}$ for all multiindices $m$  by relations~\eqref{eq:gnk_tetrahedron_relation}, \eqref{eq:gnk_far_commutativity_relation} and
\eqref{eq:gnk_order2_relation}.
\end{definition}

\begin{example} \label{ex:g32}
The group $G_{3}^{2}$ is $\langle
a,b,c\,|\,a^{2}=b^{2}=c^2=(abc)^{2}=1\rangle,$ where
$a=a_{12},b=a_{13},c=a_{23}$, the symmetric group $S_3$ on three letters $a,b,c$.

Indeed, the relation $(abc)^{2}=1$ is equivalent to the relation
$abc=cba$ because of $a^{2}=b^{2}=c^{2}=1$. This obviously yields
all the other tetrahedron relations.
\end{example}

\begin{example}
The group $G_{4}^{3}$ is isomorphic to $\langle
a,b,c,d\,|\,a^{2}=b^{2}=c^{2}=d^2=1,(abcd)^{2}=1,(acdb)^{2}=1,(adbc)^{2}=1\rangle$.
Here $a=a_{123},b=a_{124},c=a_{134},d=a_{234}.$

It is easy to check that instead of $\frac{4!}{2}=12$ relations, it
suffices to take only $\frac{3!}{2}$ relations.
\end{example}

By the {\em length} of a word we mean the number of letters in
this word, by the {\em complexity of a free $k$-braid} we mean the
minimal length of all words representing it. Such words will be
called {\em minimal representatives}. The tetrahedron relations (in
the case of free $2$-braids we call them the {\em triangle
relations}) and the far commutativity relations do not change the
complexity, and the relation $a_{m}^{2}=1$ increases or decreases
the length by $2$.

As usual in the group theory, it is natural to look for minimal
length words representing the given free $k$-braid.

If we deal with conjugacy classes of free $k$-braids, one deals with
the length of cyclic words.

The number of words of fixed length in a finite alphabet is
finite; $k$-braids and their conjugacy classes are the main
objects of the present paper.

Let us define the following two types of homomorphisms for free
braids. For each $l=1,\dots, n$, there is an  {\em index forgetting
homomorphism}\index{Homomorphism!index forgetting} $f_{l}: G_{n}^{k}\to G_{n-1}^{k-1}$; this homomorphism
takes all generators $a_{m}$ with multiindex $m$ not containing $l$
to the unit element of the group, and takes the other generators
$a_{m}$ to $a_{m'}$, where $m'=m\setminus \{l\}$; this operation is
followed by the index renumbering.

The  {\em strand-deletion homomorphism}\index{Homomorphism!strand-deletion} $d_{j}$ is defined as a
homomorphism  $G_{n}^{k}\to G_{n-1}^{k}$; it takes all generators
$a_{m}$ having multiindex containing $j$ to the unit element; after
that we renumber indices.

The free $2$-braids (called also {\em pure free braids}) were studied in~ \cite{knot_parity, MW, M3}. For free $2$-braids, the following theorem holds.

\begin{theorem}\label{theo2}
Let $b'$ be a word representing a free $2$-braid $\beta$. Then every
word $b$ which is a minimal representative of  $\beta$, is
equivalent by the triangle relations and the far commutativity
relation to some subword of the word $b'$.

Every two minimal representatives $b_{1}$ and $b_{2}$ of the same free
$2$-braid  $\beta$ are equivalent by the triangle relations and the
far commutativity relations.
\end{theorem}

Thus, for free $2$-braids, the recognition problem can be solved by
means of considering its minimal representative.

The main idea of the proof of this theorem is similar to the
classification of homotopy classes of curves in $2$-surfaces due to
Hass and Scott~\cite{HS}: in order to find a minimal representative,
one looks for ``bigon reductions''  until no longer possible, and the final
result is unique up to third Reidemeister moves for the exception of
certain special cases (multiple curves etc). For free $2$-braids
``bigon reductions'' refer to some cancellations of generators
similar generators $a_{m}$ and $a_{m}$ in good position with respect
to each other, see Fig.~\ref{bigon}, the third Reidemeister move
correspond to the triangle relations, and the far commutativity does
not change the picture at all~\cite{M3}.

Algebraically, this reduction process stems from the gradient descent algorithm in Coxeter groups.

For us, it is crucial to know that {\em when looking at a free
$2$-braid, one can see which pairs of crossings can be cancelled.}
Once we cancel all possible crossings, we get an invariant picture.

Thus, {\em we get a complete picture-valued invariant of a free
$2$-braid}.

Theorem \ref{theo2} means that this picture (complete invariant)
occurs as a sub-picture in every picture representing the same free
$2$-braid.

However, various homomorphisms  $G_{n}^{k}\to G_{n-1}^{k-1}$, whose
combination leads to homomorphisms of type $G_{n}^{k}\to
G_{n-k+2}^{2}$, allow one to construct lots of invariants of groups
$G_{n}^{k}$ valued in {\em pictures}.

In particular, these pictures allow one to get easy estimates for the
complexity of braids and corresponding dynamics.


\subsection{The Main Theorem on dynamical systems invariants}

Let ${\cal P}$ be a $k$-correct property on the space of admissible
dynamical systems with fixed initial and final states.

Let  $D$ be a pleasant stable dynamical system decribing the motion
of $n$ particles with respect to ${\cal P}$. Let us enumerate all
critical values $t$ corresponding to all multiindices for $D$, as
$t$ increases from $0$ to $1$. With $D$ we associate an element
$c(D)$ of $G_{n}^{k}$, which is equal to the product of $a_{m}$,
where $m$ are multiindices corresponding to critical values of $D$
as $t$ increases from $0$ to $1$.

\begin{theorem}\label{thm:gnk_main_theorem}
Let $D_{0}$ and $D_{1}$ be two equivalent stable pleasant dynamics
with respect to ${\cal P}$. Then $c(D_{0})=c(D_{1})$ are equal as
elements of $G_{n}^{k}$.
\end{theorem}

\begin{proof}
Let us consider an admissible deformation $D_{s}$ between $D_{0}$
and $D_{1}$. For those intervals of values $s$, where $D_{s}$ is
pleasant and stable, the word representing $c(D_{s})$, does not
change by construction. When passing through those values of $s$,
where $D_{s}$ is not pleasant or is not stable, $c(D(s))$ changes as
follows:

\begin{enumerate}
\item Let $s_{0}$ be the value of the parameter deformation, for
which the property ${\cal P}$ holds for some $(k+1)$-tuple of
indices at some time $t=t_{0}$. Note that $D_{s_0}$ is stable.

Consider the multiindex $m=(m_{1},\dots, m_{k+1})$ for which ${\cal
P}$ holds at $t=t_{0}$ for $s=s_{0}$. Let ${\tilde
m}_{j}=m\setminus m_{j},j=1,\dots, k+1$ For types $\tau_{{\tilde
m}_{j}}$, let us choose those coordinates $\zeta_{j}$, which
correspond to  the intersection $t$ at $s=s_{0}$. As $s$ changes,
these type are continuous functions with respect to $s$.

Then for small $\eps$, for $s=s_{0}+\eps$, the word $c(D_{s})$ will
contain a sequence of letters $a_{{\tilde m}_{j}}$ in a certain
order. For values  $s=s_{0}-\eps$, the word $c(D_{s})$ will contain
the same set of letters in the reverse order, see Fig.~\ref{fig:unpleasant_configurations} left. Here we have used the
fact that $D_{s}$ is stable.

\item If for some  $s=s_{0}$ we have a critical value with two different $k$-tuples $m,m'$
possessing ${\cal P}$ and $Card(m\cap m')<k-1$, then the word
$c(D_{s})$ undergoes the relation~\eqref{eq:gnk_far_commutativity_relation} as $s$ passes through $s_{0}$, see Fig.~\ref{fig:unpleasant_configurations} middle;
here we also require the stability of $D_{s_0}$.

\item If at some $s=s_{0}$, the deformation $D_{s}$ is unstable, then $c(D_{s})$
changes according to~\eqref{eq:gnk_order2_relation} as $s$ passes  through $s_{0}$, see Fig.~\ref{fig:unpleasant_configurations} right. Here we use
the fact that the deformation is admissible.

\end{enumerate}

\begin{figure}[h]
\begin{center}
 \includegraphics[width =0.8\textwidth]{ 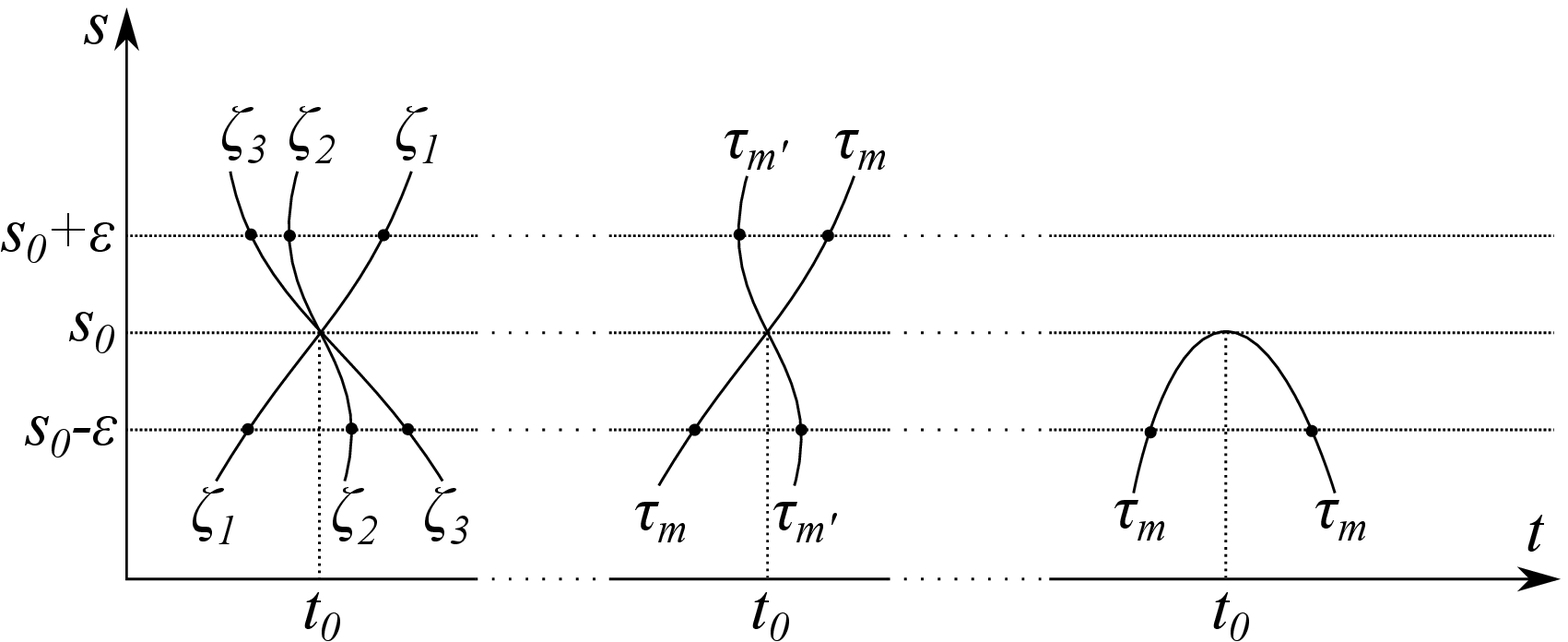}
\end{center}
\caption{Unpleasant and unstable configurations.}\label{fig:unpleasant_configurations}
\end{figure}
\end{proof}

\begin{figure}
\centering\includegraphics[width=200pt]{ 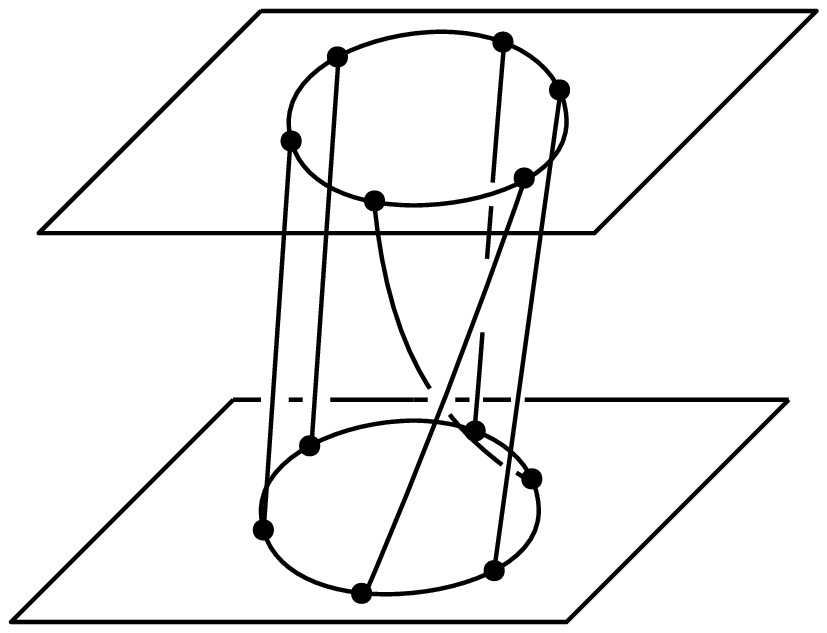}
\caption{Dynamics of points, where the initial and the final configurations are points uniformly distributed
along the unit circle}
\label{fig:exa_config}
\end{figure}

Let us now pass to our main example, the classical braid group. Here
distinct points on the plane are particles. We can require that
their initial and final positions are uniformly distributed
along the unit circle centered at $0$, see Fig.~\ref{fig:exa_config}.

For ${\cal P}$, we take the property to belong to the same line.
This property is, certainly, $3$-good. Every motion of points where
the initial state and the final state are fixed, can be approximated
by a motion where no more than $3$ points belong to the same
straight line at once, and the set of moments where three points
belong to the same line, is finite, moreover, no more than one set
of $3$ points belong to the same line simultaneously. This means
that this dynamical system is pleasant.

Finally, the correctedness of ${\cal P}$ means that if we take two
isotopic braids in general position (in our terminology: two
pleasant dynamical systems connected by a deformation), then by a
small perturbation we can get an admissible deformation for which
the following holds. There are only finitely many values of the
parameter $s$ with four points on the same line or two triples of
points on the same line at the same moment; moreover, for each such
$s$ only one such case occurs exactly for one value of $t$.

In this example, as well as in the sequel, the properties of being
{\em pleasant} and {\em correct} are based on the fact that every
two {\em general position} states can be connected by a curve passing
through states of codimension $1$ (simplest generation) finitely
many times, and every two paths with fixed endpoints, which can be
connected by a deformation, can be connected by a general position
deformation where points of codimensions $1$ and $2$ occur, the
latter happen only finitely many times.

In particular, the most complicated condition saying that the set of
some $(k+1)$ particles satisfies the property ${\cal P}$, the
corresponding derivatives are all distinct, is also a general
position argument. For example, assuming that some $4$ points belong
to the same horizontal line (event of codimension $2$), we may
require there is no coincidence of any further parameters (we avoid
events of codimension $3$).

From the definition of the invariant $c$, one easily gets the
following
\begin{theorem}
Let $D$ be a dynamical system corresponding to a classical braid.
Then the number of horizontal trisecants (see p.~\pageref{def:horizontal_trisecant}) of the braid $D$ is not
smaller than the complexity of the free $3$-braid $\beta=c(D)$.
\end{theorem}

Analogously, various geometrical properties of dynamical systems can
be analysed by looking at complexities of corresponding groups of
free $k$-braids, if one can define a $k$-correct property for these
dynamics, which lead to invariants valued in free $k$-braids.

Let us now collect some situations where the above methods can be
applied.

\begin{enumerate}\label{exa:good_properties}

\item An evident invariant of closed pure classical braids is the conjugacy class of the group $G_{n}^{3}$.
To pass from arbitrary braids to pure braids, one can take some
power of the braid in question.

\item Note that the most important partial case for $k=2$ is the classical Reidemeister braid theory.
Indeed, for a set of points on the plane $Oxy$, we can take for
${\cal P}$ the property that the $y$-coordinates of points coincide.
Then, considering a braid as motion of distinct points in the plane
$z=1-t$ as $t$ changes from $0$ to $1$, we get a set of curves in
space whose projection to $Oyz$ will have intersections exactly in
the case when the property ${\cal P}$ holds (see Fig.~\ref{fig:bisecants}). The additional
information coming from the $x$ coordinate, leads one to the
classical braid theory.

\begin{figure}[h]
\begin{center}
 \includegraphics[width =0.35\textwidth]{ 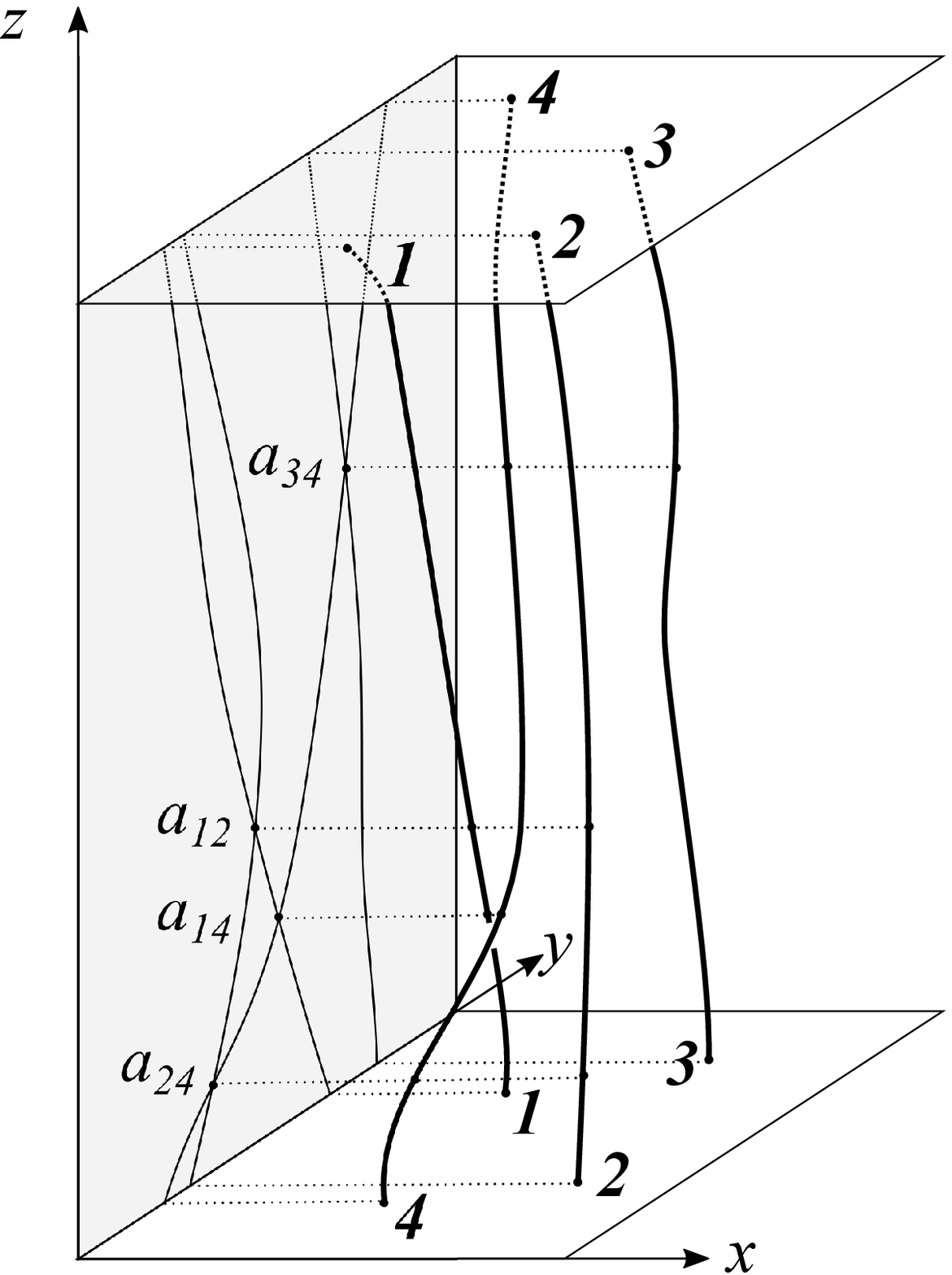}
\end{center}
\caption{Braid diagram.}\label{fig:bisecants}
\end{figure}

\item For classical braids, one can construct invariants for  $k=4$ in a way similar to $k=3$.
In this case, we again take ordered sets of $n$ points on the plane,
and the property ${\cal P}$ means that the set of points belongs to
the same circle or straight line; for three distinct points this
property always holds, and the circle/straight line is unique.

\item With practically no changes this theory can be used for the study of
weavings \cite{Weaves}, collections of projective lines in
$\R{}P^{3}$ considered up to isotopy. Here ${\cal P}$ is the
property of a set of points to belong to the same projective line.
The main difficulty here is that the general position deformation
may contain three lines having infinitely many common horizontal
trisecants. Another difficulty occurs when one of our lines becomes
horizontal; this leads to some additional relations to our groups
$G_{n}^{3}$, which are easy to handle.

\item In the case of points on a $2$-sphere we can define ${\cal P}$ to be the property of
points to belong to the same geodesic. This theory works with an
additional restriction which forbids antipodal points. Some
constraints should be imposed in the case of $2$-dimensional
Riemannian manifolds: for the space of all dynamical systems, we
should impose the restrictions which allow one to detect the
geodesic passing through two points in a way such that if two
geodesics chosen for $a,b$ and for $b,c$ coincide, then the same
geodesic should be chosen for $a,c$.

\item We can use this paradigm to study configurations of balls in a space. A dynamic of balls induces a dynamic of their centers, so any $k$-good property assigns a word in some group $G_n^k$ to the dynamic. The difference between dynamics of balls and ones of a set of points is that the centers of balls can not move too close to each other.

\item In the case of $n$ non-intersecting projective $m$-dimensional planes in $\R{}P^{m+2}$
considered up to isotopy, the theory works as well. Here, in order
to define the dynamical system, we take a one-parameter family of
projective hyperplanes in general position, for particles we take
$(m-1)$-dimensional planes which appear as intersections of the
initial planes with the hyperplane.

The properties of being good, correct etc. follow from the fact that
in general position ``particles'' have a unique same secant line (in
a way similar to projective lines, one should allow the projective
planes not to be straight).
For $m=1$ one should take $k=3$, in the general case one takes
$k=2m+1$.


\item This theory can be applied to the study of fundamental groups of various discriminant
spaces, if such spaces can be defined by several equalities and
subsets of these equalities can be thought of as property ${\cal P}$.

\item The case of classical knots, unlike classical links, is a bit more
complicated: it can be considered as a dynamical system, where the
number of particles is not constant, but it is rather allowed for
two particles to be born from one point or to be mutually
contracted.

The difficulty here is that knots do not possess a group structure,
thus, we don't have a natural order on the set of particles.
Nevertheless, it is possible to construct a map from classical knots
to {\em free $3$-knots} (or free $4$-knots) and study them in a way
similar to free $3$-braids (free $4$-braids).
\end{enumerate}

\begin{remark}
In the case of sets of points in space of dimension $3$, the
property of some points to belong to a $2$-hyperplane (or
higher-dimensional hyperplane) is not correct. Indeed, if three points
belong to the same line, then whatever fourth point we add to them,
the four points will belong to the same hyperplane, thus we can get
various multiindices of $4$ points corresponding to the same moment.

The ``triviality'' of such theory taken without any additional
constraints has the simple description that the configuration space
of sets of points in $\R^{3}$ has trivial fundamental group.
\end{remark}

\begin{remark}
In some cases we can reduce one $k$-good properties to others. For example, consider the $4$-good properties $\mathcal P_1$: ``points in the plane lie on the same circle (or the same line)'', and $\mathcal P_2$: ``points in the $3$-dimensional space lie on the same $2$-plane''.

Let $f\colon \R^2\to \R^3$ be the inverse stereographic projection:
$$
f(u,v) = \left(\frac{2u}{u^2+v^2+1}, \frac{2v}{u^2+v^2+1}, \frac{u^2+v^2-1}{u^2+v^2+1}\right)
$$

The image of the map $f$ is the unit sphere without a point, and the image of any circle or line is a circle on the sphere, i.e. an intersection of the sphere and a $2$-plane. Thus, a set $X\subset \R^2$ contains four points which lie on the same circle or line, if and only if the set $f(X)=\{f(x)\}_{x\in X}\subset \R^3$ contains four points which lie on the same $2$-plane. In other words, the property $\mathcal P_1$ can be reduced to the property $\mathcal P_2$.

We can summarise the reasonings above in the following statement.

\begin{proposition}\label{prop:Gn4_good_property_reduction}
Let $C_n(\R^2)$ be the configuration space of ordered sets consisting of $n$ distinct points in $\R^2$, and $C_n'(\R^3)$ be the configuration space of ordered sets consisting of $n$ distinct points in $\R^3$ such that there are no three points which lie on the same line. Then the following diagram is commutative
$$
\xymatrix{
  \pi_1(C_n(\R^2)) \ar[rr]^{f_*} \ar[dr]_{c_{\mathcal P_1}}
                &  &    \pi_1(C_n'(\R^3)) \ar[dl]^{c_{\mathcal P_2}}    \\
                & G_n^4                }
$$
where $c_{\mathcal P_i}$ is the homomorphism of Theorem~\ref{thm:gnk_main_theorem} determined by the property $\mathcal P_i$, $i=1,2$.
\end{proposition}

A construction above can be used for reduction of other $k$-good properties to the basic $k$-property ``points lie on the same hyperplane in $\R^{k-1}$''.

\end{remark}

\subsection{Pictures}
\label{sec:gnk_pictures}

\begin{figure}[h!]
\begin{minipage}{.4\textwidth}
\begin{center}
 \includegraphics[width =.75\textwidth]{ 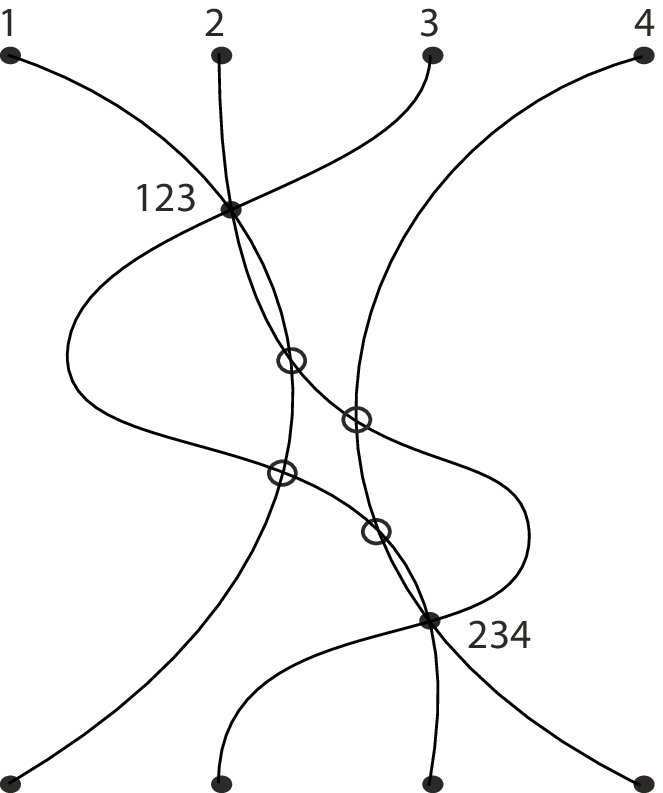}
\end{center}
 \caption{The $3$-braid $a_{234}a_{123}$}\label{3braid}
 \end{minipage}
 \begin{minipage}{.6\textwidth}
 \begin{center}
 \includegraphics[width =.9\textwidth]{ 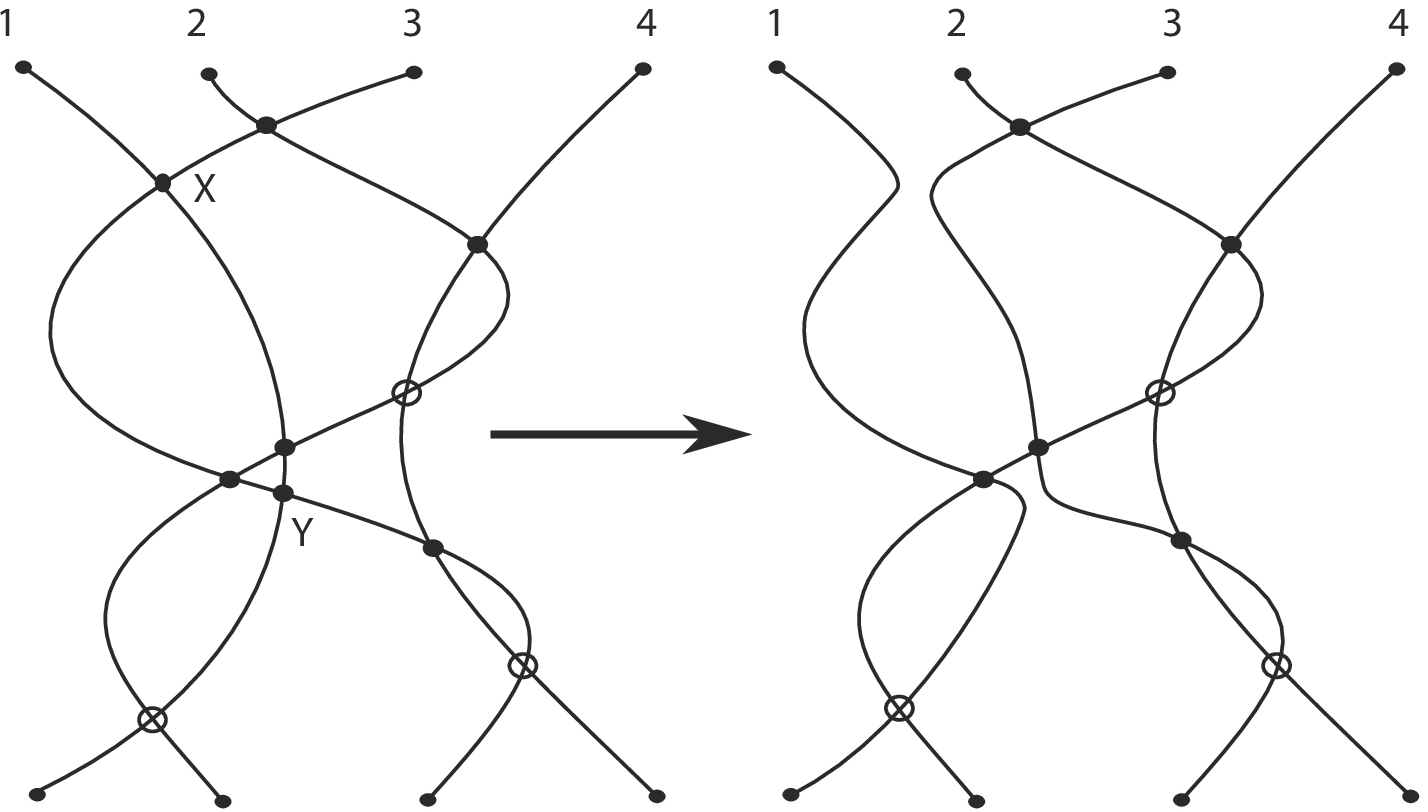}
 \end{center}
 \caption{The bigon reduction} \label{bigon}
 \end{minipage}
 \end{figure}

The free $k$-braids can be depicted by strands connecting points $(1,0),\dots, (n,0)$ to points $(1,1),\dots, (n,1)$; every strand connects $(i,0)$ to $(i,1)$; its projection to the second coordinate is a homeomorphism. We mark crossings corresponding $a_{i_1,\dots,i_k}$ by a solid dot where strands $\#i_1,\dots,\#i_k$ intersect transversally. All other crossings on the plane are artefacts of planar drawing; they do not correspond to any generator of the group, and they are encircled. An example for a 3-free braid is shown in~Fig.~\ref{3braid}.


The clue for the recognition of free $2$-braids is the {\em bigon
reduction} shown in Fig. \ref{bigon}. Here we reduce the bigon whose vertices vertices are $X,Y$.

The graph which appears after all possible bigon reductions (with opposite edge structure at vertices induced by the embedding of the diagram in the plane) is a complete invariant of the free $2$-braid (see also~\cite{M3}). \\

As we mentioned in the introduction, there are nice ``strand forgetting'' and ``strand deleting'' mappings from $G_n^k$ to $G_{n-1}^{k}$ and $G_{n-1}^{k-1}$. They are depicted in Fig.~\ref{(n,k)_(n-1,k)_(n-1,k-1)} and to some extent resemble the formula ${n \choose k} = {n -1 \choose k} + {n-1 \choose k-1}$. In the figure, in the left-hand side we have a depiction of an element of the $G_n^k$ group: a braid with $k$-crossings. In the right-hand side, we have two braids.

The upper one is obtained by deleting the $n$-th strand. All crossings incident to the deleted strand become $(k-1)$-crossings, and all other crossings are forgotten (are replaced with ${k-1 \choose 2}$ virtual crossings each; in the figure we just denote that by transforming a crossing into a collection of non-intersecting strands, not drawing the virtual crossings to keep the figure size manageable). Hence we get an element of $G_{n-1}^{k-1}$.

The bottom braid is obtained by ``strand forgetting'' mapping. Again, we eliminate the $n$-th strand, but this time we forget all crossings incident to the eliminated strand, and leave all other crossings untouched. In this manner we obtain an element of $G_{n-1}^k$.

For example, the forgetting of the last strand maps the element $\beta=a_{123}a_{234}\in G_4^3$ to the element $a_{123}\in G_3^3$, and the deleting of the last strand maps $\beta$ to the element $a_{23}\in G_3^2$.

The deep implications of the existence of these mappings are to be studied in a separate publication.

\begin{figure}[h!]
\begin{center}
 \includegraphics[width = 13cm]{ 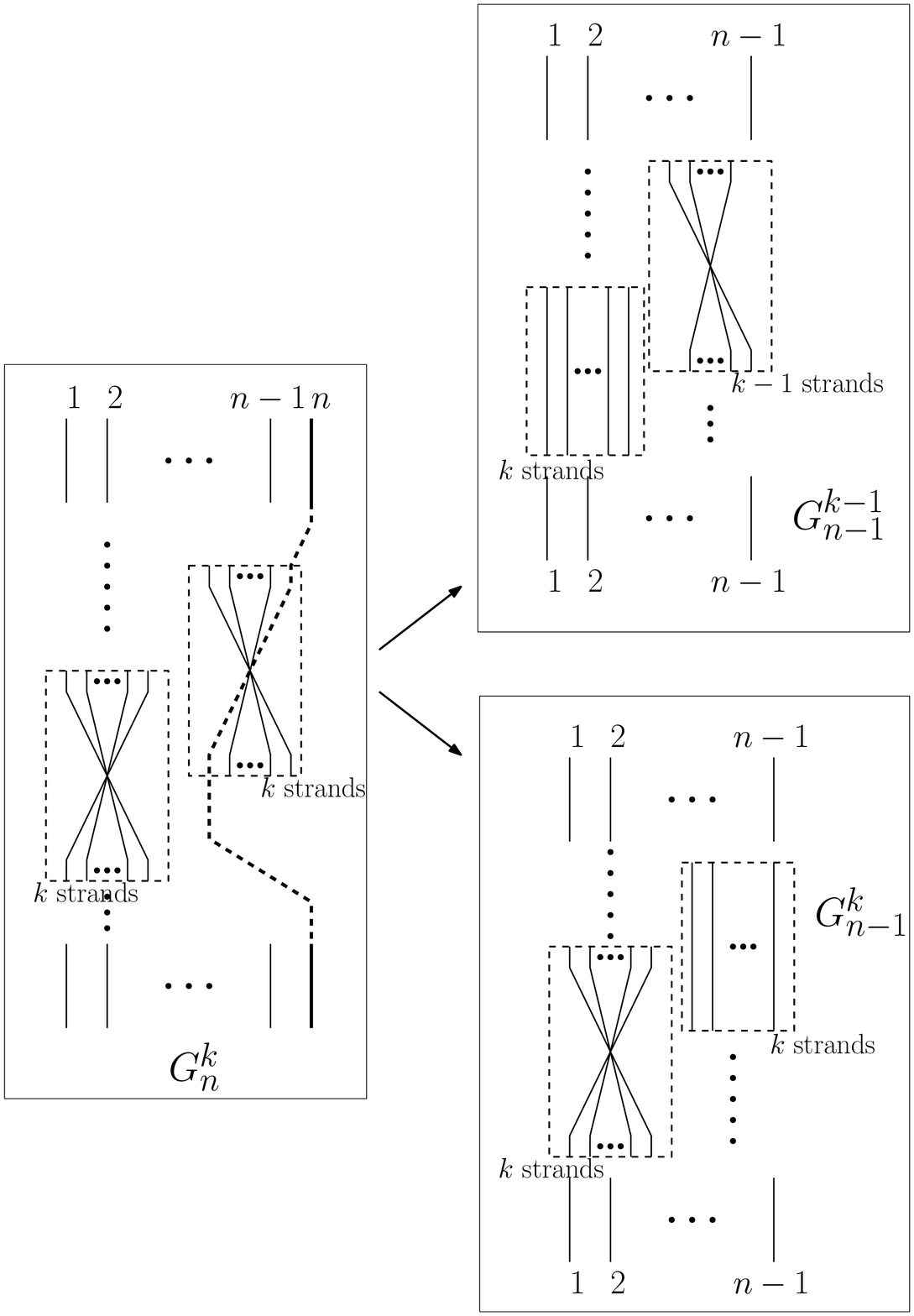}
\end{center}
\caption{Maps from $G_{n}^{k}$ to $G_{n-1}^{k}$ and $G_{n-1}^{k-1}$}\label{(n,k)_(n-1,k)_(n-1,k-1)}
\end{figure} 

%% file: homomorphisms.tex
In the present section, we continue to work with groups $G_{n}^{k}$, which generalise classical braid groups and other groups in a very broad sense.

Perhaps, the easiest highly non-trivial groups are free products of cyclic groups (finite or infinite). In
these groups, the word problem and the conjugacy problem are solved extremely easily: we
just contract a generator with its opposite while possible (in a word or in a cyclic
word) and stop when it is impossible.

This ``gradient descent'' algorithm allows one not only to solve the word problem and conjugacy problem but
also to show that ``if a word  $w$ is irreducible then it is contained in any word $w'$
equivalent to it.'' For example, $abcab$ is contained in $abaa^{-1}cbb^{-1}ab$ in the
free product $\Z*\Z*\Z=\langle a,b,c\rangle$.

Similar phenomena arose in Parity Theory of knot diagrams in low-dimensional topology
discovered by V.O.~Manturov: \\

{\em If a diagram $K$ is odd and irreducible, then it realizes itself as a subdiagram of any other diagram $K'$ equivalent to it}~\cite{knot_parity}. \\

Here ``irreducibility'' is similar to group-theoretic irreducibility and ``oddness''
means that all crossings of the diagram are odd or non-trivial in some sense, see Section~\ref{sec:free_knots} for a short overview of the parity thery.

In fact, parity theory allows one to endow each crossing of a diagram with a powerful diagram-like information,
so that if two crossings (resp., two letters $a$ and $a^{-1}$) are contracted then they have the same pictures.
Thus, a crossing possessing a non-trivial picture takes responsibility for the non-triviality of the whole diagram
(braid, word).

How are these two simple phenomena related to each other?

In the present section, we study classical braid groups. However, we can change the point of view of what we call ``a crossing''. There is a natural presentation of the braid group where generators correspond to horizontal trisecants (see Fig.~\ref{fig:trisecants}). Furthermore, there is a similar presentation where generators correspond to horizontal planes having four points lying on the same circle.

These two approaches lead to the groups $G_{n}^{3}$ and $G_{n}^{4}$ which share many nice properties with virtual braid groups and free groups. In particular, for these groups, each crossing contains very powerful information.

In Section~\ref{sect:hom_free_group}, we coarsen this information and restrict ourselves to invariants of classical braids which appear as the image of the map from groups $G_{n}^{k}$ to free groups.

One immediate advantage of this approach is that we can give some obvious estimates for
various complexities of braids which are easy to calculate.
%
%

\subsection{Homomorphism of pure braids into $G^3_n$}\label{sect:hom_G3n}

We recall (see Definition~\ref{def:gnk}, page~\pageref{def:gnk}) that the {\em $k$-free braid group on $n$ strands}\index{$G_n^k$} $G_{n}^k$ is
generated by elements $a_m$, $m$ is a $k$-element subset of $\{1,2,\dots,n\}$ with the following relations:
\begin{description}
\item[(1)] $a_m^2=1$ for each $m$.
\item[(2)] $a_m a_{m'}=a_{m'}a_m $ for any $k$-element subsets $m, m'$ such that $Card(m\cap m')\leqslant k-2$. This relation is called {\em far commutativity}.\index{Far commutativity}
\item[(3)] $(a_{m_1}a_{m_2}\dots a_{m_{k+1}})^2=1$ for each $(k+1)$-element subset $M=\{
i_1,i_2,\dots, i_{k+1}\}\subset \{1,2,\dots,n\}$ where $m_l=M\setminus\{i_l\},
l=1,\dots,k+1$. This relation is called {\em tetrahedron relation}.\index{Tetrahedron relation}
\end{description}
Groups $G_n^k$ appear naturally as groups describing dynamical systems of $n$
particles in some ``general position''; generators of $G_n^k$ correspond to codimension
$1$ degeneracy, and relations corresponds to codimension $2$ degeneracy which occurs when
performing some generic transformation between two general position dynamical systems.
Dynamical system leading to $G_n^3$ and $G_n^4$ are described in Section~\ref{chap:gnk_definition}, page~\pageref{chap:gnk_definition}; here we
describe the corresponding homomorphisms explicitly.

Generators of $G_n^3$ correspond to configurations in the evolution of dynamical systems
which contains a {\em trisecant}, that is three particles lying in one line. Generators
of $G_n^4$ correspond to configurations which includes four particles that lie in one
circle.

Let $PB_n$ be the pure $n$-strand braid group. It can be presented with the set of
generators $b_{ij}, 1\leqslant i<j \leqslant n,$ and the set of relations~\cite{Bard}
\begin{align}
b_{ij}b_{kl}=b_{kl}b_{ij},\quad  i<j<k<l \mbox{ or } i<k<l<j,\\
b_{ij}b_{ik}b_{jk} = b_{ik}b_{jk}b_{ij} = b_{jk}b_{ij}b_{ik}, \quad i<j<k,\\
b_{jl}b_{kl}b_{ik}b_{jk}=b_{jl}b_{kl}b_{ik}b_{jk},\quad i<j<k<l.
\end{align}

The following proposition is well known.
\begin{proposition}
The center $Z(PB_n)$ of the group $PB_n$ is isomorphic to $\mathbb
Z$.
\end{proposition}

For each different indices $i,j$, $1\leqslant i,j\leqslant n$, we consider the element
$c_{i,j}$ in the group $G_{n}^3$ to be the product
$$c_{i,j}=\prod_{k=j+1}^n a_{i,j,k}\cdot \prod_{k=1}^{j-1} a_{i,j,k}.$$

\begin{proposition}\label{prop_hom_PBn_Gn3}
The correspondence
$$b_{ij}\mapsto c_{i,i+1}^{-1}\dots c_{i,j-1}^{-1}c_{i,j}^2 c_{i,j-1}\dots c_{i,i+1}, i<j,$$
defines a homomorphism $\phi_n\colon PB_n\to G_n^3$.
\end{proposition}

\begin{proof}
Consider the configuration of $n$ points $z_k=e^{2\pi ik/n}, k=1,\dots,n$, in the plane
$\mathbb R^2=\mathbb C$ where the points lie on the unit circle $C=\{z\in\mathbb
C\,|\,|z|=1\}$. Pure braids can be considered as dynamical systems whose initial and
final states coincide with $z = \{z_{1}, \cdots, z_{n}\}$. We can assume that the initial state coincide with the
configuration considered above. Then by Theorem~\ref{thm:gnk_main_theorem} there is a homomorphism
$\phi_n\colon PB_n\to G_{n}^3$ and we need only describe explicitly the images of the
generators of the group $PB_n$.

For any $i<j$ the pure braid $b_{ij}$ can be presented as the following dynamical system:
\begin{enumerate}
\item the point $i$ moves along the inner side of the circle $C$, passes point $i+1, i+2,\dots, j-1$ and land on the circle before the point $j$ (Fig.~\ref{fig:bij_moves} upper left);
\item the point $j$ moves over the point $i$ (Fig.~\ref{fig:bij_moves} upper right);
\item the point $i$ returns to its initial position over the points $j,j-1,\dots, i+1$ (Fig.~\ref{fig:bij_moves} lower left);
\item the point $j$ returns to its position (Fig.~\ref{fig:bij_moves} lower right).
\end{enumerate}

 \begin{figure}
  \centering
  \begin{tabular}{cc}
    \includegraphics[width=0.3\textwidth]{ 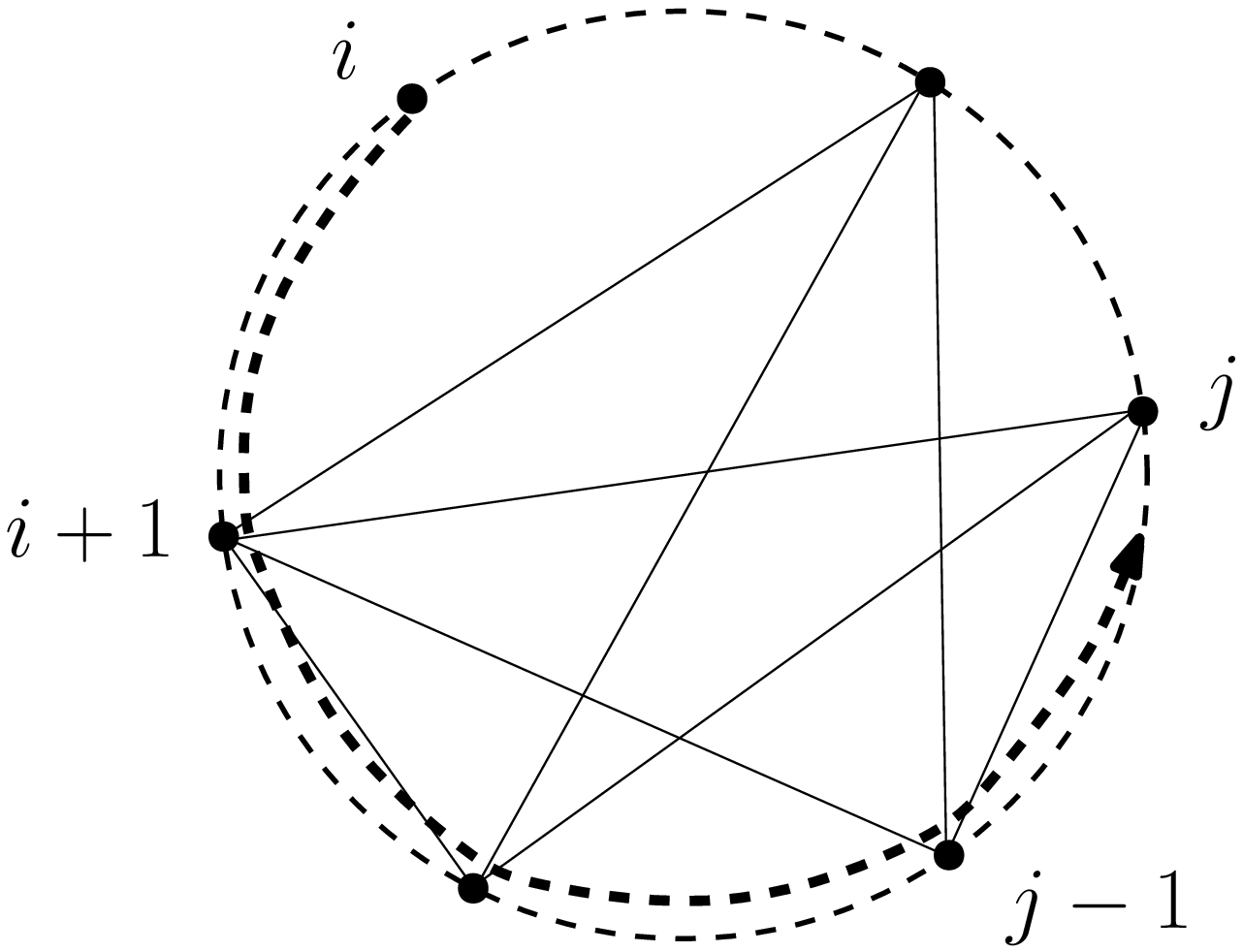} &
    \includegraphics[width=0.3\textwidth]{ 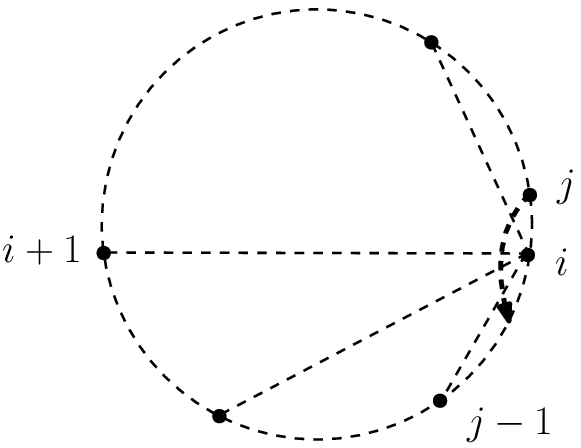} \\
    \includegraphics[width=0.3\textwidth]{ 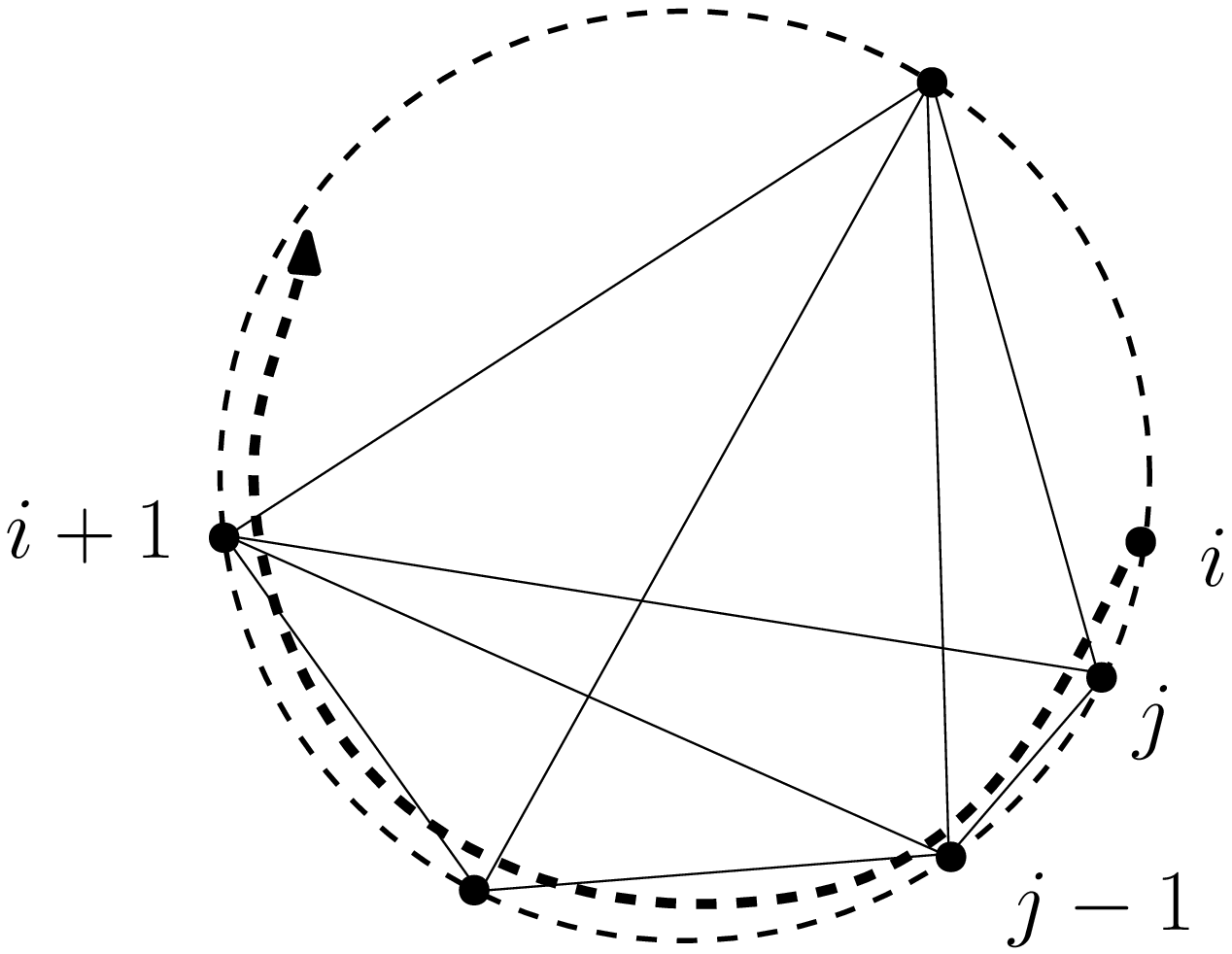} &
    \includegraphics[width=0.3\textwidth]{ 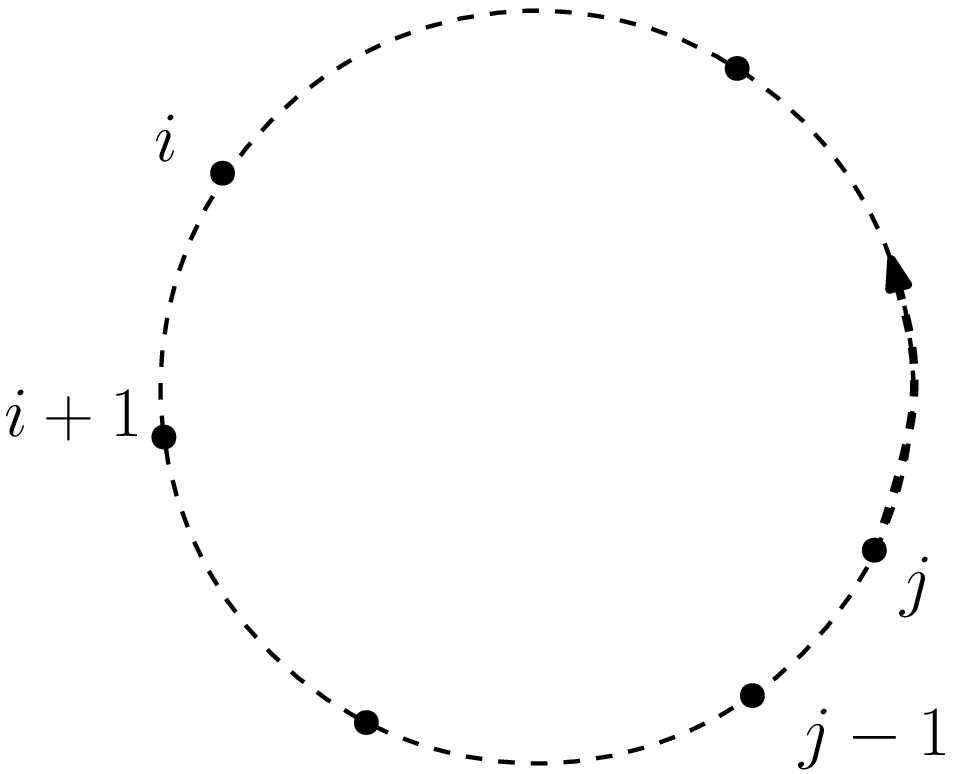}
  \end{tabular}
  \caption{Dynamical system corresponding to $b_{ij}$}\label{fig:bij_moves}
 \end{figure}

As we check all the situations in the dynamical systems where three points lie on the same line,
and write down these situations as letters in a word of the group $G_n^3$ we get exactly the element
$$c_{i,i+1}^{-1}\dots c_{i,j-1}^{-1}c_{i,j}^2 c_{i,j-1}\dots c_{i,i+1}.$$
\end{proof}

Let $\widetilde G=\Z_2 * \Z_2 * \Z_2$ and let $a_1, a_2, a_3$ be the generators of $\wG$.
Let $\wG_{even}$ be the subgroup in $\wG$ which consists of all words of even length. The $2\pi$ rotation of the whole set of points around the origin meets no trisecants, so the map $PB_3\to G_n^3$ has a kernel.

Let us prove that in the case of $3$ strings it is the {\em only} kernel of our map.

\begin{theorem}
There is an isomorphism $PB_3/Z(PB_3)\to \wG_{even}$.
\end{theorem}

\begin{proof}
The quotient group $PB_3/Z(PB_3)$ can be identified with the subgroup $H$ in $PB_3$ that
consists of the braids with unlinked strands $1$ and $2$. In other words,
the subgroup $H$ is the kernel of the homomorphism $PB_3\to PB_2$ that removes the last
strand. For any braid in $H$ we can straighten its strands $1$ and $2$. Looking at such a
braid as a dynamical system we shall see a family of states where the particles $1$ and
$2$ are fixed and the particle $3$ moves.

The points $1$ and $2$ split the line, on which they lie, into three intervals. Let us denote the
unbounded interval with the end $1$ by $a_1$, the unbounded interval with the end $2$ by
$a_2$, and the interval between $1$ and $2$ by $a_3$ (see Fig.~\ref{fig:PB3}).

 \begin{figure}
  \centering
    \includegraphics[width=0.5\textwidth]{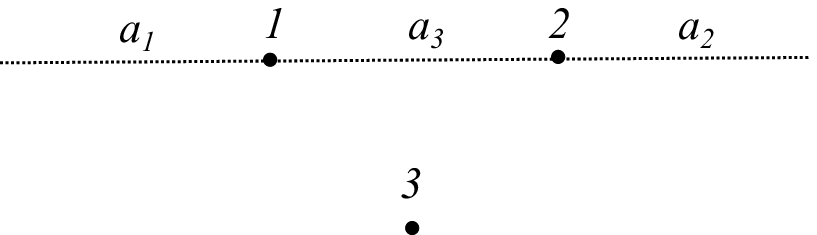}
  \caption{Initial state for a pure braid with $3$ strands}\label{fig:PB3}
 \end{figure}

We assign a word in letters $a_1,a_2,a_3$ to any motion of the point $3$ as follows. We
start with the empty word. Every time the particle $3$ crosses the line $12$ we append
the letter which corresponds to the interval the point $3$ crosses. If the point $3$
returns to its initial position, then the word has even number of letters. This
construction defines a homomorphism $\phi\colon H\to \wG_{even}$.

On the other hand, if we have an even word in letters $a_1, a_2, a_3$ then we can define a
motion of the point $3$ up to isotopy in $\mathbb R^2\setminus \{1,2\}$. Thus, we have a
well-defined map from even words to $H$ that induces a homomorphism $\psi\colon
\wG_{even}\to H\simeq PB_3/Z(PB_3)$.

The homomorphisms $\phi$ and $\psi$ defined in this manner are mutually inverse.
\end{proof}

The crucial observation, which allows us to prove that this map is an isomorphism, is that
we can restore the dynamics from the word in the case of $3$ strands. Since relations in
$G_n^3$ correspond to relations in the braid group (even though the $G_n^3$ group has additional relations), this suffices to prove the isomorphism.

When $n>4$, not all words in $G^3_n$ correspond to the dynamical systems. Thus, the study of mappings $PB_n/Z(PB_n)\to G^3_n$ requires additional techniques.

\subsection{Homomorphism of pure braids into $G^4_n$}\label{sect:hom_G4n}

In the present section, we describe an analogous mapping $PB_n\to G^4_n$; here points
$z_1,z_2,\dots,z_n$ on the plane are in general position: no four points of them belong
to the same circle (or line); codimension $1$ degeneracies will correspond to generators
of $G^4_n$, where at some moment exactly one quadruple of points belongs to the same
circle (line), and relations correspond to the case of more complicated singularities.

Let $a_{\{i,j,k,l\}}, 1\leqslant i,j,k,l\leqslant n$, be the generators of the group
$G^4_n, n>4$.

Let $1\leqslant i<j\leqslant n$. Consider the elements
\begin{align}\label{eq:d_elements}
c^{I}_{ij}&=\prod_{p=2}^{j-1}\prod_{q=1}^{p-1}a_{\{i,j,p,q\}},\\
c^{II}_{ij}&=\prod_{p=1}^{j-1}\prod_{q=1}^{n-j}a_{\{i,j-p,j,j+q\}},\\
c^{III}_{ij}&=\prod_{p=1}^{n-j+1}\prod_{q=0}^{n-p+1}a_{\{i,j,n-p,n-q\}},\\
c_{ij}&=c^{II}_{ij} c^{I}_{ij} c^{III}_{ij}.
\end{align}

\begin{proposition}\label{prop:hom_G4n}
The correspondence

\begin{equation}\label{eq:maps_to_G4n}
b_{ij}\mapsto c_{i,i+1}\dots c_{i,j-1}c_{i,j}^2 c_{i,j-1}^{-1}\dots c_{i,i+1}^{-1},\quad i<j,
\end{equation}
defines a homomorphism $\phi_n\colon PB_n\to G_n^4$.
\end{proposition}

In order to construct this map explicitly, we have to indicate the initial state in the configuration space of $n$-tuple of points. The initial state with all points lying on the circle does not work; so, we shall use the parabola instead.

Let $\Gamma = \{ (t,t^2)\,|\, t\in\R\}\subset \R^2$ be the graph of the function $y=x^2$. Consider a rapidly increasing sequence of positive numbers $t_1,t_2, \dots, t_n$ (precise conditions on the sequence growth may be calculated explicitly but here we just state that suitable sequence exists) and denote the points $(t_i,t_i^2)\in\Gamma$ by $P_i$.

Pure braids can be considered as dynamical systems whose initial and final states coincide and we assume that the initial state is the configuration $\mathcal P=\{P_1,P_2,\dots,P_n\}$. Then by Theorem~\ref{thm:gnk_main_theorem} there is a homomorphism $\psi_n\colon PB_n\to G_{n}^4$. We need to describe explicitly the images of the generators of the group $PB_n$.

\begin{figure}[t!]
  \begin{center}
  \begin{tabular}{cc}
    \includegraphics[width=0.2\textwidth]{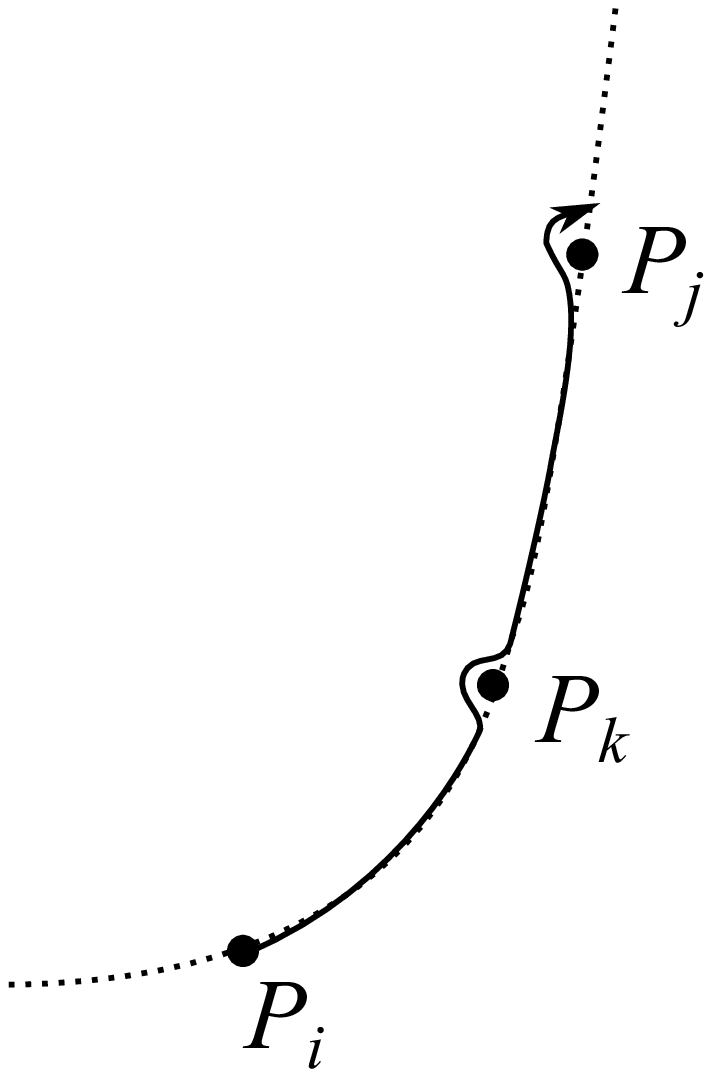} &
    \includegraphics[width=0.2\textwidth]{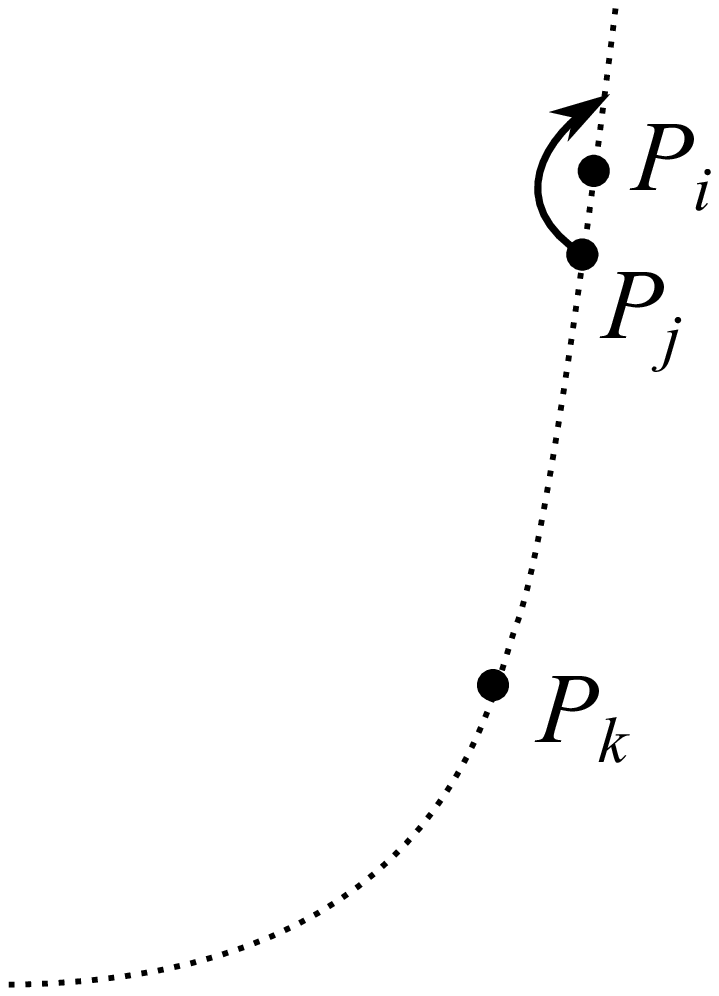} \\
    \includegraphics[width=0.2\textwidth]{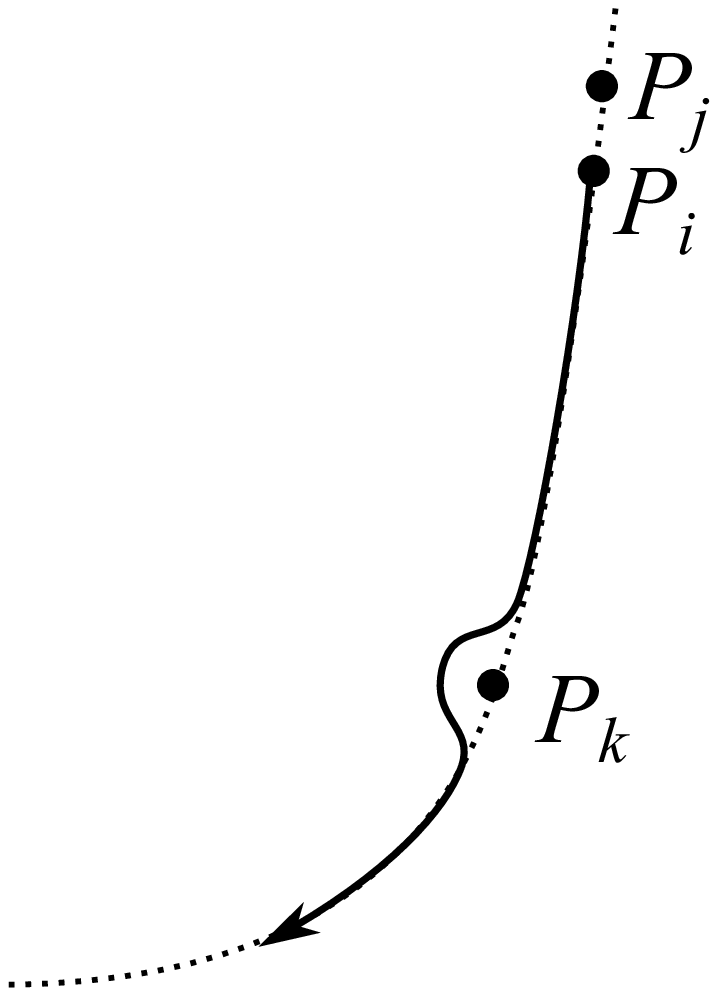} &
    \includegraphics[width=0.2\textwidth]{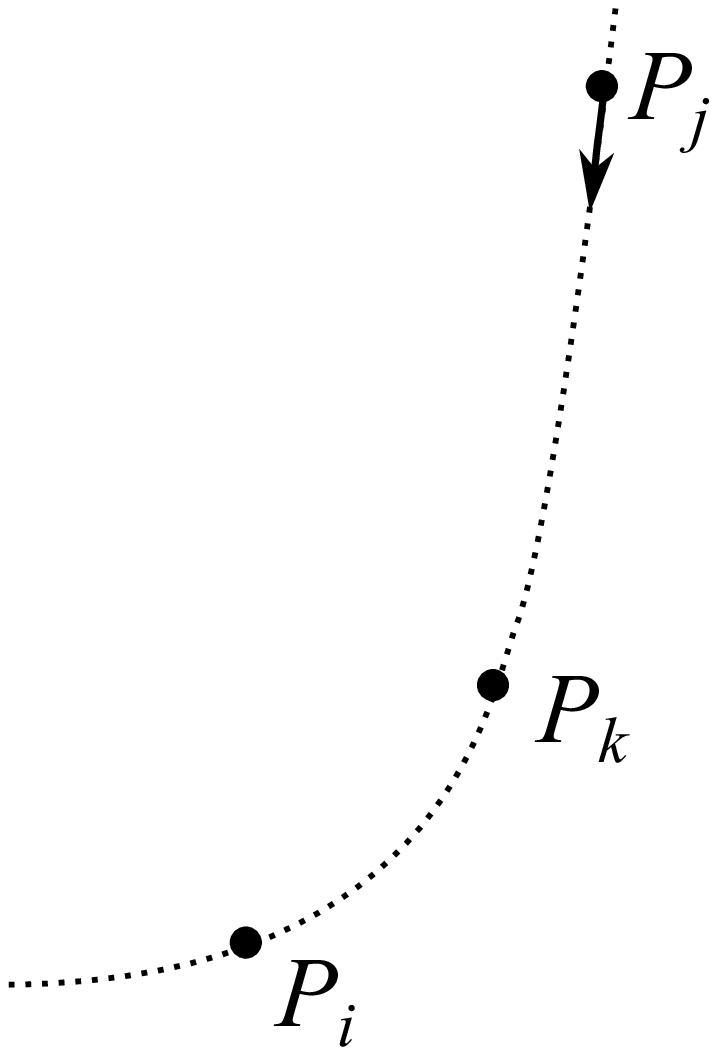}
  \end{tabular}
  \end{center}
  \caption{Dynamical system which corresponds to $b_{ij}$}\label{fig:G4_bij_moves}
 \end{figure}

For any $i<j$ the pure braid $b_{ij}$ can be presented as the following dynamical system: the point $P_{i}$ moves along the graph $\Gamma$ in the following manner:
\begin{enumerate}
\item the point $P_{i}$ moves along the graph $\Gamma$ and passes points $$P_{i+1}, P_{i+2},\dots, P_j$$ from above (Fig.~\ref{fig:G4_bij_moves} upper left);
\item the point $P_{j}$ moves from above the point $P_{i}$ (Fig.~\ref{fig:G4_bij_moves} upper right);
\item the point $P_{i}$ moves to its initial position from above the points $$P_{j-1},\dots, P_{i+1}$$ (Fig.~\ref{fig:G4_bij_moves} lower left);
\item the point $P_{j}$ returns to its position (Fig.~\ref{fig:G4_bij_moves} lower right).
\end{enumerate}

\subsection{Homomorphism into a free group}\label{sect:hom_free_group}

Let $\bar{n}=\{1,2,\dots,n\}$ be the set of indices.

Let $H^k_n\subset G^k_n$ be the subgroup whose elements are given by the {\em even}
words, that are words of $G^k_n$, which include any generator $a_m, m\subset\bar{n}, |m|=k$, evenly
many times (note that each relation in the group $G_n^k$ is of even length, hence for every even word its reduced form is also even, so we may consider {\em all} even words, not just the reduced ones). We construct a homomorphism $\phi$ from the subgroup $H^k_n$ to the free
product of $2^{(k-1)(n-k)}$ copies of the group $\mathbb Z_2$.

Roughly speaking, any letter $a_m$ in $G_n^k$ will get a collection of ``colours'' coming
from the interaction of indices from $m$ with other indices. These ``colours'' will
remain unchanged when performing the far commutativity relation $a_ma_{m'}=a_{m'}a_m$ and the tetrahedron
relation $\dots a_m \ldots = \dots a_m \dots$; moreover, if $a_m^2$ appears inside the word, both letters $a_m$ will have the same colour.

Let $a_m$ (for $m \subset \bar{n}, |m|=k$) be a generator of the group $G^k_n$. Without loss
of generality we can suppose that $m=\{1,2,\dots,k\}$. Assume that $p\in\bar{n}\setminus m$. For any $1\leqslant i\leqslant k$ consider the set
$m[i]=m\setminus\{i\}\cup\{p\}$. Define a homomorphism $\psi_p\colon G^k_n\to
\mathbb Z_2^{\oplus k-1}$ by the formulas $\psi_p(a_{m[i]})=e_i, 1\leqslant i\leqslant
k-1, \psi_p(a_{m[k]})=\sum_{i=1}^{k-1}e_i$ and $\psi_p(a_{m'})=0$ for the other
$m'\subset\bar{n}, |m'|=k$. Here $e_1,\dots, e_{k-1}$ denote the basis elements of the
group $\mathbb Z_2^{\oplus k-1}$. The homomorphism $\psi_p$ is well defined since the
relations in $G^k_n$ are generated by even words.

Consider the homomorphism $\psi = \bigoplus_{p\not\in m}\psi_p$, $\psi\colon G^k_n\to Z$, where the group $Z=\mathbb Z_2^{\oplus (k-1)(n-k)}$.
Note that $H^k_n\subset\ker\psi$.

Let $H=\mathbb Z_2^{\ast |Z|}$ be the free product of $2^{(k-1)(n-k)}$ copies of the group $\mathbb Z_2$
and the exponents of its elements are indiced with the set $Z$. Let $f_x, x\in Z,$ be the generators of the group $H$.
Define the action of the group $G^k_n$ on the set $Z\times H$ by the formula
$$a_{m'}\cdot (x,y)=\left\{\begin{array}{cl}(x,f_xy) & m'=m, \\ (x+\psi(a_{m'}),y) & m'\neq m.\end{array}\right.
$$

This action is well defined. Indeed, $a_{m'}^2\cdot (x,y)=(x,y)$ for any $m'\in\bar{n}, |m'|=k$. On the other hand, for any tuple $\tilde m=(i_1,i_2,\dots, i_{k+1})$ such
that $m\not\subset \tilde m$ one has $\left(\prod_{l=1}^{k+1}a_{\tilde
m\setminus\{i_l\}}\right)^2\cdot (x,y)=(x,y)$ since the generators act trivially here. If
$m\subset \tilde m$, then $\left(\prod_{l=1}^{k+1}a_{\tilde
m\setminus\{i_l\}}\right)^2=w_1 a_m w_2 a_m w_3$, where the word $w_2$ and the word
$w_1w_3$ are products of generators $m[i], i=1,\dots,k$ and
\begin{multline*}
\left(\prod_{l=1}^{k+1}a_{\tilde m\setminus\{i_l\}}\right)^2\cdot (x,y)=w_1 a_m w_2 a_m w_3 \cdot (x,y) =\\
 w_1 a_m w_2 a_m \cdot (x+\psi(w_3),y) =
  w_1 a_m w_2 \cdot (x+\psi(w_3),f_{x+\psi(w_3)}y)=\\ w_1 a_m\cdot (x+\psi(w_3),f_{x+\psi(w_3)}y)=
  w_1\cdot(x+\psi(w_3),f_{x+\psi(w_3)}^2y) =\\ w_1\cdot(x+\psi(w_3),y) = (x+\psi(w_1)+\psi(w_3),y)=(x,y).
\end{multline*}
The fourth and the last equalities follow from the equality $\psi(w_1)+\psi(w_3)=\psi(w_1w_3)=\psi(w_2)=\sum_{i=1}^k\psi(a_{m[i]})=0$.

For any element $g\in G^k_n$ define $\phi_x(g)$ from the relation $$g\cdot
(x,1)=(x+\psi(g),\phi_x(g)1)$$ and let $\phi(g)=\phi_0(g)$. Then $g\cdot
(x,y)=(x+\psi(g),\phi_x(g)y)$ for any $(x,y)\in Z\times H$. If $g\in H^k_n$ then
$\psi(g)=0$ and $g\cdot (0,y)=(0,\phi(g)y)$. Hence, for any $g_1,g_2\in H^k_n$ one has
$(g_1g_2)\cdot(0,1)= g_1\cdot(g_2\cdot(0,1))=g_1\cdot(0,\phi(g_2)) =
(0,\phi(g_1)\phi(g_2))$. On the other hand, $(g_1g_2)\cdot(0,1)=(0,\phi(g_1g_2))$. Thus,
$\phi\colon H^k_n\to H$ is a homomorphism.

For any element $x$ in $H$, let $c(x)$ denote the complexity of $x$, that is the length
of the irreducible representative word of $x$.

Let $b\in PB_n$ be a classical pure braid. By Proposition \ref{prop_hom_PBn_Gn3} there exists a homomorphism $\phi_n$ from $PB_n$ into the group $G_n^3$. Let $\beta = \phi_n(b)$ be the image of $b$ under this homomorphism. Note that the word $\beta$ belongs to $H^3_n$. Then for any $m\subset \bar{n}, |m|=3$, the map $\phi$ defined above is applicable to $\beta$. The geometrical complexity of the braid $b$ can be estimated by complexity of the element $\phi(\beta)$.

\begin{proposition}
The number of horizontal trisecants of the braid $b$ is not less than $c(\phi(\beta))$.
\end{proposition}
The proof of the statement follows from the definition of the maps $\phi_n$ and $\phi$.

Using the homomorphism $PB_n\to G^4_n$, we get an analogous estimation for the number of ``circular quadrisecants'' of the braid.

\subsection{Free groups and crossing numbers}\label{sect:unknotting_number}

The homomorphism described above allows one to estimate the complexity for the braid
group $BP_n$ by using $G^3_n$ and $G^4_n$. These complexities have a geometrical
meaning as the estimates of the number of horizontal trisecants ($G^3_n$) and   ``circular
quadrisecants'' ($G^4_n$).

Let $G=\Z_{2}^{*3}$ be the free product of three copies of $\Z_{2}$ with generators $a,b,c$ respectively.
A typical example of the word in this group is

\begin{equation}\label{eq:word_in_free_group}
w=abcbabca.
\end{equation}


Note that the word~\eqref{eq:word_in_free_group} is irreducible, so, every word equivalent to it contains it as a subword.
In particular, this means that {\em every word $w'$ equivalent to $w$ contains at least $8$ letters.}

The above mentioned complexity is similar to a ``crossing number'', though, crossings are treated in a non-canonical way.

In~\cite{knot_parity}, it is proved that for free knots (which are knot theoretic analogues of the
group $G^2_n$) if an irreducible knot diagram $K$ is complicated enough then it is
contained (as a smoothing) in any knot diagram equivalent to it.

In~\cite{M3}, similar statements are proved for $G^2_n$.

Let us now treat one more complicated issue, the unknotting number.\index{Unknotting number}

To make the issue simpler, let us start with the toy model.
Let $G'=\Z^{*3}$ be the free product of the three copies of $\Z$ with generators $a,b,c$.

Assume that we are allowed to perform the operation of switching the sign $a\longleftrightarrow a^{-1}$,
$b\longleftrightarrow b^{-1}$, $c\longleftrightarrow c^{-1}$ along with the usual reduction
of opposite generators $aa^{-1},bb^{-1},cc^{-1},a^{-1}a,\dots$.

Given a word $v$, how many switches do we need to get the word representing the trivial
element from $v$? How to estimate this number from below?

For the word $w=abcbabca$ 
the answer is ``infinity''. It is impossible to get $1$ from
$w$ in $G'$ because $w$ represents a non-trivial element of $G$, and the operations
$a\longleftrightarrow a^{-1}$ do not change the element of $G$ written in generators
$a,b,c$.

Now, let the element of $G'$ be given by the word $w'=a^{4}b^{2}c^{4}b^{-4}$.

This word $w'$ is trivial in $G$, however, if we look at exponents of $a,b,c$, we see
that they are $+4,-2,+4$. Thus, the number of switchings is bounded below by $\frac
{1}{2}(|4|+|2|+|4|)=5$, and one can find how to make the word $w'$ trivial in $G'$
with five switchings.

For classical braids, a crossing switching corresponds to one turn of one string of a
braid around another: consider a crossing represented by a generator $\sigma_i$; the crossing switch must transform this generator into $\sigma_i^{-1}$. Adding a full turn in the correct direction (represented by $(\sigma_i^{\pm 1})^2$) does just that: $\sigma_i\to\sigma_i(\sigma_i^{-1})^2=\sigma_i^{-1}$, see Fig.~\ref{crossing_switch}.

 \begin{figure}
  \centering
    \includegraphics[width=0.75\textwidth]{ 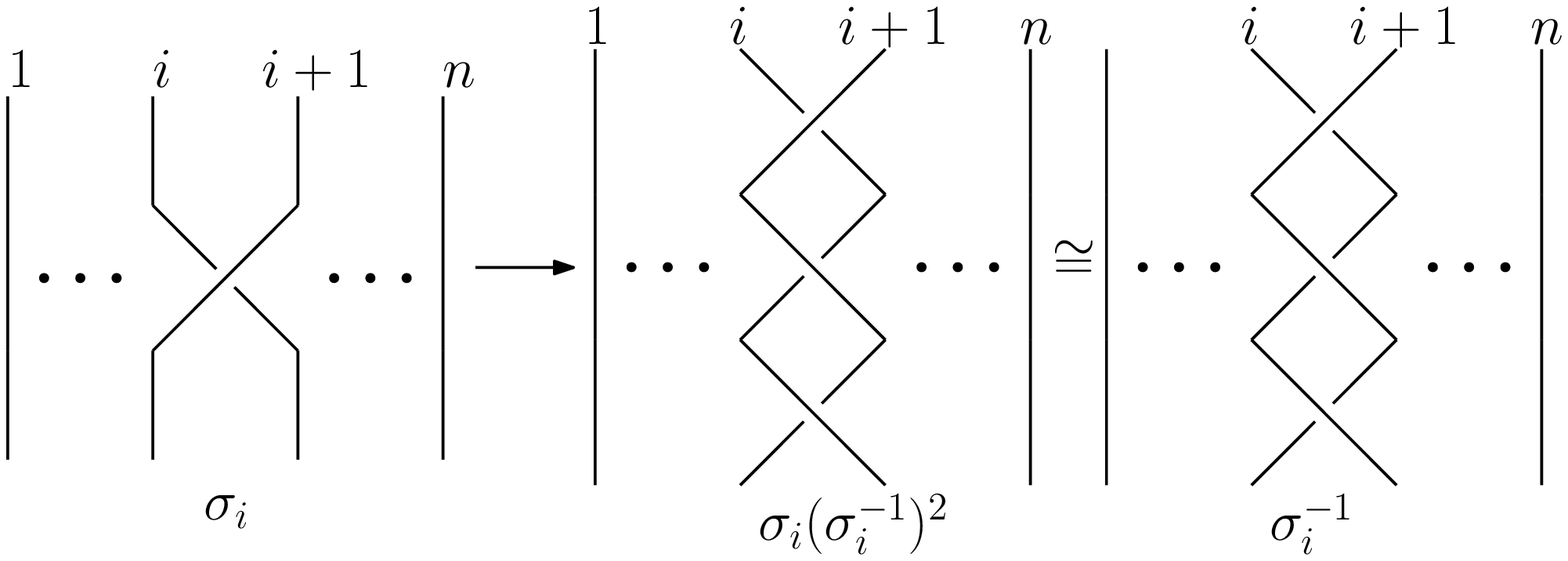}
  \caption{Interpretation of a crossing switch}
  \label{fig:crossing_switch}
 \end{figure}

If $i$ and $j$ are the numbers of two strings of some braid $B$
then the word in $G^k_n$, $k=3,4$, which corresponds to the dynamical system describing
the full turn of the string $i$ around the string $j$, looks like $c_{ij}^2$, where
$c_{ij}=\prod_{m\supset\{i,j\}}a_m$ and the product is taken over the subsets
$m\subset\{1,2,\dots, n\}$, which have $k$ elements and contain $i$ and $j$, given in some
order. Note that the word $c_{ij}^2$ is even.

Fix a $k$-element subset $m\subset\{1,2,\dots, n\}$ containing $\{i,j\}$. Consider the
homomorphism $\phi\colon H^k_n\to H$ determined by $m$, see Section~\ref{sect:hom_free_group} above. The word  $c_{ij}^2$ yields the
factor $f_xf_{x+z_{ij}}\in H$, where the index $x\in Z$ depends on the positition of
$c_{ij}^2$ inside the word $\beta$ of $G^k_n$ that corresponds to the braid $B$, and the
element $z_{ij}=\sum_{m'\,:\, m'\supset\{i,j\}, Card(m\cap m')=k-1}\psi(m')\in Z$ depends
only on $i$, $j$, and $m$ and does not depend on the order in the product $c_{ij}$.
Hence, an addition of a full turn of the string $i$ around the string $j$ corresponds to
substitution of the subword $1=f_x^2$ with the subword $f_xf_{x+z_{ij}}$, that is a
crossing switch corresponds to the switch of elements $f_x$ and $f_{x+z_{ij}}$ in the
image of the braid under the homomorphism $\phi$. Thus, the following statement holds.

\begin{proposition}\label{prop:unknotting_number_estimation}
The unknotting number of a braid $B$ is estimated below by the number of switches
$f_x\mapsto f_{x+z_{ij}}$ which are necessary to make the word $\phi(B)\in H$ trivial.
\end{proposition}

The number of switches in the Proposition~\ref{prop:unknotting_number_estimation} is not
an explicit characteristic of a word in the group $H$, but one can give several rough
estimates for it which can be computed straightforwardly. For example, consider the
following construction. Let $\pi\colon H\to H/[H,H]=\Z_2[Z]$ be the natural projection; let $Z_0$
be the subgroup in $Z$ generated by the elements $z_{ij}, \{i,j\}\subset m$. For any
element $\xi = \sum_{z\in Z}\xi_z z\in \Z_2[Z]$ and any $z\in Z$ consider the number
$c_z(\xi)=Card\left(\{z_0\in Z_0\,|\, \xi_{z+z_0}\ne 0\}\right)$, and let
$c(\xi)=\max_{z\in Z} c_z(\xi)$.

Let $\omega\in H$ be an arbitrary word and $\xi=\pi(\omega)$. Any switch $f_x\mapsto
f_{x+z_{ij}}$ in $\omega$  corresponds to a switch of $\xi$. Note that these switches map
any element $z\in Z$ into the class $z+Z_0\subset Z$. After a switch, two summands of
$\xi$, that correspond to generators in $z+Z_0$, can annihilate. Thus, removing all the
summands of $\xi$ by the generators in $z+Z_0$ takes at least $\frac 12 c_z(\xi)$
switches. Thus, we get the following estimate for the unknotting number.

\begin{proposition}\label{prop:unknotting_number_estimation_rough}
The unknotting number of a braid $B$ is estimated below by the number $\frac 12
c\left(\pi(\phi(B))\right)$.
\end{proposition}

\begin{figure}
  \centering
    \includegraphics[width=0.25\textwidth]{ 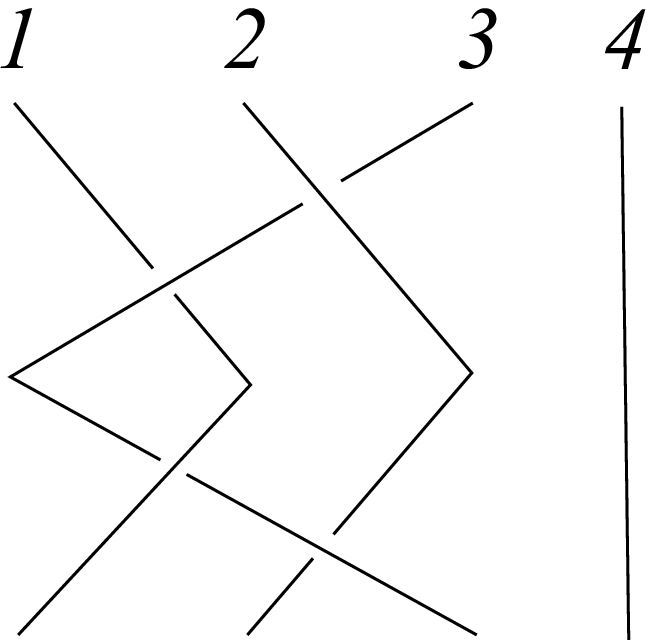}
  \caption{Braid with the unknotting number equal to $2$}\label{fig:braid_example}
 \end{figure}

Consider the braid $B$ in Fig.~\ref{fig:braid_example}. Its unknotting number is $2$. The braid corresponds (see Section~\ref{sect:hom_G3n})to the element
$$\beta=a_{123}a_{234}a_{123}a_{134}a_{123}a_{134}a_{123}a_{234}\in G^3_4.$$
Let $m=\{1,2,3\}$. Then $z_{12}=\psi(a_{124})=e_1+e_2\in Z=\Z_2^{\oplus 2}$, $z_{13}=\psi(a_{134})=e_2$, $z_{23}=\psi(a_{234})=e_1$.

The image of the element $\beta$ is $\phi(\beta)=f_{0}f_{e_1}f_{e_1+e_2}f_{e_1}$. The word $\phi(\beta)$ can not be trivialized with one switch, but it can
be made trivial with two switches:
$$
f_{0}f_{e_1}f_{e_1+e_2}f_{e_1}\stackrel{z_{13}}{\longrightarrow}f_{0}f_{e_1}f_{e_1}f_{e_1}=f_{0}f_{e_1}\stackrel{z_{23}}{\longrightarrow}f_{0}f_{e_0}=1.
$$
Note that the Proposition~\ref{prop:unknotting_number_estimation_rough} gives the
estimate $$\frac 12 c(\pi(\phi(\beta)))=\frac 12 c(f_{0}+f_{e_1+e_2})=1.$$

%% file: generalisations.tex
In the classical braid groups, a special role is played by {\em Brunnian braids on $n$ strands}.\index{Braid!Brunnian} A Brunnian braid is a braid such that each braid obtained by omitting one strand is trivial for every strand. Brunnian braids can be defined algebraically as follow.

Let us define a mapping $p_{m}\colon PB_{n+1} \rightarrow PB_{n}$ \label{map:pm}by
\begin{center}
$p_{m}(b_{ij})  = \left\{
\begin{array}{cc} 
    1 & \text{if}~ j= m, \\
        b_{ij} & \text{if}~ i,j<m, \\
       b_{i(j-1)} & \text{if}~ i<m, j>m,\\
       b_{(i-1)(j-1)} & \text{if}~ i, j>m. \\
\end{array}\right.$
\end{center}
Roughly speaking, this mapping deletes one strand from pure braids on $n$ strands.
\begin{figure}[h!]
\begin{center}
 \includegraphics[width =10cm]{ 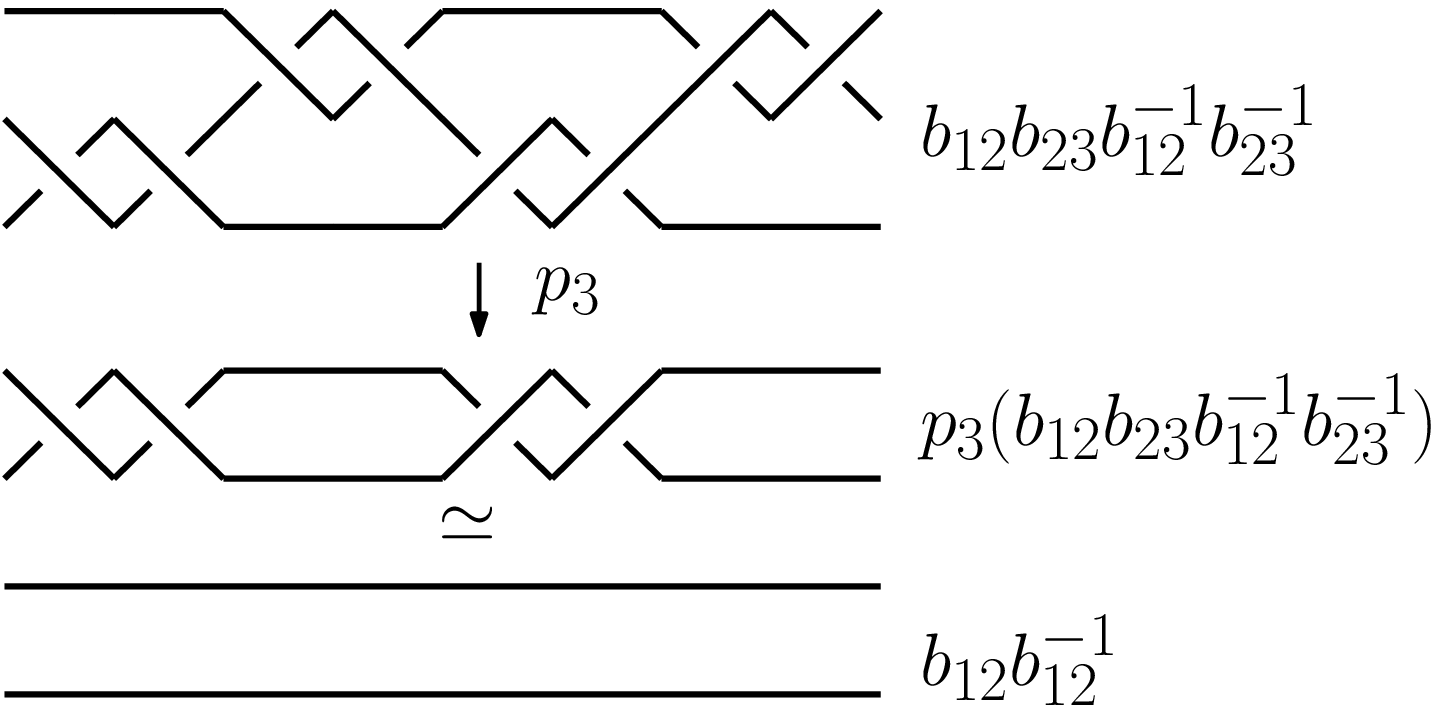}

\end{center}

 \caption{Diagrams of $b_{12}b_{13}b_{12}^{-1}b_{13}^{-1}$ and $p(b_{12}b_{13}b_{12}^{-1}b_{13}^{-1}) = b_{12}b_{12}^{-1}$ in $PB_{3}$}\label{exa-projection}
\end{figure}

\begin{definition}
{\em An $n$-strand Brunnian braid} is a pure braid $\beta$ on $n$ strands such that $p_{m}(\beta) = 1$ for every $m \in \{1,2,\dots, n\}$.
\end{definition}
Note that if an element in $PB_{n}$ is of the form $[\cdots [[b_{i_{1}j_{1}},b_{i_{2}j_{2}}],b_{i_{3}j_{3}}],\cdots],b_{i_{m}j_{m}}]$ and for every $i \in \{1, \dots, n\}$ there is $k$ such that $i \in \{i_{k},j_{k}\}$, then it is Brunnian, for example, Fig.~\ref{exa-projection}.

We can define Brunnian elements in groups $G_{n}^{k}$ and we will see that the image of a Brunnian braid is Brunnian in $G_{n}^{3}$ by a homomorphism from $PB_{n}$ to $G_{n}^{3}$.

Let us define a mapping $q_{m}\colon G_{n+1}^{3} \rightarrow G_{n}^{3}$ \label{map:qm} by

\begin{center}
$q_{m}(a_{ijk})  = \left\{
\begin{array}{cc} 
    1 & \text{if}~ m\in\{i,j,k\}, \\
       a_{ijk} & \text{if}~ i,j,k<m, \\
      a_{ij(k-1)} & \text{if}~ i,j<m, k>m, \\
       a_{i(j-1)(k-1)} & \text{if}~ i<m, j, k>m,\\
       a_{(i-1)(j-1)(k-1)} & \text{if}~ i, j, k>m.\\
   \end{array}\right.$
   \end{center}

\begin{definition}
An element $\beta$ from $G_{n+1}^{3}$ is called {\em Brunnian} if $q_{m}(\beta) = 1$ for every $m \in \{1,2,\cdots, n\}$.
\end{definition}

\begin{lemma}\label{thm_pres_brun}
For a Brunnian braid $\beta \in PB_{n}$, $\phi_{n}(\beta)$ is a Brunnian in $G_{n}^{3}$, where $\phi_n\colon PB_n\to G_n^3$ is the group homomorphism defined in Proposition~\ref{prop_hom_PBn_Gn3}.
\end{lemma}
To prove the above lemma, firstly we prove the following lemma.

\begin{lemma}
For $i,j \in  \{1,\dots ,n \}$ and $i<j$,
\begin{center}
$q_{n}(c_{i,j})  = \left\{
\begin{array}{cc} 
    1 & \text{if}~ j = n, \\
       c_{i,j} & \text{if}~j \neq n. \\
     \end{array}\right.$
   \end{center}
\end{lemma}

\begin{proof}
If $j=n$, then $$q_{n}(c_{i,n}) = q_{n}(  \prod_{k=1} a_{ink})= \prod_{k=1}q_{n}(a_{ink}) =1.$$
If $i,j \neq n$, then $$q_{n}(c_{i,j}) = q_{n}(  \prod_{k=j+1} a_{ijk} \prod_{k=1}^{j-1} a_{ijk})=$$  $$=\prod_{k=j+1} q_{n}(a_{ijk}) \prod_{k=1}^{j-1}q_{n}(a_{ijk}) = \prod_{k=j+1} a_{ijk} \prod_{k=1}^{j-1} a_{ijk}= c_{i,j}.$$\\
\end{proof}

Analogously we can show that 
\begin{center}
$q_{m}(c_{i,j})  = \left\{
\begin{array}{cc} 
    1 & \text{if}~ j= m ~\text{or}~ i=m, \\
        c_{i,j} & \text{if}~ i,j<m,\\
      c_{i,(j-1)} & \text{if}~ i<m, j>m,\\
       c_{(i-1),(j-1)} & \text{if}~ i, j>m.\\
   \end{array}\right.$
   \end{center}

\begin{proof}[Proof of Lemma~\ref{thm_pres_brun}]
It is sufficient to show that $q_{m} \circ \phi_{n} = \phi_{n-1} \circ p_{m}$, because if $p_{m}(\beta) =1$, then $q_{m} \circ \phi_{n}(\beta) = \phi_{n-1} \circ p_{m}(\beta) = \phi_{n-1}(1)=1$. For $m=n$ and $b_{ij} \in PB_{n}$, if $j=n$, then $\phi_{n-1} \circ p_{n}(b_{in}) = 1$ and

\begin{eqnarray*}
q_{n} \circ \phi_{n} (b_{in}) &=& q_{n}((c_{i,i+1})^{-1}(c_{i,i+2})^{-1} \cdots (c_{i,n-1})^{-1} (c_{i,n})^{2} c_{i,n-1} \cdots c_{i,i+2} c_{i,i+1}) \\
 &=& (c_{i,i+1})^{-1}(c_{i,i+2})^{-1} \cdots (c_{i,n-1})^{-1} q_{n}((c_{i,n})^{2}) c_{i,n-1} \cdots c_{i,i+2} c_{i,i+1} \\ 
 &=&  (c_{i,i+1})^{-1}(c_{i,i+2})^{-1} \cdots (c_{i,n-1})^{-1} c_{i,n-1} \cdots c_{i,i+2} c_{i,i+1} =1.
\end{eqnarray*}
 If $i,j \neq n$, then 
 \begin{eqnarray*}
  \phi_{n-1} \circ p_{n}(b_{ij}) &=& \phi_{n-1}(b_{ij}) \\ 
  &=&  (c_{i,i+1})^{-1}(c_{i,i+2})^{-1} \cdots (c_{i,j-1})^{-1} (c_{i,j})^{2} c_{i,j-1} \cdots c_{i,i+2} c_{i,i+1},
 \end{eqnarray*}

 and
 \begin{eqnarray*}
 q_{n} \circ \phi_{n} (b_{ij}) &=& q_{n}((c_{i,i+1})^{-1}(c_{i,i+2})^{-1} \cdots (c_{i,j-1})^{-1} (c_{i,j})^{2} c_{i,j-1} \cdots c_{i,i+2} c_{i,i+1})\\
 &=& (c_{i,i+1})^{-1}(c_{i,i+2})^{-1} \cdots (c_{i,j-1})^{-1} (c_{i,j})^{2} c_{i,j-1} \cdots c_{i,i+2} c_{i,i+1}.
\end{eqnarray*}

Analogously, it is easy to show that $q_{m} \circ \phi_{n} = \phi_{n-1} \circ p_{m}$. 
\end{proof}

Recall that in Section~\ref{sect:hom_free_group} a group homomorphism $\phi \colon H_{n}^{k} \to \mathbb{Z}_{2}^{2^{(k-1)(n-k)}}$ is defined. Since $Im(\phi_{n}) \subset H_{n}^{3} \subset G_{n}^{3}$ for $\phi_{n}\colon PB_{n} \to G_{n}^{3}$, the composition $\phi \circ \phi_{n}\colon PB_{n} \to  H_{n}^{3} \subset G_{n}^{3} \to \mathbb{Z}_{2}^{2^{2(n-3)}}$ gives an invariant for pure braids. For simplicity we call it {\em the MN invariant}.

The aim of this section is to show that the MN invariants constructed in Section~\ref{sect:hom_free_group} fail to recognize the non-triviality of Brunnian braids (Lemma~\ref{failrecognizebrunnian}). We can enhance these invariants by using the structure of $G_{n}^{k}$; in fact, we shall make only one step allowing us to recognize the commutator, see example \ref{recog_brun}.  But in principle, it is possible to go on enhancing the invariants coming from $G_{n}^{k}$ (even from $G_{n}^{3}$) to get invariants which recognize the non-triviality of commutators of arbitrary lengths: $[[[[[b_{12},b_{13}],b_{14}],b_{15},...]$. More precisely, the group $G_{n}^{3}$ itself recognizes the non-triviality of such braids, and the corresponding invariants can be derived as maps from $G_{n}^{3}$ to free products of $\mathbb{Z}_{2}$ (Example~\ref{recog_brun}). 

Firstly, we reformulate the definition of mappings $\phi\colon G_{n}^{3} \rightarrow \mathbb{Z}_{2}^{2^{2(n-3)}}$, which is defined in Section~\ref{sect:hom_free_group} and construct new ones. In the present section we denote the mapping from $G_{n}^{3}$ to $\mathbb{Z}_{2}^{2^{2(n-3)}}$ by $\phi_{(i,j,k)}$ instead of $\phi$.


 
Let us start reformulating the definition of homomorphisms $\phi_{(i,j,k)}\colon G_{n}^{3} \rightarrow \mathbb{Z}_{2}^{2^{2(n-3)}}$. \label{map:phi-ijk}

Let $\beta \in G_{n}^{3}$ be an even element. For each $c=a_{ijk}$ of $\beta$ and for $ l \in \{1,2,\cdots, n\} \backslash \{i,j,k\}$, define $i_{c}(l)$ by 
$$i_{c}(l) = (N_{jkl}+N_{ijl}, N_{ikl}+N_{ijl}) \in \mathbb{Z}_{2} \times \mathbb{Z}_{2},$$
where $N_{ikl}$ is the number of $a_{ikl}$ from the start of $\beta$ to the crossing $c$. Note that $i_{c}$ can be considered as a map from $\{1,2,\cdots, n \} \backslash \{i,j,k\} $ to $\mathbb{Z}_{2} \times \mathbb{Z}_{2}$. Fix $i,j,k\in \{1,\dots, n\}$. Let $\{c_{1}, \cdots, c_{m}\}$ be the set of $a_{ijk}$ such that for each $s , t \in \{1, 2, \cdots, m\}$, $ s< t$ if and only if we meet $c_{s}$ earlier than $c_{t}$ in $\beta$. Define a group $F_{n}^{3}$ by the group presentation generated by $\{ \sigma ~|~ \sigma\colon \{1,2,\cdots n\} \backslash \{i,j,k\} \rightarrow \mathbb{Z}_{2} \times \mathbb{Z}_{2} \}$ with relations $\{ \sigma^{2} = 1\}$. Note that $i_{c}$ is a mapping from $\{1,2,\dots n\} \backslash \{i,j,k\}$ to $\mathbb{Z}_{2} \times \mathbb{Z}_{2}$ and $i_{c}$ is in  $F_{n}^{3}$. It is easy to show that $F_{n}^{3}$ is isomorphic to the free product of $2^{2(n-3)}$ copies of $\mathbb{Z}_{2}$. Define a word $\phi_{(i,j,k)}(\beta)$ in $F_{n}^{3}$ for $\beta$ by $\phi_{(i,j,k)}(\beta) = i_{c_{1}}i_{c_{2}} \cdots i_{c_{m}}$.

\begin{lemma}\label{failrecognizebrunnian}
For a Brunnian $\beta \in PB_{n}$, $\phi_{(i,j,k)}(\phi_{n}(\beta)) =1$.
\end{lemma}

\begin{proof}
It is sufficient to show that $\phi_{(i,j,k)}(\phi_{n}(b_{lm})) =1$ for $|\{l,m \} \cap \{i,j,k\}|<2$.
For a Brunnian braid $\beta \in PB_{n}$ and for $l \not\in \{i,j,k\}$, let $\beta_{l}$ be a braid obtained by omitting $b_{lm}$ from $\beta$. Since $\phi_{(i,j,k)}(\phi_{n}(b_{lm})) =1$ for $|\{l,m \} \cap \{i,j,k\}|<2$,
$$\phi_{(i,j,k)}(\phi_{n}(\beta)) = \phi_{(i,j,k)}(\phi_{n}(\beta_{l})).$$
Since $\beta$ is Brunnian, $\beta_{l}$ is trivial and $\phi_{(i,j,k)}(\phi_{n}(\beta)) = \phi_{(i,j,k)}(\phi_{n}(\beta_{l})) =1.$ Now we will show that the statement is true. By the definition of $\phi_{n}$,
$$\phi_{n}(b_{ij}) = (c_{i,i+1})^{-1}(c_{i,i+2})^{-1} \cdots (c_{i,j-1})^{-1} (c_{i,j})^{2} c_{i,j-1} \cdots c_{i,i+2} c_{i,i+1}.$$  
Note that for $|\{l,m\} \cap \{i,j,k\}| \neq 2$, $c_{l,m}$ contains no $a_{ijk}$.  Without loss of generality we may assume that $i<j<k$ and $l<m$. There are 8 subcases:
\begin{enumerate}
\item $\{ l, m \} \cap \{i,j,k\} = \emptyset$,
\item $ l \not\in \{i,j,k\}, m \in  \{i,j,k\}$, 
\item  $l=i, i<m<j$,
\item $l=i, j<m<k$,
\item $l=i, k<m<n$,
\item $l=j, j<m<k$,
\item  $l=j, k<m<n$,
\item   $l=k, k<m<n$. 
\end{enumerate}
If $\{ l, m \} \cap \{i,j,k\} = \emptyset$, then $\phi_{n}(b_{ij})$ has no $c_{i,j}$, $c_{i,k}$ and $c_{j,k}$. Hence $\phi_{(i,j,k)}(\phi_{n}(b_{lm})) =1$. Analogously  $\phi_{(i,j,k)}(\phi_{n}(b_{lm})) =1$ in the cases of (2), (3), (6) and (8). \\
If $l=i, j<m<k$, then
$$\phi_{n}(b_{im}) = (c_{i,i+1})^{-1} \cdots (c_{i,j})^{-1} \cdots (c_{i,j-1})^{-1} (c_{i,m})^{2} c_{i,j-1} \cdots c_{i,j} \cdots c_{i,i+1}$$ 
has just two $a_{ijk}$, say $c_{1}=a_{ijk}$ in $(c_{i,j})^{-1}$ and  $c_{2}=a_{ijk}$ in $c_{i,j}$, respectively. Since the number of each $a_{stu}$ between $c_{1}$ and $c_{2}$ is even for every $s,t,u \in \{1,\dots n\}$, $i_{c_{1}} = i_{c_{2}}$. Therefore $\phi_{(i,j,k)}(\phi_{n}(b_{lm})) =1$. Analogously $\phi_{(i,j,k)}(\phi_{n}(b_{lm})) =1$ in the case of (7). \\
If $l=i, k<m<n$, 
$$\phi_{n}(b_{im}) = (c_{i,i+1})^{-1} \cdots (c_{i,j})^{-1} \cdots (c_{i,k})^{-1} \cdots (c_{i,m})^{2} \cdots c_{i,k} \cdots c_{i,j} \cdots c_{i,i+1},$$ 
and it has four $a_{ijk}$, say $c_{1} = a_{ijk}$ in $(c_{i,j})^{-1}$, $c_{2} = a_{ijk}$ in $(c_{i,k})^{-1}$, $c_{3} = a_{ijk}$ in $c_{i,k}$ and $c_{4} = a_{ijk}$ in $c_{i,j}$, respectively. Since the number of each $a_{stu}$ between $c_{2}$ and $c_{3}$ is even for every $s,t,u \in \{1,\dots n\}$, $i_{c_{2}} = i_{c_{3}}$. Similarly, $i_{c_{1}} = i_{c_{4}}$. Therefore $\phi_{(i,j,k)}(\phi_{n}(b_{im})) = i_{c_{1}} i_{c_{2}} i_{c_{3}} i_{c_{4}}=1$. The proof is complete. 
\end{proof}

By the above Lemma the MN-invariant for $G_{n}^{3}$ does not recognize the non-triviality of Brunnian braids in $PB_{n}$. In the next paragraph to make the ``indices'' stronger, we will use the parity for elements of the group $G_{n}^{2}$.

\subsection{Groups $G_{n}^{2}$ with parity and points}

When mapping $G_{n}^{3}$ to $G_{n-1}^{2}$, we lose a lot of information. But it turns out that we can save some of that information in the form of parity.

\begin{definition}
For a positive integer $n>2$, let us define {\em the group $G_n^2$ with parity} $G_{n,\mathcal{P}}^{2}$\index{$G_n^k$!with parity} as the group presentation generated by $ \{ a_{\{i,j\}}^{\epsilon} ~|~ \{i,j\} \subset \{1, \dots, n\}, i < j, ~\epsilon \in \{0,1\}  \}$ subject to the following relations:
\begin{center}\begin{enumerate}
\item $(a_{\{i,j\}}^{\epsilon})^{2} = 1$, $\epsilon \in \{0,1\}$ and $i,j \in \{1,\dots, n\}$,
\item $a_{\{i,j\}}^{\epsilon_{\{i,j\}}}a_{\{k,l\}}^{\epsilon_{\{k,l\}}} = a_{kl}^{\epsilon_{\{k,l\}}}a_{\{i,j\}}^{\epsilon_{\{i,j\}}}$ for $\{i,j\}\cap\{k,l\}=\emptyset$,
\item $a_{\{i,j\}}^{\epsilon_{\{i,j\}}}a_{\{i,k\}}^{\epsilon_{\{i,k\}}}a_{\{j,k\}}^{\epsilon_{\{j,k\}}} = a_{\{j,k\}}^{\epsilon_{\{j,k\}}}a_{\{i,k\}}^{\epsilon_{\{i,k\}}}a_{\{i,j\}}^{\epsilon_{\{i,j\}}}$, for distinct $i,j,k$, where $\epsilon_{\{i,j\}}+\epsilon_{\{i,k\}}+\epsilon_{\{j,k\}} \equiv 0$ mod $2$.
\end{enumerate}
\end{center}
We denote $a_{\{i,j\}}^{\epsilon}$ by $a_{ij}^{\epsilon}$.
\end{definition}

\begin{figure}[h!]
\begin{center}
 \includegraphics[width =8cm]{ 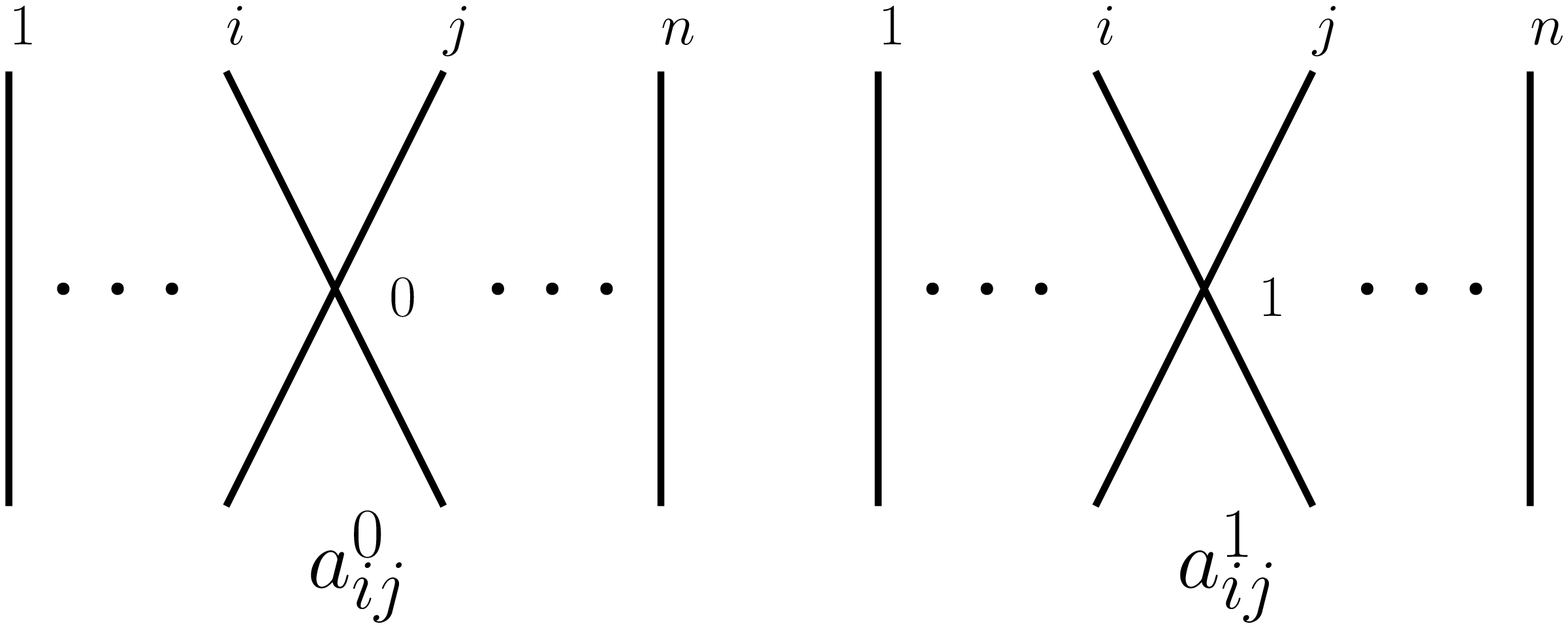}

\end{center}

 \caption{Geometrical meanings of $a_{ij}^{0}$ and $a_{ij}^{1}$ }\label{crossing_parity}
\end{figure}

 We call $a_{ij}^{0}$($a_{ij}^{1}$) {\it an even generator} ({\it an odd generator}). Geometrically, $a_{ij}^{0}$ and $a_{ij}^{1}$ correspond to crossings labeled by $0$ and $1$, respectively, see Fig.~\ref{crossing_parity}.

\begin{figure}[h!]
\begin{center}
 \includegraphics[width =6cm]{ 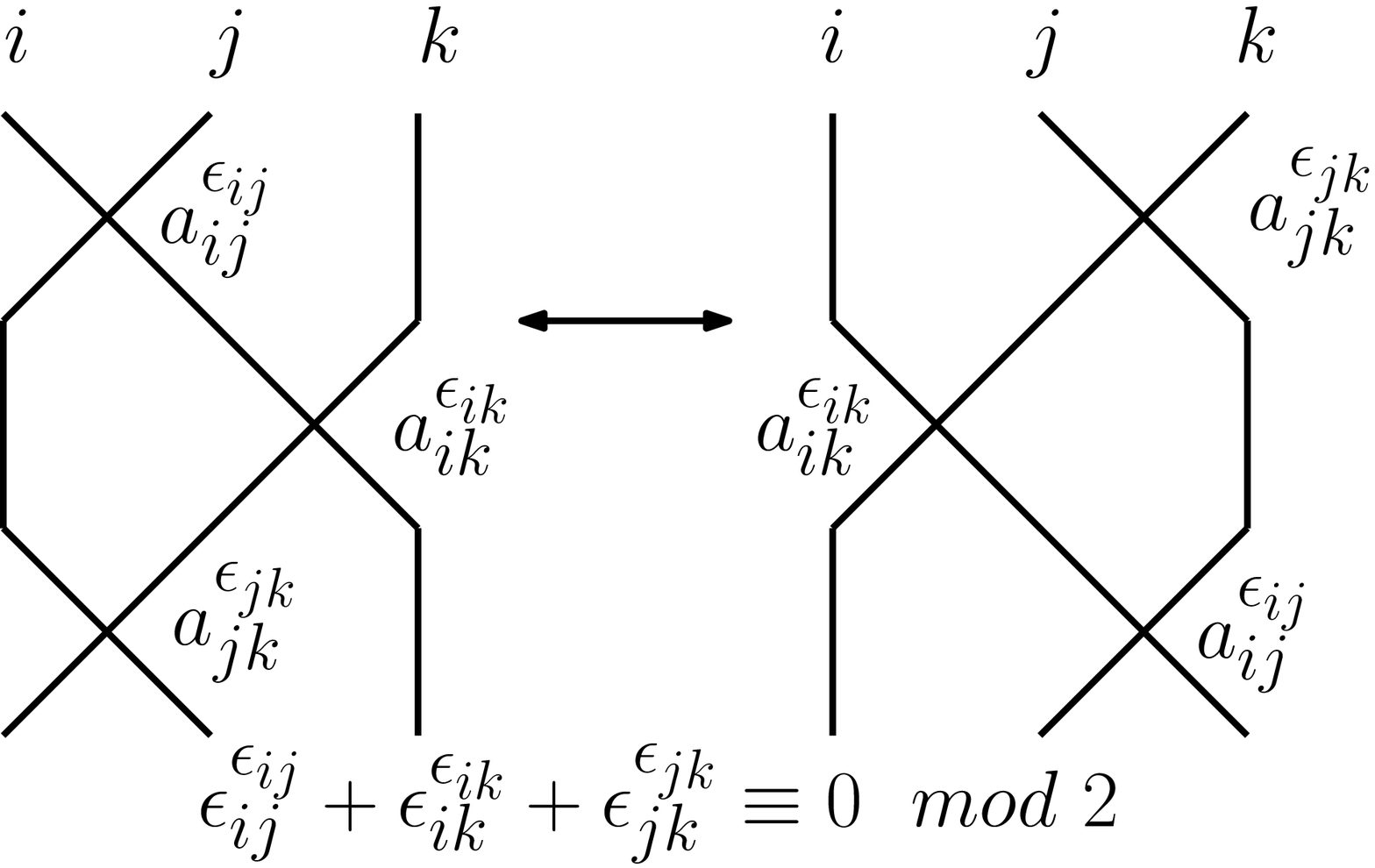}

\end{center}

 \caption{Geometrical meanings of $a_{\{i,j\}}^{\epsilon_{\{i,j\}}}a_{\{i,k\}}^{\epsilon_{\{i,k\}}}a_{\{j,k\}}^{\epsilon_{\{j,k\}}} = a_{\{j,k\}}^{\epsilon_{\{j,k\}}}a_{\{i,k\}}^{\epsilon_{\{i,k\}}}a_{\{i,j\}}^{\epsilon_{\{i,j\}}}$ where $\epsilon_{\{i,j\}}+\epsilon_{\{i,k\}}+\epsilon_{\{j,k\}} \equiv 0$ mod $2$ }\label{RM3_parity}
\end{figure}

The relations of $G_{n,\mathcal{P}}^{2}$ are closely related to the axioms of parity. Especially, the third relation ``$a_{\{i,j\}}^{\epsilon_{\{i,j\}}}a_{\{i,k\}}^{\epsilon_{\{i,k\}}}a_{\{j,k\}}^{\epsilon_{\{j,k\}}} = a_{\{j,k\}}^{\epsilon_{\{j,k\}}}a_{\{i,k\}}^{\epsilon_{\{i,k\}}}a_{\{i,j\}}^{\epsilon_{\{i,j\}}}$, for distinct $i,j,k$, where $\epsilon_{\{i,j\}}+\epsilon_{\{i,k\}}+\epsilon_{\{j,k\}} \equiv 0$ mod $2$'' corresponds to the third Reidemeister move. Recall that the third Reidemeister move is allowed if the sum of parities of three crossings contained in it is equal to zero modulo $2$ see Fig.~\ref{RM3_parity}.

 If a word $\beta$ in $G_{n,\mathcal{P}}^{2}$ has no odd generators, then we call $\beta$ {\it an even word} (or {\it an even element}). 
$G_{n}^{2}$ can be considered as a subgroup of $G_{n,\mathcal{P}}^{2}$ by a homomorphism $\iota$ from $G_{n}^{2}$ to $G_{n,\mathcal{P}}^{2}$ defined by $\iota(a_{ij}) = a_{ij}^{0}$.\label{map:iota} 
Moreover, the following statement can be proved.

\begin{lemma}\label{emb_to_parity}
The homomorphism $\iota : G_{n}^{2} \rightarrow G_{n,\mathcal{P}}^{2}$ is a monomorphism.
\end{lemma}

\begin{proof}
To prove this statement, consider projection map $pr : G_{n,\mathcal{P}}^{2} \rightarrow G_{n}^{2}$ defined by 
\begin{center}
$pr(a_{ij}^{\epsilon})  = \left\{
\begin{array}{cc} 
     a_{ij} & \text{if}~\epsilon = 0, \\
       1
 &  \text{if}~\epsilon= 1. \\
   \end{array}\right.$
   \end{center}
It is easy to show that $pr$ is well-defined function. Let $\beta$ and $\beta'$ be two words in $G_{n}^{2}$ such that $\iota(\beta)= \iota(\beta')$. By the definition of $\iota$, $pr(\iota(\beta)) =\beta$ and  $pr(\iota(\beta')) =\beta'$. That is, $\beta=\beta'$ in $G_{n}^{2}$, therefore, the proof is completed.

\end{proof}

This lemma may be reformulated in the form of the following

\begin{corollary}
If two words $\beta$ and $\beta'$ in $G_{n}^{2}$ are equivalent in $G_{n,\mathcal{P}}^{2}$, then they are equivalent in $G_{n}^{2}$.
\end{corollary}

Now, we define a group $G_{n,\mathcal{D}}^{2}$ and we call it {\em the group $G_{n}^{2}$ with points.}\index{$G_n^k$!with points}
\begin{definition}\label{dfn:Gn2withpoint}
For a positive integer $n>2$, define $G_{n,\mathcal{D}}^{2}$ by the group presentation generated by  $\{ a_{\{i,j\}} ~|~ \{i,j\} \subset \{1, \dots, n\}, i < j\}$ and $\{ \tau_{i} ~|~ i \in \{1,\dots, n\} \}$ with the following relations;
\begin{enumerate}
\item $a_{\{i,j\}}^{2} = 1$ for $\{i,j\} \subset \{1, \dots, n\}, i < j$,
\item $a_{\{i,j\}}a_{\{k,l\}} = a_{\{k,l\}}a_{\{i,j\}}$ for distinct $i,j,k,l \in \{1, \dots, n\}$,
\item $a_{\{i,j\}}a_{\{i,k\}}a_{\{j,k\}} = a_{\{j,k\}}a_{\{i,k\}}a_{\{i,j\}}$ for distinct $i,j,k \in \{1, \dots, n\}$
\item $\tau_{i}^{2} = 1$ for $i \in \{1, \dots, n\}$,
\item $\tau_{i}\tau_{j} = \tau_{j}\tau_{i}$ for $i,j\in \{1, \dots, n\}$,
\item $\tau_{i}\tau_{j}a_{ij}\tau_{j}\tau_{i} = a_{ij}$ for $i,j \in \{1, \dots, n\}$,
\item $a_{\{i,j\}}\tau_{k} = \tau_{k}a_{\{i,j\}}$ for distinct $i,j,k \in \{1, \dots, n\}$.
\end{enumerate}
We denote $a_{\{i,j\}}$ by $a_{ij}$.
\end{definition}

\begin{figure}[h!]
\begin{center}
 \includegraphics[width =8cm]{ 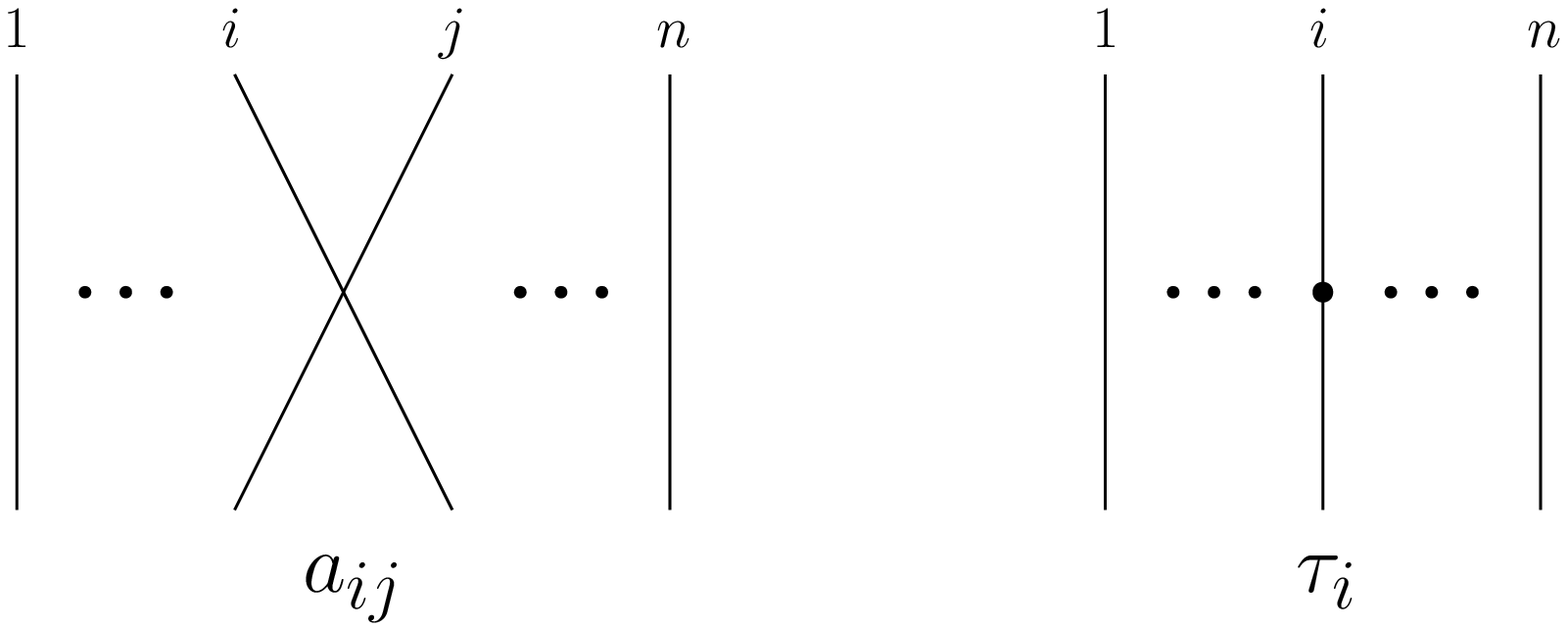}

\end{center}

 \caption{Geometrical meanings of $a_{ij}$ and $\tau_{i}$}\label{exa-cro_point}
\end{figure}

We call $\tau_{i}$ {\it a generator for a point on $i$-th component} or simply, {\it a point on  $i$-th component.} Geometrically, $a_{ij}$ is corresponded to a 4-valent vertex and $\tau_{i}$ is corresponded to a point on the $i$-th strand of the free braid, see Fig.~\ref{exa-cro_point}.

\subsection{Connection between $G_{n,\mathcal{P}}^{2}$ and $G_{n,\mathcal{D}}^{2}$ }

\begin{figure}[h!]
\begin{center}
 \includegraphics[width =10cm]{ 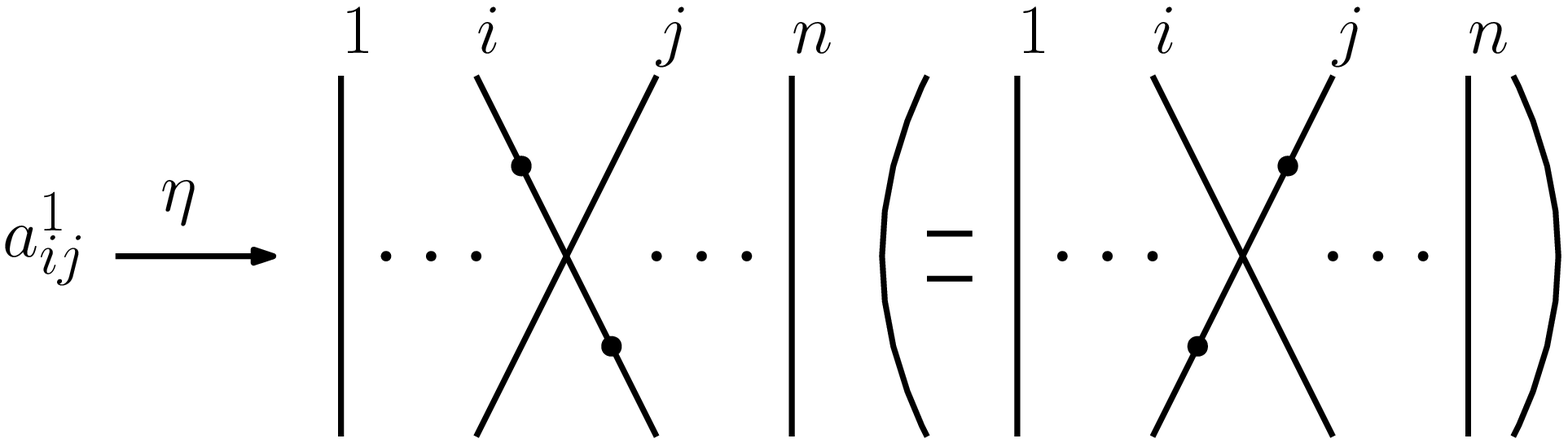}

\end{center}

 \caption{Image of $\eta(a_{ij}^{1})$}\label{imageofeta}
\end{figure}

\begin{figure}[h!]
\begin{center}
 \includegraphics[width =10cm]{ 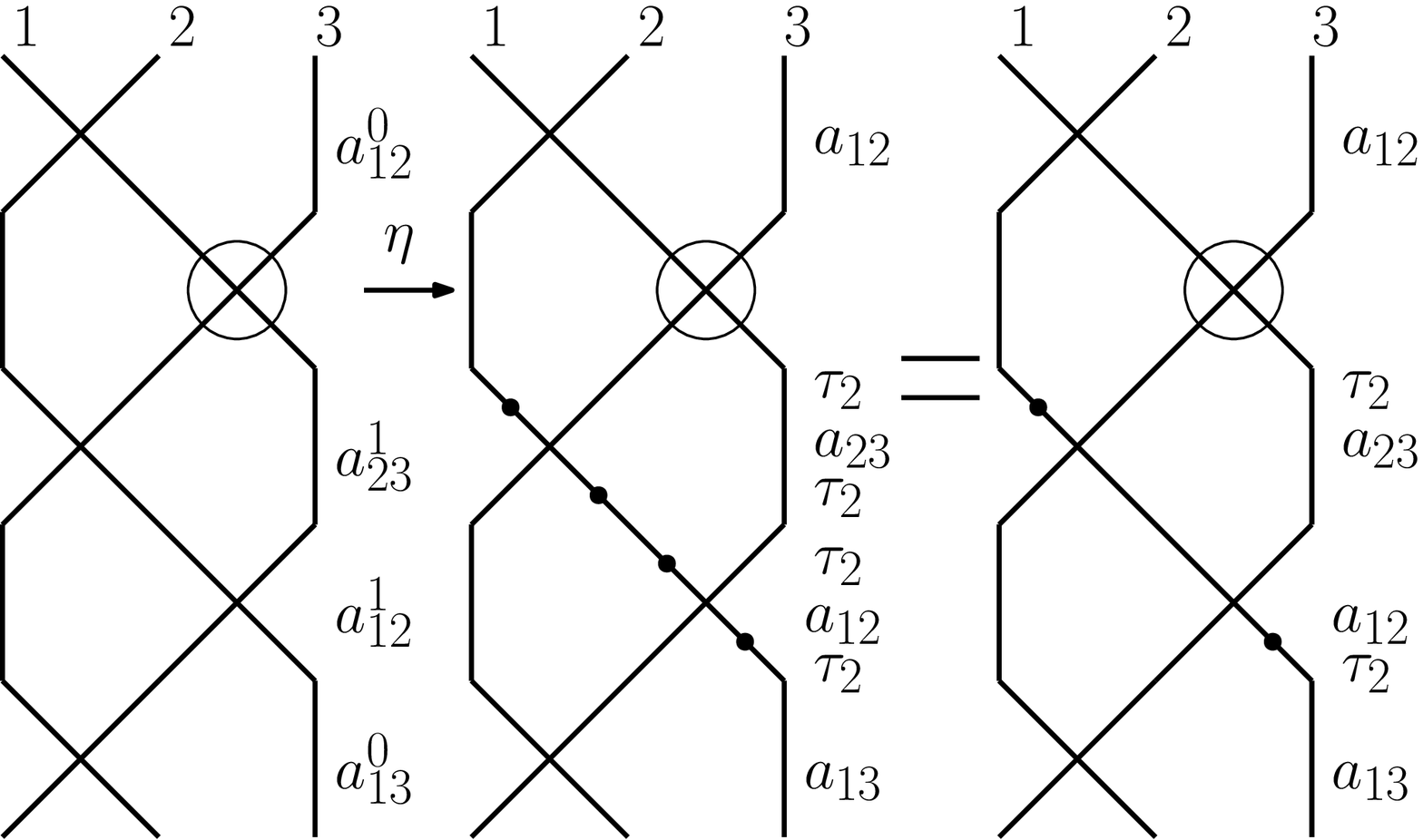}

\end{center}

 \caption{$\eta(a_{12}^{0}a_{23}^{1}a_{12}^{1}a_{13}^{0}) = a_{12}\tau_{2}a_{23}\tau_{2}\tau_{2}a_{12}\tau_{2}a_{13} = a_{12}\tau_{2}a_{23}a_{12}\tau_{2}a_{13}$}\label{exa_eta}
\end{figure}

In fact, two groups $G_{n,\mathcal{P}}^{2}$ and $G_{n,\mathcal{D}}^{2}$ are closely related to each other.
Define a homomorphism $\eta$ from $G_{n,\mathcal{P}}^{2}$ to $G_{n,\mathcal{D}}^{2}$ by
\begin{center}
$\eta(a_{ij}^{\epsilon}) = \left\{
\begin{array}{cc} 
     a_{ij} & \text{if}~\epsilon =0, \\
    \tau_{i}a_{ij}\tau_{i}
 &  \text{if}~\epsilon =1. \\
   \end{array}\right.$
\end{center} \label{map:eta}
Geometrically, the image of $a_{ij}^{1}$ is a crossing between $i$-th and $j$-th strands of the braid with two points just before and after of the crossing on $i$-th strand, see Fig.~\ref{imageofeta}. 
Note that, since we can get $\tau_{i}a_{ij}\tau_{i}=\tau_{j}a_{ij}\tau_{j}$ from the relation $\tau_{i}\tau_{j}a_{ij}\tau_{j}\tau_{i} =a_{ij}$, two diagrams with points in Fig.~\ref{imageofeta} are equivalent. On the other hand, the number of $\tau_{i}$ and $\tau_{j}$ before a given crossing $\eta(a_{ij}^{\epsilon})$ is equal to $\epsilon$ modulo $2$, because every image of $a_{kl}^{\epsilon}$ before $\eta(a_{ij}^{\epsilon})$ has two points or no points. For example, see Fig.~\ref{exa_eta}.

\begin{lemma}
The map $\eta$ from $G_{n,\mathcal{P}}^{2}$ to $G_{n,\mathcal{D}}^{2}$ is well defined.
\end{lemma}

\begin{proof}
To show that $\eta$ is well defined, it is enough to show that every relation is preserved by $\eta$. For relations $(a_{ij}^{\epsilon})^{2} = 1$, if $\epsilon =0$, then $\eta((a_{ij}^{0})^{2}) = a_{ij}^{2} = 1$. If $\epsilon =1$, then  $$\eta((a_{ij}^{1})^{2}) = \tau_{i}a_{ij}\tau_{i}\tau_{i}a_{ij}\tau_{i} = \tau_{i}a_{ij}a_{ij}\tau_{i} =  \tau_{i}\tau_{i} = 1,$$
which can be presented geometrically as Fig.~\ref{proof-lem-A2}.
\begin{figure}[h!]
\begin{center}
 \includegraphics[width =120mm]{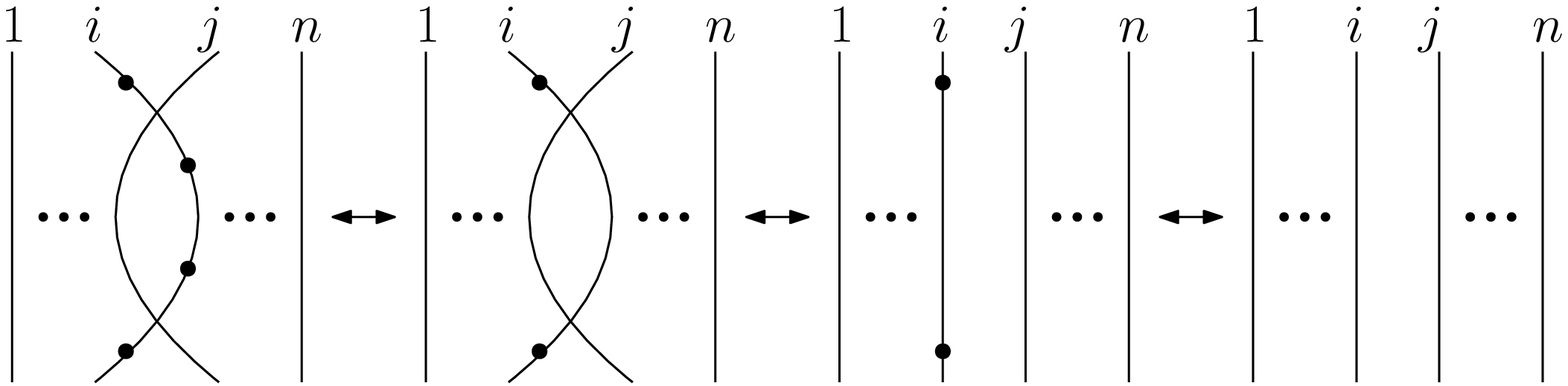}

\end{center}

 \caption{$\eta((a_{ij}^{1})^{2}) =1$}\label{proof-lem-A2}
\end{figure}
For relations $a_{ij}^{\epsilon_{ij}}a_{kl}^{\epsilon_{kl}} = a_{kl}^{\epsilon_{kl}}a_{ij}^{\epsilon_{ij}}$ with distinct $i,j,k,l$, since $i$,$j$,$k$, and $l$ are different indices, clearly the commutativity holds. For relations $a_{ij}^{\epsilon_{ij}}a_{ik}^{\epsilon_{ik}}a_{jk}^{\epsilon_{jk}} = a_{jk}^{\epsilon_{jk}}a_{ik}^{\epsilon_{ik}}a_{ij}^{\epsilon_{ij}}$, where $\epsilon_{ij}+\epsilon_{ik}+\epsilon_{jk}=0$ mod $2$, there are only two cases: all $\epsilon$'s are equal to $0$ or only two of them are equal to $1$. If every $\epsilon$ is $0$, then 
$$\eta(a_{ij}^{0}a_{ik}^{0}a_{jk}^{0}) = a_{ij}a_{ik}a_{jk} = a_{jk}a_{ik}a_{ij} = \eta( a_{jk}^{0}a_{ik}^{0}a_{ij}^{0}).$$ 
Suppose that only two of them are equal to $1$, say $\epsilon_{ij}=1$, $\epsilon_{ik}=1$ and $\epsilon_{jk}=0$. Then 

\begin{center}
$\eta(a_{ij}^{1}a_{ik}^{1}a_{jk}^{0}) = \tau_{i}a_{ij} \tau_{i} \tau_{i}a_{ik}\tau_{i}a_{jk} =  \tau_{i}a_{ij}a_{ik}\tau_{i}a_{jk} =  \tau_{i}a_{ij}a_{ik}a_{jk}\tau_{i} = \tau_{i}a_{jk}a_{ik}a_{ij}\tau_{i} = \tau_{i}a_{jk}\tau_{i} \tau_{i}a_{ik}\tau_{i}a_{ij} =\eta(a_{jk}^{0}a_{ik}^{1}a_{ij}^{1}),$
\end{center}
 which can be presented geomerically as Fig.~\ref{proof-lem-A3}.
\begin{figure}[h!]
\begin{center}
 \includegraphics[width =12cm]{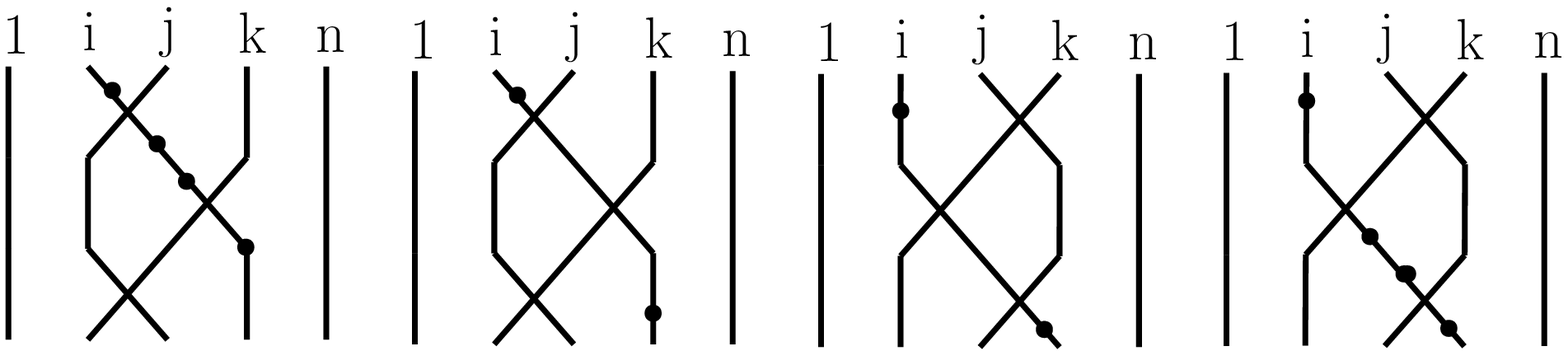}

\end{center}

 \caption{$\eta(a_{ij}^{1}a_{ik}^{1}a_{jk}^{0}) = \eta(a_{jk}^{1}a_{ik}^{1}a_{ij}^{0})$}\label{proof-lem-A3}
\end{figure}

 Therefore, $\eta$ is well-defined. 

\end{proof}

Let $H_{n,\mathcal{D}}^{2} = \{\beta \in G_{n,\mathcal{D}}^{2}~|~ N_{i}(\beta) = 0~\text{mod}~2 ~ \text{for all}~ i \in \{1,\dots, n\} \}$, where $N_{i}(\beta)$ is the number of $\tau_{i}$ in $\beta$. Note that, by the definition of $\eta$, the image of $G_{n,\mathcal{P}}^{2}$ by $\eta$ goes into the set $H_{n,\mathcal{D}}^{2}$. Since every relation of $G_{n,\mathcal{D}}^{2}$ preserves the number of $\tau$'s modulo $2$, $H_{n,\mathcal{D}}^{2}$ is a subgroup of $G_{n,\mathcal{D}}^{2}$. 

\begin{lemma}\label{lem:dot-to-parity}
Let $\beta = T_{0} a_{i_{1}j_{1}}T_{1}\cdots T_{k-1}a_{i_{k}j_{k}}T_{k} \cdots T_{m-1}a_{i_{m}j_{m}}T_{m} \in H_{n,\mathcal{D}}^{2}$, where $T_{k}$ is a product of $\tau$'s. Define a function $\chi : H_{n,\mathcal{D}}^{2} \rightarrow G_{n,\mathcal{P}}^{2}$  by 
$$ \chi(\beta) = a_{i_{1}j_{1}}^{\epsilon_{1}} \cdots a_{i_{k}j_{k}}^{\epsilon_{k}} \cdots a_{i_{m}j_{m}}^{\epsilon_{m}},$$
where $N_{i_{k}}$ is the number of $\tau_{i_{k}}$ in $T_{0}, \cdots, T_{k-1}$ and $\epsilon_{k} = N_{i_{k}} + N_{j_{k}}$. Then $\chi$ is well-defined.
\end{lemma}

\begin{proof}
Consider $\beta = T_{0} a_{i_{1}j_{1}}T_{1}\cdots T_{k-1}a_{i_{k}j_{k}}T_{k} \cdots T_{m-1}a_{i_{m}j_{m}}T_{m} \in H_{n,\mathcal{D}}^{2}$.
Assume that $a_{i_{k}j_{k}}$ is contained in one of relations, which can be applied to $\beta$. If $a_{i_{k}j_{k}}$ is contained in the relation (6) $\tau_{i}\tau_{j}a_{ij}\tau_{j}\tau_{i} = a_{ij}$ from Definition~\ref{dfn:Gn2withpoint}, then the number of  $\tau_{i_{k}}$ and $\tau_{j_{k}}$ in $T_{k-1}$ is preserved modulo $2$ and then the sum of the numbers of $\tau_{i_{k}}$ and $\tau_{j_{k}}$ in $T_{0}, \dots, T_{k-1}$ does not change modulo $2$. Therefore $\epsilon_{k} = N_{i_{k}} + N_{j_{k}}$ is even. If $a_{i_{k}j_{k}}$ is contained in the relation (7) $\tau_{k}a_{ij} = a_{ij}\tau_{k}$ from Definition~\ref{dfn:Gn2withpoint}, then the number of $\tau_{i_{k}}$ and $\tau_{j_{k}}$ in $T_{k-1}$ is preserved, since in relation (7) $i,j,k$ are different respectively. Therefore the number of $\tau_{i_{k}}$ and $\tau_{j_{k}}$ in $T_{0}, \dots, T_{k-1}$ does not change, and $\epsilon_{k} = N_{i_{k}} + N_{j_{k}}$ is preserved modulo $2$. The relation  (2) $a_{ij}a_{kl} = a_{kl}a_{ij}$ from Definition~\ref{dfn:Gn2withpoint} is preserved along $\chi$ by relation $a_{ij}^{\epsilon_{1}}a_{kl}^{\epsilon_{2}} = a_{kl}^{\epsilon_{2}}a_{ij}^{\epsilon_{1}}$ for different $i,j,k,l$. Suppose that $a_{i_{k}j_{k}}$ appears in the relation (1) $a_{ij}^{2}=1$  from Definition~\ref{dfn:Gn2withpoint}. Then $i_{k} =i_{k+1}$, $j_{k} = j_{k+1}$ and there are no $\tau$'s in $T_{k}$. Therefore $\epsilon_{k}= \epsilon_{k+1}$ and the relation is preserved by  $(a_{ij}^{\epsilon})^{2}=1$. Suppose that $a_{i_{k}j_{k}}$ appears in the relation (3) $a_{ij}a_{ik}a_{jk}=a_{jk}a_{ik}a_{ij}$ from Definition~\ref{dfn:Gn2withpoint}. Suppose that $a_{i_{k}j_{k}}a_{i_{k+1}j_{k+1}}a_{i_{k+2}j_{k+2}}=a_{i_{k+2}j_{k+2}}a_{i_{k+1}j_{k+1}}a_{i_{k}j_{k}}$ and $i_{k} = i_{k+1}$, $j_{k} = i_{k+2}$ and $j_{k+1} = j_{k+2}$. It is sufficient to show that $\epsilon_{k}+\epsilon_{k+1}+\epsilon_{k+2} = 0$ mod $2$. 
Then there are no $\tau$'s in $T_{k}$,$T_{k+1}$,$T_{k+2}$, and hence $N_{i_{k}} = N_{i_{k+1}}$, $N_{j_{k}} = N_{i_{k+2}}$, $N_{j_{k+1}} = N_{j_{k+2}}$. Therefore 
 \begin{center}
$\epsilon_{k}+\epsilon_{k+1}+\epsilon_{k+2} = N_{i_{k}}+N_{j_{k}}+N_{i_{k+1}}+N_{j_{k+1}}+N_{i_{k+2}}+N_{j_{k+2}} = 0 ~\text{mod}~2.$
\end{center}
and the proof is completed.

\end{proof}

\begin{theorem}\label{thm_isoGH}
The image of $G^2_{n,{\mathcal P}}$ under $\eta$ is isomorphic to $H_{n,\mathcal{D}}^{2}$.
\end{theorem}

\begin{proof}
We will show that $\chi \circ \eta = 1_{G_{n,\mathcal{P}}^{2}}$ and $\eta \circ \chi = 1_{H_{n,\mathcal{D}}^{2}}$. To show that $\chi \circ \eta = 1_{G_{n,\mathcal{P}}^{2}}$ let $\beta = a_{i_{1}j_{1}}^{\epsilon_{1}}  \cdots a_{i_{m}j_{m}}^{\epsilon_{m}} \in G_{n,\mathcal{P}}^{2}$. Then 
\begin{center}
$\chi(\eta(\beta)) = \chi(\tau_{i_{1}}^{\epsilon_{1}}a_{i_{1}j_{1}} \tau_{i_{1}}^{\epsilon_{1}} \cdots \tau_{i_{m}}^{\epsilon_{m}}a_{i_{m}j_{m}}\tau_{i_{m}}^{\epsilon_{m}}) = a_{i_{1}j_{1}}^{\theta_{1}}  \cdots a_{i_{m}j_{m}}^{\theta_{m}}$.
\end{center}
For each $k$, the number of all $\tau$'s which are appeared from $a_{i_{1}j_{1}}^{\epsilon_{1}}  \cdots a_{i_{k-1}j_{k-1}}^{\epsilon_{m}}$ is even. Since $\theta_{k}$ is equal to the number of all $\tau_{i_{k}}$ and $\tau_{j_{k}}$ modulo $2$, $\theta_{k} = \epsilon_{k}$. Therefore $\chi \circ \eta = 1_{G_{n,\mathcal{P}}^{2}}$.

Now we will show that $\eta \circ \chi = 1_{H_{n,\mathcal{D}}^{2}}$. Let 

$$\beta = T_{0} a_{i_{1}j_{1}}T_{1}\cdots T_{k-1}a_{i_{k}j_{k}}T_{k} \cdots T_{m-1}a_{i_{m}j_{m}}T_{m} \in H_{n,\mathcal{D}}^{2}.$$ 
Firstly we will show that $\beta$ is equivalent to an element in the form 
$$\tau_{i_{1}}^{\epsilon_{1}}a_{i_{1}j_{1}} \tau_{i_{1}}^{\epsilon_{1}} \cdots \tau_{i_{m}}^{\epsilon_{m}}a_{i_{m}j_{m}}\tau_{i_{m}}^{\epsilon_{m}}.$$ 
By the relations $\tau_{i}\tau_{j} = \tau_{j}\tau_{i}$ and $\tau_{k}a_{ij} = a_{ij}\tau_{k}$, we can assume that $\beta$ has the form $$T_{0} a_{i_{1}j_{1}}T_{1}\cdots T_{k-1}a_{i_{k}j_{k}}T_{k} \cdots T_{m-1}a_{i_{m}j_{m}}T_{m}$$
such that $T_{k}$ only has $\tau_{i_{k}}$, $\tau_{j_{k}}$, $\tau_{i_{k+1}}$ and $\tau_{j_{k+1}}$. In point of $a_{i_{k}j_{k}}$, it is a product of the following; 
$$ A_{i_{k}j_{k}} = \tau_{i_{k}}^{\theta^{f}_{i_{k}}}\tau_{j_{k}}^{\theta^{f}_{j_{k}}}a_{i_{k}j_{k}}\tau_{i_{k}}^{\theta^{b}_{i_{k}}}\tau_{j_{k}}^{\theta^{b}_{j_{k}}}.$$

Now we claim that there is an element $\beta'$ equivalent to $\beta$ in the form of
$$\beta' = A'_{i_{1}j_{1}}A'_{i_{2}j_{2}}\cdots A'_{i_{m}j_{m}},$$
where $$A'_{i_{k}j_{k}} = \tau_{i_{k}}^{\theta^{f}_{i_{k}}}\tau_{j_{k}}^{\theta^{f}_{j_{k}}}a_{i_{k}j_{k}}\tau_{i_{k}}^{\theta^{b}_{i_{k}}}\tau_{j_{k}}^{\theta^{b}_{j_{k}}},~ \theta^{f}_{i_{k}} = \theta^{b}_{i_{k}}~\text{and}~\theta^{f}_{j_{k}}=\theta^{b}_{j_{k}}.$$ 
 We call $A'_{i_{k}j_{k}}$ {\em a part of $\beta$ in the standard form} and $\beta'$ {\em a word in the standard form}.
Suppose that $\beta$ is not in the standard form. Let $k_{0} \in \{1,\dots, m\}$ be a index such that $A_{i_{k_{0}}j_{k_{0}}}$ is the first part of $\beta$ in the nonstandard form, say, $\theta_{i_{k_{0}}}^{f} \neq \theta_{i_{k_{0}}}^{b}$. Since $\beta \in H_{n,\mathcal{D}}^{2}$, the number of $\tau_{i}$ is even and there exists a part $A_{i_{l}j_{l}} = \tau_{i_{l}}^{\theta'^{f}_{i_{l}}} \tau_{j_{l}}^{\theta^{f}_{j_{l}}}a_{i_{l}j_{l}}\tau_{i_{l}}^{\theta'^{b}_{i_{l}}}\tau_{j_{l}}^{\theta^{b}_{j_{l}}}$ of $\beta$  such that $k_{0}<l$ and $i_{l} = i_{k_{0}}$(or $j_{l} = i_{k_{0}}$) of the following form:
  \begin{eqnarray*}
\beta &=& \cdots A_{i_{k_{0}}j_{k_{0}}} \cdots A_{i_{l}j_{l}} \cdots.
       \end{eqnarray*}
Let $l := min\{s \in \{1,\dots, m\}~|~i_{s} = i_{k_{0}}~ \text{or}~j_{s} = i_{k_{0}}, k_{0}<s \leq m\}.$ Assume that $i_{l} = i_{k_{0}}$.  Now we shall deform the part $A_{i_{k_{0}}j_{k_{0}}}$ in the form of $$A'_{i_{k_{0}}j_{k_{0}}} = \tau_{i_{k_{0}}}^{\theta'^{f}_{i_{k_{0}}}}\tau_{j_{k_{0}}}^{\theta^{f}_{j_{k_{0}}}}a_{i_{k_{0}}j_{k_{0}}}\tau_{i_{k_{0}}}^{\theta^{b}_{i_{k_{0}}}}\tau_{i_{k_{0}}}^{\theta'^{b}_{i_{k_{0}}}},$$ where $\theta'^{f}_{i_{k_{0}}} = \theta'^{b}_{i_{k_{0}}}$.

If $\theta_{i_{k_{0}}}^{f}=0$ and $\theta_{i_{k_{0}}}^{b}=1$, then $A_{i_{k_{0}}j_{k_{0}}}$ can be deformed to $A'_{i_{k_{0}}j_{k_{0}}}= \tau_{j_{k_{0}}}^{\theta^{f}_{j_{k_{0}}}}a_{i_{k_{0}}j_{k_{0}}}\tau_{j_{k_{0}}}^{\theta^{b}_{j_{k_{0}}}}$ as follow:
  \begin{eqnarray*}
\beta &=& \cdots A_{i_{k_{0}}j_{k_{0}}} \cdots A_{i_{l}j_{l}} \cdots\\
  &=& \cdots  \tau_{j_{k_{0}}}^{\theta^{f}_{j_{k_{0}}}}a_{i_{k_{0}}j_{k_{0}}} \tau_{i_{k_{0}}} \tau_{j_{k_{0}}}^{\theta^{b}_{j_{k_{0}}}} \cdots \tau_{i_{l}}^{\theta^{f}_{i_{l}}}\tau_{j_{l}}^{\theta^{f}_{j_{l}}}a_{i_{l}j_{l}}\tau_{i_{l}}^{\theta^{b}_{i_{l}}}\tau_{i_{l}}^{\theta^{b}_{i_{l}}} \cdots \\
  &=&  \cdots  \tau_{j_{k_{0}}}^{\theta^{f}_{j_{k_{0}}}}a_{i_{k_{0}}j_{k_{0}}}\tau_{j_{k_{0}}}^{\theta^{b}_{j_{k_{0}}}}  \tau_{i_{k_{0}}} \cdots \tau_{i_{l}}^{\theta^{f}_{i_{l}}}\tau_{j_{l}}^{\theta^{f}_{j_{l}}}a_{i_{l}j_{l}}\tau_{i_{l}}^{\theta^{b}_{i_{l}}}\tau_{i_{l}}^{\theta^{b}_{i_{l}}} \cdots \\
   &=&  \cdots  A'_{i_{k_{0}}j_{k_{0}}}   \cdots  \tau_{i_{k_{0}}}A_{i_{l}j_{l}}\cdots .
       \end{eqnarray*}             
Assume that $\theta_{i_{k_{0}}}^{f}=1$ and $\theta_{i_{k_{0}}}^{b}=0$. Note that in $\beta$ there are no $A_{i_{k_{0}}j_{k}}$ between $A_{i_{k_{0}}j_{k_{0}}}$ and $A_{i_{l}j_{l}}$ because $l := min\{s \in \{1,\dots, m\}~|~i_{s} = i_{k_{0}}~ \text{or}~j_{s} = i_{k_{0}}, k_{0}<s \leq m\}.$
If $\theta^{f}_{i_{l}} =1$,
then by applying relations $a_{ij}\tau_{k} = \tau_{k}a_{ij}$ and $\tau_{i}\tau_{j} = \tau_{j}\tau_{i}$ the generator $\tau_{i_{l}}$ can be moved near to $a_{i_{k_{0}}j_{k_{0}}}$ as follows:

\begin{eqnarray*}
 \beta 
        &=& \cdots  \tau_{j_{k_{0}}}^{\theta^{f}_{j_{k_{0}}}} \tau_{i_{k_{0}}} a_{i_{k_{0}}j_{k_{0}}}\tau_{j_{k_{0}}}^{\theta^{b}_{j_{k_{0}}}} \cdots A_{**} \cdots \tau_{i_{l}}\tau_{j_{l}}^{\theta^{f}_{j_{l}}}a_{i_{l}j_{l}}\tau_{i_{l}}^{\theta^{b}_{i_{l}}}\tau_{j_{l}}^{\theta^{b}_{j_{l}}} \cdots\\
 &=&  \tau_{j_{k_{0}}}^{\theta^{f}_{j_{k_{0}}}} \tau_{i_{k_{0}}} a_{i_{k_{0}}j_{k_{0}}}\tau_{j_{k_{0}}}^{\theta^{b}_{j_{k_{0}}}}\tau_{i_{l}} \cdots A_{**} \cdots \tau_{j_{l}}^{\theta^{f}_{j_{l}}}a_{i_{l}j_{l}}\tau_{i_{l}}^{\theta^{b}_{i_{l}}}\tau_{j_{l}}^{\theta^{b}_{j_{l}}} \cdots\\
 &=& A'_{i_{k_{0}}j_{k_{0}}} \cdots A_{**} \cdots \tau_{j_{l}}^{\theta^{f}_{j_{l}}}a_{i_{l}j_{l}}\tau_{i_{l}}^{\theta^{b}_{i_{l}}}\tau_{j_{l}}^{\theta^{b}_{j_{l}}} \cdots.
       \end{eqnarray*}
       
If $\theta^{f}_{i_{l}} =0$, then by applying the relation $\tau_{i}^{2} =1$ two of the generator $\tau_{i_{k_{0}}} = \tau_{i_{l}}$ can be added and by applying the relations $a_{ij}\tau_{k} = \tau_{k}a_{ij}$ and $\tau_{i}\tau_{j} = \tau_{j}\tau_{i}$ the generator $\tau_{i_{l}}$ can be moved near to $a_{i_{k_{0}}j_{k_{0}}}$ as follows:

\begin{eqnarray*}
 \beta 
        &=& \cdots  \tau_{j_{k_{0}}}^{\theta^{f}_{j_{k_{0}}}} \tau_{i_{k_{0}}} a_{i_{k_{0}}j_{k_{0}}}\tau_{j_{k_{0}}}^{\theta^{b}_{j_{k_{0}}}} \cdots A_{**} \cdots \tau_{j_{l}}^{\theta^{f}_{j_{l}}}a_{i_{l}j_{l}}\tau_{i_{l}}^{\theta^{b}_{i_{l}}}\tau_{j_{l}}^{\theta^{b}_{j_{l}}} \cdots\\
 &=&  \tau_{j_{k_{0}}}^{\theta^{f}_{j_{k_{0}}}} \tau_{i_{k_{0}}} a_{i_{k_{0}}j_{k_{0}}}\tau_{j_{k_{0}}}^{\theta^{b}_{j_{k_{0}}}} \cdots A_{**} \cdots \tau_{i_{l}}\tau_{i_{l}}\tau_{j_{l}}^{\theta^{f}_{j_{l}}}a_{i_{l}j_{l}}\tau_{i_{l}}^{\theta^{b}_{i_{l}}}\tau_{j_{l}}^{\theta^{b}_{j_{l}}} \cdots\\
  &=&  \tau_{j_{k_{0}}}^{\theta^{f}_{j_{k_{0}}}} \tau_{i_{k_{0}}} a_{i_{k_{0}}j_{k_{0}}}\tau_{j_{k_{0}}}^{\theta^{b}_{j_{k_{0}}}}\tau_{i_{l}} \cdots A_{**} \cdots \tau_{i_{l}}\tau_{j_{l}}^{\theta^{f}_{j_{l}}}a_{i_{l}j_{l}}\tau_{i_{l}}^{\theta^{b}_{i_{l}}}\tau_{j_{l}}^{\theta^{b}_{j_{l}}} \cdots\\
 &=& A'_{i_{k_{0}}j_{k_{0}}} \cdots A_{**} \cdots \tau_{i_{l}}\tau_{j_{l}}^{\theta^{f}_{j_{l}}}a_{i_{l}j_{l}}\tau_{i_{l}}^{\theta^{b}_{i_{l}}}\tau_{j_{l}}^{\theta^{b}_{j_{l}}} \cdots.
       \end{eqnarray*}
Since the parts $A_{i_{k}j_{k}}$, $k< k_{0}$, of $\beta$ are not changed, we can deform $\beta$ to the word $\beta'$ in the standard form inductively.
  
 Now assume that $\beta = \prod_{k=1}^{m} = \tau_{i_{k}}^{\theta_{i_{k}}}\tau_{j_{k}}^{\theta_{j_{k}}}a_{i_{k}j_{k}}\tau_{i_{k}}^{\theta_{i_{k}}}\tau_{j_{k}}^{\theta_{j_{k}}}$ is in the standard form. Then 
\begin{eqnarray*}
\eta \circ \chi(\beta) &=& \eta \circ \chi( \prod_{k=1}^{m} = \tau_{i_{k}}^{\theta_{i_{k}}}\tau_{j_{k}}^{\theta_{j_{k}}}a_{i_{k}j_{k}}\tau_{i_{k}}^{\theta_{i_{k}}}\tau_{j_{k}}^{\theta_{j_{k}}})\\
 &=& \eta(\prod_{k=1}^{m} a_{i_{k}j_{k}}^{\theta_{i_{k}j_{k}}}) = \prod_{k=1}^{m}(\tau_{i_{k}}^{\theta_{i_{k}j_{k}}}a_{i_{k}j_{k}}\tau_{i_{k}}^{\theta_{i_{k}j_{k}}}),
\end{eqnarray*}  
where $\theta_{i_{k}j_{k}} = \theta_{i_{k}}+\theta_{j_{k}}$ mod $2$. If $\theta_{i_{k}}= \theta_{j_{k}}=1$, by the relation $\tau_{i}\tau_{j}a_{ij}\tau_{j}\tau_{i} =a_{ij}$, $\beta$ and $\eta \circ \chi(\beta)$ are same elements in $G_{n,\mathcal{D}}^{2}$. If $\theta_{i_{k}}=1$ and $ \theta_{j_{k}}=0$ (or $\theta_{i_{k}}=0$ and $ \theta_{j_{k}}=1$) by the relation $\tau_{i}\tau_{j}a_{ij}\tau_{j}\tau_{i} =a_{ij}$ (in other words, $\tau_{i}a_{ij}\tau_{i}=\tau_{j}a_{ij}\tau_{j}$), $\beta$ and $\eta \circ \chi(\beta)$ are same elements in $G_{n,\mathcal{D}}^{2}$ and hence $\eta \circ \chi = 1_{G_{n,\mathcal{D}}^{2}}$.
\end{proof}

\subsection{Connection between $G_{n,\mathcal{D}}^{2}$ and $G_{n+1}^{2}$}

We recall (Fig.~\ref{(n,k)_(n-1,k)_(n-1,k-1)} and Section~\ref{sec:gnk}) that maps from $G_{n}^{k}$ to $G_{n-1}^{k}$ and from $G_{n}^{k}$ to $G_{n-1}^{k-1}$ are constructed by deleting one subindex. 
Here we justify these maps and extended to the case of mappings from $G_{n}^{2}$ to $G_{n-1,\mathcal{D}}^{2}$. We obtained the group monomorphisms $\iota \colon G_{n}^{2} \rightarrow G_{n,\mathcal{P}}^{2}$ and $\eta \colon G_{n,\mathcal{P}}^{2} \rightarrow G_{n,\mathcal{D}}^{2}$. To show that they are group monomorphisms, we found the right inverses for them. In this section we consider relations between $G_{n,\mathcal{D}}^{2}$ and $G_{n+1}^{2}$.
 
Let us define a homomorphism $\omega$ from $G_{n+1}^{2}$ to $G_{n,\mathcal{D}}^{2}$ by 
\begin{center}
$\omega(a_{ij}) = \left\{
\begin{array}{cc} 
     a_{ij} & \text{if}~n+1 \not\in \{i,j\}, \\
    \tau_{i}
 &  \text{if}~ j= n+1, \\
 \tau_{j} &  \text{if}~ i= n+1. \\
   \end{array}\right.$
\end{center}

To this end, we prove the following:
\begin{lemma}\label{func(n+1)topoint}
The mapping $\omega$ from $G_{n+1}^{2}$ to $G_{n,\mathcal{D}}^{2}$ is well-defined.
\end{lemma}

\begin{proof}
It is enough to show that every relation for $G_{n+1}^{2}$ is preserved by $\omega$. If every index is different with $n+1$, it is clear.
Consider the relation $a_{ij}a_{kl}=a_{kl}a_{ij}$ for distinct $i,j,k,l \in \{1,\dots, n+1\}$. If one of $i,j,k$ and $l$ is $n+1$, say $i=n+1$, then 
$$\omega(a_{(n+1)j}a_{kl}) = \tau_{j}a_{kl}= a_{kl}\tau_{j} =\omega(a_{kl}a_{(n+1)j}).$$

 For the relations $a_{ij}a_{ik}a_{jk} = a_{jk}a_{ik}a_{ij}$, if one of $i,j$ and $k$  is $n+1$, say $i=n+1$, then 
 $$\omega( a_{ij}a_{ik}a_{jk}) = \tau_{j}\tau_{k}a_{jk} = a_{jk}\tau_{k}\tau_{j} = \kappa(a_{jk}a_{ik}a_{ij}).$$
Clearly, the relation $a_{ij}^{2}=1$ is preserved along $\omega$ and the statement is proved. 
\end{proof}

\begin{figure}[h!]
\begin{center}
 \includegraphics[width =8cm]{ 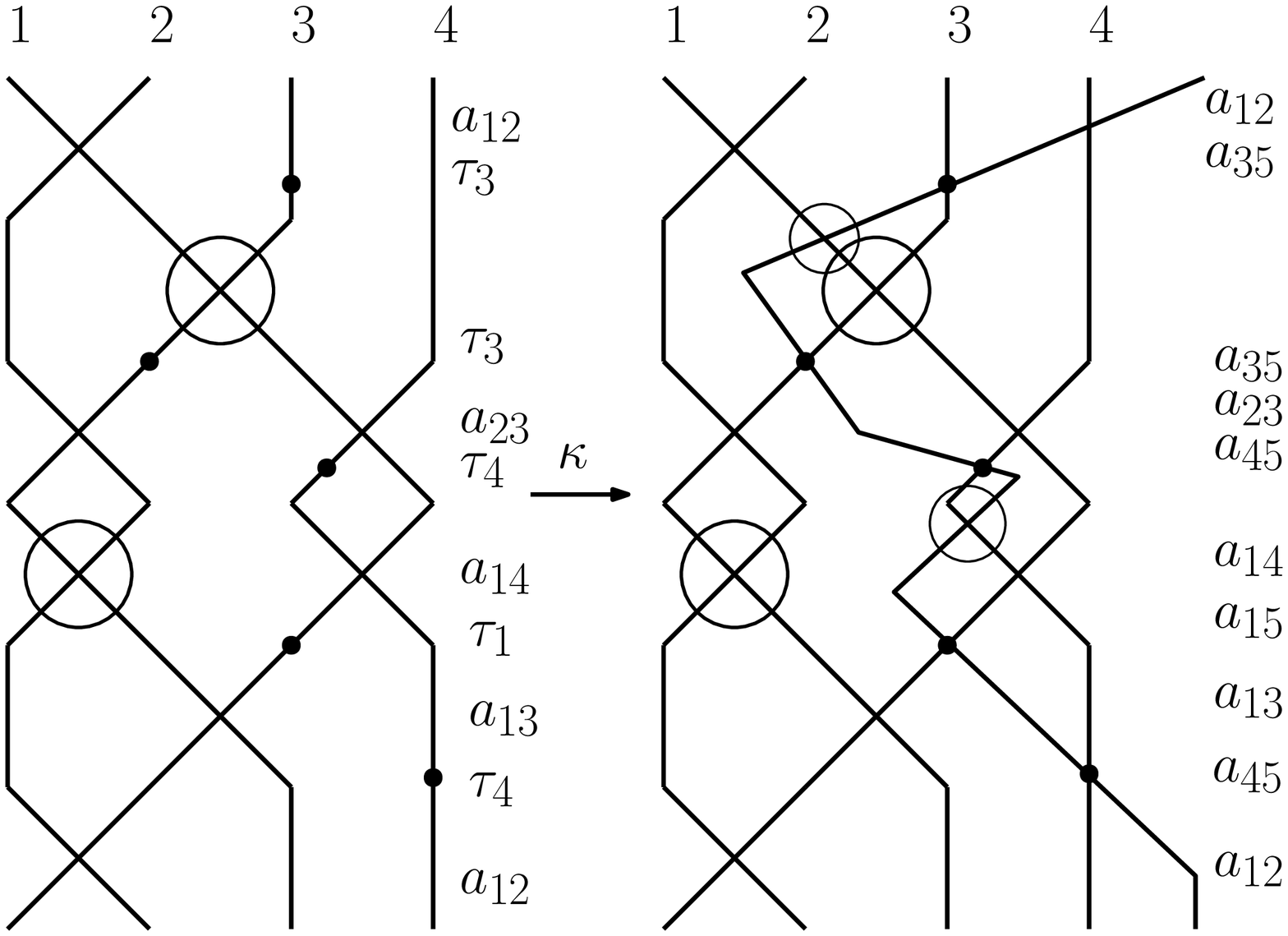}

\end{center}

 \caption{Example of an image of the mapping $\kappa$}\label{exa_kappa}
\end{figure}

Now, let us define the inverse mapping. We define $\kappa$ from $G_{n,\mathcal{D}}^{2}$ to $G_{n+1}^{2}$ by $\kappa(a_{ij}) = a_{ij}$ and $\kappa(\tau_{i}) = a_{i(n+1)}$, for example, see Fig.~\ref{exa_kappa}. \label{map:kappa}

 But this mapping is not well-defined because the relation $\tau_{i}\tau_{j} = \tau_{j}\tau_{i}$ in $G_{n,\mathcal{D}}^{2}$ goes to $a_{i(n+1)}a_{j(n+1)}= a_{j(n+1)}a_{i(n+1)}$ but $a_{i(n+1)}$ and $a_{j(n+1)}$ do not commute. However, we can modify the mapping and get the following lemma.
\begin{lemma}
The mapping $\kappa \colon G_{n,\mathcal{D}}^{2} \rightarrow G_{n+1}^{2}/ \langle a_{i(n+1)}a_{j(n+1)}= a_{j(n+1)}a_{i(n+1)}\rangle$, defined by $\kappa(a_{ij}) = a_{ij}$ and $\kappa(\tau_{i}) = a_{i(n+1)}$, is well-defined.
\end{lemma}

\begin{proof}
It is sufficient to show that every relation for $G_{n,\mathcal{D}}^{2}$ goes to an identity in $G_{n+1}^{2}/ \langle a_{i(n+1)}a_{j(n+1)}= a_{j(n+1)}a_{i(n+1)}\rangle$. For relations (1), (2), (3) and (4) from Definition~\ref{dfn:Gn2withpoint}, it is clear. Consider the relation $\tau_{i}\tau_{j} = \tau_{j}\tau_{i}$. Then
$$\kappa(\tau_{i}\tau_{j}) =  a_{i(n+1)}a_{j(n+1)}= a_{j(n+1)}a_{i(n+1)} = \kappa(\tau_{j}\tau_{i}).$$ 
The relation (6) $\tau_{i}\tau_{j}a_{ij}\tau_{j}\tau_{i} = a_{ij}$ from Definition~\ref{dfn:Gn2withpoint} can be rewritten as $\tau_{i}\tau_{j}a_{ij} = a_{ij}\tau_{j}\tau_{i}$. Then 
$$\kappa(\tau_{i}\tau_{j}a_{ij}) = a_{i(n+1)}a_{j(n+1)}a_{ij}=a_{ij}a_{j(n+1)}a_{i(n+1)} = \kappa(a_{ij}\tau_{j}\tau_{i}).$$
 Finally, the relations (5) and (7) from Definition~\ref{dfn:Gn2withpoint} are preserved by the relation $a_{i(n+1)}a_{j(n+1)}= a_{j(n+1)}a_{i(n+1)}$ and the proof is completed. 
\end{proof}

\begin{figure}[h!]
\begin{minipage}{.4\textwidth}
\begin{center}
 \includegraphics[width =.85\textwidth]{ 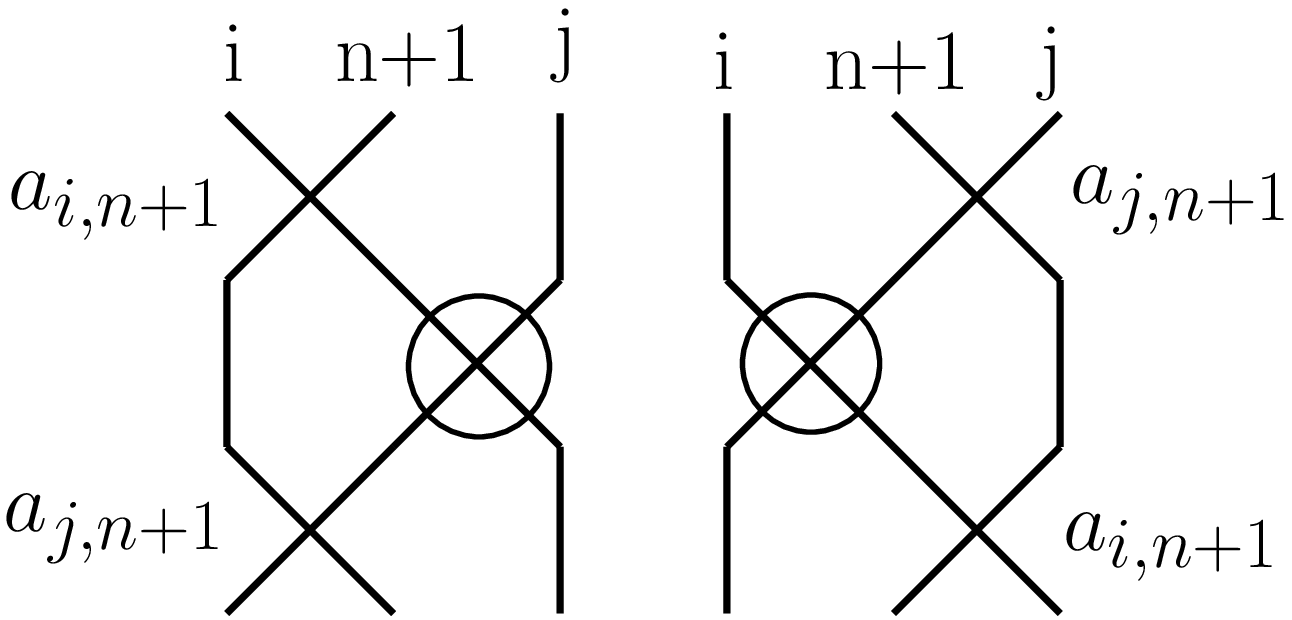}
\end{center}
\vspace{-0.4cm}
 \caption{Forbidden move}\label{forbidden_comp}
 \end{minipage}
 \hspace{3mm}
 \begin{minipage}{.6\textwidth}
 \begin{center}
 \includegraphics[width = 1\textwidth]{ 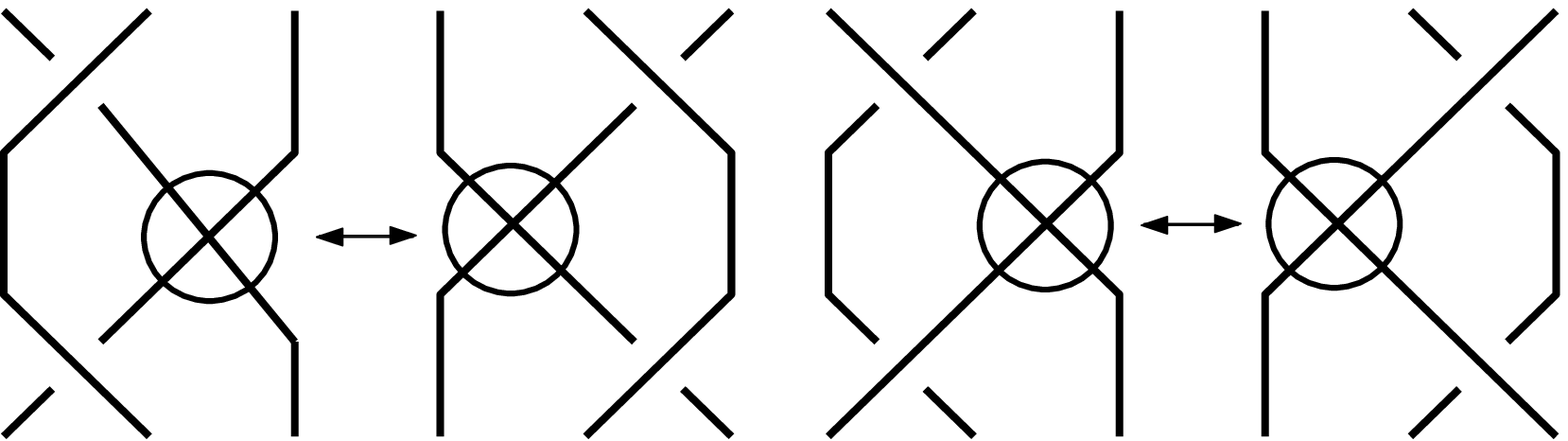}

\end{center}

 \caption{Forbidden moves for virtual links}\label{forbidden_virt}
 \end{minipage}
\end{figure}

Geometrically the relations $a_{i(n+1)}a_{j(n+1)} = a_{j(n+1)}a_{i(n+1)}$ mean that two consecutive classical crossings pass a virtual crossing, see Fig.~\ref{forbidden_comp}. Generally, this move in this case, is called {\it a forbidden move by $(n+1)$-th strand.}

In the virtual knot theory {\it forbidden moves}\index{Forbidden moves} consists of two moves in Fig.~\ref{forbidden_virt}. 

The following proposition is well known (see \cite{NelsonUnknotting, Kanenobu}):
\begin{proposition}
Every virtual knot is equivalent to the trivial knot under the virtual Reidemeister moves and the forbidden moves.
\end{proposition}
It follows that the forbidden moves for free knots unknot free knots. But in the case of (free or virtual) links, it is not true.

\begin{corollary}\label{Gnd2_Gn+12}
The groups $G_{n,\mathcal{D}}^{2}$ and $G_{n+1}^{2}/ \langle a_{i(n+1)}a_{j(n+1)}= a_{j(n+1)}a_{i(n+1)}\rangle$ are isomorphic.
\end{corollary}

\begin{proof}
We will show that $\omega \circ \kappa = 1_{G_{n,\mathcal{D}}^{2}}$ and $\kappa \circ \omega = 1_{G_{n+1}^{2}/ \langle a_{i(n+1)}a_{j(n+1)}= a_{j(n+1)}a_{i(n+1)}\rangle}$. For generators $a_{ij}$ in $G_{n+1}^{2}/ \langle a_{i(n+1)}a_{j(n+1)}= a_{j(n+1)}a_{i(n+1)}\rangle$, 
\begin{center}
$\omega(a_{ij}) =  \left\{
\begin{array}{cc} 
     a_{ij} & \text{if}~n+1 \in \{i,j\}, \\
    \tau_{i}
 &  \text{if}~ j= n+1, \\
 \tau_{j} &  \text{if}~ i= n+1. \\
   \end{array}\right.$
\end{center}
By definition of $\kappa$, $\omega(\kappa(a_{ij}))= a_{ij}$. Hence  $\omega \circ \kappa = 1_{G_{n,\mathcal{D}}^{2}}$. Analogously, it can be shown that $\kappa \circ \omega = 1_{G_{n+1}^{2}/ \langle a_{i(n+1)}a_{j(n+1)}= a_{j(n+1)}a_{i(n+1)}\rangle}$.
\end{proof}

Now define a mapping $\omega_{m} : G_{n+1}^{2} \rightarrow G_{n,\mathcal{D}}^{2}$ \label{map:omega} by
\begin{center}
$\omega_{m}(a_{ij})  = \left\{
\begin{array}{cc} 
    \tau_{i} & \text{if}~ j= n+1,~ i<m, \\
     \tau_{i-1} & \text{if}~ j= n+1, ~i>m. \\
      a_{ij} & \text{if}~ i,j<m, ~i,j \neq n+1 \\
       a_{i(j-1)} & \text{if}~ i<m, j>m, ~i,j \neq n+1\\
       a_{(i-1)(j-1)} & \text{if}~ i, j>m,~ i,j \neq n+1\\
   \end{array}\right.$
   \end{center} 
Analogously $\omega_{m}$ is well-defined (see Lemma~\ref{func(n+1)topoint}). Then, for each braid $\beta$ on $(n+1)$ strands such that the number of crossings between $i$-th and $m$-th strands is even for each index $i \neq m$, that is, $\beta \in \omega_{m}^{-1}(H_{n,\mathcal{D}}^{2})$, and a braid with parity can be obtained.
\begin{example}
Let  $\beta = a_{12}a_{23}a_{13}a_{23}a_{13}a_{23}a_{12}a_{23}$ in $G_{3}^{2}$. Note that $\omega_{m}(\beta) \in H_{2,\mathcal{D}}^{2}$ for each index $m$. Then
\begin{itemize}
\item $\chi(\omega_{1}(\beta))=\chi(\tau_{1}a_{12}\tau_{2}a_{12}\tau_{2}a_{12}\tau_{1}a_{12}) = a_{12}^{1}a_{12}^{0}a_{12}^{1}a_{12}^{0} \neq 1,$
\item $\chi(\omega_{2}(\beta)) =\chi(\tau_{1}\tau_{2}a_{12}\tau_{2}a_{12}\tau_{2}\tau_{1}\tau_{2})=\chi(a_{12}\tau_{1}a_{12}\tau_{1})= a_{12}^{0}a_{12}^{1}  \neq 1 ,$
\item $\chi(\omega_{3}(\beta)) =\chi(a_{12}\tau_{2}\tau_{1}\tau_{2}\tau_{1}\tau_{2}a_{12}\tau_{2})= \chi(a_{12}\tau_{2}a_{12}\tau_{2})= a_{12}^{0}a_{12}^{1}  \neq 1 .$
\end{itemize}
Since $\psi_{m}$ is well defined, $\beta$ is not trivial in $G_{3}^{2}$.
\end{example}

To avoid confusion let us gather mappings, which are considered in the present section, without fomulations:

\begin{itemize}
\item $p_{m} : PB_{n+1} \rightarrow PB_{n}$ in page \pageref{map:pm}.
\item $\phi_{n} : PB_{n} \rightarrow H_{n}^{3} \subset G_{n}^{3}$ in Proposition~\ref{prop_hom_PBn_Gn3} in page \pageref{prop_hom_PBn_Gn3}.
\item $q_{m} : G_{n+1}^{3} \rightarrow G_{n}^{3}$ in page \pageref{map:qm}.
\item $\phi_{(i,j,k)} : G_{n}^{3} \rightarrow \mathbb{Z}_{2}^{2^{2(n-3)}}$ in page \pageref{map:phi-ijk}.
\item $\iota : G_{n}^{2} \rightarrow G_{n,\mathcal{P}}^{2}$ in page \pageref{map:iota}.
\item $pr : G_{n,\mathcal{P}}^{2} \rightarrow G_{n}^{2}$ in the proof of Lemma \ref{emb_to_parity} in page \pageref{emb_to_parity}.
\item $\eta : G_{n,\mathcal{P}}^{2} \rightarrow G_{n,\mathcal{D}}^{2}$ in page \pageref{map:eta}.
\item $\chi : H_{n,\mathcal{D}}^{2}  \rightarrow G_{n,\mathcal{P}}^{2}$ in Lemma \ref{lem:dot-to-parity} in page \pageref{lem:dot-to-parity}.
\item $\kappa : G_{n,\mathcal{D}}^{2} \rightarrow G_{n+1}^{2}/ \langle a_{i(n+1)}a_{j(n+1)}= a_{j(n+1)}a_{i(n+1)}\rangle$ in page \pageref{map:kappa}.
\item $\omega_{m}\colon G_{n+1}^{2} \rightarrow G_{n, \mathcal{D}}^{2}$ in page \pageref{map:omega}.
\end{itemize}

\subsection[Parity for $G_{n}^{2}$ and invariant of pure braids]{Parity for $G_{n}^{2}$ and invariant of pure braids valued in a free product of $\mathbb{Z}_{2}$}

Let $\beta \in G_{n,\mathcal{P}}^{2}$. For a fixed pair $i,j \in \{1,\dots, n\}$ and for $k \in \{1,\dots, n\} \backslash \{i,j\}$, define $i_{a_{ij}^{\epsilon}}^{\mathcal{P}}(k)$ for each $a_{ij}^{\epsilon}$ in $\beta$ by 

\begin{wrapfigure}{R}{0.4\textwidth}
\begin{center}
 \includegraphics[width = 0.35\textwidth]{ 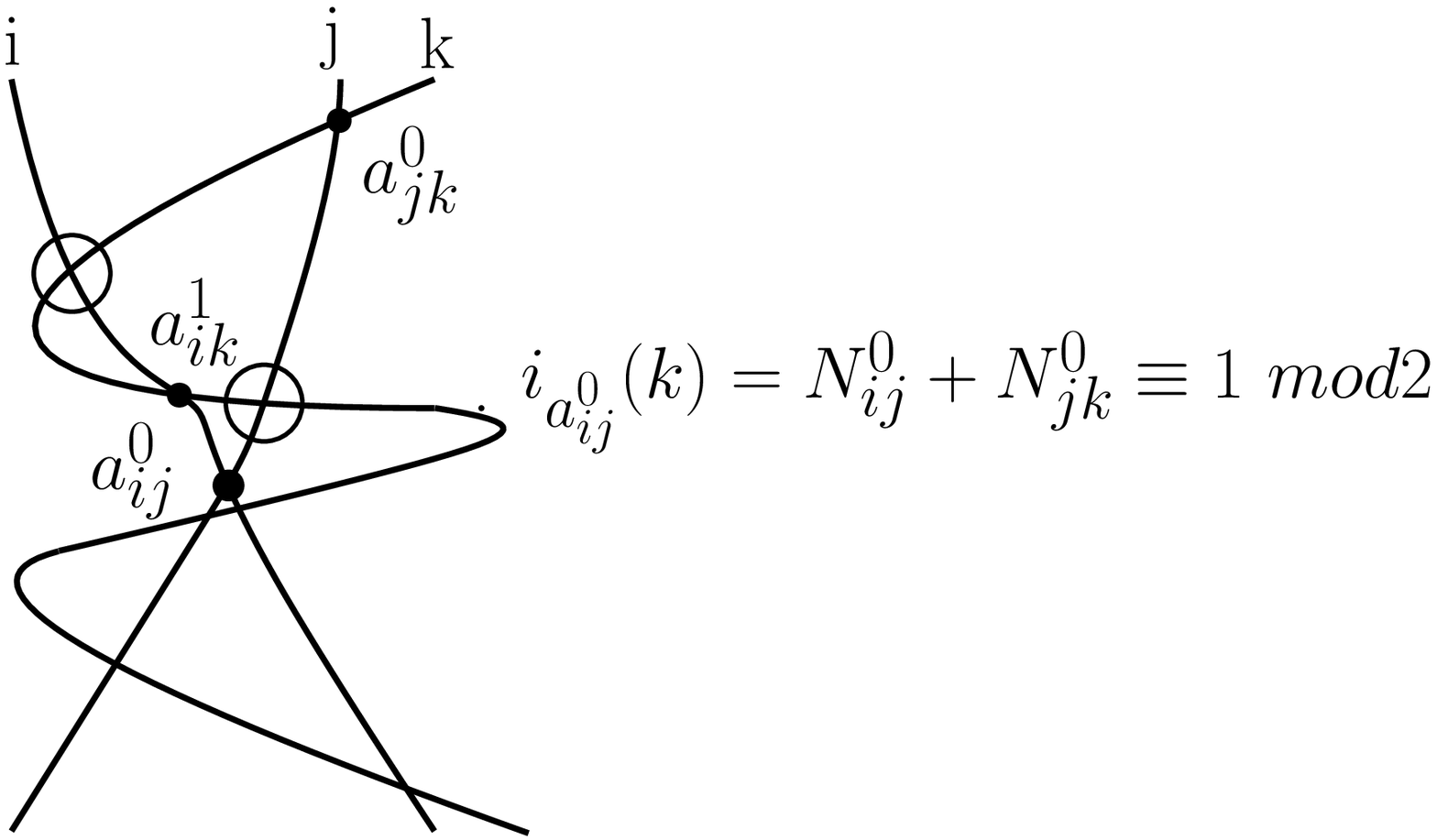}
\end{center}
\vspace{-0.5cm}
 \caption{The value of $i^{\mathcal{P}}_{a_{ij}}(k)$}\label{exa_ic_parity}
\end{wrapfigure}

\begin{center}
$ i_{a_{ij}^{\epsilon}}^{\mathcal{P}}(k) = \left \{
\begin{array}{cc} 
    N_{ik}^{0}+N_{jk}^{0} ~ mod ~2 & \text{if}~\epsilon =0, \\
    N_{ik}^{0}+N_{jk}^{1} ~ mod~ 2 &  \text{if}~\epsilon =1, \\
   \end{array}\right. $
\end{center}
where $N_{ik}^{\epsilon}$ is the number of $a_{ik}^{\epsilon}$, which appears before $a_{ij}^{\epsilon}$, for example, see Fig.~\ref{exa_ic_parity}.

Let $\{c_{1}, \cdots, c_{m}\}$ be the ordered set of $a_{ij}^{\epsilon_{m}}$'s in $\beta$ such that the order agrees with the order of position of $a_{ij}^{\epsilon}$'s. Define $w^{\mathcal{P}}_{\{i,j\}} \colon G_{n,\mathcal{P}}^{2} \rightarrow F_{n}^{2}$ by $w^{\mathcal{P}}_{\{i,j\}}(\beta) =\prod_{s=1}^{m} i_{c_{s}}^{\mathcal{P}} $. Here, the superscript $p$ means ``parity''.

\begin{lemma}\label{invGnp2}
$w^{\mathcal{P}}_{\{i,j\}} $ is well defined.
\end{lemma}

\begin{proof}
It suffices to show that the image of $w^{\mathcal{P}}_{\{i,j\}}$ does not change when relations of $G_{n,\mathcal{P}}^{2}$ are applied to $\beta$. In the cases of $(a_{ij}^{\epsilon})^{2}=1$ and $a_{ij}^{\epsilon_{1}}a_{kl}^{\epsilon_{2}}= a_{kl}^{\epsilon_{2}}a_{ij}^{\epsilon_{1}}$, it is easy. For relations $a_{ij}^{\epsilon_{ij}}a_{ik}^{\epsilon_{ik}}a_{jk}^{\epsilon_{jk}} = a_{jk}^{\epsilon_{jk}}a_{ik}^{\epsilon_{ik}}a_{ij}^{\epsilon_{ij}}$, where $\epsilon_{ij}+\epsilon_{ik}+\epsilon_{jk} \equiv 0$ mod $2$, suppose that the relation is not applied on $c_{s}$. Then the number of $a_{ik}^{\epsilon}$ and $a_{jk}^{\epsilon}$ before $c_{s}$ remains the same and then $w^{\mathcal{P}}_{\{i,j\}}(\beta)$ does not change. Suppose that $c_{s}$ is in the applied relation, say $c_{s}^{\epsilon_{1}}a_{il}^{\epsilon_{2}}a_{jl}^{\epsilon_{3}} = a_{jl}^{\epsilon_{3}}a_{il}^{\epsilon_{2}}c_{s}^{\epsilon_{1}}$, where ${\epsilon_{1}}+{\epsilon_{2}}+{\epsilon_{3}} \equiv 0$ mod $2$. If $l \neq k$, then the number of $a_{ik}^{\epsilon}$ and $a_{jk}^{\epsilon}$ before $c_{s}$ is not changed. Suppose that $l=k$. If $\epsilon_{1}=0$, then $\epsilon_{2} = \epsilon_{3} =0 $ or  $\epsilon_{2} = \epsilon_{3} =1$. Then the sum of the number of $a_{ik}^{0}$ and $a_{jk}^{0}$ remains the same modulo $2$. If $\epsilon_{1}=1$, then  $\epsilon_{2} = 1,\epsilon_{3} =0 $ or  $\epsilon_{2} = 0$,$\epsilon_{3} =1$. Then $i_{k}^{\mathcal{P}}(c_{s}) = N_{ik}^{0}+N_{jk}^{1} $ is not changed modulo $2$ and $w^{\mathcal{P}}_{\{i,j\}}(\beta)$ is not changed and the proof is completed.
\end{proof}

\begin{figure}
\begin{center}
 \includegraphics[height=120pt]{ 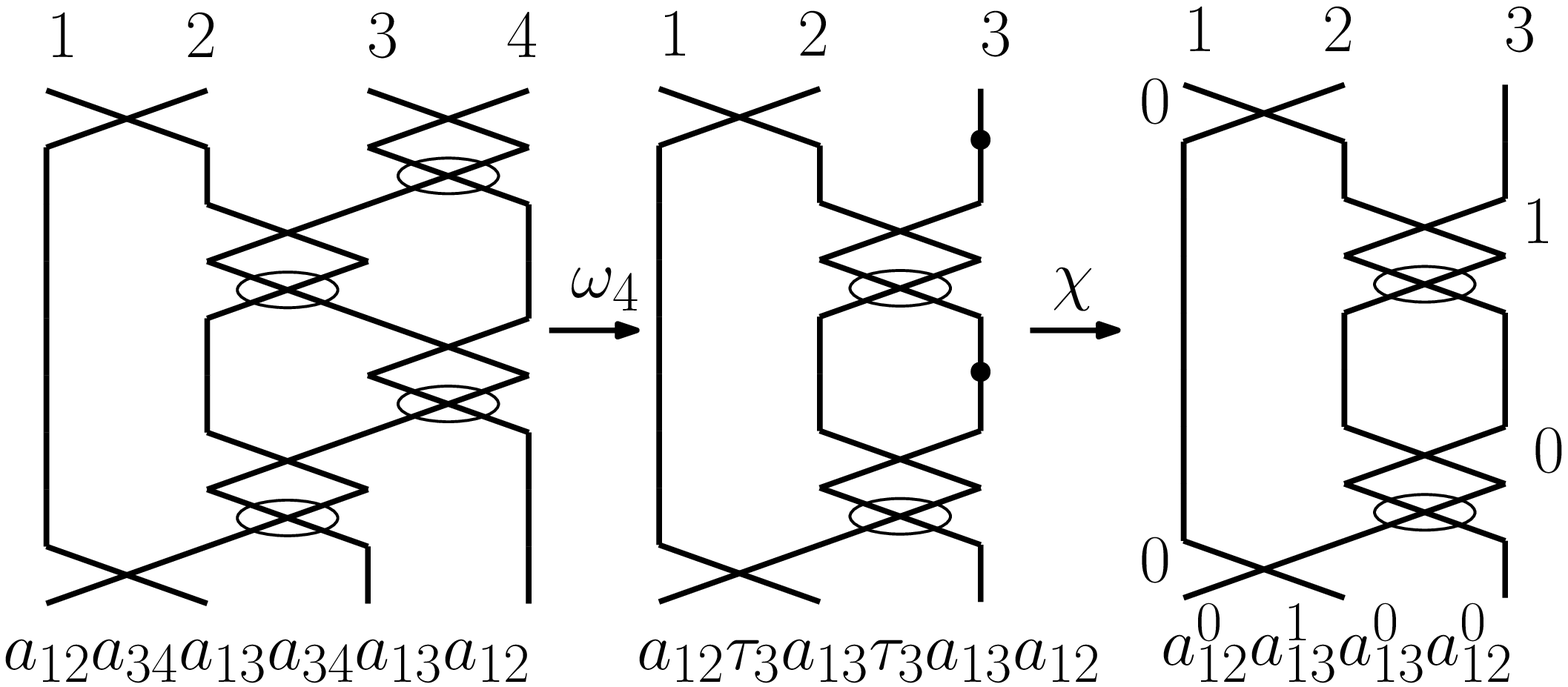}
\end{center}
 \caption{$\beta = a_{12}a_{34}a_{13}a_{34}a_{13}a_{12} \in G_{4}^{2}$ and  $\phi_{\{1,2\}}^{4} = 01 \in F_{3}^{2}$ }\label{exa_parity_MN}
\end{figure}

Let us consider the mapping $\phi_{\{i,j\}}^{l} \colon G_{n+1}^{2} \rightarrow F_{n}^{2}$ defined by $\phi_{\{i,j\}} ^{l} = w_{\{i,j\}} ^{\mathcal{P}} \circ \chi \circ \omega_{l}$.
For example, for $\beta = a_{12}a_{34}a_{13}a_{34}a_{13}a_{12}$, we obtain $\chi \circ \omega_{4}(\beta) = a_{12}^{0}a_{13}^{1}a_{13}^{0}a_{12}^{0}$. Then $\phi_{\{1,2\}}^{4}(\beta) = w_{\{1,2\}}^{\mathcal{P}} \circ \chi \circ \omega_{4}(\beta) = 01 \in F_{3}^{2}$, see Fig.~\ref{exa_parity_MN}.

\begin{corollary}
$\phi_{\{i,j\}}^{l}$ is an invariant for $\beta \in G_{n+1}^{2}$. 
\end{corollary}

\begin{proof}
Since $w_{l}$ and $\chi$ are homomorphisms, by Lemma~\ref{invGnp2}, $\phi_{\{i,j\}}^{l}$ is an invariant for $\beta \in G_{n+1}^{2}$.
\end{proof}

\begin{example}
Let $X =  a_{12}a_{13}a_{12}a_{13}$ and $Y = a_{23}a_{35}a_{23}a_{35}$ and 
$\beta = [X,Y]$ in $G_{5}^{2}$. Let us show now that $\beta$ is not trivial. To this end, we consider the element $w_{\{1,2\}}^{5}$ where the parity is obtained from $5$-th strand. For a pair $(1,2)$, the value $\phi(\beta)$ is trivial, where $\phi \colon H_{n}^{2} \to \mathbb{Z}^{2^{n-2}}$ is the mapping defined in Section~\ref{sect:hom_free_group}, because $Y$ is an even word and $Y$ contains no $a_{12}$. But $\phi_{\{1,2\}}^{5}(\beta)$ is not trivial. Now we calculate it. Firstly, we obtain 
$$\chi \circ \omega_{5}(\beta) = a_{12}^{0}a_{13}^{0}a_{12}^{0}a_{13}^{0}a_{23}^{0}a_{23}^{1}a_{13}^{0}a_{12}^{0}a_{13}^{0}a_{12}^{0}a_{23}^{1}a_{23}^{0}.$$ Let $c_{1} = a_{12}^{0}$, $c_{2} = a_{12}^{0}$,  $c_{3} = a_{12}^{0}$,  $c_{4} = a_{12}^{0}$ such that the order of $c_{i}$ agrees with the order of $a_{12}^{\epsilon}$ in $\beta$. Then $i_{c_{1}}^{\mathcal{P}}(3) = 0$, $i_{c_{2}}^{\mathcal{P}}(3) = 1$, $i_{c_{3}}^{\mathcal{P}}(3) = 0$, $i_{c_{4}}^{\mathcal{P}}(3) = 1$ and $i_{c_{s}}(4)= 0$. Then $\phi_{\{1,2\}}^{5}(\beta) = \zeta_{(0,0)}\zeta_{(1,0)}\zeta_{(0,0)}\zeta_{(1,0)}$ and it cannot be canceled in $F_{4}^{2}$, where $\zeta_{a,b}$ is defined by $\zeta_{a,b}(3)=a$ and $\zeta_{a,b}(4)=b$. 
\end{example}

This invariant can be used for $\beta \in G_{n+1}^{3}$ by the homomorphism $r_{m} : G_{n+1}^{3} \rightarrow G_{n}^{2}$ defined by

$$ r_{m}(a_{ijk}) = \left \{
\begin{array}{cc} 
   a_{ij} & \text{if}~k=m,i,j<m, \\
   a_{i(j-1)}  & \text{if}~k=m,i<m,j>m, \\
   a_{(i-1)(j-1)}  & \text{if}~k=m,i>m,j>m, \\
    1 &  \text{if}~i,j,k \neq m. \\
   \end{array} \right. $$

\begin{example}
Let 
\begin{eqnarray*}
\beta  &=&a_{124}a_{123}a_{135}a_{134}a_{124}a_{134}a_{135}a_{123}a_{134}a_{135}a_{134}a_{123}a_{135}a_{134}a_{124}a_{134}a_{135}a_{123}\\
&&a_{124}a_{134}a_{135}a_{134},
\end{eqnarray*}

in $G_{5}^{3}$. Then
$$\beta_{1} = r_{1}(\beta) =a_{24}a_{23}a_{35}a_{34}a_{24}a_{34}a_{35}a_{23}a_{34}a_{35}a_{34}a_{23}a_{35}a_{34}a_{24}a_{34}a_{35}a_{23}a_{24}a_{34}a_{35}a_{34}.$$
For index $5$, \\
$$\chi \circ \omega_{5}(\beta_{1}) 
 = a_{24}^{0}a_{23}^{0}a_{34}^{1}a_{24}^{0}a_{34}^{1}a_{23}^{0}a_{34}^{0}a_{34}^{1}a_{23}^{1}a_{34}^{0}a_{24}^{0}a_{34}^{0}a_{23}^{1}a_{24}^{0}a_{34}^{1}a_{34}^{0}.$$
Let $c_{1} = a_{24}^{0}$, $c_{2} = a_{24}^{0}$, $c_{3} = a_{24}^{0}$ and $c_{4} = a_{24}^{0}$ such that the order of $c_{i}$ agrees with the order of $a_{24}^{\epsilon}$ in $\chi \circ \omega_{5}(\beta_{1}) $.
Then \begin{itemize}
\item $i_{c_{1}}(3)= N_{23}^{0}+N_{34}^{0} = 0+0 =0, mod~2$
\item $i_{c_{2}}(3)= N_{23}^{0}+N_{34}^{0}=1+0 =1,mod~2$
\item $i_{c_{3}}(3)= N_{23}^{0}+N_{34}^{0}=2+2 =0,mod~2$
\item $i_{c_{4}}(3)= N_{23}^{0}+N_{34}^{0} = 2+3 =1.mod~2$
\end{itemize}
Therefore $\phi_{\{2,4\}}^{5}(\beta_{1}) = (0)(1)(0)(1) \neq 1$ and hence $\beta$ is not trivial in $G_{5}^{3}$. 
\end{example}

\begin{example}\label{recog_brun} Let
\begin{eqnarray*}
\beta  &=& [[[b_{12},b_{14}],b_{16}] ,[b_{13},b_{15}]] \\ 
&=& b_{12} b_{14} b_{12}^{-1} b_{14}^{-1} b_{16} b_{14} b_{12} b_{14}^{-1} b_{12}^{-1} b_{16}^{-1}
 b_{13} b_{15} b_{13}^{-1} b_{15}^{-1}b_{16} b_{12} b_{14} b_{12}^{-1} b_{14}^{-1} b_{16}^{-1} b_{14} b_{12} b_{14}^{-1} \\
 &&b_{12}^{-1}b_{15} b_{13} b_{15}^{-1} b_{13}^{-1} \in PB_{6}.
\end{eqnarray*}

  Note that for each $k \in \{ 1,2,3,4,5,6\}$, $p_{k}(\beta) =1$, that is, $\beta$ is Brunnian in $PB_{n}$. Then\\ \\
 $\chi \circ \omega_{6}(r_{1}(\phi_{6}(\beta))) = a_{23}^{0}a_{24}^{0}a_{25}^{0}a_{23}^{1}a_{35}^{1}a_{34}^{1} a_{45}^{0}a_{24}^{0}a_{34}^{0}a_{45}^{1}a_{24}^{1}a_{35}^{1}a_{23}^{1}a_{25}^{0}a_{24}^{0}a_{23}^{0}\\
a_{25}^{1}a_{24}^{1}a_{24}^{0}a_{25}^{1}a_{23}^{0}a_{24}^{0}a_{25}^{0}a_{23}^{1}a_{35}^{1}a_{24}^{1}a_{45}^{1}a_{34}^{0}a_{24}^{0}a_{45}^{0}a_{34}^{1}a_{35}^{1}a_{23}^{1}a_{25}^{0}a_{24}^{0}a_{23}^{0}a_{25}^{1}\\
a_{24}^{1}a_{35}^{1}a_{24}^{1}a_{35}^{1}a_{25}^{0}a_{25}^{1}a_{35}^{1}a_{24}^{1}a_{35}^{1}a_{23}^{1}a_{34}^{0}a_{35}^{0}a_{23}^{0}a_{34}^{1}a_{24}^{1}a_{45}^{1}a_{25}^{0}a_{35}^{0}a_{45}^{0}a_{25}^{1}a_{35}^{1}\\
a_{24}^{1}a_{34}^{1} a_{23}^{0}a_{35}^{0}a_{24}^{0}a_{35}^{0}a_{25}^{1}a_{45}^{0}a_{35}^{1}a_{25}^{0}a_{45}^{1}a_{24}^{1}a_{35}^{0}a_{24}^{1}a_{25}^{1}a_{23}^{1}a_{24}^{0}a_{25}^{0}a_{23}^{0}a_{35}^{0}a_{34}^{0}\\
a_{45}^{0}a_{24}^{0}a_{34}^{1}a_{45}^{1}a_{24}^{1}a_{35}^{0}a_{23}^{0}a_{25}^{0}a_{24}^{0}a_{23}^{1}a_{25}^{1}a_{24}^{1}a_{35}^{0}a_{24}^{1}a_{45}^{1}a_{34}^{1}a_{24}^{0}a_{45}^{0}a_{34}^{0}a_{24}^{1}a_{35}^{0}\\
a_{25}^{1}a_{25}^{1}a_{35}^{0}a_{24}^{1}a_{34}^{0}a_{45}^{0}a_{24}^{0}a_{34}^{1}a_{45}^{1}a_{24}^{1}a_{35}^{0}a_{24}^{1}a_{25}^{1}a_{23}^{1}a_{24}^{0}a_{25}^{0}a_{23}^{0}a_{35}^{0}a_{24}^{1}a_{45}^{1}a_{34}^{1}\\
a_{24}^{0}a_{45}^{0}a_{34}^{0}a_{35}^{0}a_{23}^{0}a_{25}^{0}a_{24}^{0}a_{23}^{1}a_{25}^{1}a_{24}^{1}a_{35}^{0}a_{24}^{1}a_{45}^{1}a_{25}^{0}a_{35}^{1}a_{45}^{0}a_{25}^{1}a_{35}^{0}a_{24}^{1}a_{35}^{0}a_{23}^{0}\\
a_{34}^{1}a_{35}^{1}a_{23}^{1}a_{34}^{0}a_{24}^{1}a_{35}^{0}a_{25}^{1}a_{45}^{0}a_{35}^{1}a_{25}^{0}a_{45}^{1}a_{24}^{1}a_{34}^{0}a_{23}^{1}a_{24}^{1}a_{25}^{1}$\\ \\
and there are 40 $a_{24}^{\epsilon}$. We obtain that 
$$\phi_{\{2,4\}}^{6} = \zeta_{(0,0)}\zeta_{(0,1)}\zeta_{(1,1)}\zeta_{(0,0)}\zeta_{(0,1)}\zeta_{(1,1)}\zeta_{(0,0)}\zeta_{(0,1)}\zeta_{(0,0)}\zeta_{(0,1)}\zeta_{(1,1)}\zeta_{(0,1)},$$
where $\zeta_{a,b}$ is defined by $\zeta_{a,b}(3)=a$ and $\zeta_{a,b}(5)=b$ and $\phi_{\{2,4\}}^{\mathcal{P}}$ is not trivial in $F_{5}^{2}$. 

\end{example}

%% file: realisation.tex
The aim of the present section is to construct spaces where $G_{n}^{k}$ (or a finite subgroup of it) acts faithfully. To that end we construct a space, the fundamental group of which is commensurable with the group $G_n^k$. This construction is useful for the solution of the word problem in those groups.

\subsection{Realisation of the groups $G_{k+1}^k$}
\label{sec:realisation_n=k+1}

In the present section we tackle the problem of realisation of the groups $G_n^k$ for the special case $n=k+1$ (the general case is handled in the subsection~\ref{subsec:reliasation_Gk+1k}).

More precisely, we shall prove the following
\begin{theorem}
There is a subgroup ${\tilde G}_{k+1}^{k}$
of the group $G_{k+1}^{k}$ of index $2^{k-1}$ which
is isomorphic to $\pi_{1}({\tilde C}'_{k+1}(\R{}P^{k-1}))$.\label{th0}
\end{theorem}

The space ${\tilde C}'_{n}$ and the group ${\tilde G}_{k+1}^{k}$ will be defined later in this section.

The simplest case of the above Theorem is

\begin{theorem}\label{thm:realisation_G43}
The group ${\tilde G}_{4}^{3}$
(which is a finite index subgroup of $G_4^3$) is isomorphic to $\pi_{1}(FBr_{4}(\R P^{2}))$,
the $4$-strand braid group on $\R{}P^{2}$ with two points fixed.\label{th1}
\end{theorem}

That theorem gives a solution to the word problem in the group $G_4^3$. We shall prove the theorem in Section~\ref{sect:braid_from_Gkk-1}.

\subsubsection{Preliminary definitions}

We begin with the definition of spaces $C'_{n}(\R{}P^{k-1})$
and maps from the corresponding fundamental groups to
the groups $G_{n}^{k}$.

Let us fix a pair of natural numbers $n>k$. A point in $C'_{n}(\R{}P^{k-1})$
is an ordered set of $n$ pairwise distinct points in $\R{}P^{k-1}$ such that any $(k-1)$
of them are in general position. Thus, for instance, if $k=3$,
then the only condition is that these points are pairwise distinct. For $k=4$
 for points $x_{1},\dots, x_{n}$ in $\R{}P^{3}$
 we impose the condition that no three of them belong
 to the same line (though some four are allowed to belong to the same plane),
 and for  $k=5$ a point in $C'_{n}(\R{}P^{4})$ is a set of ordered $n$ points in  $\R{}P^{4}$,
 with no four of them belonging to the same $2$-plane. \\

Let us use the projective coordinates $(p_{1}:p_{2}:\cdots: p_{k})$ in $\R{}P^{k-1}$. Recall that $p_1$ through $p_k$ are not all zero, and $\vec{v} \sim \vec{w}$ if and only if $v = \lambda w$ for some $\lambda \in {\mathbb R}$. For a point $y\in \R{}P^{k-1}$ we shall denote by $p_i(y)$ the $i-$th projective coordinate of $y$.

Now let us fix the following  $k-1$ points in general position, $y_{1},y_{2}, \dots, y_{k-1}\in \R{}P^{k-1}$, where $p_{i}(y_{j})=\delta_{i}^{j}$. Let us define the subspace ${\tilde C}'_{n}(\R{}P^{k-1})$ taking those $n$-tuples of points $x_{1},\dots, x_{n}\in \R{}P^{k-1}$ for which $x_{i}=y_{i}$ for $i=1,\dots, k-1$. In other words, we take $k-1$ ``standard'' unit points, and append their set with $n-k+1$ arbitrary ones.

 We say that a point $x\in C'_{n}(\R^{k-1})$  is {\em singular}, if the set of points $x=(x_{1},\dots, x_{n})$, corresponding to $x$, contains some subset of $k$ points lying on the same $(k-2)$-plane. Let us fix two non-singular points $x,x'\in C'_{n}(\R{}P^{k-1}).$

We shall consider smooth paths $\gamma_{x,x'}: [0,1]\to C'_{n}(\R^{k-1})$ connecting the points $x$ and $x'$.
For each such path there are values $t$ for which $\gamma_{x,x'}(t)$ is not in general position (some $k$ of them belong to the same $(k-2)$-plane). We call these values $t\in [0,1]$ {\em singular}.

On the path $\gamma$, we say that the moment $t$ of passing through the singular point $\xi$, corresponding to the set  $\xi_{i_1},\dots, \xi_{i_k}$, is {\em transverse} (or {\em stable})\index{Path!transverse}\index{Path!stable} if for any sufficiently small perturbation ${\tilde \gamma}$ of the path  $\gamma$, in the neighbourhood of the moment $t$ there exists exactly one time moment $t'$ corresponding to some set of points $\xi_{j_1},\dots, \xi_{j_k}$ not in general position.

\begin{definition}
We say that a path is {\em good and transverse}\index{Path!good} \index{Path!transverse} if the following holds:

\begin{enumerate}

\item
The set of singular values $t$ is finite;

\item For every singular value  $t=t_{l}$ corresponding to $n$ points representing $\gamma_{x,x'}(t_{l})$, there exists only one subset of $k$ points belonging to a $(k-2)$-plane;

\item Each singular value is  {\em transverse}.

\end{enumerate}
\end{definition}

\begin{definition}
We say that the path without singular values is {\em void}.
\end{definition}

Note that by definition a void path is good and transverse.

We shall call good and transverse paths {\em braids}. We say that a braid whose ends $x,x'$ coincide with respect to the order, is {\em pure}. We say that two  braids $\gamma,\gamma'$ with endpoints $x,x'$ are {\em isotopic} if there exists a continuous family  $\gamma^{s}_{x,x'},s \in [0,1]$ of smooth paths with fixed ends such that  $\gamma^{0}_{x,x'}=\gamma,\gamma^{1}_{x,x'}=\gamma'$. By a small perturbation, any path can be made good and transverse (if endpoints are generic, we may also require that the endpoints remain fixed).

%
%

\subsubsection{The realisability of $G_{k+1}^k$}\label{subsec:reliasation_Gk+1k}

The main idea of the proof of Theorem \ref{th0} is to associate with every word in $G_{k+1}^{k}$ a braid in ${\tilde C}'_{k+1}(\R{}P^{k-1})$.

Let us start with the main construction from \cite{HigherGnk}.

\begin{figure}
\centering
\includegraphics[width=0.8\textwidth]{ 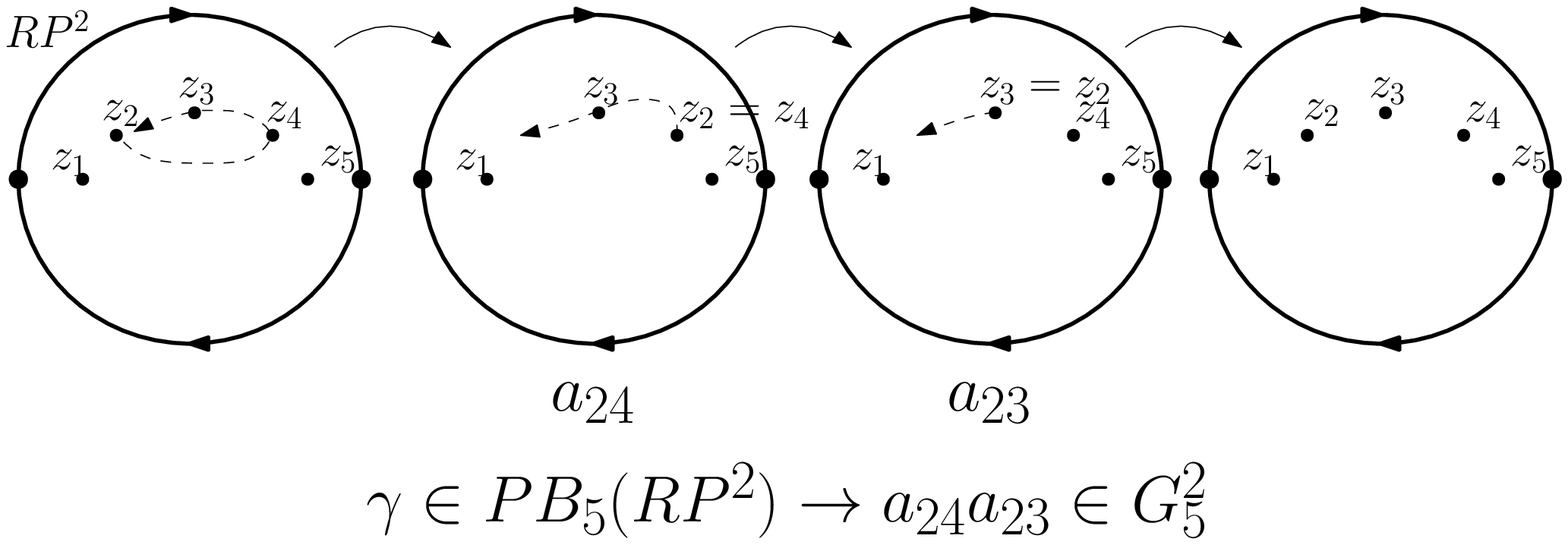}
\caption{An example of associating an element of the group $G_n^k$ with a pure braid in a projective space}
\label{fig:exa-PBnRP2}
\end{figure}

With each good and stable braid from $PB_{n}(\R{}P^{2})$ we associate an element of the group  $G_{n}^{k}$ as follows. We enumerate all singular values of our path $0<t_{1}<\dots <t_{l}<1$ (we assume that $0$ and $1$ are not singular). For each singular value $t_{p}$ we have a set $m_{p}$ of $k$ indices corresponding to the numbers of points which are not in general position. With this value we associate the letter  $a_{m_{p}}$. With the whole path  $\gamma$ (braid) we associate the product $f(\gamma)=a_{m_{1}}\dots a_{m_{l}}$. And example is shown in Fig~\ref{fig:exa-PBnRP2}.

\begin{theorem}\cite{HigherGnk}
The map $f$ takes isotopic braids to equal elements of the group $G_{n}^{k}$.
For pure braids, the map $f$ is a homomorphism  $f\colon\pi_{1}(C'_{n}(\R{}P^{2}))\to G_{n}^{3}$.
\label{thgn3}
\end{theorem}

Now we claim that \\

{\em Any word from  $G_{k+1}^{k}$ can be realised by a path of the above form.} \\

\begin{wrapfigure}{r}{0.4\textwidth}
\centering\includegraphics[width=0.25\textwidth]{ 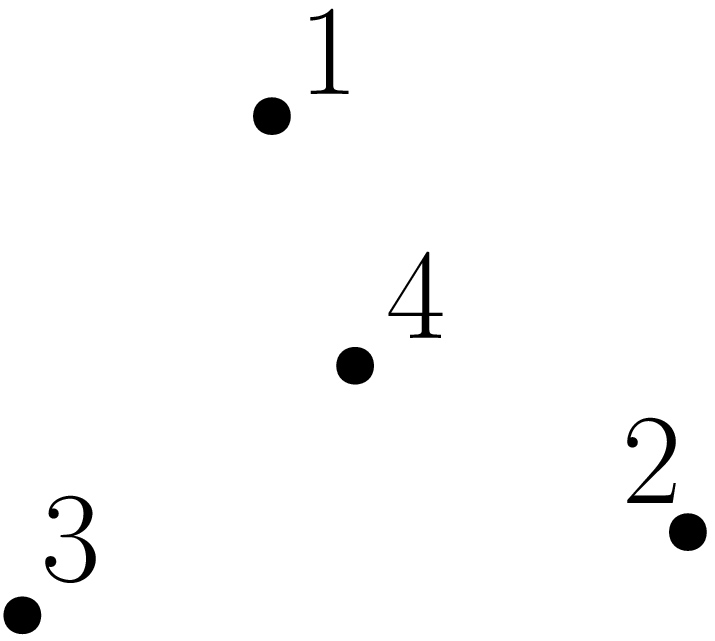}
\caption{The ``locked'' position for the move $a_{123}$}
\label{locked}
\end{wrapfigure}

Note that if we replace $\R{}P^{k-1}$ with $\R^{k-1}$, the corresponding statement will fail.
Namely, starting with the configuration of four points, $x_{i},\;i=1,\dots, 4$, where
$x_{1},x_{2},x_{3}$ form a triangle and $x_{4}$ lies inside this triangle, see Fig.~\ref{locked}, we see that any path starting from this configuration will lead to a word starting from $a_{124},a_{134}, a_{234}$ but not from $a_{123}$. In some sense the point $4$ is ``locked'' and the roles of the points are not the same.

The following well-known theorem (see, e.g., \cite{Wu}) plays a crucial role in the construction

\begin{theorem}
For any two sets of $k+1$ points in general position in $\R{}P^{k-1}$,
$(x_{1},\dots, x_{k+1})$ and $(y_{1},\dots, y_{k+1})$ there is an action of
$PGL(k, \R)$ taking all $x_{i}$ to $y_{i}$.\label{wulem}
\end{theorem}

For us this  will mean that there is no difference between all possible ``non-degenerate starting positions'' for $k+1$ points in $\R{}P^{k-1}$.

We shall consider paths in ${\tilde C'}_{k+1}(\R{}P^{k-1})$ starting and ending in two points from this set (possibly, the same). 

We shall denote homogeneous coordinates in $\R{}P^{k-1}$ by $(p_{1}:\cdots: p_{k})$ in
contrast to points (which we denote by $(x_{1},\dots, x_{k+1})$).

\subsubsection{Constructing a braid from a word in $G_{k+1}^{k}$}\label{sect:braid_from_Gkk-1}

Our main goal is to construct a braid from a word. To this end, we need a base point for the braid. For the sake of convenience, we shall use not one, but rather $2^{k-1}$ reference points. For the first $k$ points $y_{1}=(1:0:\cdots:0),\dots, y_{k}=(0:\cdots:0:1)$ fixed, we will have $2^{k-1}$ possibilities for the choice of the last point. Namely, let us consider all possible strings of length $k$ of $\pm 1$ with the last coordinate $+1$: $$(1:1:\cdots: 1:1),(1:\cdots: 1:-1:1),\cdots, (-1:-1:\cdots: -1:1)$$ with $p_{k} = +1$. We shall denote these points by $y_{\vec{\epsilon}}$ where $\vec{\epsilon}=(\pm 1,\dots ,\pm 1)$ records the first $(k-1)$ coordinates of the point (the last coordinate always is fixed to 1).

Now, for each string $\vec{\epsilon}=(\pm 1,\dots ,\pm 1)$ of length $k-1$, we set $z_{\vec{\epsilon}}=(y_{1},y_{2},\dots, y_{k},y_{s})$.

The following lemma is evident:
\begin{lemma}
For every point $z\in \R{}P^{k-1}$ with projective coordinates
$(p_{1}(z):\cdots: p_{k-1}(z):1)$, let ${\tilde z}=(\sign(p_{1}(z)):\sign(p_{2}(z)):
\cdots: \sign(p_{k-1}(z)):1)$. Then there is a path
between $(y_{1},\dots, y_{k},z)$ and $(y_{1},\dots, y_{k}, {\tilde z})$
in ${\tilde C}'_{k+1}(\R{}P^{k})$ with the first points $y_{1},\dots, y_{k}$
fixed, and the corresponding path in ${\tilde C}'_{k+1}$ is void.
\label{lemmaB}
\end{lemma}

\begin{proof}
Indeed, it suffices just to connect $z$ to ${\tilde z}$ by the shortest geodesic.
\end{proof}

From this we easily deduce the following
\begin{lemma}
Every  point $y\in {\tilde C}'_{k+1}(\R{}P^{k-1})$ can be connected by a void
path to some
$(y_{1},\dots, y_{k}, y_{\vec{\epsilon}})$ for some $\vec{\epsilon}$.\label{lemBB}
\end{lemma}

\begin{proof}
Indeed, the void path can be constructed in two steps. At the first step, we
construct a path which moves both $y_{k}$ and $y_{k+1}$, so that $y_{k}$
becomes $(0:\cdots :0:1)$, and at the second step, we use Lemma \ref{lemmaB}.
To realise the first step, we just use linear maps which keep the hyperplane
$p_{k}=0$ fixed.
\end{proof}

The lemma below shows that the path mentioned in Lemma \ref{lemBB} is unique up to homotopy.

\begin{lemma}
Let $\gamma$ be a closed path in ${\tilde C}'_{k+1}(\R{}P^{k-1})$ such that
the word $f(\gamma)\in G_{k+1}^k$ is empty. Then $\gamma$ is homotopic to the trivial braid.
\label{lmA}
\end{lemma}

\begin{proof}

In ${\tilde C}'_{k+1}$, we deal with the motion of points, where all but $x_{k},x_{k+1}$ are fixed.

Consider the projective hyperplane ${\cal P}_{1}$ passing through $x_{1},\dots, x_{k-1}$ given by the equation $p_{k}=0$. We know that none of the points $x_{k}, x_{k+1}$ is allowed to belong to ${\cal P}_{1}$. Hence, we may fix the last coordinate $p_{k}(x_{k})=p_{k}(x_{k+1})=1$.

Now, we may pass to the affine coordinates of these two points (still to be denoted by $p_{1},\dots, p_{k}$). The condition $\forall i=1,\dots, k-1: p_{j}(x_{k})\neq p_{j}(x_{k+1})$ follows from the fact that the points $x_{1},\dots, {\hat{x}}_{j}, \dots, x_{k+1}$ are generic. This means that $\forall i=1,\dots, k-1$ the sign of $p_{i}(x_{k})-p_{i}(x_{k+1})$ remains fixed.

The motion of points $x_{k},x_{k+1}$ is determined by their coordinates $p_{1},\dots, p_{j}$, and since their signs are fixed, the configuration space for this motions is simply connected. This means that the loop $\gamma$ is described by a loop in a two-dimensional simply connected space.
\end{proof}

Our next strategy is as follows. Having a word in $G_{k+1}^{k}$, we shall associate with this word a path in ${\tilde C}'_{k+1}(\R{}P^{k-1})$. After each letter, we shall get to  $(y_{1},\dots ,y_{k},y_{\vec{\epsilon}})$ for some $\vec{\epsilon}$.

Let us start from $(y_{1},\dots, y_{k},y_{1,\dots ,1})$.

After making the final step, one can calculate the coordinate of the $(k+1)$-th points. They will be governed by Lemma~\ref{governed} (see ahead). As we shall see later, those words we have to check for the solution of the word problem in $G_{k+1}^{k}$, will lead us to closed paths, i.e., pure braids.

Let us be more detailed.

\begin{lemma}
Let a non-singular set of points ${\bf y} = \{y_1, \dots, y_n\}$ in $\R{}P^{k-1}$ be given. Then for every set of $k$ numbers $i_{1},i_{2},\dots, i_{k}\in \bar{n}$, there exists a path $y_{i_{1}\dots i_{k}}(t) = {\bf y}(t) = \{y_1(t),\dots,y_n(t)\}$ in $C'_{n}(\R{}P^{k-1})$, having $y(0)=y(1)={\bf y}$ as the starting point and the final point and with only one singular moment $t_0$ corresponding to the numbers $i_{1},\dots, i_{k}$ in the sense that exactly the points $y_{i_1}(t_0),\dots,y_{i_k}(t_0)$ are not in general position. Moreover, we may assume that at this moment all points except the $i_1-$th one, are fixed during the path.

Moreover, the set of paths possessing this property is connected: any other path ${\tilde y}(t)$, possessing the above properties, is homotopic to $y(t)$ in this class.
 \label{lm2}
\end{lemma}

\begin{proof}
Indeed, for the first statement of the Lemma,
it suffices to construct a path for some initial position of points and
then apply Theorem \ref{wulem}.

For the second statement, let us take two different paths $\gamma_{1}$ and
$\gamma_{2}$ satisfying the conditions of the Lemma. By a small perturbation, we may assume that for both of them,
 $t=\frac{1}{2}$ is a singular moment with the same position of $y_{i_{1}}$.

We can contract the loop formed by $\gamma_{1}|_{t\in [\frac{1}{2},1]}$
and the inverse of $\gamma_{2}|_{t\in [\frac{1}{2},1]}$ by using Lemma \ref{lmA} as
this is a small perturbation of a void braid. We are left with
$\gamma_{1}|_{t\in [0,\frac{1}{2}]}$ and the inverse of $\gamma_{2}|_{t\in [0,\frac{1}{2}]}$
which is contractible by Lemma \ref{lmA} again.

\end{proof}

\begin{remark}
Note that in the above lemma, we deal with the space $C'_{n}(\R{}P^{k-1})$ (whose points are ordered sets of $n$ pairwise distinct points in $\R{}P^{k-1}$ such that any $(k-1)$ of them are in general position), and not with ${\tilde C}'_{n}(\R{}P^{k-1})$ (where the first $k-1$ points are fixed to be the ``standard'' unit ones). On the other hand, we may always choose $i_{1}\in \{k,k+1\}$; hence, the path in question can be chosen in ${\tilde C}'(\R{}P^{k-1})$.
\end{remark}

For every subset
 $m\subset \bar{n}, Card(m)=k+1$ we can create a path
 $p_{m}$ starting from any of the base points listed
 above and ending at the corresponding basepoints.

We construct our path step-by step by applying Lemma \ref{lm2} and returning
to some of base points by using Lemma \ref{lemBB}.

From \cite{HigherGnk}, we can easily get the following
\begin{lemma}
Let $i_{1},\dots, i_{k+1}$ be some permutation of $1,\dots, k+1$. Then the concatenation
of paths $p_{i_{1}i_{2}\cdots i_{k}}p_{i_{1}i_{3}i_{4}\cdots i_{k+1}}\cdots
p_{i_{2}i_{3}\cdots i_{k}}$ \\ is homotopic to the concatenation of paths in the inverse
order $$p_{i_{2}i_{3}\cdots i_{k}}\cdots p_{i_{1}i_{3}i_{4}\cdots i_{k+1}}p_{i_{1}i_{2}\cdots i_{k}}.$$ \label{quadrisec}
\end{lemma}

\begin{proof}
Indeed, in \cite{HigherGnk}, some homotopy corresponding to the above mentioned relation
corresponding to {\em some} permuation is discussed. However, since all basepoints are
similar to each other as discussed above, we can transform the homotopy from \cite{HigherGnk}
to the homotopy for any permutation.
\end{proof}

\begin{lemma}
For the path starting from the point $(y_{1},\cdots, y_{k},y_{s})$ constructed as in
Lemma \ref{lm2} for the set of indices $j$,
we get to the point $(y_{1},\cdots, y_{k},y_{s'})$ such that:
\begin{enumerate}
\item if $j=1,\dots, k$, then $s'$ differs from $s$ only in coordinate $a_{j}$;
\item if $j=k+1$, all coordinates of $s'$ differ from those coordinates of $s$ by sign.
\end{enumerate}
\label{governed}
\end{lemma}

Denote the map from words in $G_{k+1}^{k}$ to paths between basepoints by $g$.

By construction, we see that for every word $w$ we have $f(g(w))=w\in G_{k+1}^{k}$.

Now, we define the group ${\tilde G}_{k+1}^{k}$ as the subgroup of $G_{k+1}$ which is
taken by $g$ to {\em braids}, i.e., to those paths with coinciding initial and final
points in the configurational space. From Lemma \ref{governed}, we see that this is a subgroup of index $(k-1)$: there are exactly
$(k-1)$ coordinates. \\

Let us pass to the proof of Theorem \ref{th1}.
Our next goal is to see that equal words can originate only from homotopic paths.

To this end, we shall first show that the map $f$
from Theorem \ref{thgn3}
is an isomorphism for $n=k+1$. To perform this goal, we should construct the
inverse map $g:{\tilde G}_{k+1}^{k}\to \pi_{1}({\tilde C'}_{k+1}(\R{}P(k-1)))$.

Note that for $k=3$ we deal with
the pure braids $
PB_{4}(\R{}P^{2}). $

Let us fix a point $x\in C'_{4}(\R{}P^{2})$.
With each generator $a_{m}, m\subset \bar{n},$ $Card(m)=k$ we associate a path  $g(m)=y_{m}(t)$,
described in Lemma \ref{lm2}. This path is not a braid: we can return to any
of the $2^{k-1}$ base points. However, once we take the concatenation of paths correspoding
to ${\tilde G}_{k+1}^{k}$, we get a braid.

By definition of the map $f$, we have $f(g(a_{m}))=a_{m}$.
Thus, we have chosen that the map $f$ is a surjection.

Now, let us prove that the kernel of the map $f$ is trivial. Indeed,
assume there is a pure braid $\gamma$ such that $f(\gamma)=1\in G_{k+1}^{k}$.
We assume that $\gamma$ is good and stable. If this path has $l$ critical points, then
we have the word corresponding to it  $a_{m_1}\cdots a_{m_l}\in G_{k+1}^{k}$.

Let us perform the transformation $f(\gamma)\to 1$ by applying the relations of $G_{k+1}^{k}$
to it and transforming the path $\gamma$ respectively. For each relation of the sort $a_{m}a_{m}=1$ for a generator $a_{m}$ of the group
$G_{k+1}^{k}$, we see that the path $\gamma$ contains two segments whose concatenation
is homotopic to the trivial loop (as follows from the second statement of Lemma \ref{lm2}).

Whenever we have a relation of length $2k+2$ in the group $G_{k+1}$, we use the Lemma
\ref{quadrisec} to perform the homotopy of the loops.

 Thus, we have proved that if the word $f(\gamma)$ corresponding to a braid $\gamma\in G_{k+1}^{k}$
is equal to $1$ in $G_{k+1}^{k}$ then the braid $\gamma$ is isotopic to a braid $\gamma'$
such that the word corresponding to it is empty.
Now, by Lemma \ref{lmA}, this braid is homotopic to the trivial braid. This ends the proof of Theorem~\ref{thm:realisation_G43}.

\subsubsection{The group $H_{k}$ and the algebraic lemma}

The aim of the present section is to reduce the word problem in $G_{k+1}^{k}$ to the word problem
in a certain subgroup of it, denoted by $H_{k}$.\index{$H_k$}

Let us rename all generators of $G_{k+1}^{k}$ lexicographically:
$$b_{1}=a_{1,2,\dots,k},\cdots, b_{k+1}=a_{2,3,\dots, k+1}.$$

Let $H_{k}$ be the subgroup of $G_{k+1}^{k}$ consisting of all elements $x\in G_{k+1}^{k}$ that can be
represented by words with no occurrences of the last letter $b_{k+1}$.

Our task is to understand whether a word in $G_{k+1}^{k}$ represents an element in $H_{k}$.
To this end, we recall the map from \cite{MN}. Consider the group $F_{k-1}=\Z_{2}^{*2^{k-1}}=\langle
c_{m}|c_{m}^{2}=1\rangle$, where all generators $c_{m}$ are indexed by $(k-1)$-element strings of $0$
and $1$ with only relations being that the square of each generator is equal to $1$.
We shall construct a map\footnote{This map becomes a homomorphism when restricted to a
finite index subgroup.} from $G_{k+1}^{k}$ to $F_{k-1}$ as follows.

Take a word $w$ in generators of $G_{k+1}^{k}$ and list all occurrences of the last letter
$b_{k+1}=a_{2,\cdots, k+1}$ in this word. With any such occurrence we first associate the string of indices $0,1$
of length $k$. The $j$-th index is the number of letters $b_{j}$ preceding this occurrence of $b_{k+1}$ modulo $2$.
Thus, we get a string of length $k$ for each occurrence.

Let us consider ``opposite pairs'' of strings $(x_{1},\cdots, x_{k})\sim(x_{1}+1,\cdots, x_{k}+1)$ as equal.
Now, we may think that the last ($k$-th) element of our string is always $0$, so, we can restrict ourselves
with $(x_{1},\cdots, x_{k-1},0)$. Such a string of length $k-1$ is called the {\em index}\index{Index} of
the occurrence of $b_{k+1}$.

Having this done, we associate with each occurrence of $b_{k+1}$ having index $m$ the generator $c_{m}$ of $F_{k-1}$.
With the word $w$, we associate the word $f(w)$ equal to the product of all generators $c_{m}$ in order.

In \cite{MN}, the following Lemma is proved:
\begin{lemma}
The map $f\colon G_{k+1}^{k}\to F_{k-1}$ is well defined.
\end{lemma}

Now, let us prove the following crucial
\begin{lemma}
If $f(w)=1$ then $w\in H_{k}$.
\end{lemma}

In other words, the group $F_{k-1}$ yields the only obstruction for an element from $G_{k+1}^{k}$ to have a presentation
with no occurrence of the last letter.

\begin{proof}
Let $w$ be a word such that $f(w)=1$. If $f(w)$ is empty, then there is nothing to prove.
Otherwise $w$ contains two ``adjacent''
occurrences of the same index. This means that
$w=A b_{k+1} B b_{k+1} C$, where $A$ and $C$ are some words, and $B$ contains no occurrences of $b_{k+1}$
and the number of occurrences of $b_{1},b_{2},\cdots, b_{k}$ in $B$ are of the same parity.

Our aim is to show that $w$ is equal to a word with smaller number of $b_{k+1}$ in $G_{k+1}^{k}$. Then we will be able to induct
on the number of $b_{k+1}$ until we get a word without $b_{k+1}$.

Thus, it suffices for us to show that
$b_{k+1}Bb_{k+1}$ is equal to a word from $H_{k}$. We shall induct on the length of $B$.
Without loss of generality, we may assume that $B$ is reduced, i.e., it does not contain
adjacent $b_{j}b_{j}$.

Let us read the word $B$ from the very beginning $B=b_{i_1}b_{i_2}\cdots$ If all elements
$i_{1},i_{2},\cdots$ are distinct, then, having in mind that the number of occurrences of all generators
in $B$ should be of the same parity, we conclude that $b_{k+1}B= B^{-1} b_{k+1}$, hence $b_{k+1}Bb_{k+1}=B^{-1}b_{k+1}b_{k+1}=
B^{-1}$
is a word without occurrences of $b_{k+1}$.

Now assume ${i_1}={i_p}$ (the situation when the first repetition is for ${i_j}={i_p},1<j<p$
is handled in the same way). Then we have
$b_{k+1}B=b_{k+1}b_{i_1}\cdots b_{i_{p-1}} b_{i_1} B'$. Now we collect all indices distinct from
$i_{1},\cdots, i_{p-1},{k+1}$ and write the word $P$ containing exactly one generator for each of these indices
(the order does not matter). Then the word $W = P b_{k+1} b_{i_{1}}\cdots b_{i_{p-1}}$ contains
any letter exactly once and we can invert the word $W$ as follows: $W^{-1}=b_{i_{p-1}}\cdots b_{i_{1}}b_{k+1}P^{-1}$.
Thus, $b_{k+1}B=P^{-1}(Pb_{k+1}b_{i_1}\cdots b_{i_{p-1}})b_{i_1}B'=P^{-1}b_{i_{p-1}}\cdots b_{i_1}b_{k+1}P^{-1}b_{i_1}B'$.

We know that the letters in $P$ (hence, those in $P^{-1}$) do not contain $b_{i_1}$. Thus,
the word $P^{-1} b_{i_1}$ consists of distinct letters. Now we perform the same trick:
create the word $Q=b_{i_2}b_{i_3}\cdots b_{i_{p-1}}$ consisting of remaining letters from $\{1,\dots, k\}$,
we get:

$$ P^{-1}b_{i_{p-1}}\cdots b_{i_1}b_{k+1}P^{-1}b_{i_1}B'$$
$$=P^{-1}b_{i_{p-1}}\cdots b_{i_1}QQ^{-1}b_{k+1}P^{-1}b_{i_1}B'$$
$$=P^{-1}b_{i_{p-1}}\cdots b_{i_1}Qb_{i_1}Pb_{k+1}QB'.$$

Thus, we have moved $b_{k+1}$ to the right and between the two occurrences
of the letter
$b_{k+1}$, we replaced $b_{i_1}\cdots, b_{i_{p-1}}b_{i_1}$ with
just $b_{i_2}\cdots b_{i_{p-1}}$, thus, we have shortened the distance
between the two adjacent occurrences of $b_{k+1}$.

Arguing as above, we will finally cancel these two letters $b_{k+1}$ and perform the
induction step.

\end{proof}

\begin{theorem}
In the group $G_{k+1}^{k}$, the problem whether $x\in G_{k+1}^{k}$ belongs
to $H_{k}$, is solvable.\label{kl}
\end{theorem}

\begin{proof}
Indeed, having a word $w$ in generators of $G_{k+1}^{k}$, we can look at the image of this
word by the map $f$. If it is not equal to $1$, then, from \cite{MN},
it follows that $w$ is non-trivial, otherwise we can construct a word ${\tilde w}$ in $H_{k}$
equal to $w$ in $G_{k+1}^{k}$.
\end{proof}

\begin{remark}
A shorter proof of Theorem \ref{kl} based on the same ideas was communicated to the authors by A.A.~Klyachko.
We take the subgroup $K_{k}$ of $G_{k+1}^{k}$ generated by products
$B_{\sigma}=b_{\sigma_{1}}\cdots b_{\sigma_{k}}$ for all permutations $\sigma$
of $k$ indices. This group $K_{k}$ contains the commutator of $H_{k}$ and is a normal
subgroup in $G_{k+1}^{k}$.

Moreover, the quotient group $G_{k+1}^{k}/ K_{k}$ is naturally isomorphic to
the free product $(H_{k}/K_{k})*\langle b_{k+1}\rangle$. Hence, the problem whether an element
of $G_{k+1}^{k}/ K_{k}$ belongs to $H_{k}/ K_{k}$ is trivially solved,
which solves the initial problem because of the normality of $K_{k}$ in $G_{k+1}^{k}$.
\end{remark}

Certainly, to be able to solve the word problem in $H_{k}$, one needs to know a presentation for
$H_{k}$. It is natural to take $b_{1},\cdots, b_{k}$ for generators of $H_{k}$. Obviously,
they satisfy the relations $b_{j}^{2}=1$ for every $j$.

To understand the remaining relations for different $k$, we shall need geometrical arguments. \\

We have completely constructed the isomorphism between the (finite index subgroup) of
the group $G_{k+1}^{k}$ and a fundamental group of some configuration space.

This completely solves the word problem in $G_{4}^{3}$ for braid groups in projective spaces are very
well studied, see, e.g., \cite{JvB}. The same can be said about the conjugacy problem in ${\tilde G}_{4}^{3}$.

Besides, we have seen that the word problem for the case of general $G_{k+1}^{k}$
can be reduced to the case of $H_{k}$.

Let us mention that it was proved by A.B.~Karpov (unpublished)
that the only relations in $H_{3}$ are $a^{2}=b^{2}=c^{2}$ which also follows
from the geometrical techniques of the present section.

The main open question which remains unsolved is how to construct a configuration
space which can realise $G_{n}^{k}$ for $n>k+1$. Even in the case of $G_{4}^{2}$
this problem seems very attractive though the word and conjugacy problems for $G_{n}^{2}$ can
be solved by algebraic methods, see \cite{Coxeter}.

The word problem and the partial case of conjugacy problem in $G_{4}^{3}$ was first
solved by A.B.~Karpov, see \cite{upcoming_book} for details.

It would be very interesting to compare the approach of $G_{n}^{k}$ with various generalisations
of braid groups, e.g., Manin-Schechtmann groups \cite{ManinSchechtmann,KapranovVoevodsky}.\index{Group!Manin-Schechtmann}

\subsection{Realisation of $G_{n}^{k}$ for $n\neq k+1$} \label{subsec:reliasation_Gnk}

In the present section, we shall tackle the question, {\em how to represent all words from $G_{n}^{k}$
by paths in some configuration space in $n>k+1$}.

The answer will be: for whatever $k$, we shall deal with configurations of $n$ points in $\R{}P^{n-2}$,
and the condition will be: $k$ points are not in the general position.

An obvious objection by the reader should be: {\em this condition is not of codimension $1$.}

Certainly, this is true, and to handle the problem, let us be more detailed.

\subsubsection{A simple partial case} \label{sec:realisation_1}

We begin with a simple example: the group $G_{4}^{2}$.

By $\tilde{C}(\R{}P^{2},4)$ (for shortness, $\tilde{C}_{4}^{2}$) we denote the following configuration space.
A point of ${\tilde C}_{4}^{2}$ is a collection of four points $z_{1},z_{2},z_{3},z_{4}\in \R{}P^{2}$
(not necessarily pairwise distinct) together with six selected (projective) lines $l_{ij}$,
passing through $z_{i},z_{j}$. Here all $l_{ij}$ belong to $\R{}P^{1}$.

Note that there is a natural ``line forgetting map'' $f:{\tilde C}_{4}^{2}\to (\R{}P^{2})^{4}$
which takes $(z_{i}, l_{ij})$ to $(z_{i})$.

Note that the restriction of this map to the subspace where the four points are pairwise distinct
is a bijection. If there is exactly one pair of coinciding points, say, $z_{1}=z_{2}$,
this map is a fibration with fibre $\R{}P^{1}$.

Now, let us make a restriction to the above configuration space:
by ${\tilde C}' (\R P^{2},4)$ (for shortness, $\tilde{C'}_{4}^{2}$) we
denote those points $(z_{i}, \cdots, l_{ij},\cdots)$ of $\tilde{C}_{4}^{2}$
where no three distinct points $z_{i},z_{j},z_{k}$ are collinear.

A point in $\tilde{C'}(\R{}P^{2},4)$ with $z_{i}=z_{j}$ for some $i\neq j$
will be called {\em non-generic}. Note that the set of all non-generic points in
$\tilde{C'}(\R{}P^{2},4)$ has codimension $2$.

Indeed, let us calculate the dimensions explicitly.

Fix two generic points $p,p^{*}\in {\tilde C}'(\R{}P^{2},4)$. A {\em generic path} from $p$ to $p^{*}$ is a smooth path in the space $\tilde{C'}_n^k$ connecting the points $p$ and $p^*$ and passing through a finite number of non-generic points.

With any generic path $\gamma$ from $p$ to $p^{*}$ we associate a word in letters
$a_{12},a_{13},a_{14},a_{23},a_{24},a_{34}$ as follows: when passing through
a point where $z_{i}$ coincides with $z_{j}$, we write down $a_{ij}$.

Thus we get a word $w(\gamma)$ in $a_{ij}$ and the corresponding element of $G_{4}^{2}$.

\begin{theorem} \label{thm:path_and_word}
The above element of $G_{4}^{2}$ depends only on the homotopy class of the path $\gamma$.

In particular, if $p=p^{*}$, we get a map from $\pi_{1}(\tilde{C'}_{4}^{2})$ to $G_{4}^{2}$.
\end{theorem}

 \begin{figure}
 \label{fig:g42sketch}
 \centering
    \includegraphics[width=0.8\textwidth]{ 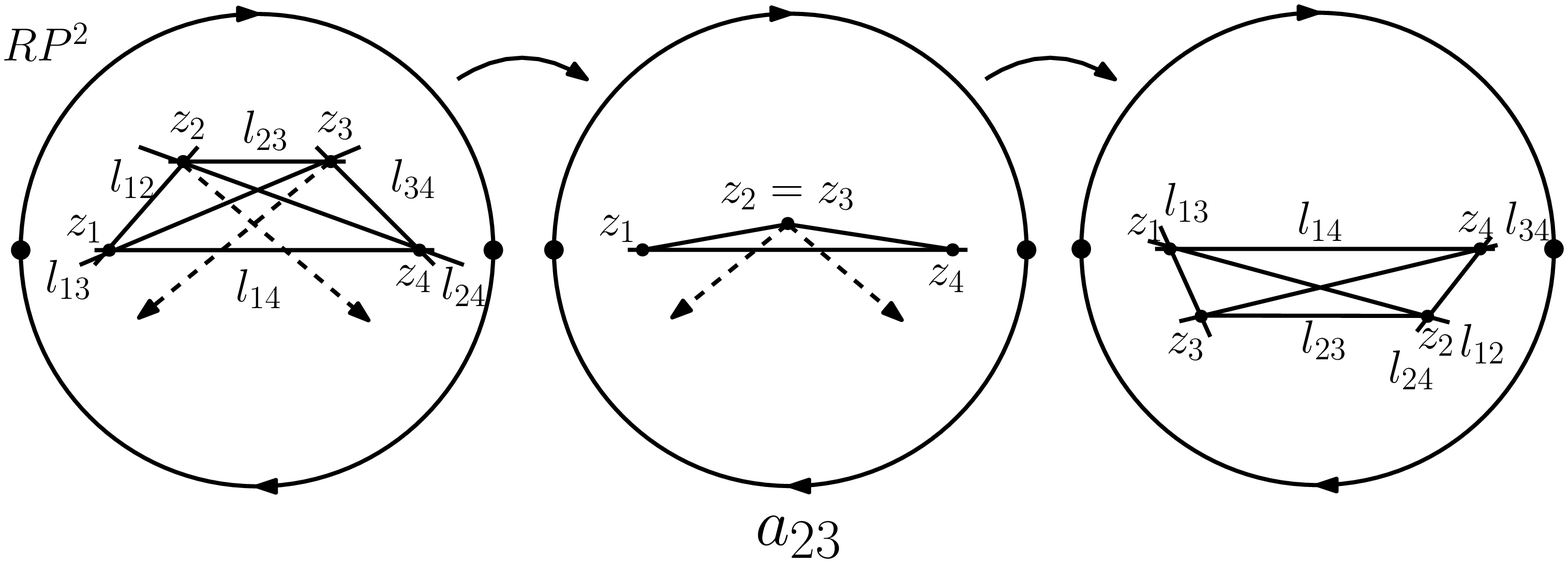}
  \caption{A step in recovering a path from a word in $G_4^2$}
 \end{figure}

Certainly, we have already had a map from homotopy groups of ``some'' configuration
spaces to $G_{n}^{k}$. The main advantage here is that there is an ``inverse'' operation: having a word in letters of $G_{4}^{2}$, we can restore a certain path $\gamma$ which gives rise to such a word.

Let $w$ be such a word. We put some four generic points $z_{1},z_{2},z_{3},z_{4}$ in $\R{}P^{2}$. The point $p$ will be $(z_{1},z_{2},z_{3},z_{4},l_{12},l_{13},l_{14},l_{23},l_{24},l_{34})$ (note that $l_{ij}$ are automatically restored from $z_{i}$).

Take the first letter $w_{1}=a_{ij}$ of $w$.

In the path $\gamma$ corresponding to $a_{ij}$ all points $z_{j}$ except $z_{1}$ will be fixed. Connect $z_{1}$ with $z_{2}$ by shortest geodesic; let $z'_{1}$ be a point such that $z_{1}z'_{1}$ is a geodesic segment containing $z_{2}$. In the path $\gamma$ the point $z_{1}$ will go through $z_{2}$ to $z'_{1}$; $l_{12}$ will stay the same as an element of $\R{}P^{1}$; the motion of other $l_{1j},j\neq 2$ is uniquely restored by the motion of $z_{j}$. After that we get to a generic point $(z'_{1},z_{2},z_{3},z_{4},l'_{12},\cdots)$.

Take the second letter $w_{2}=a_{kl}$ of the word $w$ and proceed as above. A sketch of this process is depicted in Fig.~\ref{fig:g42sketch}.

By construction, we get the following
\begin{statement}
$w(\gamma_{w})=w$.
\end{statement}

\subsubsection{General construction}

Fix $n>k+1$. We will follow the same pattern as used above in the case of $G_4^2$.

Consider the configuration space $\tilde{C}_n^k$. A point of the space $\tilde{C}_n^k$ is a collection of $n$ points $z_1, \dots, z_n \in \R{}P^{n-2}$ such that any $k-1$ of them are in general position, together with a collection of $(k-1)$-planes $P_{i_1 \dots i_k}$ passing through the points $z_{i_1}, \dots, z_{i_k}$ for each $k-$set $(i_1,\dots, i_k), \, i_j\in\{1,\dots,n\}$. Furthermore, we define a restriction of this space: by $\tilde{C'}_n^k$ we denote the set of such points $(z_i, \dots, P_{i_1 \dots i_k}, \dots)$ of $\tilde{C}_n^k$ that if some $k+1$ points $z_{i_1},\dots,z_{i_{k+1}}$ are not in general position then some $k$ among them are not in general position.

A point $z$ in $\tilde{C'}_n^k$ such that some $k$ points $z_{i_1},\dots,z_{i_k}$ are not in general position will be called {\em non-generic}. In the same manner as in the case of the group $G_4^2$ one can define a {\em generic path}\index{Path!generic} in the space $\tilde{C'}_n^k$ connecting two generic points $p, p^*$ and associate with it a word in letters $a_{m_i}$ where $m_i$ denote the $k-$subsets of the the set $\bar{n}$ and a letter $a_m, m=(i_1,\dots,i_k)$ is written when passing through a point where $z_{i_1},\dots,z_{i_{k}}$ are not in general position. Thus we obtain a word $w(\gamma)$ in letters $a_{m_i}$ and the corresponding element of $G_n^k$.

\begin{theorem}
	The above element of $G_n^k$ depends only on the homotopy class of the path $\gamma$.
	
	In particular, if $p=p^*$, we get a map from $\pi_1(\tilde{C'}_n^k)$ to the group $G_n^k$.
\end{theorem}

This theorem may be proved by checking the possible transformations of the words, induced by homotopies of the path $\gamma$. For example, if we homotop the path by going through a non-generic point $z$ (``puncture a plane'', in essence) and then return, going through the same plane (that is, through a non-generic point $z'$, not necessary coinciding with the point $z$ but such that the same set of $z_{i_!},\dots,z_{i_k}$ are not in general position, we obtain a subword $a_m a_m$ which is trivial due to the relations in the group $G_n^k$ and thus may be removed.

Just as in Section \ref{sec:realisation_1}, Theorem \ref{thm:path_and_word}, the advantage here is that we have an ``inverse'' operation: given a word in letters of $G_n^k$ we can restore a path $\gamma$ which gives rise to such a word. To be more precise, the following holds:

\begin{proposition}
	For every word $w\in G_n^k$ there exists a path $\gamma_w$ in the space $\tilde{C'}_n^k$ such that $w(\gamma_w)=w$. Moreover, for a given word the path $\gamma_w$ may be constructed algorithmically.
\end{proposition}

\subsection{The $G_{n}^{k}$-complex}

The aim of the present subsection is to construct a complex whose fundamental group
is commensurable with the group $G_{n}^{k}$.\index{$G_n^k$-complex} \\

Recall that two groups are called {\em commensurable} if they have isomorphic finite index subgroups.

We start with the toy model $G_{3}^{2}$. The group has ${3\choose 2} =3$ generators
$a_{12},a_{13},a_{23}$ which are all involutions, and all the relations (except $a_{ij}^{2}=1$)
can be represented by words of length six: $a_{ij}a_{ik}a_{jk}a_{ij}a_{ik}a_{jk}=1$.

Let $\{0,1\}^{3}$ be the $0$-skeleton of the $3$-dimensional unit cube in the space with coordinates $x,y,z$ which we associate with the pairs $(12), (13), (23)$ in the following way: $x\to(12),y\to(13),z\to(23)$.

Now consider the $1$-skeleton $\Sigma^1_{3}$ of the $3$-dimensional cube with the above mentioned vertices. Let us take the homotopy classes of paths beginning at the vertex $(0,0,0)$ in this 1-dimensional complex. They are in correspondence with the elements of the group $$\mathbb{Z}_{2}*\mathbb{Z}_{2}*\mathbb{Z}_{2}=\langle a_{12},a_{13},a_{23}|a_{ij}^{2}=1\rangle.$$

Note that in this construction the group relations are given by {\em homotopically trivial loops}, not {\em any} loops returning to origin. In particular, a quadrilateral $(0,0,0)\to(0,0,1)\to(1,0,1)\to(1,0,0)\to(0,0,0)$ which corresponds to the word $w=a_{23}a_{12}a_{23}a_{12}$ does not yield a relation (and hence does not mean that $w=1$ in the group) because the corresponding loop is not trivial in the 1-complex. That means that the only relations provided by the complex are of the form $a_{ij}^2=1$ coming from loops going along an edge and back again. That gives us exactly the group $\mathbb{Z}_{2}*\mathbb{Z}_{2}*\mathbb{Z}_{2}$.

The words in the alphabet $\{a_{12}, a_{13}, a_{23}\}$ are naturally split into eight classes defined with respect to the parities of the number of occurrences of the letters $a_{12}$, $a_{13}$, and $a_{23}$. These classes correspond to eight vertices of the cube: if we start at the origin $(0,0,0)$, we can get to all eight vertices going along some paths. Two paths (group elements) ending in the same vertex belong to the same class.

\begin{figure}
\centering\includegraphics[width=300pt]{ 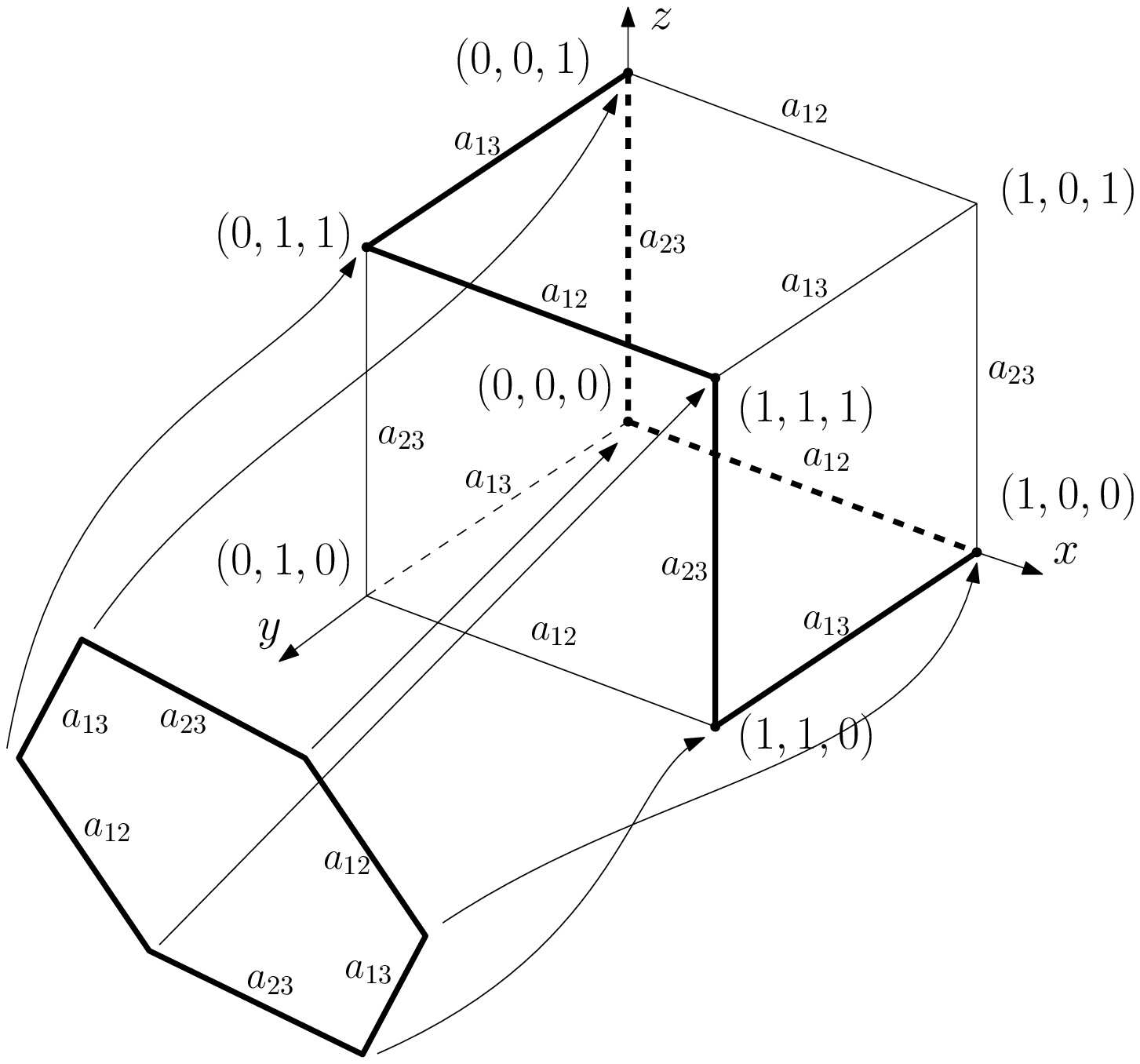}
\caption[A cell corresponding to the word $a_{12}a_{13}a_{23}a_{12}a_{13}a_{23}$]{Black edges of the cube form the boundary of the hexagon corresponding to the word $a_{12}a_{13}a_{23}a_{12}a_{13}a_{23}$}
\label{fig:gnk_cube}
\end{figure}

%

In order to get to the group $G_3^2$, we need to add the relations, so we should paste $2$-cells into the 1-complex $\Sigma_1$ --- those 2-cells will make some loops trivial, and hence produce the relations. The cells we are going to paste are hexagons. Each cell goes through exactly six out of eight vertices of the cube. For example, the relation $a_{12}a_{13}a_{23}a_{12}a_{13}a_{23}=1$ is represented by a loop starting from $(0,0,0)$ and going through the vertices $(1,0,0),(1,1,0),(1,1,1),(0,1,1),(0,0,1)$ and misses the vertices $(0,1,0)$ and $(1,0,1)$, see Fig.~\ref{fig:gnk_cube}. We paste four hexagons (one for each pair of opposite vertices). Therefore, every vertex is incident to exactly three hexagons. Fig.~\ref{fig:gnk_cube} illustrates how we paste one hexagon, the other three are pasted in the same manner.

Thanks to the $a_{ij}^2=1$ relations, all cyclic permutations of a loop yield the same relations, as well as the reversal of the loop. Hence when recording the hexagons we may begin with arbitrary vertex and choose any direction of walking around its boundary. Taking that into consideration, the four hexagons we add are of the following form:
$$(0,0,0)(1,0,0)(1,1,0)(1,1,1)(0,1,1)(0,0,1);$$
$$(0,0,0)(0,1,0)(0,1,1)(1,1,1)(1,0,1)(1,0,0);$$
$$(0,0,0)(0,0,1)(1,0,1)(1,1,1)(1,1,0)(0,1,0);$$
$$(1,0,0)(1,1,0)(0,1,0)(0,1,1)(0,0,1)(1,0,1).$$
For each of the first three hexagons, its boundary is a loop beginning at $(0,0,0)$, so they naturally give the following relations:
$$a_{12}a_{13}a_{23}a_{12}a_{13}a_{23}=1;$$
$$a_{13}a_{12}a_{23}a_{13}a_{12}a_{23}=1;$$
$$a_{23}a_{12}a_{13}a_{23}a_{12}a_{13}=1.$$
To get a loop from the last one, we can add a little ``tail'': add the letter $a_{12}$ to the beginning and the end of the word represented by the boundary of the last hexagon. We get the relation $a_{12}a_{13}a_{12}a_{23}a_{13}a_{12}a_{23}a_{12}=1$ which is equivalent to the relation
$$a_{13}a_{12}a_{23}a_{13}a_{12}a_{23}=1.$$
Analysing these relations we see that the second and the fourth ones coincide, the first one and the third one differ by a cyclic permutation (and hence are equivalent). Finally, the first and the second one differ by inversion and cyclic permutation, hence are equivalent as well. That means that all four hexagons give one and the same relation --- $a_{12}a_{13}a_{23}a_{12}a_{13}a_{23}=1,$ the tetrahedral relation of the group $G_3^2$, which yields all other tetrahedral relations in this group (see Example~\ref{ex:g32}).

In this manner we have constructed a $2$-complex which we denote by $\Sigma^2_{3}$.

\begin{statement}
Elements of the group $G_{3}^{2}$ are in bijection with homotopy types of paths in $\Sigma_{3}^{2}$ starting from $(0,0,0)$.
\end{statement}

Now, we are ready to do the same for the group $G_{n}^{k}$.

We begin with a ${n\choose k}$-dimensional unit cube; its coordinates will correspond to $k$-tuples of indices from $\{1,\dots, n\}$. In the same manner as before, we take its 1-skeleton and  paste $2k$-gons corresponding to relations of the tetrahedron type, but in the general case we also paste the quadrilaterals corresponding to far commutativity relations (note that when $n=k+1$, there is no far commutativity). We denote the resulting complex by $\Sigma_n^k$ (note, that it is obtained from a ${n\choose k}-$dimensional cube, not from an $n-$dimensional one).

\begin{statement}
Elements of the group $G_{n}^{k}$ are in bijection with homotopy types of paths in $\Sigma_{n}^{k}$ starting from $(0,0,0)$.
\end{statement}

%% file: manifolds_invariants.tex
In the main principle formulated in Section \ref{chap:gnk_definition}, we spoke about the dynamics of $n$ {\em particles} in a certain topological space.

In all examples considered in the previous sections, all particles were {\em points} in certain topological spaces which led us to various {\em configuration spaces}. In the present section, we shall deal with {\em submanifolds} instead of {\em points} playing roles of particles, hence, with {\em moduli spaces} of such submanifolds.

One of the goals of the present section is {\em to construct invariants of arbitrary manifolds}. The groups $G_{n}^{k}$ allow one to study {\em braids} of certain spaces, however, the moduli space for the $G_{n}^{k}$ method is quite restrictive. When general submanifolds are used, there are many more degrees of freedom.

We shall start with the example of projective duality: when we look at hypersurfaces in projective spaces, projective duality helps indicate the types of submanifolds that are required.

\subsection{Projective duality}

Consider the main examples from Sections \ref{chap:gnk_definition} and \ref{section:realisation} --- configurations of points in a linear or projective space.

Assume we deal with a collection of $n$ points in $\R{}P^{k-1}$, $2\le k$, such that any $k-1$ of them are in general position.
Our property of codimension $1$ will be:~
\begin{center}
{\em $k$ points are not generic.}
\end{center}

For example, when dealing with points in $\R{}P^{1}$, this means that some two points coincide and leads to a map from the configuration space to $G_{n}^{2}$, in the case of $\R{}P^{2}$ we deal with braids (loops in the configuration space of $n$ points in $\R{}P^2$), and we get a map to $G_{n}^{3}$~(cf. Proposition~\ref{prop_hom_PBn_Gn3}), for points in $\R{}P^{3}$ without collinear triples, we get a map to $G_{n}^{4}$, and so on.

Now, we just pass to the dual space. Namely, in the projective space of dimension $k-1$, a point $(z_{0}:\cdots: z_{k-1})$ corresponds to the hyperplane $\{(w_0:\cdots:w_{k-1})\in\R{}P^{k-1}\,|\, \sum_{i=0}^{k-1}w_iz_i=0\}$. Here we require that {\em any $k-1$ hyperplanes are generic}, that is, {\em the intersection of any $k-1$ planes is a unique point in $\R{}P^{k-1}$.} Note that a $(k-2)$-dimensional plane in ${\mathbb R}P^{k-1}$ corresponds to a $(k-1)$-dimensional vector subspace in ${\mathbb R}^k$.
Now we define our restricted moduli space $M'(\R{}P^{k-1},n)$ as the set $n$-tuples of hyperplanes such that any $k-1$ of them are generic.\label{text:moduli_space_condition} Certainly, the projective duality leads to a homeomorphism $i:M'(\R{}P^{k-1},n)\to C'(\R{}P^{k-1},n)$ which takes the hyperplane with coordinates $(z_{0}:\cdots: z_{k-1})$ to the point with the same coordinates $(z_{0}:\cdots: z_{k-1})$ in $\R{}P^{k-1}$.

Our {\em good $k$-condition} is:~

\begin{center}
{\em some $k$ hyperplanes share a point.}
\end{center}

For example, three lines in $\R{}P^{2}$ obey this $3$-condition if they have a common intersection point.

This condition is of codimension $1$, it precisely corresponds to the condition  ``{\em some $k$ points are not generic in $C'(\R{}P^{k-1},n)$}''.

Just as in Section~\ref{sec:realisation_n=k+1}, we take all good and transverse paths in $M'(\R{}P^{k-1})$ with fixed ends and take those paths to words in $G_{n}^{k}$.

In this case the isomorphism $i:M'(\R{}P^{k-1},n)\to C'(\R{}P^{k-1})$ allows one to follow Theorem \ref{thm:gnk_main_theorem}. Namely, we get

\begin{theorem}
Having fixed a generic point in $M'(\R{}P^{k-1})$, we have a well-defined map
from $\pi_1(M'(\R{}P^{k-1}), *)$ to elements of $G_{n}^{k}$.
\end{theorem}

If we analyse the case of $n=5$ and $k=4$, we see the exact meaning of the tetrahedron equation.

Namely, for $k=4$ we look for hyperplanes (i.e., $2$-dimensional planes in $\R{}P^{3}$), such that each three of them intersect in a point. In this case, the letter $a_{ijkl}$ corresponds exactly to a collection of four hyperplanes which pass through a point.


If we look at what happens in the neighbourhood of such a point, we shall see ``inversion of a tetrahedon''. 

If we slightly perturb it to get $\eps>0$ or $\eps<0$, we get two tetrahedra. 

Now, consider some fifth plane which is not parallel to any of the four above planes.
Then it passes through each of them in some order.

\begin{figure}
\centering\includegraphics[width=300pt]{ 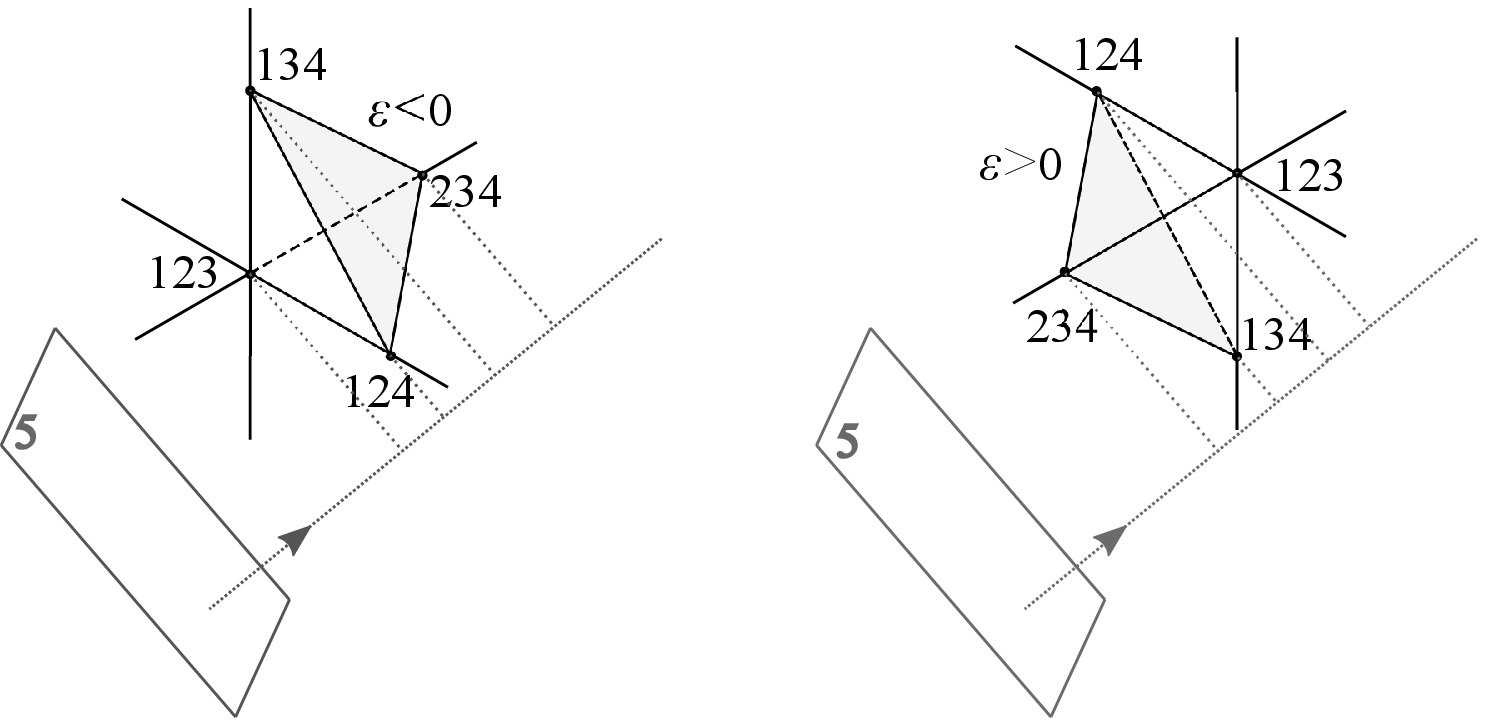}
\caption{Zamolodchikov relation}\label{fig:zamolodchikov}
\end{figure}

If in the picture $\eps<0$ we have the order $$(1,2,3),(1,2,4),(1,3,4),(2,3,4),$$
then in the picture $\eps>0$ we have the opposite order, see Fig~\ref{fig:zamolodchikov},
$$(2,3,4),(1,3,4),(1,2,4),(1,2,3).$$

This leads us exactly to the Zamolodchikov relation\index{Zamolodchikov relation} $$a_{1234}a_{1235}a_{1245}a_{1235}a_{2345}=a_{2345}a_{1345}a_{1245}a_{1235}a_{1234}.$$

\subsection{Embedded hypersurfaces}

The aim of the present section is to ``copy'' the above pattern to
the general case of a moduli space $M'(M,\Sigma_{1},\dots, \Sigma_{n})$,
where $M$ is a certain manifold of dimension $k-1$, and $\Sigma_{i}, i=1,\dots,n,$
are manifolds of dimension $k-2$. The elements of the moduli space are embeddings of $\Sigma_{i}$ in $M$.

Copying the $k-$good condition 
from the previous subsection (on page~\pageref{text:moduli_space_condition}), we require that under the above embedding, any $k-1$ submaninfolds $\Sigma_{i_1},\dots,
\Sigma_{i_{k-1}}$ should intersect transversely in exactly one point.

Then, using the main Theorem \ref{thm:gnk_main_theorem}, we get the following

\begin{theorem}
Let $M^{k-1}$ be a closed manifold of dimension $k-1$ and let
$\Sigma_{1},\cdots, \Sigma_{n}$ be a fixed collection of manifolds of dimension $k-2$.

Consider the moduli space $M'(M^{k-1},\Sigma_{1},\dots, \Sigma_{n})$ of embeddings in general position (in the sense above) of $\Sigma_{i}, i=1,\dots,n,$  in $M$  ($M'(M,\Sigma_{i})$, for short).
Then for each connected component $\#_{j}{\mathcal M}(M,\Sigma_{i})$ there is a well defined map from the fundamental group of this connected component to the group $G_{n}^{k}$.

Hence we get a collection of subgroups in $G_{n}^{k}$ which can serve as an invariant of the entire manifold $M$. \label{thm:manifolds_invariant_main}
\end{theorem}

\subsubsection{Examples}

If we deal with case when $M=S^{k-1}$ and all $\Sigma^{i}$ are homeomorphic
to $S^{i-1}$, we get precisely the subgroup in $G_{n}^{k}$ corresponding
to ``realisable'' elements, i.e., braids in the sense of~\cite{HigherGnk}.

\begin{example}
If $M=\R{}P^{k-1}, \Sigma_{1}=\cdots = \Sigma_{k}=\R{}P^{k-2}$, we get exactly
the image of the braid group of the projective space in  $G_{n}^{k}$.
\end{example}

\begin{example}
If $k=3$ (i.e., the dimension of the above total manifold is $2$), the above
construction works for any manifold except the sphere.
\end{example}

\begin{figure}
\begin{minipage}{.45\textwidth}
\centering\includegraphics[width=.45\textwidth]{ 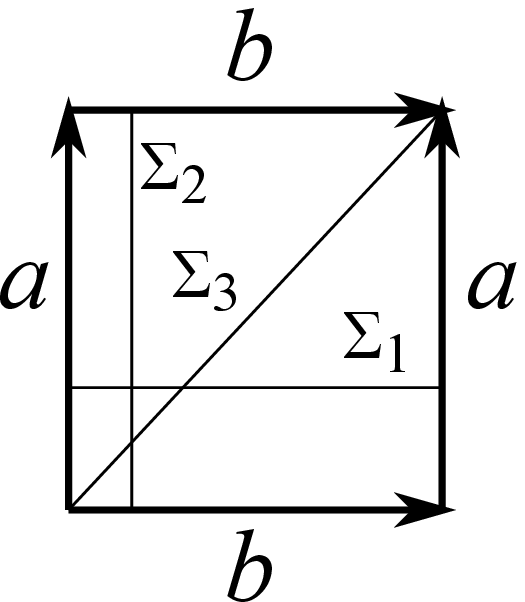}
\caption{Submanifolds in the torus}\label{fig:torus_gnk}
 \end{minipage}
 \hspace{5mm}
 \begin{minipage}{.45\textwidth}
 \centering\includegraphics[width=.45\textwidth]{ 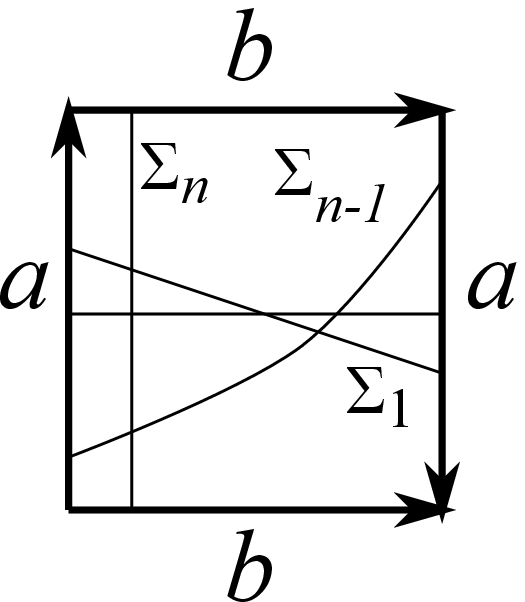}
\caption{Submanifolds in the Klein bottle}\label{fig:klein_gnk}
  \end{minipage}
\end{figure}

Namely, for the torus we have $M=T^2$ and $\Sigma_1=\dots=\Sigma_n=S^1$. Thanks to the gradient flow we can think of the submanifolds $\Sigma_i$ as simple geodesics. Let $(u_i,v_i)\in\Z^2$ be the coordinates of $\Sigma_i$ in $H_1(M,\Z)$ with the basis $\Sigma_1$, $\Sigma_2$, i.e. $\Sigma_i = u_i\Sigma_1+v_i\Sigma_2$ in $H_1(M,\Z)$. The condition that any two submanifolds $\Sigma_i$ have exactly one intersection point implies that their intersection index is $\pm 1$, hence, $|u_i|= 1, |v_i|=1$ for any $i\ge 3$. For the same reason any two submanifolds $\Sigma_i$, $\Sigma_j$ can not have the same homological coordinates. And we can not have simultaneously two submanifolds with coordinates $(\pm1,\pm1)$. Thus, $n\le 3$ and up to isomorphisms the homological classes of the submanifolds in the chosen basis are $\Sigma_1=(1,0), \Sigma_2=(0,1), \Sigma_3=(1,1)$, see Fig.~\ref{fig:torus_gnk}.

Let $n=3$. The fundamental group of the moduli space $(T^2,\Sigma_1,\Sigma_2,\Sigma_3)$ is generated by loops $\gamma_i,\ i=1,2,3,$ where $\gamma_i$ pushes $\Sigma_i$ in the orthogonal direction in the torus till the submanifold makes a round trip and returns to the initial position; the other submanifolds do not move. Each generator $\gamma_i, i=1,2,3,$ is mapped to the generator $a_{123}$ of the group $G^3_3$ because the moving submanifold $\Sigma_i$ passes once the intersection point of the other two.

\begin{example}

Consider the Klein bottle $K^2$ and submanifolds $\Sigma_{1},\cdots, \Sigma_{n}$ of dimension $1$ in it such that any two submanifolds intersect in exactly one point, and the intersection is transverse. We have the following situation:
\begin{enumerate}
\item the orientation-reverting cycles among the submanifolds $\Sigma_i$ must belong to the same homology class;
\item there can be at most one orientation-preserving cycle among $\Sigma_i$.
\end{enumerate}


Let $\Sigma_1,\dots,\Sigma_{n-1}$ are orientation-reversing cycles and $\Sigma_n$ are an orientation-preserving cycle, see Fig.~\ref{fig:klein_gnk}. Let $\gamma$ be the loop in the fundamental group of the moduli space given the double circular move of $\Sigma_n$ along an orientation-reversing cycle, e.g. $\Sigma_1$. We consider double move because the orientation of the $\Sigma_n$ must not change. The image of $\gamma$ in $G_n^3$ is equal to $\left(\prod_{k=1}^{(n-1)(n-2)/2}a_{i_kj_kn}\right)^2$ where $\Sigma_{i_k}\cap\Sigma_{j_k}$ is the $k$-th intersection point the submanifold $\Sigma_n$ encounters while moving along the orientation-reverting cycle. For example, for the configuration of the orientation-reverting cycles given in Fig.~\ref{fig:klein_gnk_conf}, the image of $\gamma$ equals to $\left(\prod_{i=1}^{n-2}\prod_{j=i+1}^{n-1}a_{ijn}\right)^2$. For other configurations of the orientation-reversing cycles the result will differ by conjugation.

\begin{figure}
\centering\includegraphics[width=100pt]{ 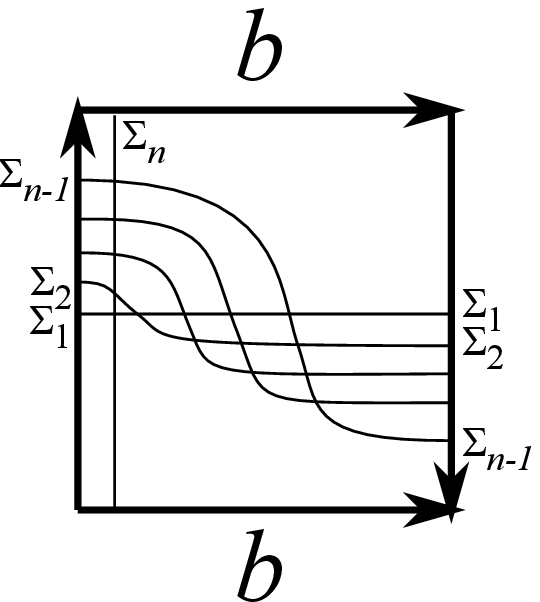}
\caption{A configuration of submanifolds in the Klein bottle}\label{fig:klein_gnk_conf}
\end{figure}
\end{example}

However, a careful reader may note that under some circumstances, the space $M'(M,\Sigma_i)$ is empty. This happens, for example, in the case when $M$ is a sphere. Then any two closed manifolds $\Sigma_{1},\Sigma_{2}$ in this sphere of codimension one have evenly many intersection which contradicts the main condition.

We shall tackle this case in the next subsection.

\subsection{Immersed hypersurfaces}

The most important case which is not covered by the main Theorem \ref{thm:manifolds_invariant_main}
is the case when the ambient manifold $M^{k}$ is the sphere. In this case,
the moduli space described above is {\em void}.

How to overcome this difficulty and to construct meaningful invariants?

Let us start with the $2$-sphere and a collection of $1$-circles in it.

When performing an isotopy, the number of intersection points between two
circles change, moreover, generically, {\em two circles never have exactly
one intersection point}. This means that we have to be more restrictive.




Note that in the case when the ambient manifold is the sphere of any
dimension $S^{n}$, and submanifolds are just spheres $S^{n-1}$ of one dimension
less, the above principle does not work at all: $H_{j}(S^{n})=0$ for
$j<n$, hence, in general position, $n$ spheres $S^{n-1}$ have zero
intersection, which means that generically the number of intersections
is {\em even}.

Below, we suggest two solutions to the above problem. The first one is to
replace the group $G_{n}^{3}$ with its abelinisation.


The other approach deals with the situation when we {\em immerse} several $S^{k-1}$
into a $S^{k}$ and for each $(k-1)$-tuple of $S^{k-1}$ we fix exactly one intersection
point and keep track of that one. We pay no attention to the other intersections.

We are sure that this approach (with minor or major modifications)
works for arbitrary manifolds and submanifolds, but here we restrict ourselves
with the case of $2$-dimensional manifolds and curves in them.

\subsection{Circles in $2$-manifolds and the group $G_{n}^{3}$}

Fix a positive integer $n$ and let us consider an oriented sphere $S^{2}$.
We shall consider the space of all embeddings of $n$ oriented $1$-dimensional circles
$S_{1},\cdots, S_{n}$ to $S^{2}$ such that every two circles $S_{i}$, $S_{j}$
have exactly two common points.
Denote this moduli space by $M_{imm}(S^{2},n)$.

We claim that there is a well-defined map from $\pi_{1}(M_{imm}(S^{2},n))$ to
the abelinisation of the group $G_{n}^{3}$.
Indeed, generically we do not have any triple intersections. When considering
a path in the moduli space, we meet codimension $1$ situations: three circles have an intersection.


As for the case of line on the plane, we have two configurations of codimension $2$:
\begin{enumerate}
\item two triple points (note that points may belong to the same triple of circles);
\item a quadruple point.
\end{enumerate}

\begin{figure}
\centering\includegraphics[width=100pt]{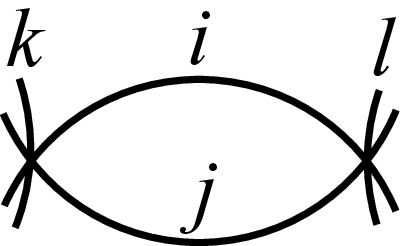}
\caption{A configuration of codimension $2$ that gives relation $a_{ikj}a_{ijl}=a_{ijl}a_{ikj}$. Note that the generators have two common indices.}
\label{fig:gnk_double_commutativity_relation}
\end{figure}

The first case corresponds to commutativity relations $a_{ijk}a_{i'j'k'}=a_{i'j'k'}a_{ijk}$ for any triples of indices, see for example Fig.~\ref{fig:gnk_double_commutativity_relation}. The second case gives a tetrahedron relation

$$a_{ijk}a_{ijl}a_{ikl}a_{jkl}=a_{jkl}a_{ikl}a_{ijl}a_{ijk}.$$

These relations
with the relation $a_{ijk}^2=1$
define a presentation for the group $G=G_{n}^{3}/[G_{n}^{3},G_{n}^{3}]\simeq\Z_2^{\oplus\frac 16 n(n-1)(n-2)}$. Thus, we get a map from $\pi_{1}(M_{imm}(S^{2},n))$ to $G$. It is easy to see the image of the map is trivial.

\begin{figure}
\centering\includegraphics[width=200pt]{ 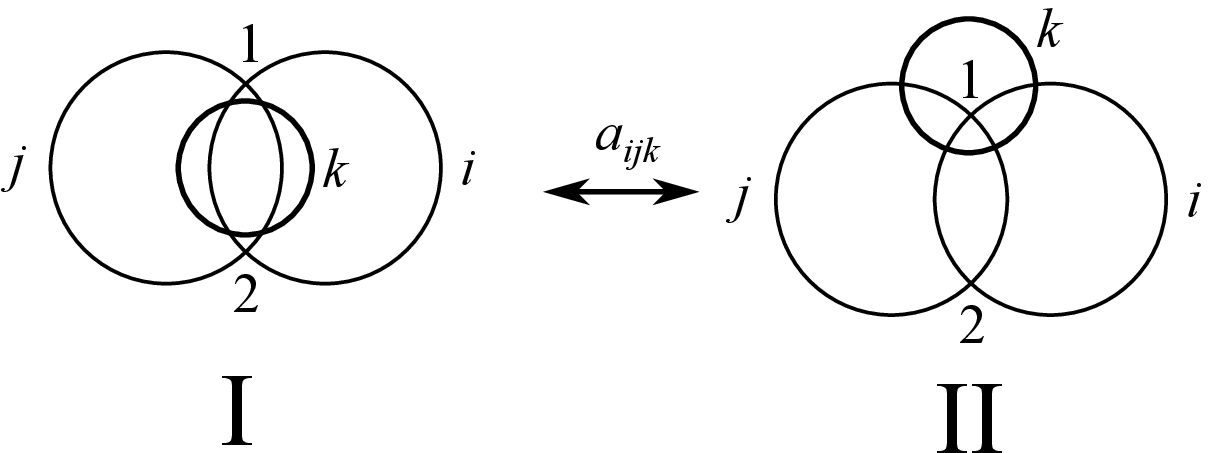}
\caption{Two types of configurations of circles $S_i, S_j, S_k$.}
\label{fig:circle_ijk_configurations}
\end{figure}

Indeed, given indices $i,j,k$, there are two types of possible configurations of the circles $S_i, S_j, S_k$ on the sphere, see Fig.~\ref{fig:circle_ijk_configurations}. In the first configuration any circle separates the intersection points of the other two, in the configuration circles do not separate the intersection points. When one switches from one configuration to the other the generator $a_{ijk}$ appears. Since for a loop in the configuration space the initial and final configurations coincide, the number of the generators $a_{ijk}$ in the word corresponding to the loop will be even.

Since the configurations of codimension $1$ and $2$ are local, the result can be generalised to the moduli space consisting of embeddings of circles into a surface $\Sigma$ such that any two circles have two intersection points in $\Sigma$.


\subsection{Immersed curves in $M^{2}$}

In this section we shall again restrict ourselves to the two-dimensional case,
namely, we consider the moduli space of immersed curves in $M^{2}$.
The situation with immersed hypersurfaces will cause much more difficulties.

Fix a $2$-surface $\Sigma$ and a natural number $n$; elements of our moduli space
are $n$ closed curves $\gamma_{1},\cdots, \gamma_{n}$ smoothly
immersed in $\Sigma$ together with a choice of one transverse intersection
for each pair of curves: $P_{ij}\in \gamma_{i}\cap \gamma_{j}$ (the points may coincide). The total
number of intersection between any two curves may be an arbitrary positive number.
We denote this space by $P(\Sigma,n)$.

\begin{theorem}
There is a well defined map from $\pi_{1}(P(\Sigma,n))$ to $G_{n}^{3}$.
\end{theorem}

\begin{proof}
We are looking for those moments when fixed intersection points for some triple of points coincide: $P_{ij}=P_{ik}=P_{ik}$ and associate the generator $a_{ijk}$ of the group $G_{n}^{3}$ to that. We pay no attention to the other intersections.

The rest of the proof follows the proof of the main Theorem~\ref{thm:gnk_main_theorem} from Section~\ref{chap:gnk_definition}.
\end{proof}


\begin{remark}
In the above setup we can drop the manifold $\Sigma$ and deal with curves
by themselves. In order to fix transverse intersection, we should deal with diagrams of
{\em $n$-component free links}, such that every two components have exactly
one intersection. We see that any equivalence of diagrams in this class, presented by a sequence of Reidemeister moves, produces a word in $G_{n}^{3}$ (Fig.~\ref{fig:free_knot_transformation}).

\begin{figure}
\centering\includegraphics[width=355pt]{ 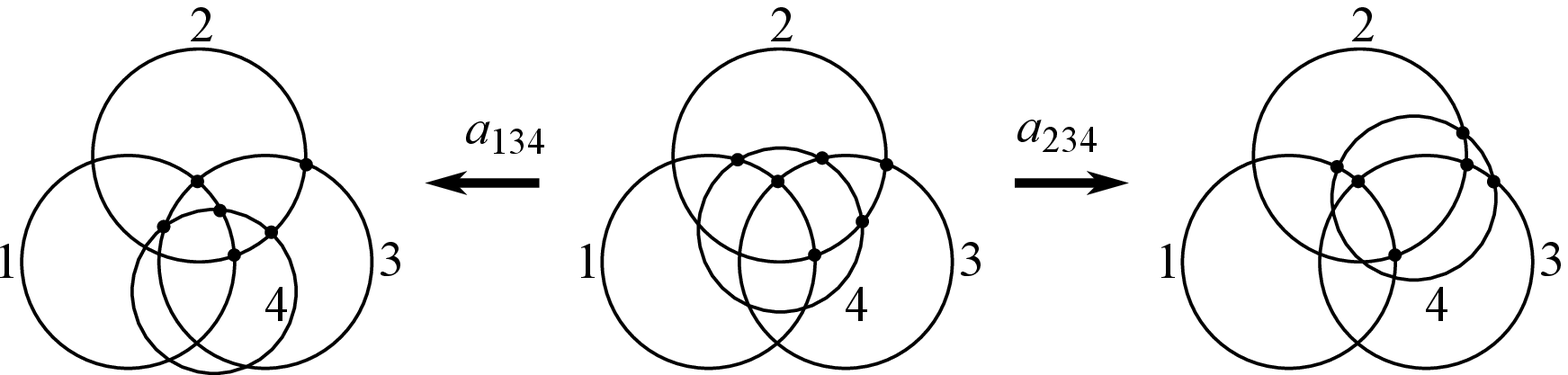}
\caption{Transformations of free links and the corresponding words in $G_n^3$. The crossings of the link are marked with dots.}
\label{fig:free_knot_transformation}
\end{figure}
\end{remark}

%% file: gamma_n4.tex
In Part~\ref{chap:gnk_definition} we defined a family of groups $G_{n}^{k}$ for two positive integers $n>k$, and formulated the following principle:

\begin{center}
{\em If dynamical systems describing a motion of n particles admit some nice codimension one property governed by exactly $k$ particles, then these dynamical systems have a topological invariant valued in $G_{n}^{k}$. }
\end{center}

The main examples coming from $G_{n}^{k}$-theory are homeomorphisms from the $n$-strand pure braid group to the groups $G_{n}^{3}$ and $G_{n}^{4}$, introduced in Section~\ref{chap:gnk_homomorphisms}. If we consider a motion of $n$ pairwise distinct points on the plane and choose the property ``some three points are collinear'', then we get a homomorphism from the pure $n$-strand braid group $PB_{n}$ to the group $G_{n}^{3}$. If we choose the property ``some four points belong to the same circle or line'', we shall get a group homomorphism from $PB_{n}$ to $G_{n}^{4}$.

In other words, in our examples we look for ``walls'' in the configuration space $C_{n}(\mathbb{R}^{2})$ where some three points are collinear (or some four points are on the same circle or line). This condition can be well defined for $C_{n}(\mathbb{R}P^{2})$.

But what would be the ``nice'' codimension one property if we try to study similar configuration spaces or braids for some other topological spaces? First, our conditions will heavily depend on the metric: the property ``three points are collinear'' is metrical. On the other hand, even having some good metric chosen, we meet other obstacles because we need to know what is a ``line''. For example, there is no unique geodesic passing through two points in the general case. For example, if you consider the property ``three points belong to one geodesic'' on a flat torus, fix two different points on the torus and try to record as a $G_n^k$-word all critical moments for a simple movement of a third point, the word would be infinite because the union of geodesics connecting the fixed two points is dense.

The detour for this problem will be as follows: we shall consider the condition ``locally'' and instead of ``global configurations''. We shall consider only Vorono\"{i} tilings or Delaunay triangulations~\cite{AKL}. From this point of view, we deal with the space of triangulations with a fixed number of vertices, where any two adjacent triangulations are related by a Pachner move~\cite{Nabutovsky}, which is closely related to the group $\Gamma_{n}^{4}$ (see Definition~\ref{dfn_Gamma}).

On the other hand, triangulations of spaces and Pachner moves are also related to Yang-Baxter maps (see \cite{Dynnikov}). Moreover, in~\cite{ChoZichertYun, HikamiInoue} a boundary-parabolic $PSL(2,\mathbb{C})$-representation of $\pi_{1}(S^{3} \backslash K)$ for a hyperbolic knot $K$ is studied by using cluster algebras and {\em flips} -- Pachner moves for 2-dimensional triangulations.

As we noticed in the previous paragraph, the group $\Gamma_{n}^{4}$ is closely related to triangulations of spaces and Pachner moves. In other words, we can consider triangulations of spaces and Pachner moves not only geometrically, but also algebraically. It can be expected to obtain invariants by means of the group $\Gamma_{n}^{4}$ not only for braids, but also for knots, which are obtained by closing braids.

\subsection{A group homomorphism from $PB_{n}$ to $\Gamma_{n}^{4}$}

\begin{definition}\label{dfn_Gamma}
The group $\Gamma_{n}^{4}$\index{Group!$\Gamma_n^4$} is the group generated by $\{ d_{(ijkl)}~|~ \{i,j,k,l\} \subset \bar{n}, |\{i,j,k,l\}| = 4\}$ subject to the following relations:

\begin{enumerate}
\item $d_{(ijkl)}^{2} = 1$ for $(i,j,k,l) \subset \bar{n}$,
\item $d_{(ijkl)}d_{(stuv)} = d_{(stuv)}d_{(ijkl)}$, for $| \{i,j,k,l\} \cap \{s,t,u,v\} | < 3$,
\item $d_{(ijkl)}d_{(ijlm)}d_{(jklm)}d_{(ijkm)}d_{(iklm)} = 1$ for distinct $i,j,k,l,m$.
\item $d_{(ijkl)}=d_{(kjil)}=d_{(ilkj)}=d_{(klij)}=d_{(jkli)}=d_{(jilk)}=d_{(lkji)}=d_{(lijk)}$ for distinct $i,j,k,l$.
\end{enumerate}
The equality among these 8 generators is given by the dihedral group acting on the indices.
\end{definition}

\begin{figure}[h!]
 \centering
\includegraphics[width = 8cm]{ 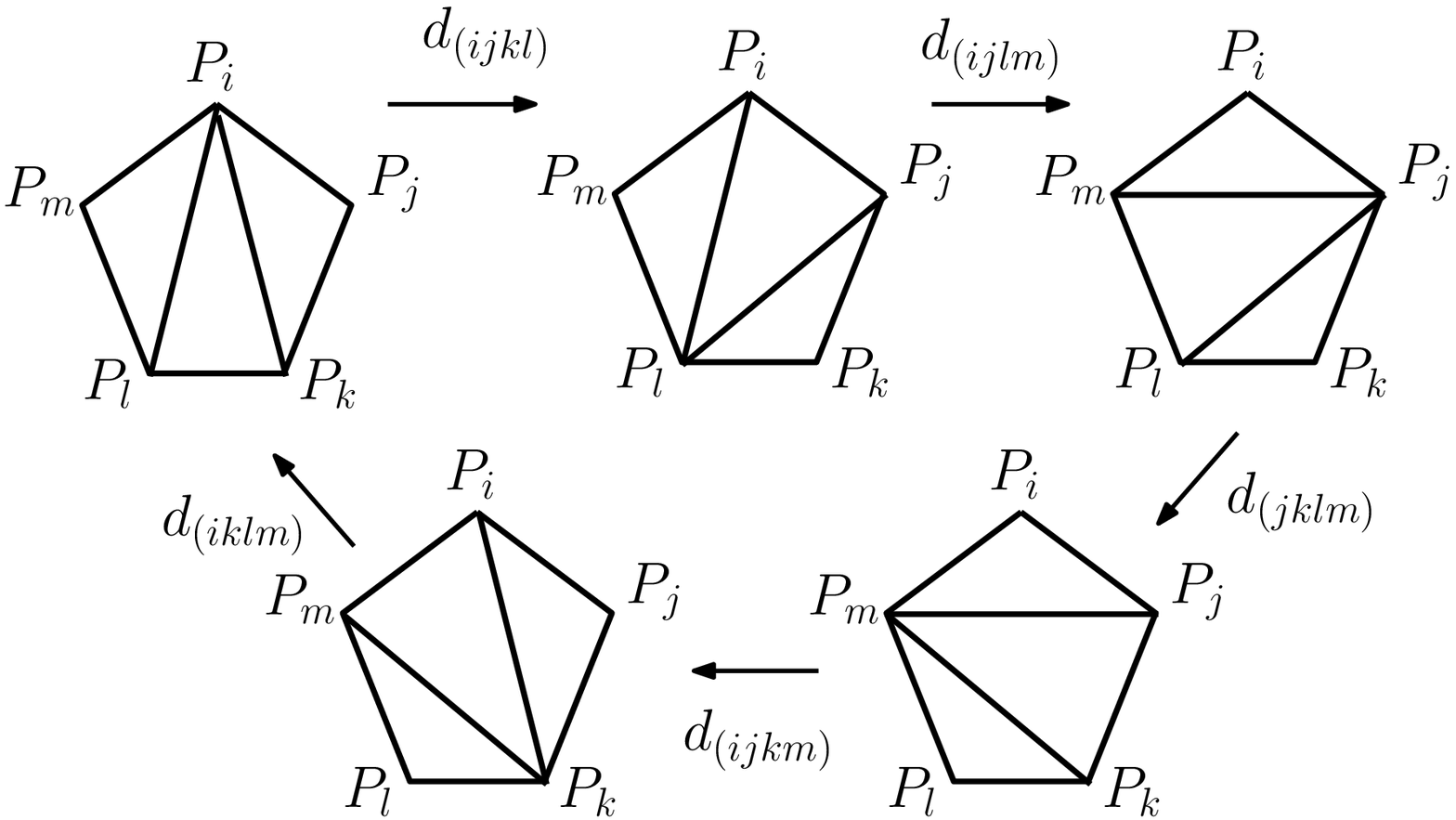}
 \caption{A sequence of five subsets of $\{P_{i},P_{j},P_{k},P_{l},P_{m}\}$, corresponding to flips on the pentagon}\label{proof_rel3}
\end{figure}

The group $\Gamma_{n}^{4}$ is naturally related to triangulations of $2$-surfaces and the Pachner moves for the two dimensional case, called ``flip'', see Fig.~\ref{fig:flip}. More precisely, a generator $d_{(ijkl)}$ of $\Gamma_{n}^{4}$ corresponds to a flip in the sequence constituting the Pentagon relation as described as Fig.~\ref{proof_rel3}, the most important relation for the group~$\Gamma_{n}^{4}$.

\begin{figure}[h!]
\label{fig:flip}
 \centering
 \includegraphics[width = 6.5cm]{ 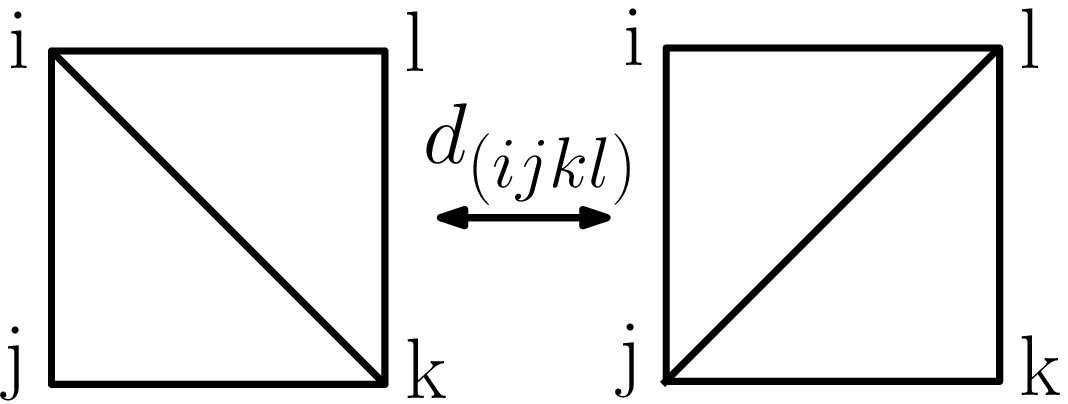}
 \caption{Flip on a rectangle $ijkl$}\label{fig:flip}
\end{figure}

More specifically, the relation $d_{(ijkl)}d_{(ijlm)}d_{(jklm)}d_{(ijkm)}d_{(iklm)} = 1$ corresponds to the sequence of flips applied to a pentagon, see Fig.~\ref{proof_rel3}.

Below we construct a group homomorphism $f_{n}$ from $PB_{n}$ to $\Gamma_{n}^{4}$. The topological background for that is very easy: we consider codimension 1 ``walls'' which correspond to generators (flips) and codimension 2 singularities correspond to relations (of the group $\Gamma$). Having this, we construct a map on the level of generators.

\begin{figure}[h!]
 \centering
\includegraphics[width = 0.4\textwidth]{ 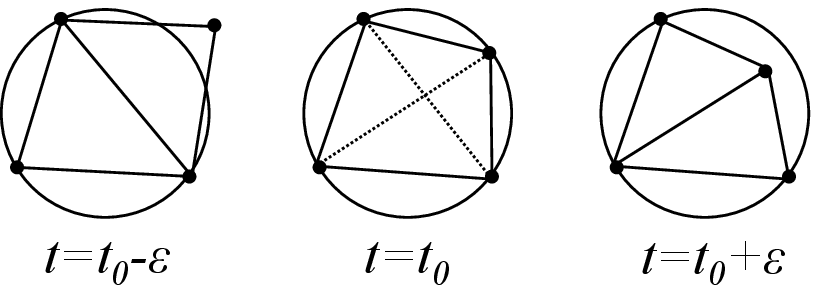}
 \caption{Flip of Delaunay triangulation at some critical moment $t_0$}\label{fig:pachner_circle}
\end{figure}

The ``walls'' in our case are built on the property ``four points on one circle which has no other points inside itself''. Indeed, given a braid it determines a dynamical system of $n$ points $P_1(t),\dots, P_n(t), t\in[0,1],$ in the plane. We look at the dynamics of the Delaunay triangularion spanned on the point set $P_i(t)$ and mark the moments when its combinatorics changes. It is when a flip occurs. Geometrically, it corresponds to a configuration with four points on one circle that contains no points inside (because Delaunay triangles are those whose circumscribed circle has no point inside), see Fig.~\ref{fig:pachner_circle}.

\begin{figure}[h]
 \centering
\includegraphics[width = 0.6\textwidth]{ 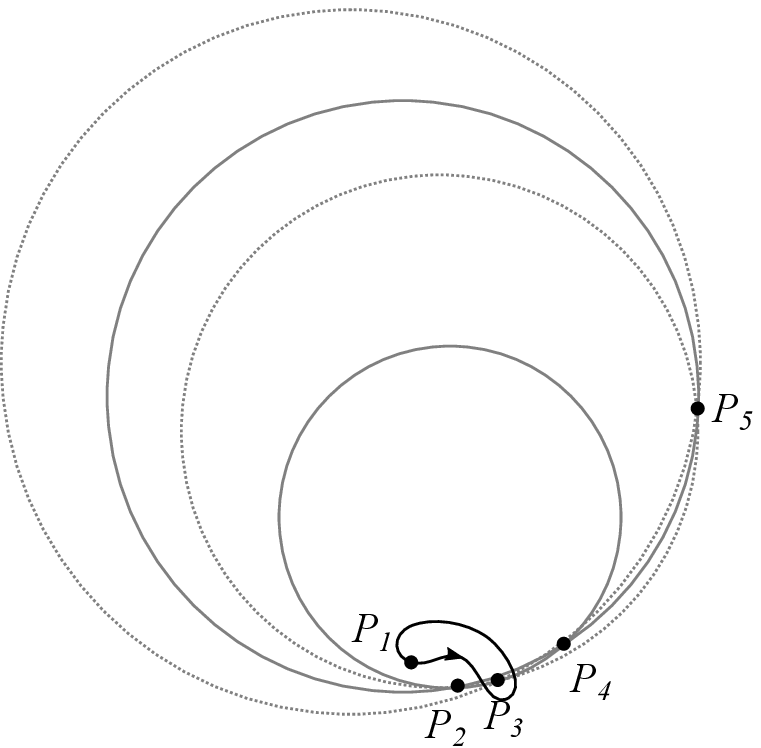}
 \caption{Dynamical system corresponding to the generator $b_{13}$}\label{fig:generator_bl3}
\end{figure}

{\bf Geometric description of the mapping from $PB_{n}$ to $\Gamma_{n}^{4}$.}~Let us consider $b_{ij}$ as the dynamical system described in Section~\ref{sect:hom_G4n}.
The homomorphism $f_{n}$ from $PB_{n}$ to $\Gamma_{n}^{4}$ can be defined as follows; for the above dynamical system for each generator $b_{ij}$, let us enumerate the moments $0<t_{1}<t_{2}< \cdots < t_{l} <1$ such that at the moment $t_{k}$ four points belong to one circle without points inside the circle. At the moment $t_{k}$, if $P_{s}, P_{t}, P_{u}, P_{v}$ are positioned on one circle as the indicated order and there are no points inside the circle, then set $\delta_{k} = d_{(stuv)}$. With the pure braid $b_{ij}$ we associate the product $f_{n}(b_{ij}) = \delta_{1}\delta_{2}\cdots \delta_{l}$.
As an example, consider the generator $b_{13}$ in $PB_5$. The corresponding dynamical system is shown in Fig.~\ref{fig:generator_bl3}: the point $P_1$ moves around $P_3$. While moving, $P_1$ goes through the circles containing other points,and mark these moments with letters $d$ indiced with the numbers of the incident points in the same cyclic order. We ignore the dotted circles $P_2P_3P_5$ and $P_3P_4P_5$ because they contains other points inside ($P_4$ and $P_2$ respectively). Thus, we get the word $d_{2145}d_{2134}d_{2314}d_{2145}$ which is the image of $b_{13}$ in the group $\Gamma_5^4$.

\begin{remark}\label{nbr_inside_points}
From the construction of base points $\{P_{1},\cdots, P_{n}\}$ it follows points $\{P_{1}, \cdots, P_{j-1}\} \cup \{P_{p+1}, \cdots P_{q-1}\}$ are placed inside the circle passing through points $P_{j},P_{p},P_{q}$ for $j<p<q$. That is, inside the circle passing through $P_{j},P_{p},P_{q}$ for $j<p<q$ there are $j-1+q-p-1$ points from $\{P_{1},\cdots, P_{n}\}$.
\end{remark}

{\bf Algebraic description of the mapping from $PB_{n}$ to $\Gamma_{n}^{4}$.}~Let us describe the mapping $f_{n} : PB_{n} \rightarrow \Gamma_{n}^{4}$ as follows:
Let us denote
\begin{center}
$d_{\{p,q, (r,s)_{s}\}}  = \left\{
\begin{array}{cc} 
     d_{(pqrs)} & \text{if}~p<q<s, \\
      d_{(prsq)} & \text{if}~p<s<q, \\
  d_{(rspq)} & \text{if}~s<p<q,\\
  d_{(qprs)} & \text{if}~q<p<s, \\
      d_{(qrsp)} & \text{if}~q<s<p, \\
  d_{(rsqp)} & \text{if}~s<q<p.
   \end{array}\right.$
   \end{center}
 \begin{remark}
Note that the generator $d_{\{p,q, (r,s)_{s}\}}$ corresponds to four points $P_{p},P_{q},P_{r},P_{s}$ such that they are placed on a circle according to the order of $p,q,s$ and the point $P_{r}$ is placed close to $P_{s}$ for the orientation $P_{r}$ to $P_{s}$ to be the counterclockwise orientation, see Fig.~\ref{exa_pts_circle}. The subscript $s$ of $(r,s)_{s}$ means that the point $P_{s}$ does not move, but the point $P_{r}$ moves turning around the point $P_{s}$ after this moment. In other words, the generator $d_{\{p,q, (r,s)_{s}\}}$ corresponds to the moment when the point $P_{r}$ is ``moving'' closely to the point $P_{s}$, turning around $P_{s}$. We would like to highlight that $d_{\{p,q, (r,s)_{s}\}} \neq d_{\{p,q, (s,r)_{s}\}}$.

\begin{figure}[h!]
 \centering
\includegraphics[width = 8cm]{ 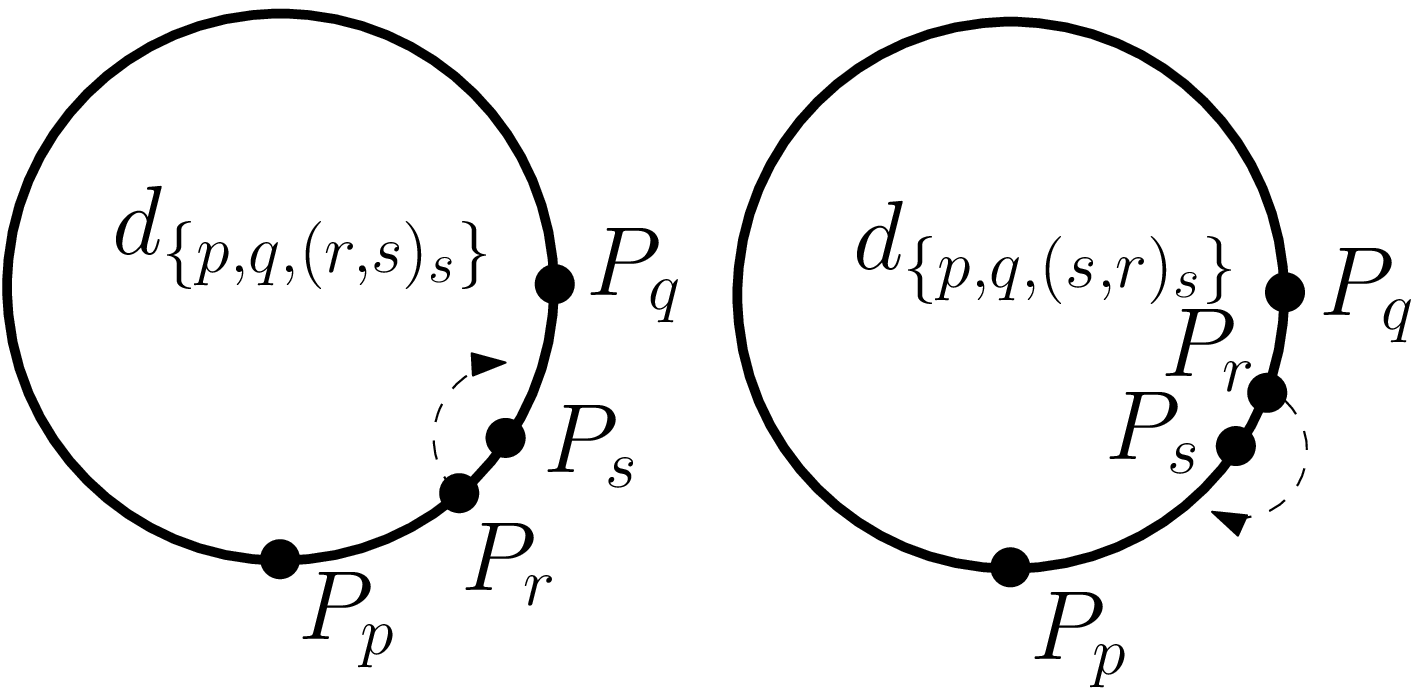}
 \caption{For $p<s<q$, $d_{\{p,q, (r,s)_{s}\}} = d_{(prsq)}$, but $d_{\{p,q, (s,r)_{s}\}} = d_{(psrq)}$.}\label{exa_pts_circle}
\end{figure}
\end{remark}

Let us denote $k \in \{p,q,r\}$ such that $min\{p,q,r\}< k <max\{p,q,r\}$ by $mid\{p,q,r\}$. Let us define $\gamma_{\{p,q,(i,j)_{j}\}}$ as follows: If $min\{p,q,j\}<i<mid\{p,q,j\}$ or $i>max\{p,q,j\}$, then

\begin{center}
\begin{equation}\label{for:gen_Gamma1}
\gamma_{\{p,q,(i,j)_{j}\}}= \left\{
\begin{array}{cc}
     d_{\{p,q,(i,j)_{j}\}}, & \text{if}~min\{p,q,j\}-mid\{p,q,j\}+max\{p,q,j\}-2=0, \\
  1,
 &   \text{if}~min\{p,q,j\}-mid\{p,q,j\}+max\{p,q,j\}-2 \neq 0. \\
   \end{array}\right.
   \end{equation}
\end{center}

If $i<min\{p,q,j\}$ or $mid\{p,q,j\}<i<max\{p,q,j\}$, then

\begin{center}
\begin{equation}\label{for:gen_Gamma2}
\gamma_{\{p,q,(i,j)_{j}\}}= \left\{
\begin{array}{cc}
     d_{\{p,q,(i,j)_{j}\}}, & \text{if}~min\{p,q,j\}-mid\{p,q,j\}+max\{p,q,j\}-2=1, \\
  1,
 &   \text{if}~min\{p,q,j\}-mid\{p,q,j\}+max\{p,q,j\}-2 \neq 1. \\
   \end{array}\right.
   \end{equation}
\end{center}

As asserted in Remark~\ref{nbr_inside_points} $min\{p,q,j\}-mid\{p,q,j\}+max\{p,q,j\}-2$ is the number of points inside the circle passing points $P_{j},P_{p},P_{q}$. Formulas (\ref{for:gen_Gamma1}) and (\ref{for:gen_Gamma2}) are describing that if four points $P_{i},P_{j},P_{p},P_{q}$ belong to one circle without points inside the circle, then we write a generator $d_{\{p,q,(i,j)_{j}\}}$, otherwise, the identity $1$.

Let $b_{ij} \in PB_{n}$, $1 \leq i<j \leq n$, be a generator. Consider the elements
\begin{eqnarray*}
\Delta_{i,(i,j)}^{I} = \prod_{p=2}^{j-1}\prod_{q=1}^{p-1}\gamma_{\{p,q,(i,j)_{j}\}},\\
\Delta_{i,(i,j)}^{II} = \prod_{p=1}^{j-1}\prod_{q=1}^{n-j}\gamma_{\{(j-p),(j+q),(i,j)_{j}\}},\\
\Delta_{i,(i,j)}^{III} = \prod_{p=1}^{n-j-1}\prod_{q=0}^{p-1}\gamma_{\{(n-p),(n-q),(i,j)_{j}\}},\\
\Delta_{i,(i,j)}=\Delta_{i,(i,j)}^{II}\Delta_{i,(i,j)}^{I}\Delta_{i,(i,j)}^{III}.
\end{eqnarray*}
Now we define $f_{n} : PB_{n} \rightarrow \Gamma_{n}^{4}$ by
$$f_{n}(b_{ij}) = \Delta_{i,(i,(i+1))}\cdots \Delta_{i,(i,(j-1))}\Delta_{i,(i,j)}\Delta_{i,(j,i)}\Delta_{i,((j-1),i)}^{-1}\cdots \Delta^{-1}_{i,((i+1),i)},$$
for $1 \leq i<j \leq n$.

For the example in Fig.~\ref{fig:generator_bl3} we have $f_5(b_{13})=\Delta_{1,(1,2)}\Delta_{1,(1,3)}\Delta_{1,(3,1)}\Delta_{1,(2,1)}^{-1}$ where
$\Delta_{1,(1,2)}=d_{2145}d_{2134}$, $\Delta_{1,(1,3)}=d_{2314}$, $\Delta_{1,(3,1)}=d_{2134}$ and $\Delta_{1,(2,1)}=d_{2145}d_{2134}$, so
$f_5(b_{13})=d_{2145}d_{2134}d_{2314}d_{2145}$ as we have seen.

\begin{theorem}\label{thm_to_gamma}
The map $f_{n} : PB_{n} \rightarrow \Gamma_{n}^{4}$, which is defined above, is a well defined homomorphism.
\end{theorem}

\begin{proof}
When we consider isotopies between two pure braids, it suffices to take into account only singularities of codimension at most two. Singularities of codimension one give rise to generators, and relations come from singularities of codimension two. Now we list the cases of singularities of codimension two explicitly.

\begin{enumerate}
\item One point moving on the plane is tangent to the circle, which passes through three points, see Fig.~\ref{proof_rel1}. This corresponds to the relation $d_{(ijkl)}^{2} =1$.

\begin{figure}[h!]
 \centering
\includegraphics[width = 8cm]{ 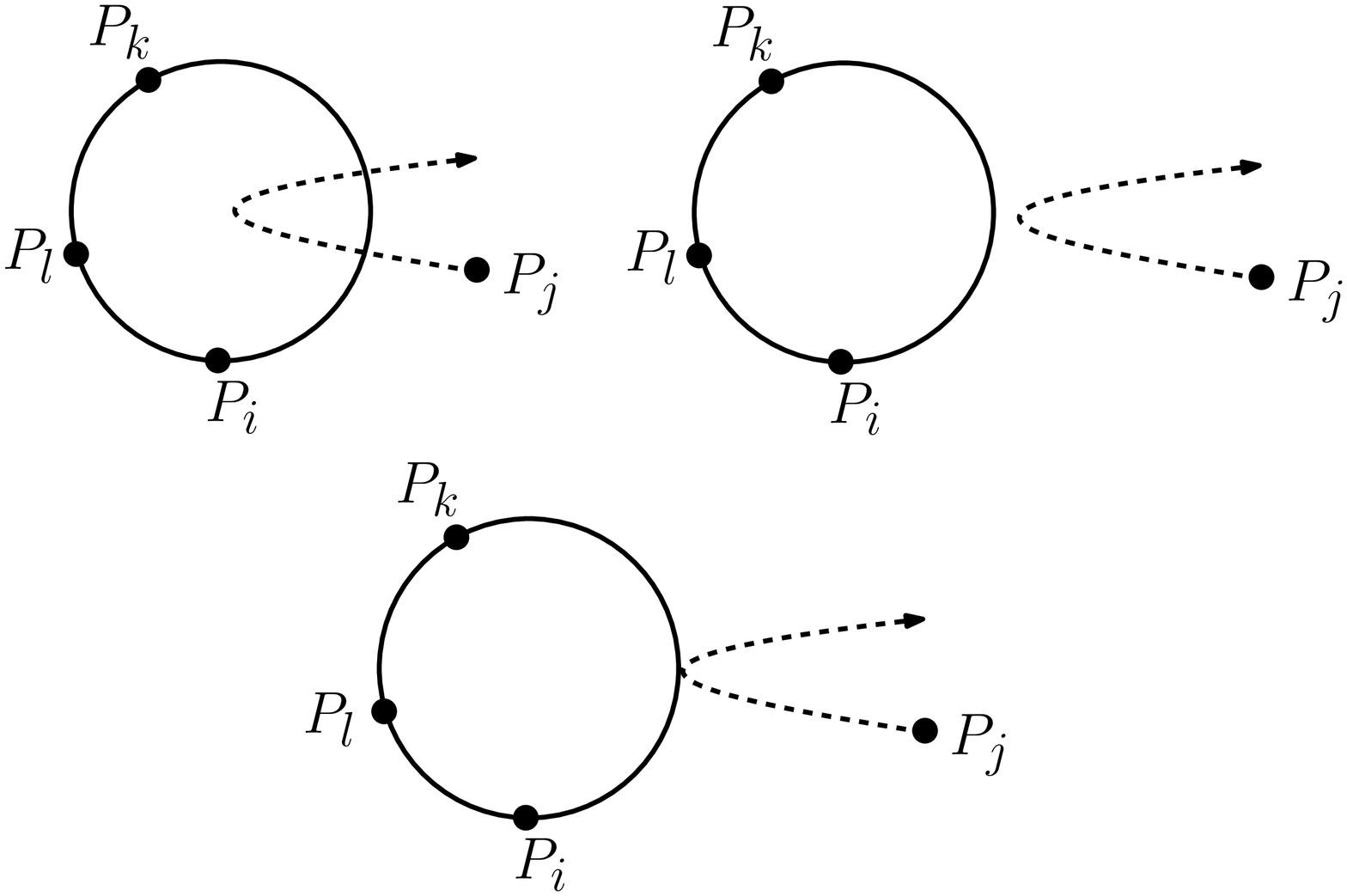}
 \caption{A point $P_{j}$ moves, being tangent to the circle, which passes through $P_{i},P_{k}$ and $P_{l}$}\label{proof_rel1}
\end{figure}

\item There are two sets $A$ and $B$ of four points each, which are on the same circles such that $|A\cap B| \leq 2$, see~Fig.~\ref{proof_rel2}. This corresponds to the relation $d_{(ijkl)}d_{(stuv)}=d_{(stuv)}d_{(ijkl)}$.

\begin{figure}[h!]
 \centering
\includegraphics[width = 5cm]{ 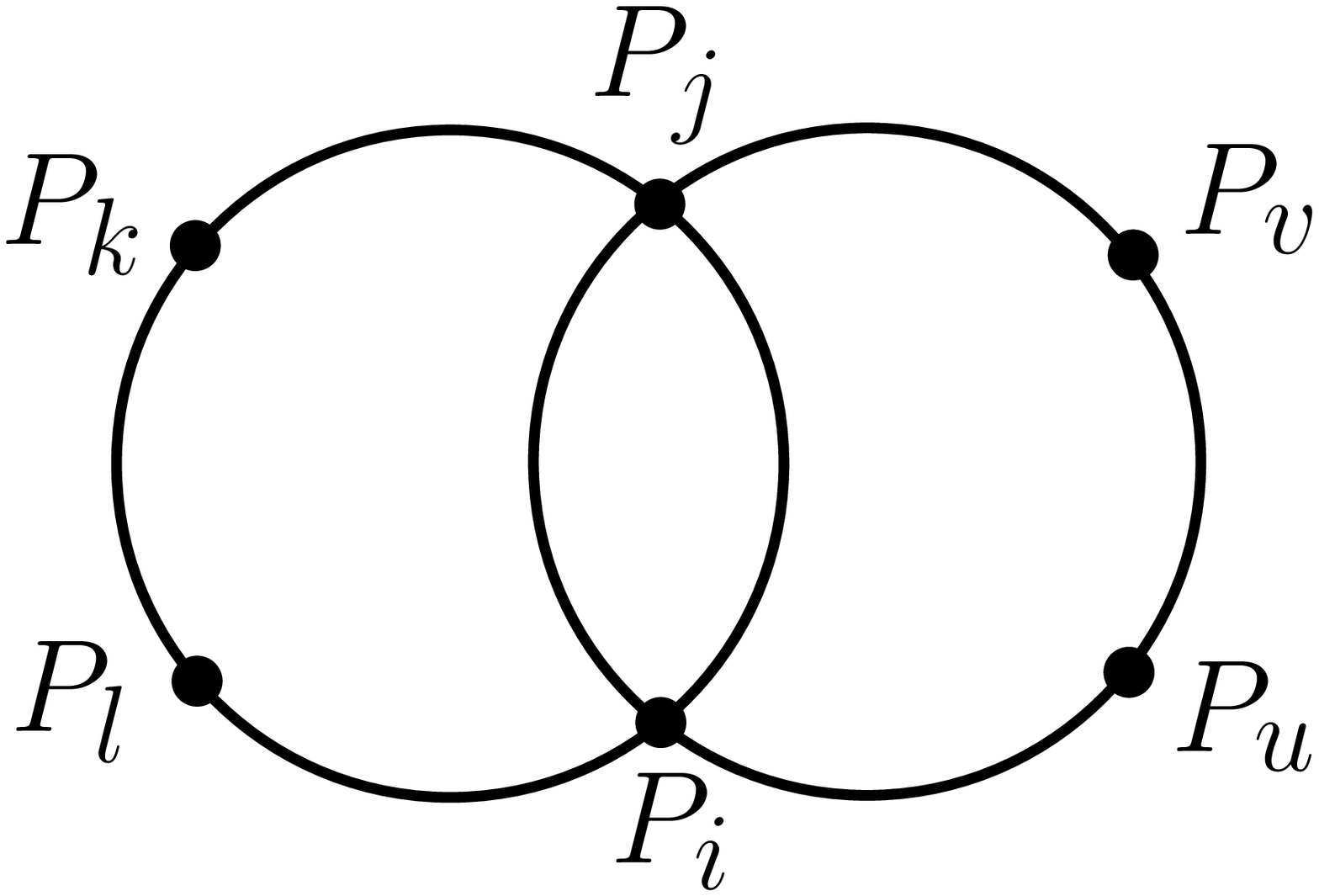}
 \caption{Two sets $A$ and $B$ of four points on the circles such that $|A\cap B| =2$}\label{proof_rel2}
\end{figure}

\item There are five points $\{P_{i},P_{j},P_{k},P_{l},P_{m}\}$ on the same circle. We obtain the sequence of five subsets of $\{P_{i},P_{j},P_{k},P_{l},P_{m}\}$ with four points on the same circle, which corresponds to the flips on the pentagon, see Fig.~\ref{proof_rel3}. This corresponds to the relation $d_{(ijkl)}d_{(ijlm)}d_{(jklm)}d_{(ijkm)}d_{(iklm)}=1$.

\end{enumerate}

\end{proof}

\subsection{A group homomorphism from $PB_{n}$ to $\Gamma_{n}^{4} \times \Gamma_{n}^{4} \times \dots \times \Gamma_{n}^{4}$}\label{sec:braids_product_Gamma}

For a positive integer $n>5$, the homomorphism $f_{n}\colon PB_{n} \rightarrow \Gamma_{n}^{4}$ can be extended to a mapping from $PB_{n}$ to the product of $([\frac{n-4}{2}]+1)$ copies of $\Gamma_{n}^{4}$ as follows (here and in the sequel $[x]$ means that $x$ is rounded down to the nearest integer). We denote the product of $(s+1)$ copies of the group $\Gamma_n^4$ by
$$\Gamma_{n,0}^{4} \times \Gamma_{n,1}^{4} \times \cdots \times \Gamma_{n,s}^{4}.$$
The idea is that we count this number of points modulo $n-4$.

Let $r=n-4$. For the dynamical system for a generator $b_{ij}$, which is described in Section~\ref{sect:hom_G4n}, let us enumerate $0<t_{1}<t_{2}< \cdots < t_{l} <1$ such that at the moment $t_{k}$ four points $P_{s}, P_{t}, P_{u}, P_{v}$ are positioned on one circle or one line.


If $P_{s}, P_{t}, P_{u}, P_{v}$ are positioned on one line $L$ or one circle $C$ at the moment $t_{k}$ in the indicated order, then there are two connected components of $\mathbb{R}^{2} \backslash L$ (or $\mathbb{R}^{2} \backslash C$). If the numbers of points contained in connected components of $\mathbb{R}^{2} \backslash L$ (or $\mathbb{R}^{2} \backslash C$) are $\alpha$ and $n-4-\alpha$ mod $r$ for some $\alpha \leq [\frac{r}{2}]$, then $t_{k}$ corresponds to $\delta_{k} = (1,\dots ,1,d_{(stuv)}, 1,\dots, 1) \in  \{1\} \times \cdots \{1\} \times \Gamma_{n,\alpha}^{4} \times \{1\} \times \cdots \{1\} \subset   \Gamma_{n,0}^{4} \times \Gamma_{n,1}^{4} \times \cdots \times \Gamma_{n,[\frac{r}{2}]}^{4}$.

With the pure braid $b_{ij}$ we associate the product $f^{r}_{n}(b_{ij}) = \delta_{1}\delta_{2}\cdots \delta_{l}$.

Algebraically this construction can be presented as follows:

\begin{eqnarray*}
D_{i,(i,j)}^{I} = \prod_{p=2}^{j-1}\prod_{q=1}^{p-1}\delta^{r}_{\{p,q,(i,j)_{j}\}},\\
D_{i,(i,j)}^{II} = \prod_{p=1}^{j-1}\prod_{q=1}^{n-j}\delta^{r}_{\{(j-p),(j+q),(i,j)_{j}\}},\\
D_{i,(i,j)}^{III} = \prod_{p=1}^{n-j+1}\prod_{q=0}^{n-p+1}\delta^{r}_{\{(n-p),(n-q),(i,j)_{j}\}},\\
\end{eqnarray*}
where
\begin{eqnarray*}
\delta^{r}_{\{p,q,(i,j)_{j}\}}&= &
    (1,\dots, d_{\{p,q,(i,j)_{j}\}},1,\dots, 1) \in  \{1\} \times \cdots \{1\} \times \Gamma_{n,\alpha}^{4} \times \{1\} \times \cdots \{1\}  \\  && \subset \Gamma_{n,0}^{4} \times \Gamma_{n,1}^{4} \times \cdots \times \Gamma_{n,[\frac{r}{2}]}^{4},
    \end{eqnarray*}
if
\begin{center}
$min\{j,p,q\}-mid\{j,p,q\}+max\{j,p,q\} -2 \equiv \left\{
\begin{array}{cc}
    \alpha~mod~r, & \text{if}~min\{p,q,j\}<i<mid\{p,q,r\}~\text{or}~i>max\{p,q,j\}, \alpha \leq [\frac{r}{2}], \\
     -\alpha~mod~r, & \text{if}~min\{p,q,j\}<i<mid\{p,q,r\}~\text{or}~i>max\{p,q,j\}, \alpha \geq [\frac{r}{2}], \\
   \alpha-1~mod~r,
 &   \text{if}~i<min\{p,q,j\} ~\text{or}~ mid\{p,q,j\}<i<max\{p,q,j\}, \alpha \leq [\frac{r}{2}],\\
 -\alpha-1 ~mod~r,
 &   \text{if}~i<min\{p,q,j\} ~\text{or}~ mid\{p,q,j\}<i<max\{p,q,j\}, \alpha \geq [\frac{r}{2}]. \\
   \end{array}\right.$
\end{center}

$$D_{i,(i,j)}=D_{i,(i,j)}^{II}D_{i,(i,j)}^{I}D_{i,(i,j)}^{III}.$$
Now we define $f^{r}_{n} : PB_{n} \rightarrow \Gamma_{n,0}^{4} \times \Gamma_{n,1}^{4} \times \cdots \times \Gamma_{n,[\frac{r}{2}]}^{4}$ by
$$f^{r}_{n}(b_{ij}) = D_{i,(i,(i+1))}\cdots D_{i,(i,(j-1))}D_{i,(i,j)}D_{i,(i,j)}D_{i,(i,(j-1))}^{-1}\cdots D^{-1}_{i,(i,(i+1))},$$
for $1 \leq i<j \leq n$.

\begin{theorem}\label{main_theorem_rgamma}
The map $f^{r}_{n} : PB_{n} \rightarrow \Gamma_{n,0}^{4} \times \Gamma_{n,1}^{4} \times \cdots \times \Gamma_{n,[\frac{r}{2}]}^{4}$ defined above is a well defined homomorphism.
\end{theorem}

\begin{proof}
This statement might be proved similarly to the proof of Theorem~\ref{thm_to_gamma}. Let us list cases of singularities of codimension two explicitly. From now on let us consider four points on one circle. For the case of four points on one line the proof is analogous.
Note that the image of $f^{r}_{n}$ from four points on one circle depends on the number of points inside the circle. Without loss of generality, the number of points inside the circle is less than or equal to the number of points outside the circle.

\begin{enumerate}
\item One point moving on the plane is tangent to the circle, which passes through three points, see the center in Fig.~\ref{proof_rel1}. Note that the number of points inside the circle does not change when the point $P_{j}$ moves. It is easy to see that the image, when one point passes through the circle twice (upper left in Fig.~\ref{proof_rel1}), is
\begin{eqnarray*}
 (1,\dots,1, d_{(ijkl)}d_{(ijkl)},1,\dots, 1) \in  \Gamma_{n,0}^{4} \times \cdots \times \Gamma_{n,\alpha}^{4} \times \cdots \times \Gamma_{n,[\frac{r}{2}]}^{4},
    \end{eqnarray*}
where the $\alpha$ is the number of points mod $r$ inside the circle, passing through $P_{i},P_{k},P_{l}$.
If the point does not pass through the circle (upper right in Fig.~\ref{proof_rel1}), then the image is
\begin{eqnarray*}
 (1,\dots, 1,\dots, 1) \in  \Gamma_{n,0}^{4} \times \cdots \times \Gamma_{n,\alpha}^{4} \times \cdots \times \Gamma_{n,[\frac{r}{2}]}^{4}.
    \end{eqnarray*}
 The equality of those two images is obtained by the relation $d_{(ijkl)}^{2} =1$.

\item There are two sets $A=\{P_{i},P_{j},P_{k},P_{l}\}$ and $B=\{P_{s},P_{t},P_{u},P_{v}\}$ of four points, which are on the same circles such that $|A\cap B| \leq 2$, see Fig.~\ref{proof_rel2}. If the numbers of points inside the circles, which pass through points $\{P_{i},P_{j},P_{k},P_{l}\}$ and $\{P_{s},P_{t},P_{u},P_{v}\}$ respectively, are the same mod $r$, then the images of them are
\begin{eqnarray*}
 (1,\dots, d_{(ijkl)}d_{(stuv)},1,\dots, 1) \in  \Gamma_{n,0}^{4} \times \cdots \times \Gamma_{n,\alpha}^{4} \times \cdots \times \Gamma_{n,[\frac{r}{2}]}^{4},
    \end{eqnarray*}
and
\begin{eqnarray*}
 (1,\dots, d_{(stuv)}d_{(ijkl)},1,\dots, 1) \in  \Gamma_{n,0}^{4} \times \cdots \times \Gamma_{n,\alpha}^{4} \times \cdots \times \Gamma_{n,[\frac{r}{2}]}^{4},
    \end{eqnarray*}
The equality of them follows from the relation $d_{(ijkl)}d_{(stuv)} = d_{(stuv)}d_{(ijkl)}$.

If the number of points inside circles, where the points $\{P_{i},P_{j},P_{k},P_{l}\}$ and $\{P_{s},P_{t},P_{u},P_{v}\}$ are positioned respectively, are different mod $r$, then the images of them are
\begin{eqnarray*}
 (1,\dots, d_{(ijkl)}, \cdots, 1 , \cdots 1)  (1,\dots,1,\dots, d_{(stuv)}, \cdots 1) \\ \in  \Gamma_{n,0}^{4} \times \cdots \times \Gamma_{n,\alpha}^{4} \times \cdots \times \Gamma_{n,\beta}^{4} \times \cdots \times \Gamma_{n,[\frac{r}{2}]}^{4},
    \end{eqnarray*}
and
 \begin{eqnarray*}
 (1,\dots,1,\dots, d_{(stuv)}, \cdots 1) (1,\dots, d_{(ijkl)}, \cdots, 1 , \cdots 1) \\ \in  \Gamma_{n,0}^{4} \times \cdots \times \Gamma_{n,\alpha}^{4} \times \cdots \times \Gamma_{n,\beta}^{4} \times \cdots \times \Gamma_{n,[\frac{r}{2}]}^{4},
    \end{eqnarray*}
where $\alpha$ and $\beta$ depend on the numbers of points inside the circles, which pass through $\{P_{i},P_{j},P_{k},P_{l}\}$ and $\{P_{s},P_{t},P_{u},P_{v}\}$, respectively. It is easy to obtain the equality of them.

\item There are five points $\{P_{i},P_{j},P_{k},P_{l},P_{m}\}$ on the same circle. We obtain 10 generators. 
Note that if there are $\alpha$ points inside a circle, where four points lie, then in the very next time, when four points among $\{P_{i},P_{j},P_{k},P_{l},P_{m}\}$ lie on the same circle, there are $\alpha+1$ points inside the circle, but in the very next step, there are $\alpha$ points inside the circle again.
More precisely, we obtain the following element;
\begin{eqnarray*}
\indent(1,\dots, d_{(ijkl)}d_{(ijlm)}d_{(jklm)}d_{(ijkm)}d_{(iklm)},d_{(ijkm)}d_{(iklm)}d_{(ijkl)}d_{(ijlm)}d_{(jklm)}, \cdots, 1) \\
\noindent \in  \Gamma_{n,0}^{4} \times \cdots \times \Gamma_{n,\alpha}^{4} \times \Gamma_{n,\alpha+1}^{4}\times \cdots \times \Gamma_{n,[\frac{r}{2}]}^{4}.
    \end{eqnarray*}

From the relation $$\indent\indent d_{(ijkl)}d_{(ijlm)}d_{(jklm)}d_{(ijkm)}d_{(iklm)} =1 =d_{(ijkm)}d_{(iklm)}d_{(ijkl)}d_{(ijlm)}d_{(jklm)}$$
of $\Gamma_{n}^{4}$ we obtain that it is equal to
\begin{eqnarray*}
\indent(1,\dots, 1,1,\dots, 1) \in  \Gamma_{n,0}^{4} \times \cdots \times \Gamma_{n,\alpha}^{4} \times \Gamma_{n,\alpha+1}^{4}\times \cdots \times \Gamma_{n,[\frac{r}{2}]}^{4}.
    \end{eqnarray*}
and the proof is completed.

\end{enumerate}

\end{proof}

\subsection{Braids in $\mathbb{R}^{3}$ and groups $\Gamma_{n}^{4}$}

In Section~\ref{section:realisation}, we studied the notion of braids for $\mathbb{R}^{3}$ and $\mathbb{R}P^{3}$. {\em A braid for $\mathbb{R}^{3}$}\index{Braid!for $\mathbb{R}^3$} (or $\mathbb{R}P^{3}$) is a path in the configuration space $\tilde{C}'_{n}(\mathbb{R}^{3})$ (or $\tilde{C}'_{n}(\mathbb{R}P^{3})$), see page~\pageref{th1}, with some restrictions. If the initial and final points of the path in $\tilde{C}'_{n}(\mathbb{R}^{3})$ coincide, then the path is called {\it a pure braid for $\mathbb{R}^{3}$ $(\mathbb{R}P^{3})$}\index{Braid!for $\mathbb{R}^3$!pure}. In the present section we will construct a group homomorphism from pure braids on $n$ strands in $\mathbb{R}^{3}$ to $\Gamma_{n}^{4}$.



We shall consider (good and transverse) pure braids on $n$ strands in $\mathbb{R}^{3}$ and construct a group homomorphism from pure braids on $n$ strands in $\mathbb{R}^{3}$ to the group $\Gamma_{n}^{4}$. Each element of $\Gamma_{n}^{4}$ corresponds to a moment when four points lie on $(4-2)$-dimensional plane in $\mathbb{R}^{4-1}$, but in the case of $\Gamma_{n}^{4}$ ``the order'' of four points on $(4-2)$-dimensional plane is very important. This order will be ignored when the group homomorphism from pure braids in $\mathbb{R}^{3}$ to $\Gamma_{n}^{4}$ will be constructed. Now we formulate more precisely how the group homomorphism from pure braids on $n$ strands in $\mathbb{R}^{3}$ to $\Gamma_{n}^{4}$ is constructed.

Let $\gamma$ be a good and transverse pure braid on $n$ strands with base point $x = (x_{1},\cdots, x_{n})$ in $\mathbb{R}^{3}$. We call $t \in [0,1]$ {\em a special singular moment of $\gamma$} if the following holds:

\begin{enumerate}
\item at the moment $t$ four points $x_{p},x_{q},x_{r},x_{s}$ are on the same plane $\Pi_{t}$,
\item the four points $x_{p},x_{q},x_{r},x_{s}$ form a convex quadrilateral on $\Pi_{t}$,
\item all of $\{x_{1},\cdots, x_{n}\} \backslash \{ x_{p},x_{q},x_{r},x_{s} \}$ are situated only in one of the connected components of $\mathbb{R}^{3} \backslash \Pi_{t}$.
\end{enumerate}

A normal vector $v_{t}$ of $\Pi_{t}$ pointing to the connected component, in which the points $$\{x_{1},\cdots, x_{n}\} \backslash \{ x_{p},x_{q},x_{r},x_{s} \}$$ are placed, is called {\it the pointing vector} at the special singular moment $t$.

\begin{figure}[h!]
 \centering
\includegraphics[width = 2.5cm]{ 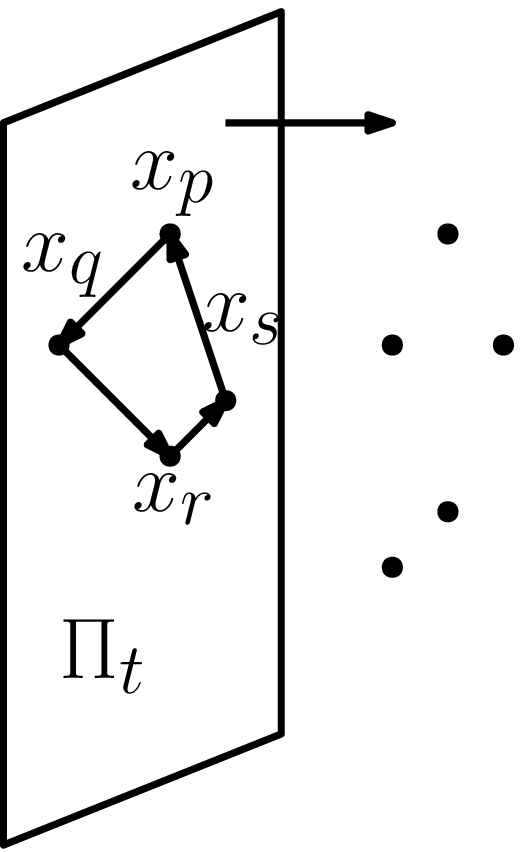}
 \caption{A special singular moment, corresponding to $d_{(pqrs)}$}\label{singular_moment}
\end{figure}

\begin{remark}
\begin{enumerate}
\item The plane $\Pi_{t}$ admits a unique orientation with respect to $v_{t}$.
\item Naturally, the quadrilateral admits the orientation with respect to $v_{t}$, see Fig.~\ref{singular_moment}.
\end{enumerate}
\end{remark}

Let us enumerate all special singular moments $0<t_{1}<\cdots< t_{l}<1$ of $\gamma$. For each $t_{s}$ by definition of good pure braids on $n$ strands there are exactly four points $\{ x_{p},x_{q},x_{r},x_{s}\}$ on plane $\Pi_{t_{s}}$. As indicated in the previous remark the convex quadrilateral on plane $\Pi_{t_{s}}$ with four vertices $\{ x_{p},x_{q},x_{r},x_{s}\}$ admits the orientation with respect to $v_{t_{s}}$. If the four points $x_{p},x_{q},x_{r},x_{s}$ are positioned in the indicated order in accordance with the orientation with respect to $v_{t_{s}}$, then we associate the moment $t_{s}$ to $d_{t_{s}}=d_{(pqrs)}$. Let us define a map $g\colon \pi_{1}(\tilde{C}_{n}'(\mathbb{R}^{3})) \rightarrow \Gamma_{n}^{4}$ by the equality $g(\gamma) = d_{t_{1}}\cdots d_{t_{l}}$.

\begin{theorem}
The map $g\colon \pi_{1}(\tilde{C}_{n}'(\mathbb{R}^{3})) \rightarrow \Gamma_{n}^{4}$ is well-defined.
\end{theorem}

\begin{proof}
We consider moments of isotopy between two paths, when the path at some moment in the isotopy between two paths is not good or not transverse. Let us list such cases explicitly.

\begin{enumerate}
\item There are four points on a plane, which disappear after a small perturbation, see Fig.~\ref{rel1_gamma}. This corresponds to the relation $d_{(pqrs)}^{2} = 1$.

\begin{figure}[h!]
	\begin{minipage}{.45\textwidth}
		\begin{center}
 			\includegraphics[width =.95\textwidth]{ 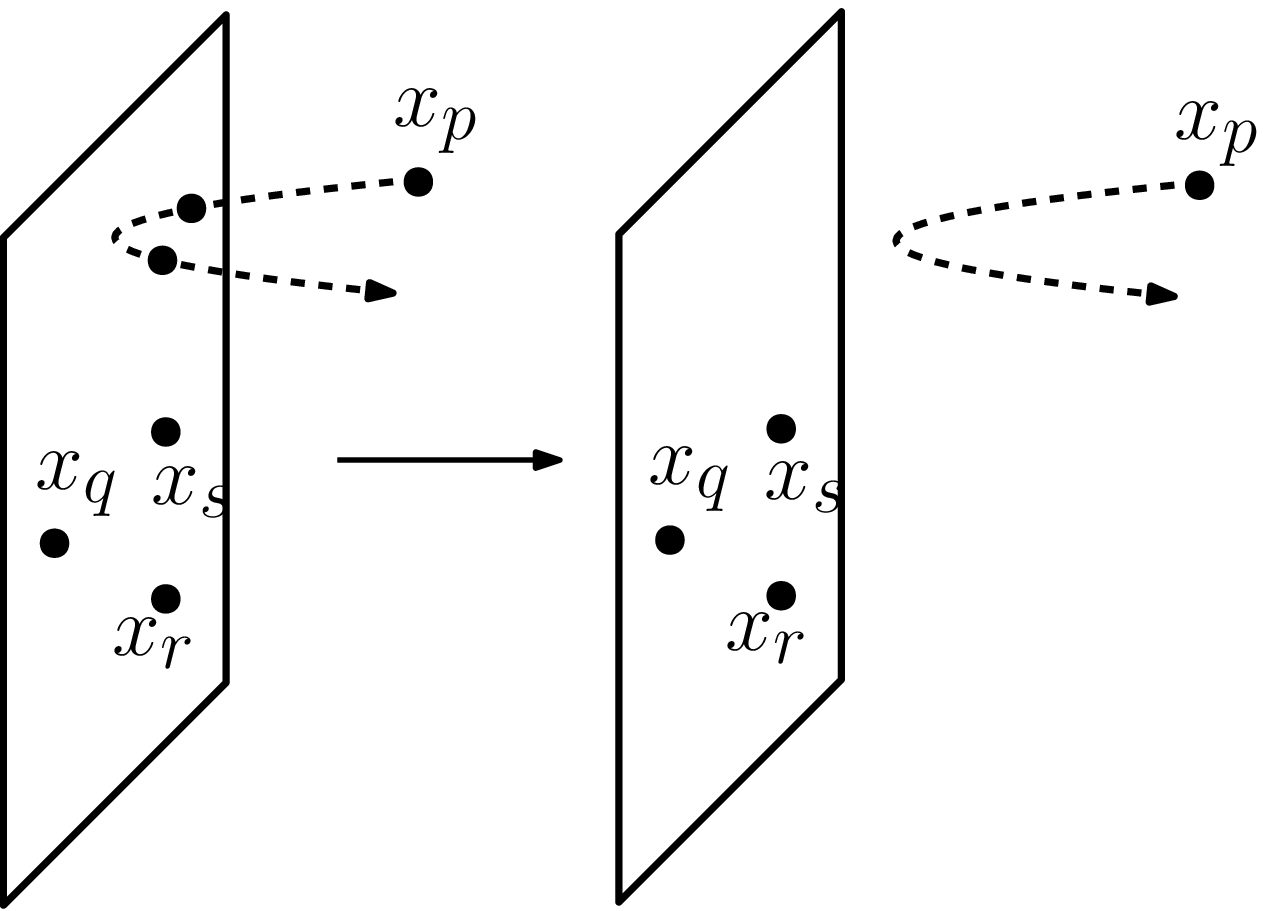}
		\end{center}
		\vspace{-0.1cm}
 		\caption{Case 1: Four points on a plane, which disappears after a small perturbation} \label{rel1_gamma}
 	\end{minipage}
 	\hspace{5mm}
 	\begin{minipage}{.45\textwidth}
 		\begin{center}
 			\includegraphics[width =.6\textwidth]{ 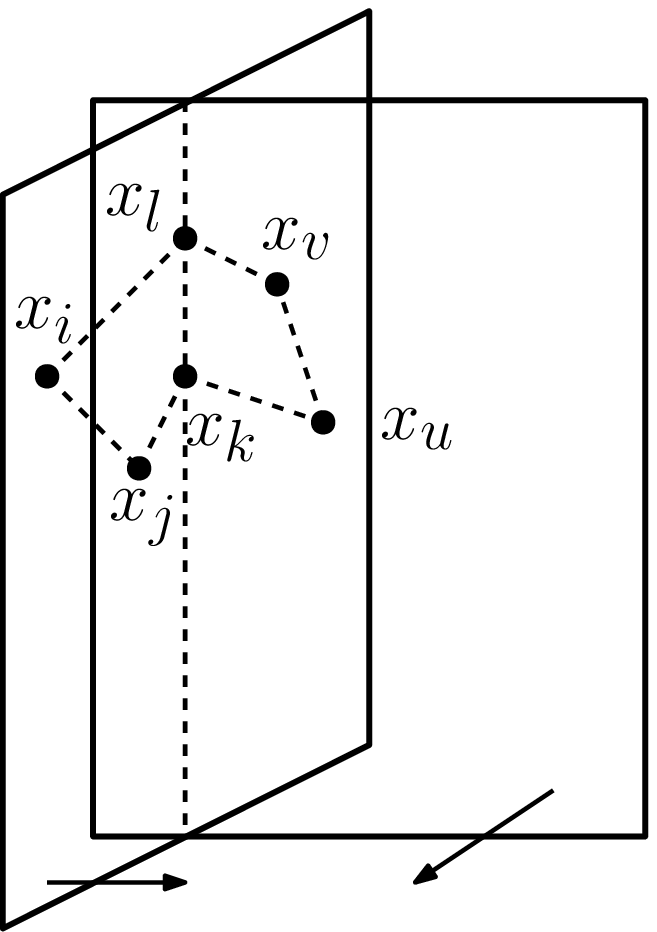}
		\end{center}
		\caption{Case 2: Two sets of four points $\{x_{i},x_{j},x_{k},x_{l}\}$ and $\{x_{l},x_{k},x_{u},x_{v}\}$ on planes at the same moment} \label{rel2_gamma}
 	\end{minipage}
 \end{figure}


\item At a moment there are two sets of four points $m$ and $m'$ with $|m \cap m'|<3$, which are placed on plane at the same moment, see Fig.~\ref{rel2_gamma}. This corresponds to the relation $d_{(ijkl)}d_{(stuv)} =d_{(stuv)} d_{(ijkl)}$.


\item At a moment five points lie on a plane. This is similar to the case of ``five points on the circle'' in the proof of Theorem~\ref{thm_to_gamma}. This corresponds to the relation $d_{(ijkl)}d_{(ijlm)}d_{(jklm)}d_{(ijkm)}d_{(iklm)} = 1$ of $\Gamma_{n}^{4}$.

\end{enumerate}

\end{proof}

\begin{figure}[h!]
 \centering
\includegraphics[width = 4cm]{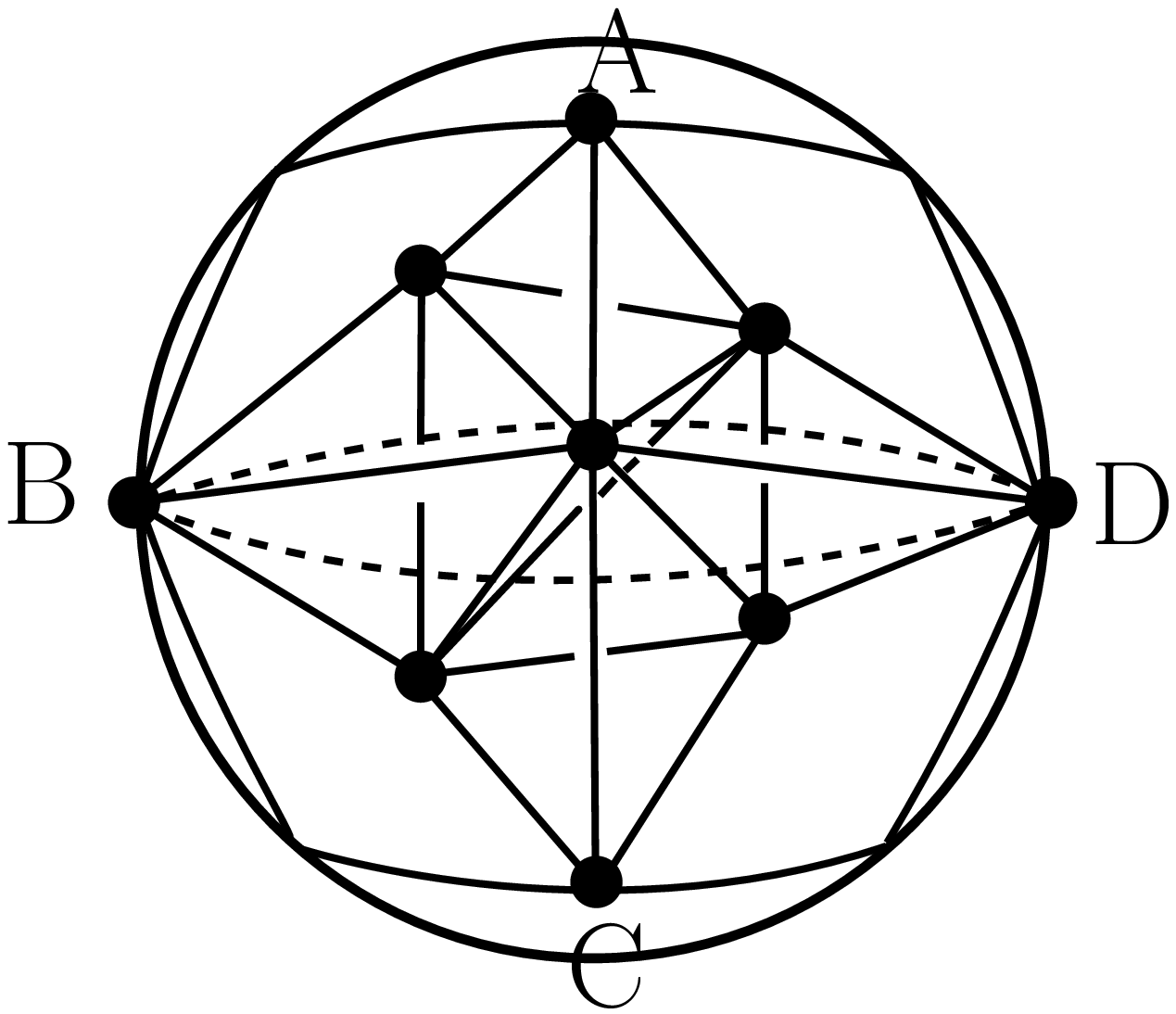}
 \caption{Triangulation of 3-disc with $\{P_{1}(t), \cdots, P_{n}(t)\}$ inside the sphere and four points $\{A,B,C,D\}$ on the sphere}\label{triangulation_R3}
\end{figure}

Let $\{P_{1}(t), \cdots, P_{n}(t)\}_{t \in [0,1]}$ be $n$ moving points in $\mathbb{R}^{3}$, corresponding to a path in $\pi_{1}(C_{n}'(\mathbb{R}^{3}))$. We may assume that the points $\{P_{1}(t), \cdots, P_{n}(t)\}_{t \in [0,1]}$ move inside a sphere with sufficiently large diameter. Let us fix four points $\{A,B,C,D\}$ on the sphere. A triangulation of the 3-ball with vertices $$\{P_{1}(t), \cdots, P_{n}(t)\} \cup \{A,B,C,D\}$$ can be obtained for each $t \in [0,1]$, for example, see Fig.~\ref{triangulation_R3}.

\begin{figure}[h!]
 \centering
\includegraphics[width = 9cm]{ 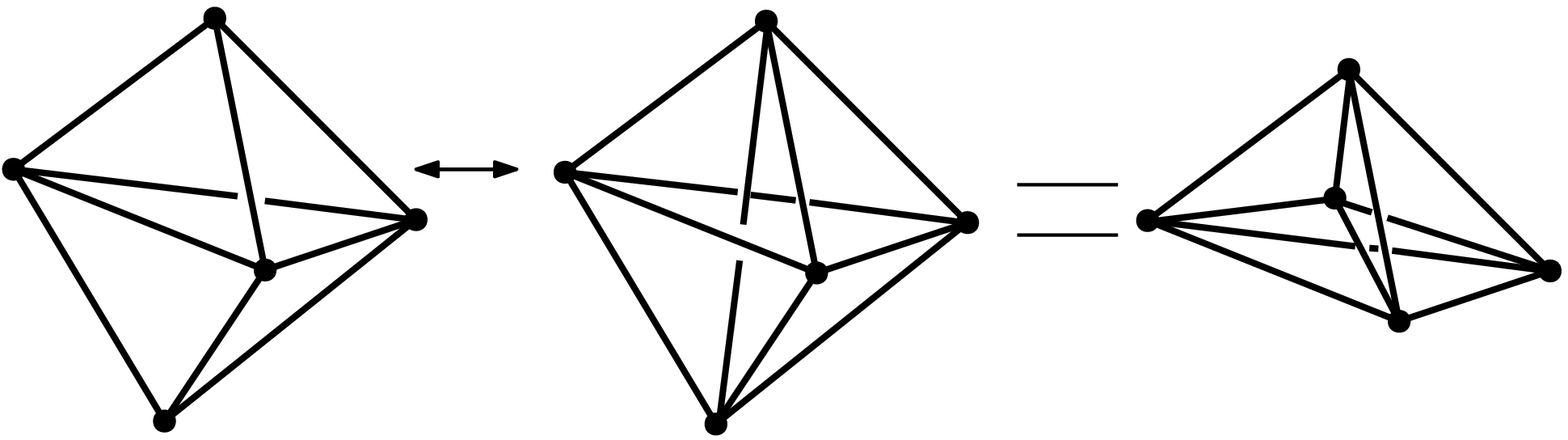}
 \caption{Pachner move and moving of a vertex of the triangulation of 3-dimensional space}\label{pachner_move}
\end{figure}

On the other hand, as described in Fig.~\ref{pachner_move}, the moving of a vertex of the triangulation can be described by applying the {\em Pachner moves} to the triangulation of 3-dimensional space.

\begin{figure}[h!]
 \centering
\includegraphics[width = 8cm]{ 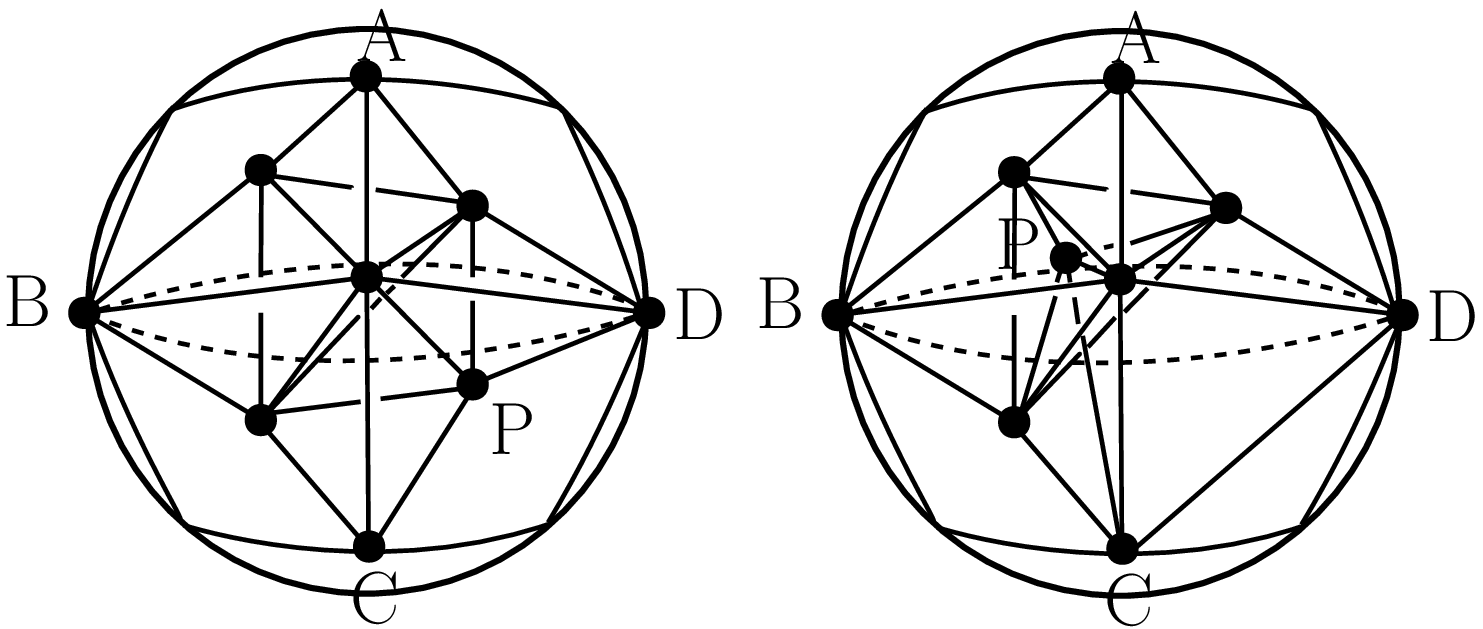}
 \caption{Applying a Pachner move to the triangulation of a 3 dimensional space and the moving of a vertex}\label{point_moving_triangulation}
\end{figure}

In other words, a path $\{P_{0}(t),\cdots, P_{n}(t)\}$ in $\pi_{1}(C_{n}'(\mathbb{R}^{3}))$ can be described by a finite sequence of Pachner moves applied to the triangulations of the sphere, see Fig.~\ref{point_moving_triangulation}.

It can be expected that by using the triangulation of a sphere with $n+4$ points and the sequence of the Pachner moves, we obtain an invariant of (pure) braids.

%% file: braid_inv_triangulation.tex
\subsection{A representation of braids via Ptolemy transformation}\label{sec:inv_braids_ptolemy}
The aim of the present section is to construct an (infinite-dimensional)
representation of braid groups for arbitrary $2$-surfaces valued in operators
on the space of rational functions over a ring. For related work see \cite{Kashaev_lecture}.

The main idea behind this representation is the operator which is associated with the flip, see Fig~\ref{fig:flip}.\index{Flip}

%

Let $\Sigma$ be a $2$-surface endowed with a metric $G$, and let $N$ be a large natural number.
We shall consider the space of all triangulations of $\Sigma$ with $N$ vertices. Some of these triangulations originate from Delaunay triangulations with vertices
in the given $N$ points (for this to happen, some length constraints should hold). When moving $N$ points around, the Delaunay triangulation\index{Delaunay triangulation} undergoes some flips. 

In the sequel, by an {\em ordered triangulation} we mean a triangulation with enumerated
vertices. Consider the set of all triangulations of $\Sigma$ with $N$ enumerated vertices. 
With any ordered triangulation having $X$ edges, we associate variables $\{x_{1}, \dots, x_{X}\}$, one to each oriented edge. An edge is usually oriented from the vertex with smaller index towards the one with larger index. If we associate a variable $a$ to an edge, we associate $-a$ to the opposite edge. Two triangulations are called {\em adjacent} if one may be obtained from the other via a single flip. Consider two adjacent ordered triangulations $T_{1}$ and $T_{2}$. $T_2$ is obtained from $T_1$ by removing a single edge and replacing it with another diagonal in the same quadrilateral (note that all other edges remain unchanged by the flip).

\begin{figure}
\centering\includegraphics[width=230pt]{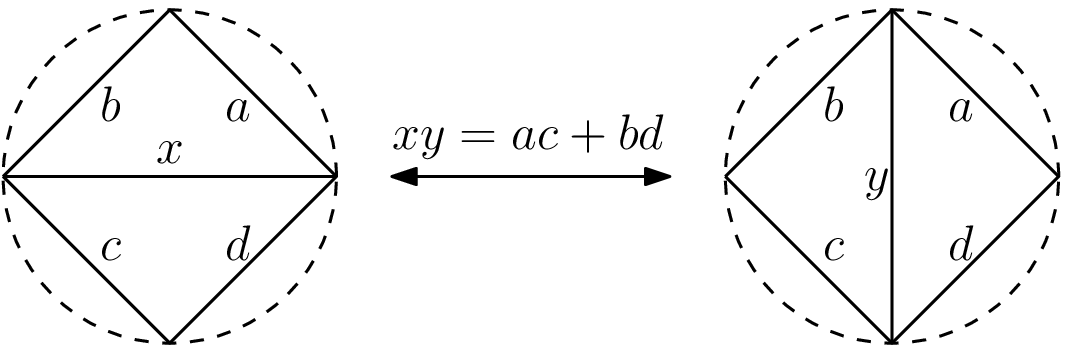}
\caption{The Ptolemy transformation}
\label{flip_ind}
\end{figure}

Now, consider the ring $\mathbb{Q} (x_{1},\dots, x_{X})$ of all rational functions in variables $\{x_{1}, \dots, x_{X}\}$. We shall relate variables as follows: If the triangulation $T_{2}$ is obtained from the triangulation $T_{1}$ by a flip as in Fig.~\ref{flip_ind} within the quadrilateral $a,b,c,d$, then the diagonal $y$ is equal to $\frac{ac+bd}{x}$. 
The variables transformation $x\to y=\frac{ac+bd}{x}$ is called a {\em Ptolemy transformation}. Note that each flip gives us a Ptolemy transformation of variables, but there are elements of $\mathbb{Q} (x_{1},\dots, x_{X})$ which cannot be obtained from sequences of Ptolemy transformations\index{Ptolemy transformation}.

Starting with a triangulation $T$, we can consider loops in the space of triangulations,
i.e., sequences of triangulations $T=T_{0}\to T_{1}\to T_{2}\to \dots \to T_{k}=T$ where
for each $i$ the triangulations $T_{i}, T_{i+1}$ differ by a flip.

Certainly, when we deal with loops, the labels of each $T_{i}$ are expressible in terms
of $T_{0}$, subsequently, for the triangulation $T_{k}=T$ we get a new set of labels.

\begin{theorem}\label{inv_braids_ptolemy}
If a loop $T=T_{0}\to \dots \to T_{k}=T$ in the space of triangulations originates from a braid
which is homotopic to the identity, then the labels of $T_{k}$ valued in $\mathbb{Q}(x_{1},\dots,x_{X})$ coincide with those of $T_{0}$.
\end{theorem}
The main ingredient in the proof is the pentagon equation. In other words, when a pentagon with two diagonals is given, if we apply Ptolemy transformations five times to quadrilaterals of the pentagon, the pentagon and its two chosen diagonals returns to the initial state with the initial labels. 

Relations of the form $a^2=1$ follow immediately from the construction: applying a flip twice to the same quadrilateral returns it to its initial state with the initial labels.

Finally, the far commutativity relation and the relations are also straightforward: the order in which we perform flips in two quadrilaterals with at most two common vertices does not change the outcome, because those flips change only labels of the diagonals, and such quadrilaterals do not have common diagonals.


\begin{example}
Let the  three points $1,2,3$ be the vertices of a triangle. Consider a point $4$ in the center of this triangle and a point $5$ which goes around $4$ in the neighbourhood of $4$. This gives rise to a sequence of triangulations given in Fig. \ref{exmpl}. We denote the triangulations appearing in the sequence by $T_0$ -- $T_6$.
\begin{figure}
\centering\includegraphics[width=330pt]{ 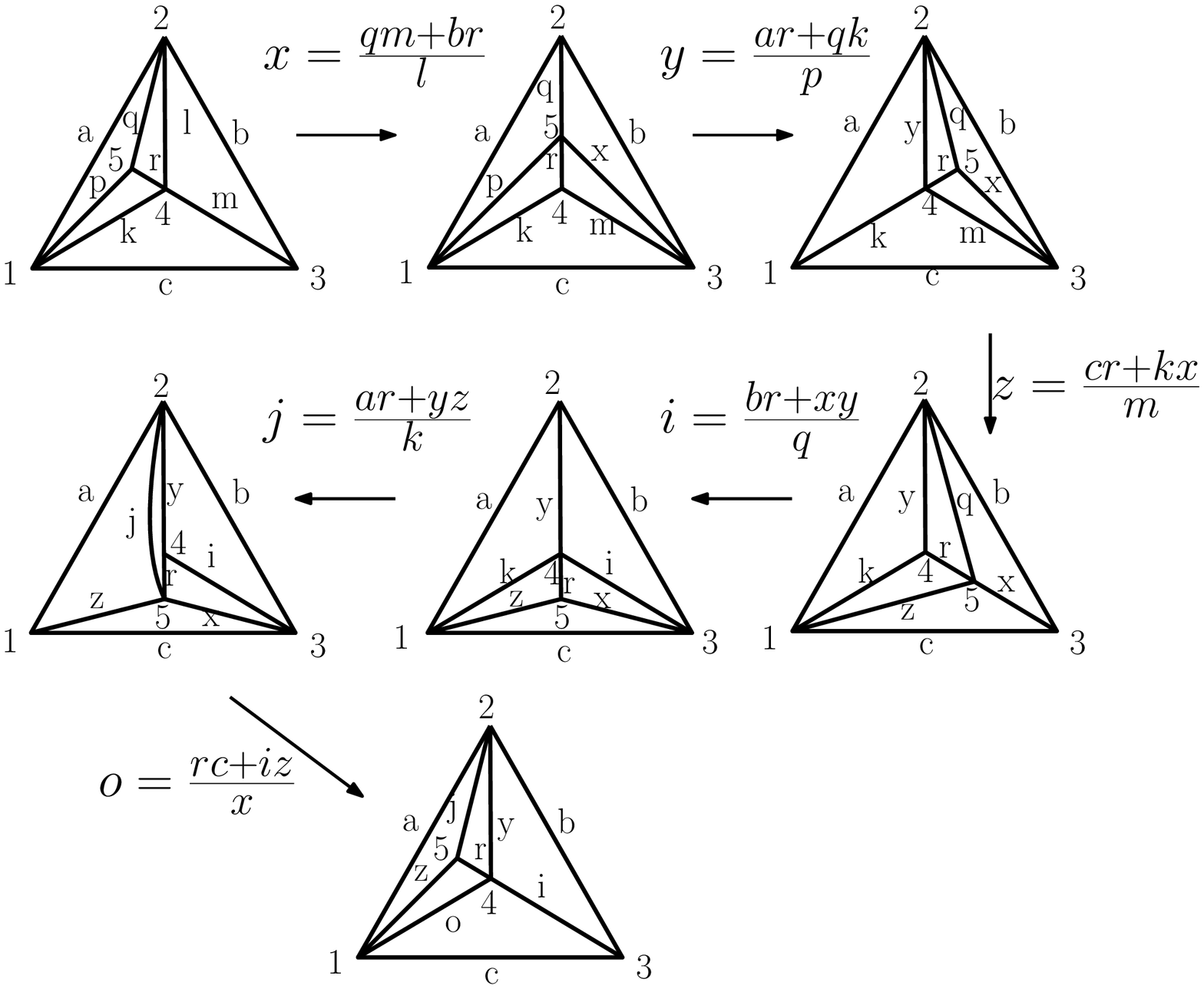}
\caption{Six triangulations corresponding to a rotation of point}
\label{exmpl}
\end{figure}
Let us calculate the labels of edges corresponding to these pictures. For the initial
one, we denote the labels by 
as shown in Fig. \ref{exmpl}.

In all subsequent figures the points $1,2,3,4$ as well as labels $a,b,c,d,k,l,m$
remain unchanged.
Besides the initial triangulation $T_{0}$ and the final one $T_{6}$ it is worth taking glance at $T_{3}$. It looks similar to $T_{0}$, however, with one important difference: the small edge denoted by $r$ is now oriented in the opposite direction. At $T_{6}$ we obtain the following labels:
$$x=\frac{qm+br}{l}, y=\frac{ar+qk}{p}, z=\frac{crl+kqm+kbr}{ml};$$
$$i=\frac{bprl+(ar+qk)(qm+br)}{pql}, j=\frac{arpml+(ar+qk)(crl+kqm+kbr)}{kpml};$$
$$o=\frac{k^{2}mq^{2}+bklpr+cl^{2}pr+cklqr+akmqr+abkr^{2}+aclr^{2}}{lmpq}.$$
From this we see that the variables corresponding to $T_{6}$ do not coincide with those
corresponding to $T_{0}$, which says that the above invariant of this pure braid
is non-trivial.
\end{example}

\subsection{A representation of braids via tropical relation}

Below another relation which is associated to flips (Fig.~\ref{fig:flip}) is presented. We would like to highlight that for the relation, which is associated with each flip, to give an invariant for braids as Theorem~\ref{inv_braids_ptolemy}, it should satisfy {\em the pentagon equation}\index{Pentagon equation} as shown in Fig~\ref{proof_rel3}.

\begin{figure}[h!]
	\begin{minipage}{.45\textwidth}
		\begin{center}
 			\includegraphics[width =.95\textwidth]{ 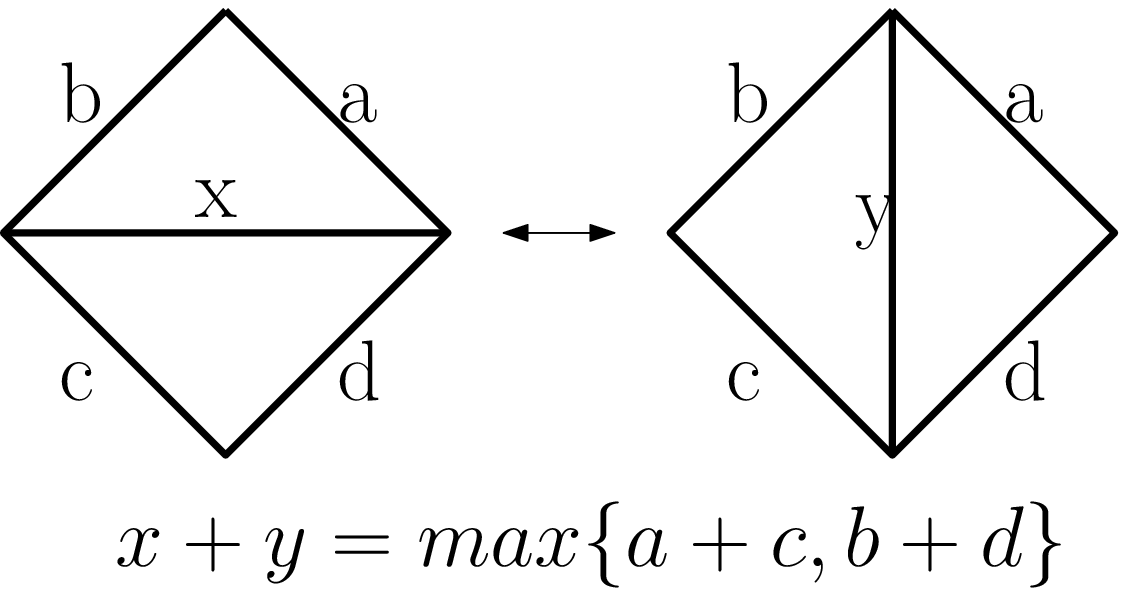}
		\end{center}
		\vspace{-0.1cm}
 		\caption{The tropical equation} \label{flip_tropical}
 	\end{minipage}
 	\hspace{5mm}
 	\begin{minipage}{.45\textwidth}
 		\begin{center}
 			\includegraphics[width =.43\textwidth]{ 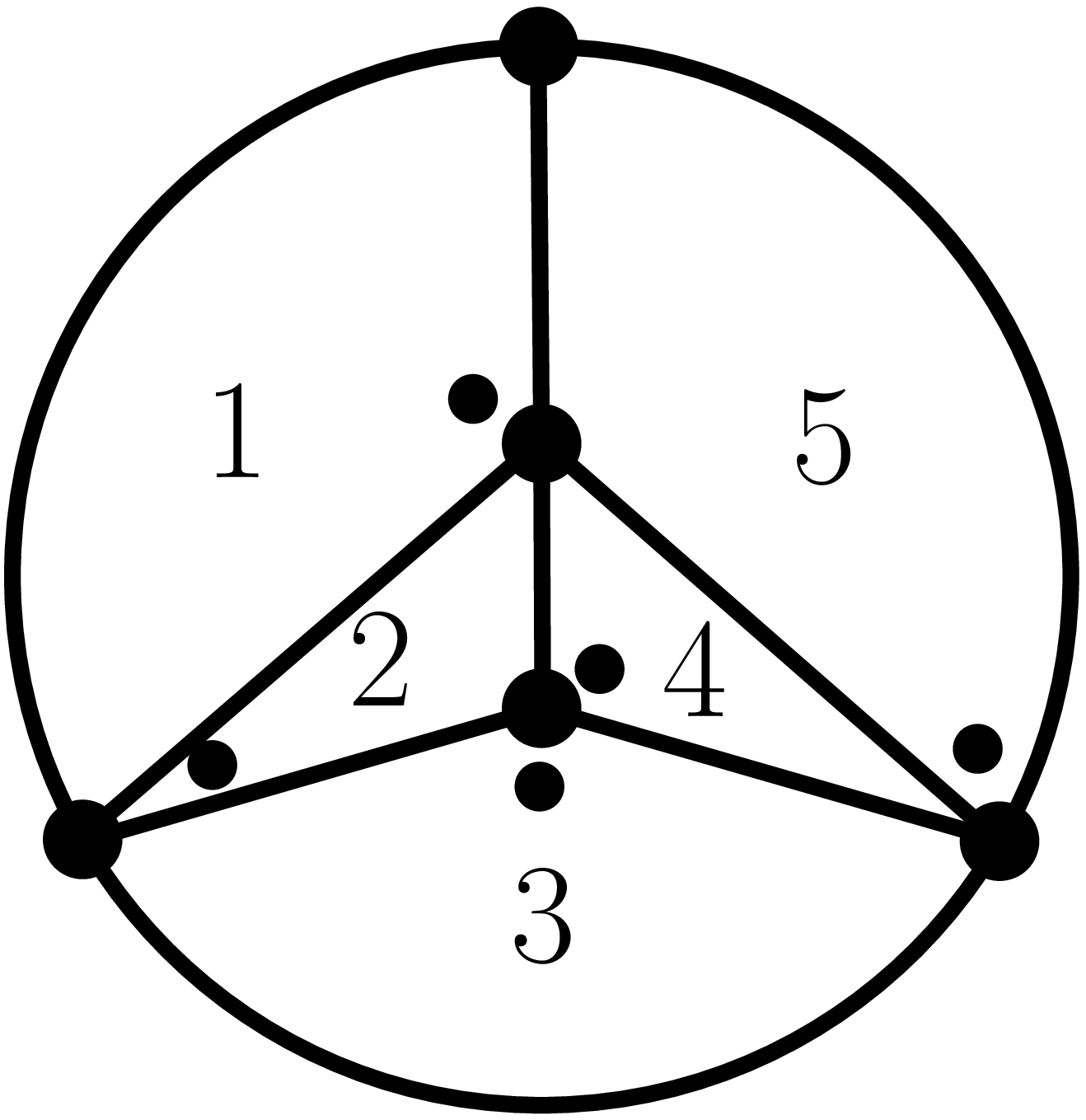}
		\end{center}
		\caption{A decorated ideal triangulation of the disc} \label{exa_decorated_triang}
 	\end{minipage}
 \end{figure}

Recall, that in {\em tropical geometry} the following operations are considered:

$$x\oplus y=\max\{x,y\},$$
$$x\otimes y=x+y,$$
where in the right-hand side we consider the usual maximum and addition operations. Note that this operation endow $\mathbb{R}\cup\{-\infty\}$ with a structure of {\em semiring}: there is no requirement of the existence of inverse elements under the addition $\oplus$. The Ptolemy relation $xy=ab+cd$, considered before, in the tropical case transforms into
\begin{equation}
x\otimes y = (a\otimes c) \oplus (b\otimes d) \Leftrightarrow x + y = max\{ a+c,b+d \},
\end{equation}
see Fig.~\ref{flip_tropical}.

\begin{lemma}
The tropical relation satisfies the far commutativity relation, the relations of the form $a^2=1$, and the pentagon relation.
\end{lemma}

This lemma holds due to the following consideration: if the prove does not use the additive inverses, it should work just like in the usual case considered before. The most difficult part is to check the pentagon relation. But we can do it ``forward'' by looking at what happens to the variables under the five steps of the process (five flips in the geometric language). This operation does not use the additive inverse and hence is done straightforwardly. The calculation verifies the claim.

The far commutativity relations again (as in the discussion of Theorem~\ref{inv_braids_ptolemy}) come from the fact that if two quadrilaterals are ``far enough'' from each other, their diagonals are all different and the variables transformations would not intervene with one another.

Finally, the $a^2=1$ may be checked by hand. Indeed, let us apply a flip twice to the left-hand side of Fig.~\ref{flip_tropical}. First we get $x\to y=\max\{a+c,b+d\} - x$. After the second flip we get $y\to \max\{a+c,b+d\}-y$. In other words, we have the transformation $$x\to \max\{a+c,b+d\}-(\max\{a+c,b+d\}-x)=x.$$
This equality conclude the proof.


\subsection{A representation of braids via decorated triangulations}
The purpose of the present section is to introduce the notion of decorated triangulation and to construct a representation of braids group, for more detail, see \cite{Kashaev_lecture}. 

\begin{definition}\cite{Kashaev_lecture}
{\em A decorated triangulation}\index{Decorated triangulation} is a triangulation where all triangles are linearly ordered and in each triangle, a distinguished corner is fixed.
\end{definition}

An example of a decorated triangulation of the disc with 5 points is described in Fig.~\ref{exa_decorated_triang}.


\begin{definition}\label{def:BAS}\cite{Kashaev_lecture}
Let $\mathcal{C} = (\mathcal{C}, \times, \{P_{X,Y}\}_{X,Y \in Ob\mathcal{C}})$ be a category of groups. A {\em basic algebraic system} (BAS) in $\mathcal{C}$ consists of a group $V \in Ob\mathcal{C}$ and two morphisms $R \in End (V)$ and $W \in End (V\times V)$ such that
\begin{enumerate}
\item $R_{i}^{3} = id_{V}$,
\item $W_{i,j}W_{j,k} = W_{j,k}W_{i,k}W_{i,j}$,
\item $R_{i}R_{j}W_{j,i}R_{i}^{-1}W_{i,j} = P_{i,j}$
\end{enumerate}
for all $i,j,k\in V$, where $R_i:=R(i), \, W_{i,j}:=W(i,j), \, P_{i,j}:=P_{V,V}(i,j)$.
\end{definition}

\begin{figure}
\centering\includegraphics[width=130pt]{ 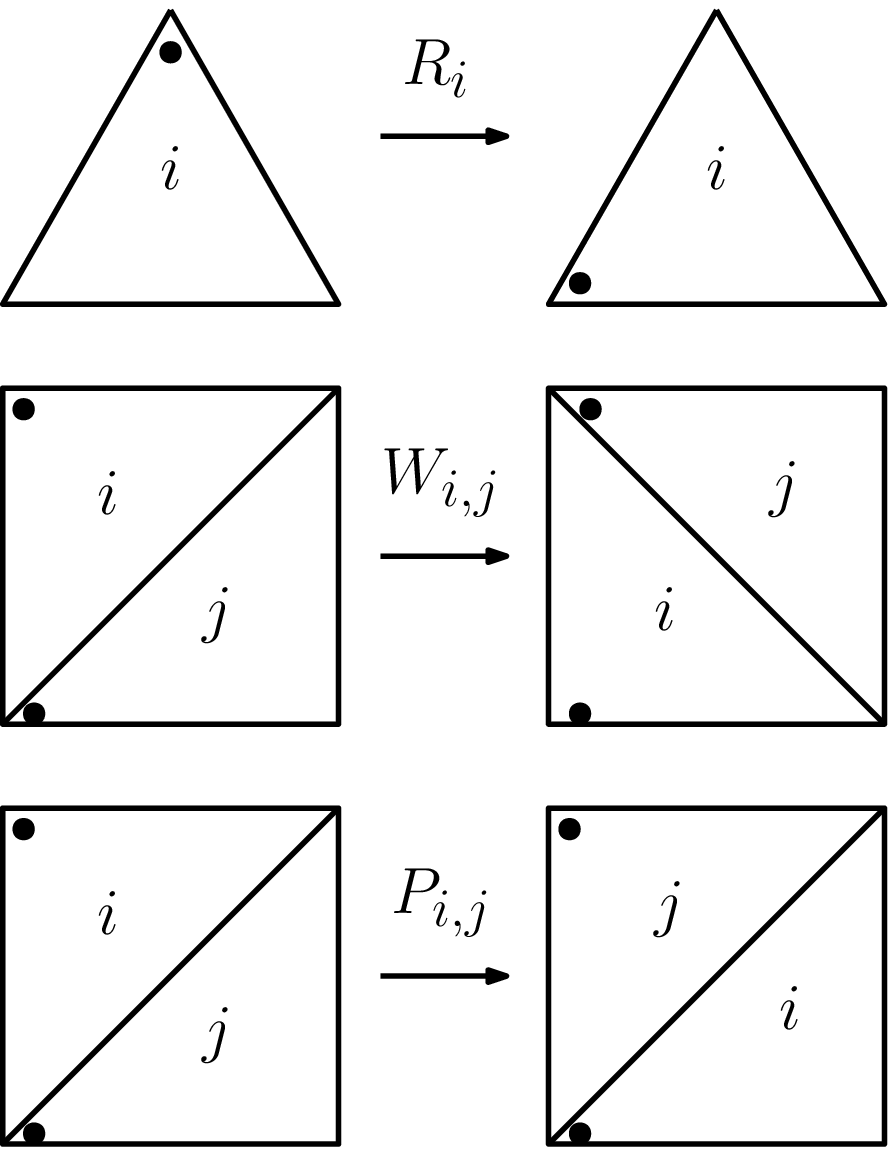}
\caption{Distinguished corner changes and the flip}
\label{corner_flip}
\end{figure}

In the context of decorated triangulations the morphisms $R_{i}\colon V \rightarrow V$ and $W_{i,j}\colon V \otimes V \rightarrow V \otimes V$ correspond to distinguished corner changes and the flip, respectively; the morphism $P$ corresponds to numbering switch, see Fig.~\ref{corner_flip}.

\begin{figure}
\centering\includegraphics[width=250pt]{ 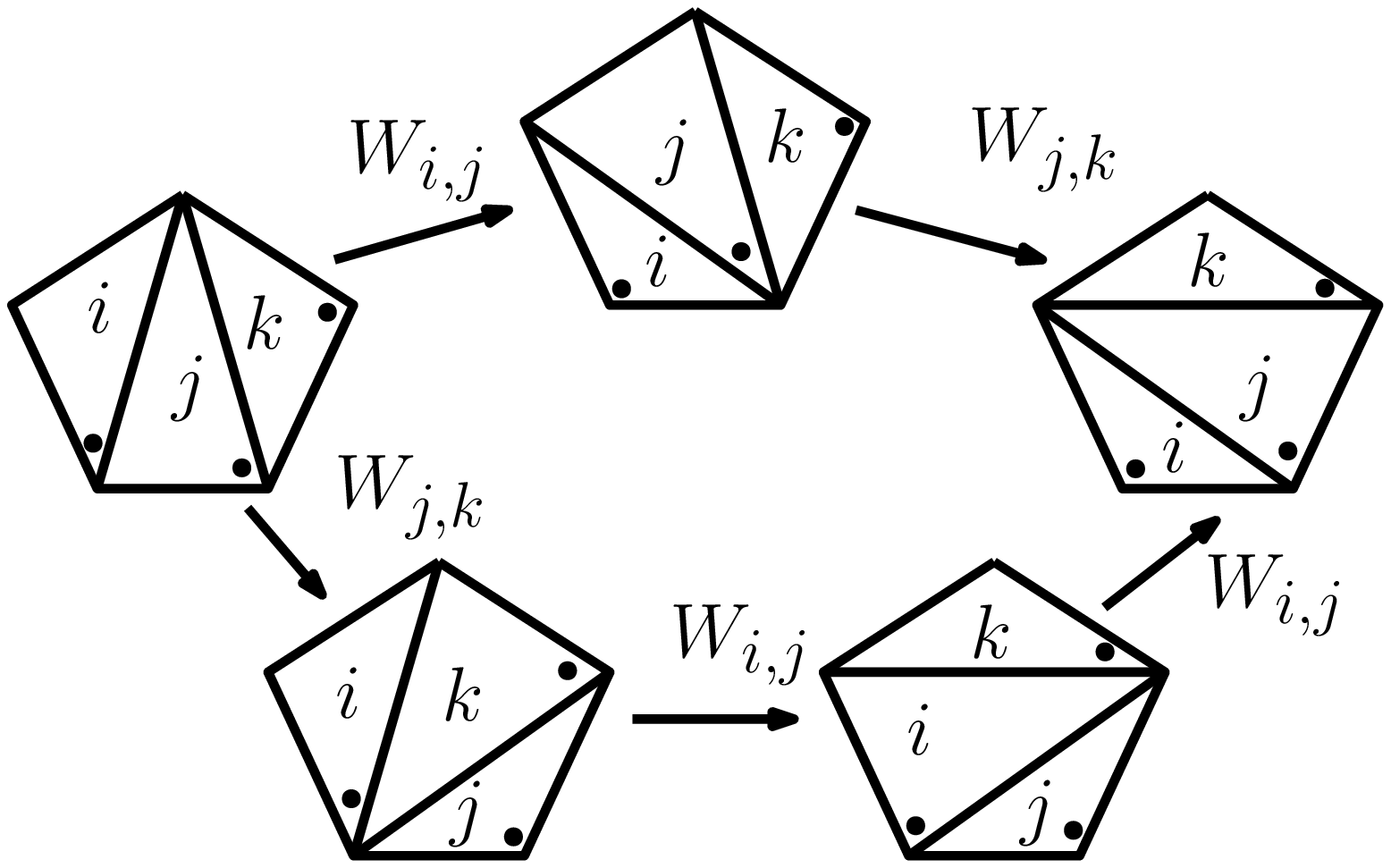}

\caption{Geometric description of the condition (2) from Definition~\ref{def:BAS}}
\label{pentagon_BAS}
\end{figure}

\begin{figure}
\centering\includegraphics[width=200pt]{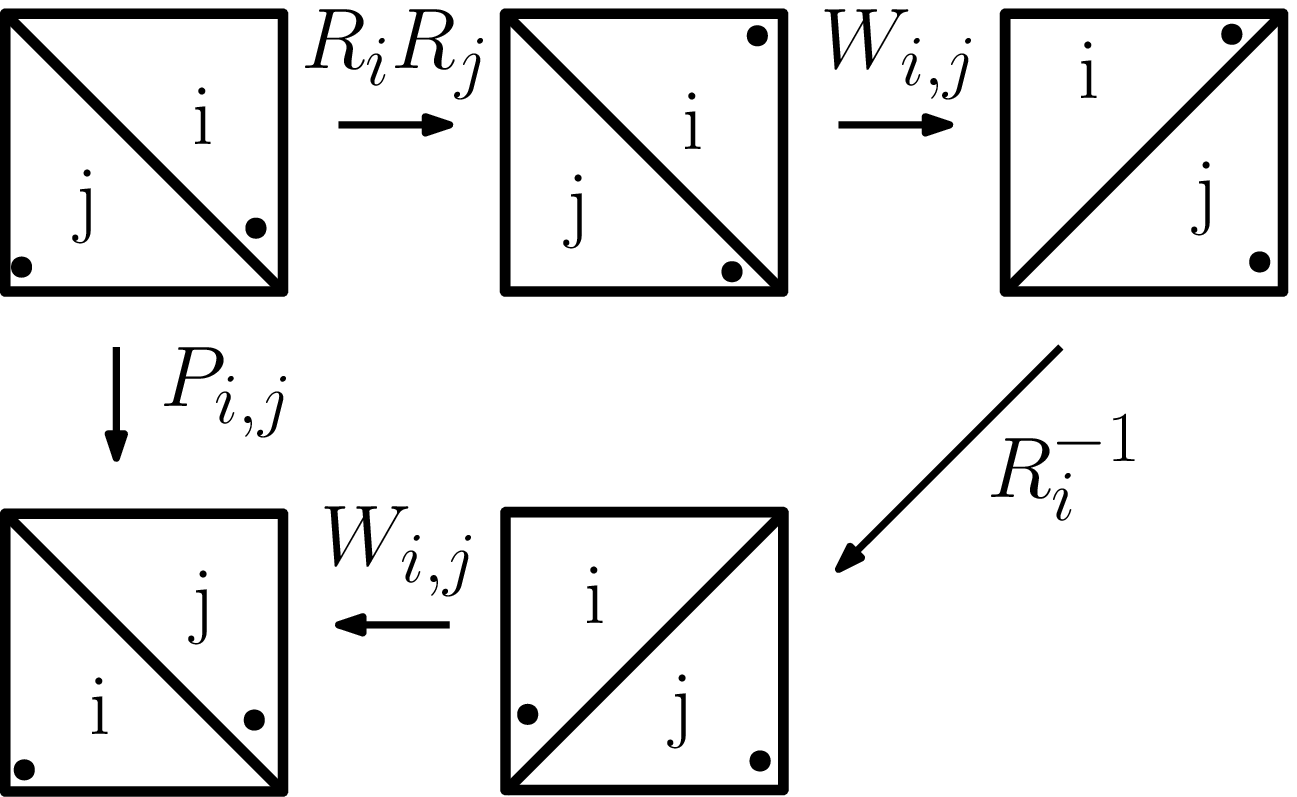}
\caption{Geometric description of the condition (3) from Definition~\ref{def:BAS}}
\label{geo_change_order_BAS}
\end{figure}

Moreover the conditions (2) and (3) from Definition~\ref{def:BAS} can be described by the pentagon relation and changes of the order of triangles, see Fig.~\ref{pentagon_BAS} and Fig.~\ref{geo_change_order_BAS}.

\begin{example}\label{exa:BAS_Ptolemy}
Let $V = \mathbb{R}_{+}^{2}$ (with the usual component-wise addition as the operation) and define $R\colon V \rightarrow V$ and $W\colon V \times V \rightarrow V \times V$ by

\begin{equation}
R(x_{1},x_{2}) = \left(\frac{1}{x_{1}},\frac{x_{2}}{x_{1}}\right), W((x_{1},x_{2}),(y_{1},y_{2})) = \left(\frac{x_{2}y_{1}}{x_{1}y_{2}+x_{2}}, \frac{y_{2}}{x_{1}y_{2}+x_{2}}\right).
\end{equation}
Then this is a BAS as in Definition~\ref{def:BAS}. We call this {\em ratio coordinates}.
\end{example}

\begin{figure}
\centering\includegraphics[width=150pt]{ 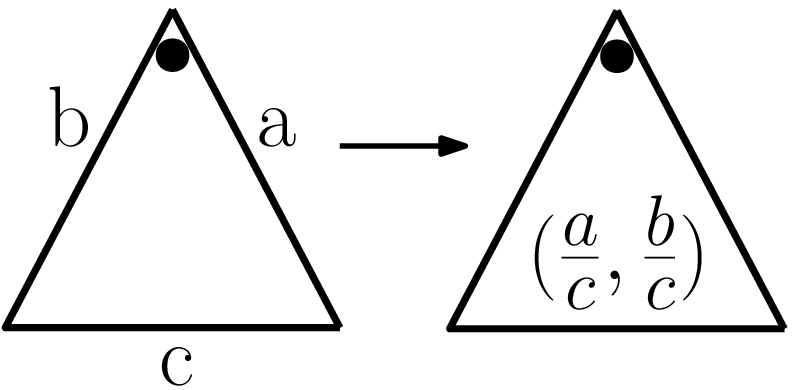}
\caption{$a,b,c \in \mathbb{Q}(x_{1}, \dots x_{X})$, which are associated to each edge of the triangulation $\tau$}
\label{ratio_coordinate}
\end{figure}

Let us consider the set of all ordered triangulations of $\Sigma$ with $N$ enumerated vertices as in Section~\ref{sec:inv_braids_ptolemy}. For each triangulation $\tau$ of $\Sigma$ we can obtain a decorated triangulation $\tilde{\tau}$ by fixing a corner for each triangle and enumerating the triangles. Now we define a map from triangles of $\tilde{\tau}$ to $\mathbb{Q}(x_{1}, \dots, x_{X}) \times \mathbb{Q}(x_{1}, \dots, x_{X})$ as described in Fig.~\ref{ratio_coordinate}. Since the ratio coordinates BAS satisfies the pentagon equation, this mapping gives us an invariant of braids.

\subsection{A representation of braids via truncated triangulation}
\label{subsec:truncated_Zickert}

\begin{figure}
\centering\includegraphics[width=150pt]{ 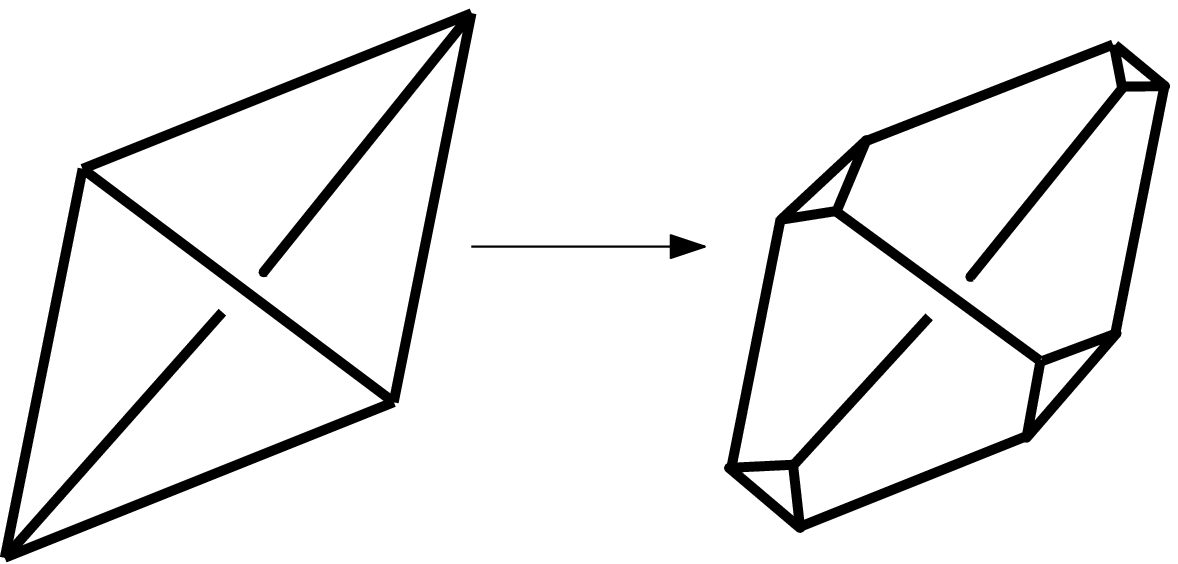}
\caption{Truncated triangulation, obtained from triangulation}
\label{truncated_tri}
\end{figure}

In this section we shall define a representation of braids by using {\em truncated triangulations}\index{Truncated triangulation}. Truncated triangulations can be obtained from usual triangulations as Fig.~\ref{truncated_tri}.

Edges of truncated triangulations can be associated with matrices in $SL(2,\mathbb{C})$ satisfying cocycle conditions. Moreover, there is 1-1 correspondence between the labelling of edges of usual triangulations (described in in Section~\ref{sec:inv_braids_ptolemy}) valued in $\mathbb{C} \backslash\{0\}$, and the labelling of edges of truncated triangulations valued in $SL(2,\mathbb{C})$, for details, see \cite{GoernerZickert}.

\begin{figure}
\centering\includegraphics[width=200pt]{ 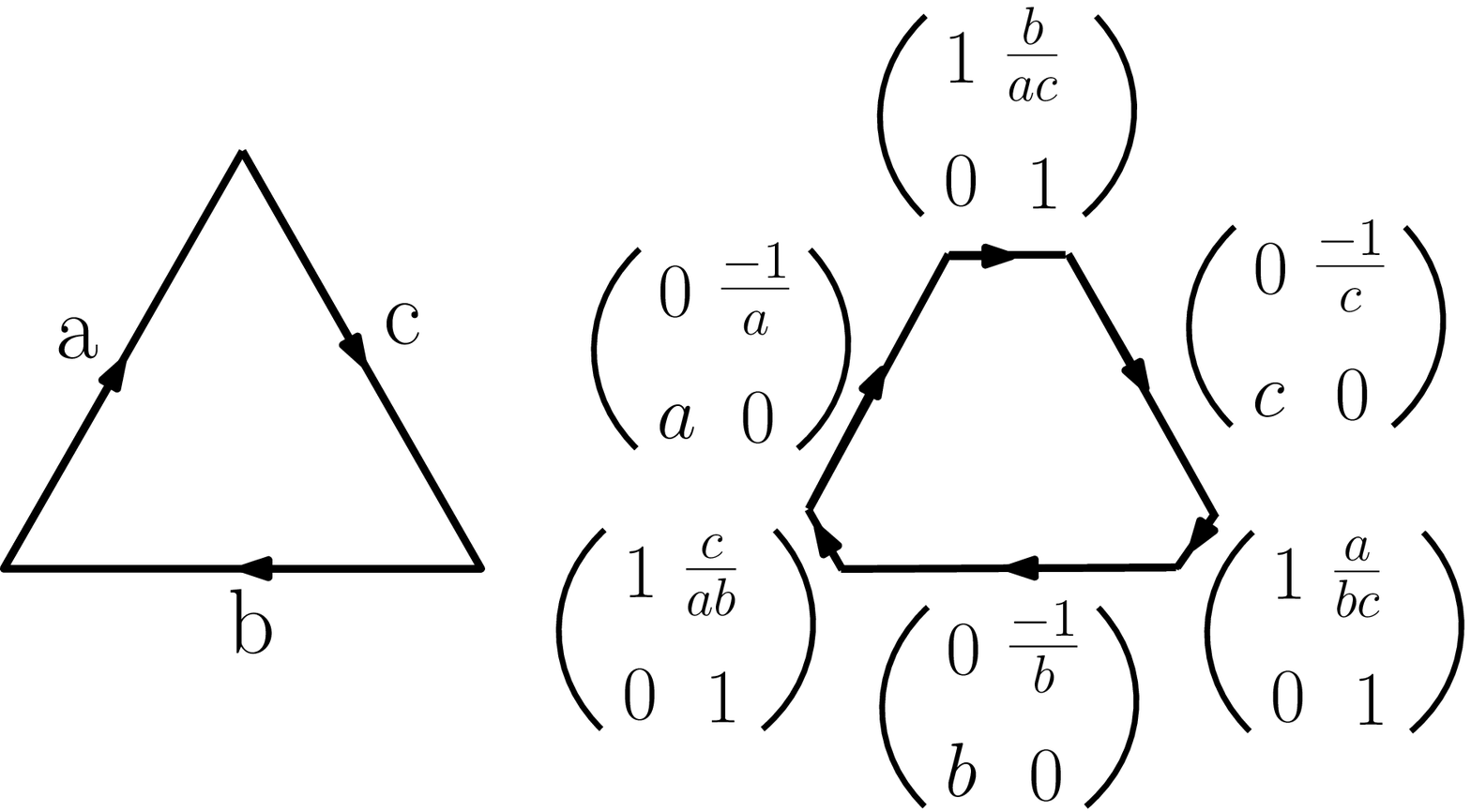}
\caption{Labeling of edges of the hexagon, obtained from triangle}
\label{labeling_hexagon}
\end{figure}

Now we associate each edges with $2\times 2$ matrices $A$ such that $det(A)=1$ for multiplications of matrices, associated edges of triangles or hexagons, to be the identity matrix according to given orientation. For the hexagon obtained from triangle we associate edges with matrices as Fig.~\ref{labeling_hexagon}. By simple calculation, we can see that the multiplication of matrices according to the given orientation is the identity matrix.

\begin{remark}
\begin{enumerate}
\item The hexagon from a triangle has {\em short and long} edges and the forms of matrices, which are associated with short and long edges, are different.
\item Each matrix, associated with an edge as in Fig.~\ref{labeling_hexagon}, is determined by labels of edges of the triangle.
\item When associating edges with matrices the orientation of edges is important.
\end{enumerate}

\end{remark}
For hexagons, which are obtained from triangles in the flip, we associate edges with matrices as Fig.~\ref{labeling_hexagon}.
\begin{figure}
\centering\includegraphics[width=250pt]{ 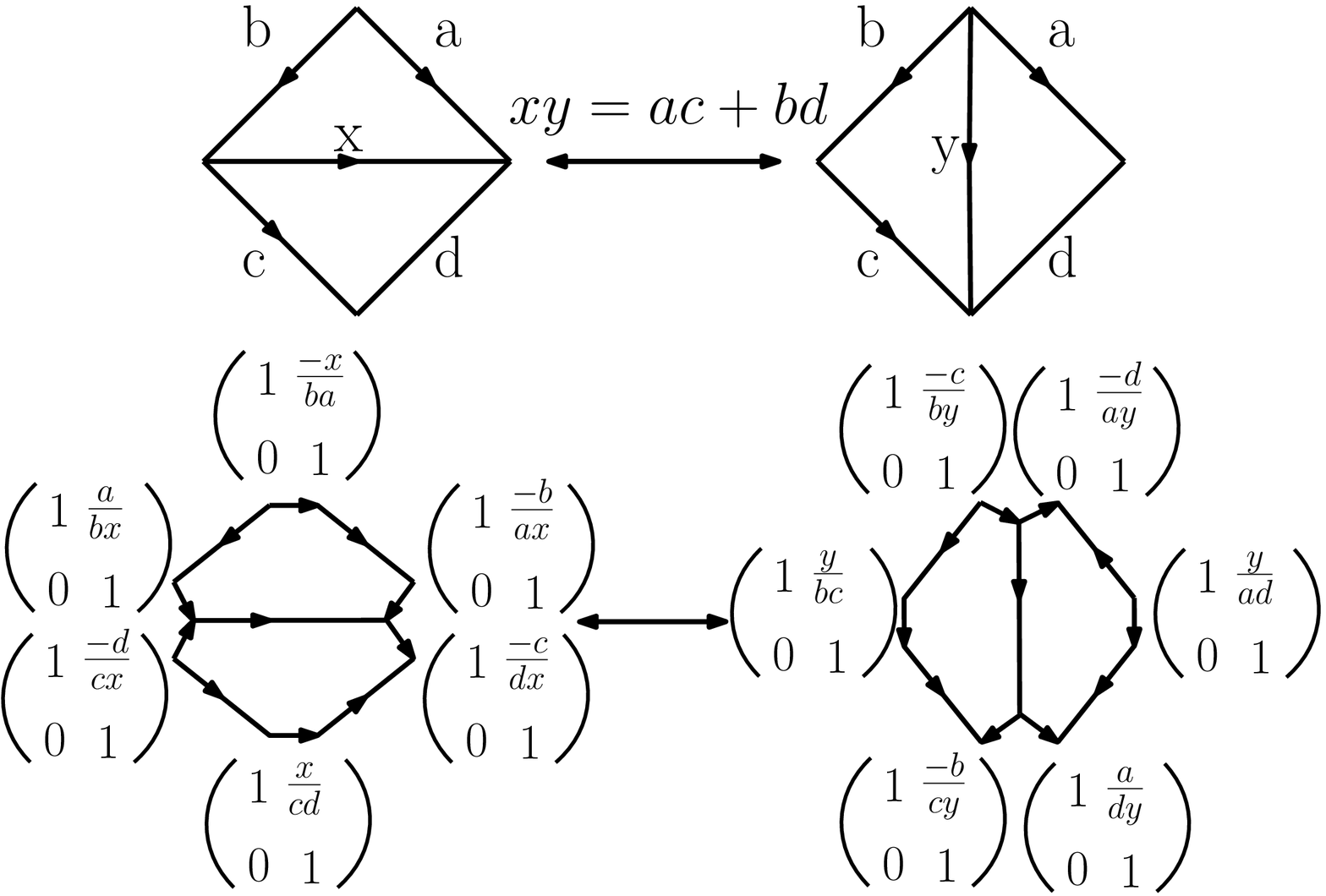}
\caption{Labelling of edges of hexagons, which are obtained from triangles in the flip}
\label{labeling_hexagons_flip}
\end{figure}
We need to assert that the matrices, attributed to short edges in the left-hand side of the figure equal the ones in the right-hand side. Taking into consideration that concatenation of two short edges yields a multiplication of the corresponding matrices, we obtain the following system of matrix equations:

\begin{align*}
\begin{pmatrix} 1 & \frac{-x}{ba} \\ 0 & 1 \end{pmatrix} &= \begin{pmatrix} 1 & \frac{-c}{by} \\ 0 & 1 \end{pmatrix} \begin{pmatrix} 1 & \frac{-d}{ay} \\ 0 & 1 \end{pmatrix}, \\
\begin{pmatrix} 1 & \frac{y}{bc} \\ 0 & 1 \end{pmatrix} &= \begin{pmatrix} 1 & \frac{a}{bx} \\ 0 & 1 \end{pmatrix} \begin{pmatrix} 1 & \frac{-d}{cx} \\ 0 & 1 \end{pmatrix}^{-1}, \\
\begin{pmatrix} 1 & \frac{-y}{ad} \\ 0 & 1 \end{pmatrix} &= \begin{pmatrix} 1 & \frac{-b}{ax} \\ 0 & 1 \end{pmatrix} \begin{pmatrix} 1 & \frac{-c}{dx} \\ 0 & 1 \end{pmatrix}, \\
\begin{pmatrix} 1 & \frac{x}{cd} \\ 0 & 1 \end{pmatrix} &= \begin{pmatrix} 1 & \frac{-b}{cy} \\ 0 & 1 \end{pmatrix}^{-1} \begin{pmatrix} 1 & \frac{a}{dy} \\ 0 & 1 \end{pmatrix}.
\end{align*}

This system is equivalent to the equation $xy = ac+bd$, which is the Ptolemy equation. It is easy to see that this labelling valued in $SL(2,\mathbb{Q}(x_{1}, \dots, x_{X}))$ satisfies the pentagon equation and the following theorem holds, which is analogous to Theorem~\ref{inv_braids_ptolemy}:

\begin{theorem}\cite{GoernerZickert}\label{inv_braids_truncated}
If a loop $T=T_{0}\to T_{k}=T$ in the space of triangulations originates from a braid
which is homotopic to the identity, then the labels valued in $SL(2,\mathbb{Q}(x_{1}, \dots, x_{X}))$ of $T_{k}$ coincide with those of $T_{0}$.
\end{theorem}

%% file: unsolved.tex
\addcontentsline{toc}{part}{Further work and unsolved problems}
\part*{Further work and unsolved problems}
 
Below, we present several unsolved problems we find the most valuable as well as some work
in progress. A much more detailed version will appear in the forthcoming book \cite{upcoming_book}.
 
\section*{A small modification of groups $\Gamma$ admits a simplicial structure}
 
It is well known that braid groups (with some enhancement called steady flow) lead to a simplicial structure which is very useful for calculating homotopy groups of spheres, \cite{Wu}.
 
The simplicial structure naturally arises when one is able to ``delete strands'' and to ``double strands''. For the case of the groups $\Gamma_{4}^{k}$ a hint how to double the strands is presented in Section \ref{subsec:truncated_Zickert}.
 
\section*{Upgrading braids to knots}
 
The groups $G_{n}^{k}$ and $\Gamma_{n}^{k}$ are a sort of braid groups. What is the ``knot counterpart'' for this group? 
 
\section*{Far triangulations}
A generic set of points on the plane (on the sphere) gives rise to a triangulation of this plane (sphere), but in fact it gives rise to {\em many triangulations}. Indeed, $n+1$ points in $S^{n}$ form a simplex if the sphere passing through these points has no other points inside it. But we can consider $(n+1)$-tuple of points containing exactly $p$ points on one side, as it was done in Section \ref{sec:braids_product_Gamma}. This will give rise to triangulations of some cell complexes, which sometimes turn out to be manifolds.
 
Thus we get a collection of cell complexes in the sphere. It would be interesting to understand their connection to Pontrjagin duality, characteristic classes, etc.
 
\section*{Embedding of braids of one type into braids of another type}
 
In Section \ref{sec:higher_gamma} we defined the manifold of triangulations in topological case and in smooth case. It can be done in situations when manifold possesses some structures. Thus, we have many braid groups, and that yields two sorts of questions.
 
Let us write $S\ge S'$ if structure $S$ yields structure $S'$ (say, $S$=smooth, $S'$= topological).
 
\begin{enumerate}
\item The question on the nature of $\pi_0$: is the number of connected components for the structure $S$ the same as for the structure $S'$?
 
\item Is it true that two $S'$-braids equivalent as $S$-braids are equivalent as $S'$-braids?
 
\end{enumerate}

\section*{One-dimensional knots in higher dimensional spaces}

Let us define a {\em 1-dimensional knot in $n-$dimensional space}\index{Knot in higher dimensions} as a curve in the $n-$dimensional space such that in every horizontal $(n-1)-$dimensional section any $n-1$ points should be in general position.

The problem is to study isotopy classes of such knots.

\section*{$G$ and $\Gamma$ structures of invariants of manifolds}

For invariants constructed via statistic sums there should be $G_{n}^{k}$ or
$\Gamma_{n}^{k}-$equivariant versions. For instance, the Viro-Turaev invariant is defined via triangulations and Matveev-Piergallini moves (see \cite{TurCurves, Weaves, mat1, mat2, pierg}). There should be a $\Gamma_{n}^{5}-$version.
 
The question is, how to define $G_{n}^{k}$ and $\Gamma_{n}^{k}-$equivariant homology groups of manifolds?

%% file: appendum.tex
\addcontentsline{toc}{part}{Appendix}
\part*{Appendix}
In the present Appendix we give a brief outline of the free knot theory and parity theory, and then present an overview of some of the most recent developments of the theory, which currently present an open field of work and study. In particular, we concentrate on the notion of {\em manifolds of triangulations}, the construction of {\em 3-free links} arising form the $G_n^k$ groups, and the study of {\em domino tilings} of surfaces, which lead both to groups similar to $\Gamma_n^4$, and to certain invariants of braids.

\section{Free knots, parity, and parity bracket}
\label{sec:free_knots}

Here we give a short description of the theory of free knots, parity, and parity bracket: here we deal only with one parity and one variant of bracket. 

\begin{figure}
\centering\includegraphics[width=300pt]{ 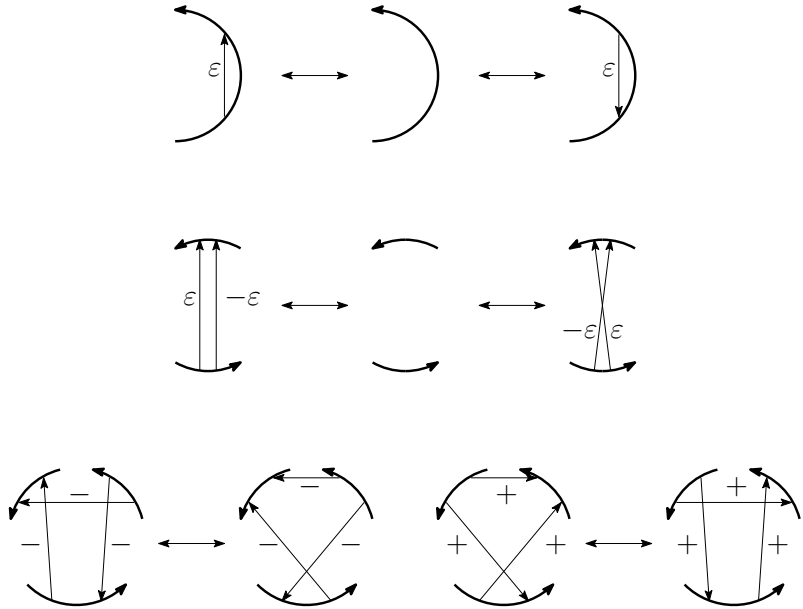}
\caption{Reidemeister moves for Gau{\ss} diagrams. Reidemeister moves for formal Gau{\ss} diagrams are obtained by forgetting the signs and arrows}
\label{fig:chord_reid}
\end{figure}

By {\em formal Gau{\ss} diagrams} we mean Gau{\ss} diagrams without sings/arrows; for Reidemeister moves on Gau{\ss} diagrams see Fig.~\ref{fig:chord_reid}. Consider the equivalence classes $$\langle\text{Formal Gauss diagrams}\rangle/\langle\text{Reidemeister moves on Gau{\ss} diagrams}\rangle.$$ This is the same as $$\langle\text{Framed $4$-graphs}\rangle/\langle\text{Reidemeister moves on framed $4$-graphs}\rangle.$$ Here by {\em framed $4$-graphs} we mean such $4$-graphs that at each crossing of the graph all four halfedges incident to the vertex are split into two pairs, edges in each pair are called {\em formally opposite}. By slightly abusing notation, we also admit circles without vertices as $4$-graphs or connected components of $4$-graphs.

\begin{figure}
\centering\includegraphics[width=300pt]{ 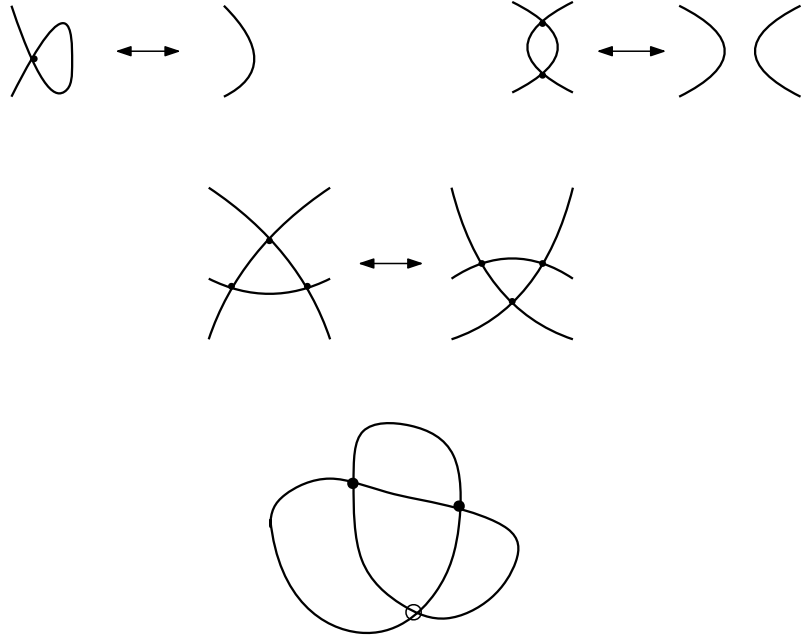}
\caption{Upper part: Reidemeister moves framed 4-graphs; framing is natural: two half-edges are opposite if they constitute a smooth curve. Lower part: an example of a framed 4-graph; a virtual vertex is encircled}
\label{fig:graph_reid}
\end{figure}

Framed $4$-graphs and moves on them are given in Fig.~\ref{fig:graph_reid}. When drawn on the plane, classical crossings are denoted by thick dots; encircled (virtual) crossings are artefacts of projection.

We say that two chords of a formal Gau{\ss} diagrams are {\em linked} if their ends alternate. A chord is {\em (Gau{\ss}) odd} if it is linked with {\em oddly many} other chords; otherwise it is called {\em (Gau{\ss}) even}. Crossings of a framed $4$-graph are named {\em odd} or {\em even} with respect to the parity of the corresponding chords.

Denote the set of Gau{\ss} diagrams by $\mathcal{D}$. The set of all crossings of the diagrams of the set $\mathcal{D}$ will be denoted by $\mathfrak{C}=\mathfrak{C}(\mathcal{D})$. To define a {\em parity}, we decorate every crossing of the set $\mathfrak{C}$ with either 0 (such crossings are called {\em even}) or 1 (such crossings are called {\em odd}). In other words, a parity is a mapping from the set $\mathfrak{C}$ to the set $\mathbb{Z}_2$ which obeys certain conditions.\index{Parity axioms} To be precise: \\

\begin{enumerate}
	\item {\bf Correspondence axiom:} for a move transforming a diagram $D$ into a diagram $D'$ and for every crossing $v$ of the diagram $D$ which has a naturally corresponding to it crossing $v'$ of the diagram $D'$ we have $p(v)=p(v')$;
	\item {\bf $\Omega_1$-axiom:} for a vertex $v$ of a loop (which is the subject of a decreasing first Reidemeister move) we have $p(v)=0$;
	\item {\bf $\Omega_2$-axiom:} for two vertices $v_1, v_2$ of a bigon (which are the subject of a decreasing second Reidemeister move) we have $p(v_1)+p(v_2)=0$;
	\item {\bf $\Omega_3$-axiom:} for three vertices $v_1, v_2, v_3$ of a triangle (which are the subject of a third Reidemeister move) we have $p(v_1)+p(v_2)+p(v_3)=0$.
\end{enumerate}

\begin{definition}
\label{def:traditional_parity}
	A {\em parity} for a knot theory is a mapping $p\colon\mathfrak{C}\to\mathbb{Z}_2$ satisfying the axioms 1--4.\index{Parity}
\end{definition}

An example of a parity is the Gau\ss{}ian parity defined above.

Now we will proceed to the construction of the {\em parity bracket} for knots. But we still need some preliminary considerations on the subject of smoothings\index{Smoothing} of a knot diagram.  \\

Let a framed 4-graph $K$ be a diagram of a free knot. Given a smoothing $s$ we denote the resulting diagram by $K_s$. Here $s$ (also called a {\em state}) is a choice of smoothing for each smoothed vertex (recall that there are two possible smoothings for each vertex). By definition the result of the smoothing of a diagram with the empty set of vertices is the diagram itself.

Let $K$ be a diagram and let $v_1,\dots,v_n$ be all its even crossings with respect to the Gau\ss{}ian parity $p_G$. The smoothing at all even crossings of the diagram is called an {\em even smoothing}\index{Smoothing!even} of $K$. There exist $2^n$ even smoothings. Moreover, if an even smoothing produces a graph with one unicursal component, then such smoothing is called a {\em 1-even smoothing}.

For a given graph we single out the following sets of smoothings:
\begin{itemize}
	\item the set $S$ of all smoothings;
	\item the set $S_{even}$ of all even smoothings;
	\item the set $S_1$ of all smoothings having one unicursal component;
	\item the set $S_{even,1}$ of all 1-even smoothings.
\end{itemize}
The elements of those sets will be denoted by $s, s_{even}, s_1$ and $s_{even,1}$, respectively.

Now consider framed 4-graphs with one unicursal component modulo the equivalence relation generated by the second Reidemeister move. Let us define the linear space $\mathfrak{G}$ as the set of $\mathbb{Z}_2$-linear combinations of such equivalence classes.

\begin{definition}
	\label{def:g_space}
	The linear space $\tilde{\mathfrak{G}}$ is the set of $\mathbb{Z}_2-$linear combinations of framed 4-graphs modulo the following equivalence relations:
	\begin{enumerate}
		\item the second Reidemeister move $\Omega_2$;
		\item the relation $K\sqcup \bigcirc = 0$, i.e. if a framed 4-graph has more than one component and at least one component is trivial, such graph is assumed to be zero in the space $\tilde{\mathfrak{G}}$.
	\end{enumerate}
\end{definition}

We have a natural map $g:\tilde{\mathfrak{G}}\to\mathfrak{G}$ which act by taking to zero all equivalence classes of framed graphs with more than one unicursal component. The map $g$ is an epimorphism of groups.

Now we have the necessary tools to define the bracket invariant.	 Let $K$ be a diagram (a framed 4-graph). Then its parity bracket $[K]$ is defined by the formula $$[K]=\sum_{s_{even,1}}K_{s_{even,1}}.$$

\begin{theorem}
	\label{th:first_bracket_invariant}
	The bracket $[\cdot]$ is an invariant of free knots.
\end{theorem}


The most striking formula for this bracket is $[K]=K$ for diagrams with all chords odd. Since the graph is odd, its bracket $[K]$ consists of only one summand --- graph $K$ itself. This graph is its own minimal representative.


It is crucial, however, to note that in the right-hand side we still have some equivalence (minor, but non-vacuous) for framed $4$-graphs: we allow to apply second Reidemeister moves. It is well known that such diagrams have a unique minimal representative which can be obtained by subsequent reductions.

So, it is important to demand that $K$ is not only {\em odd} (having all crossings odd) but also {\em irreducible}. Thus we can say that a {\em local} minimality yields {\em global} minimality: the fact that no single decreasing Reidemeister move can be applied to the diagram (the second due to irreducibility, the first and the third for the reason of the graph being odd) yields the minimality of the graph. Note that on the left-hand side of this equality, the graph $K$ is considered as a representative of a free knot, and on the right-hand side it is considered as an element from $\mathfrak{G}$, that is, a single graph. That means that an irreducibly odd graph is its own invariant.

Interesting new results on invariants of free knots valued in a free group may be found in \cite{manturov2020freegroup}.

\section{Manifolds of triangulations} \label{sec:manifold_of_triang}
The spaces of triangulations of a given manifold have been widely studied. The celebrated theorem of Pachner~\cite{pach} says that any two triangulations of a given manifold can be connected by a sequence of bistellar moves, or Pachner moves. Below we consider groups which naturally appear when considering the set of triangulations with fixed number of simplices of maximal dimension.

There are three ways of introducing this groups: the geometrical one, which depends on the metric, the topological one, and the combinatorial one. The second one can be thought of as a ``braid group'' of the manifold and, by definition, is an invariant of the topological type of manifold; in a similar way, one can construct the smooth version. We construct a canonical map from the braid group of any $d$-dimensional manifold to $\Gamma_{n}^{d+2}$ thus getting topological/smooth invariants of these manifolds.

\subsection{Geometrical manifold of triangulations}

Fix a smooth manifold $M$ of dimension $d$ with a Riemannian metric $g$ on it.
Let $n\gg d$ be a large natural number. Consider the set of all Delaunay triangulations of $M$ with $n$ vertices (if such triangulations exist). Such Delaunay triangulations are indexed by the sets of vertices $x_{1},\dots, x_{n}$, hence, they form a subset of the configuration space, which we denote by $C(M,n)$. This subset will be an open (not necessarily connected) manifold $M_{g}^{dn}$ of dimension $dn$. We call it the {\em geometrical manifold of triangulations}\index{triangulation manifold!metrical}.

The manifold of triangulations $M_{g}^{dn}$ has a natural stratification. A point of $M_{g}^{dn}$, given by a set of vertices $(x_{1},\dots, x_{n})$, is called {\em generic} if there is no sphere $S^{d-1}$ in $M$ such that exactly $d+2$ points among $x_{k}$ belong to this sphere without any points inside the sphere. The condition means that there is a unique Delaunay triangulation with the set of vertices $(x_{1},\dots, x_{n})$.

We say that a point $(x_{1},\dots, x_{n})$ is {\em of codimension $1$}, if there exists exactly one sphere with exactly $d+2$ points on it and no points inside it; for this configuration there are two possible Delaunay triangulations related to each other by  a flip (Pachner move).   Analogously, {\em codimension $2$ strata} correspond to either one sphere with $d+3$ points or two spheres containing $d+2$ points each.

In codimension $2$ this corresponds to either one point of valency five in Vorono\"{\i} tiling or two points of valency four in Vorono\"{\i} tiling.

Note that the manifold $M_{g}^{dn}$ may be not connected, i.e. there can exist non-equivalent triangulations. Denote the connected components of $M_{g}^{dn}$ by $(M_{g}^{dn})_{1},\dots, (M_{g}^{dn})_{p}$. On the other hand, if one considers the spines of some manifold, they all can be transformed into each other by Matveev--Piergallini moves~\cite{Matveev}.

\begin{definition}
The {\em geometrical $n$-strand braid groups} of the manifold $M_{g}$ are the fundamental groups
$$B_g(M_{g},n)_{j} = \pi_{1}((M_{g}^{dn})_{j}), \quad j=1,\dots, p.$$
\label{def:geometric_n_braid}
\end{definition}

\subsection{Topological manifold of triangulations}

Note that the structure of the geometrical manifold of triangulations heavily depends on the metrics of the manifold $M$. For example, if we take $M$ to be the torus glued from the square $1\times 10$, the combinatorial structure of the manifold of triangulations will differ from the combinatorial structure for the case of the manifold of triangulations for the case of the torus glued from the square $1\times 1$. The {\em topological manifold of triangulations} defined below is independent on the metric.

Consider a topological manifold $M^{d}$. We consider all Riemannian metrics $g_{\alpha}$ on this manifold. They naturally lead to stratified manifolds $M_{g_{\alpha}}^{dn}$ as described above. By a {\em generalised cell} of such a stratification we mean a connected component of the set of generic points of $M_{g_{\alpha}}^{dn}$.

Given a metric $g_\alpha$, we say that two generalised cells $C_{1}$ and $C_{2}$ in $M_{g_\alpha}^{dn}$ are {\em adjacent} it there exist two points, say, $x=(x_{1},\dots, x_{n})$ in $C_{1}$ and $x'=(x'_{1},\dots, x'_{n})$ in $C_{2}$ and a path $x_{t}=(x_{1}(t),\dots, x_{n}(t))$, such that $x_{i}(0)=x_{i}$ and $x_{i}(1)=x'(i)$ such that all points on this path are generic except for exactly one point, say, corresponding to $t=t_{0}$, which belongs to the stratum of codimension $1$.

We say that two generic strata of $M_{g_{\alpha}}^{dn}$ and of $M_{g_{\beta}}^{dn}$ are {\em equivalent} if there is a homeomorphism of $M_{g_{\alpha}}^{dn}\to M_{g_{\beta}}^{dn}$ taking one stratum to the other.

A generalised {\em $0$-cell} of the manifold of triangulations is an equivalence class of generic strata.
Analogously, we define generalised {\em $1$-cells} of the manifold of triangulations as equivalence
classes of pairs of adjacent vertices for different metrics $M_{g_{\alpha}}^{dn}$,
$M_{g_{\beta}}^{dn}$ to the pair of adjacent vertices equivalent to the initial ones.
In a similar manner, we define generalised {\em $2$-cells} as equivalence classes of discs for
metrics $M_{g_{\alpha}}^{dn}$ such that:
\begin{enumerate}
\item vertices of the disc are points in $0$-strata;
\item edges of the disc connect vertices from adjacent $0$-strata;
each edge intersects codimension $1$ set exactly in one point;
\item the cycle is spanned by a disc which intersects codimension $2$ set exactly
at one point;
\item equivalence is defined by homeomorphism taking disc to disc,
edge to an equivalent edge and vertex to an equivalent vertex and
respects the stratification.
\end{enumerate}

Thus, we get the $2$-frame of some manifold $M^{dn}_{top}$. This manifold might be disconnected.

\begin{definition}
The {\em topological $n$-strand braid groups} of the manifold $M$ are the fundamental groups
$$B_t(M,n)_{j} = \pi_{1}((M^{dn}_{top})_{j}), \quad j=1,\dots, q.$$
\label{def:topological_n_braid}
\end{definition}

\subsection{Combinatorial manifold of triangulations}

We did not actually discuss the {\em existence} of such triangulation for a given manifold and a given set of points, requesting only that the number of points is sufficiently large and points are sufficiently dense. It turns out that this condition is not always sufficient, as shows the work \cite{BDGM} by Boissonnat et al. Even for a large number of points and a Riemannian manifold with metric arbitrary close to Euclidian one, there may not exist a Delaunay triangulation.

One may impose additional restrictions on the vertex set or on the manifold itself to overcome this difficulty, but to work with the most general situation we shall consider the third notion of manifolds of triangulations: the {\em combinatorial manifold of triangulations} which we denote by $M^{dn}_{comb}$. We construct the 2-frame of the manifold $M^{dn}_{comb}$ because we are interested only in the fundamental group of it.

First, we fix $n$ points on the manifold $M^d$ and the triangulations of the manifold with vertices in those points.

The vertices (0-cells) of the manifold $M^{dn}_{comb}$ are all the triangulations of the manifold $M$ (we identify triangulations which are combinatorially equivalent). Two triangulations are connected by an edge (a 1-cell of the frame) if and only if they differ by a {\em flip} --- a Pachner move.

Finally, the 2-cells of the frame are chosen to correspond to the configurations of codimension $2$ for Delaunay triangulations. There are two types of such configurations. One corrsponds to two {\em independent} flips $\alpha,\beta$ (flips related to different edges), here we define a quadrilateral consisting of subsequent flips corresponding to $\alpha,\beta,\alpha^{-1},\beta^{-1}$. The other type corresponds to all possible simplicial polytopes with $d+3$ vertices inscribed in the unit sphere. With each polytope of such sort, we associate a relation of length $d+3$ like one that appears in the definition of $\Gamma_n^k$.

\begin{definition}
The {\em combinatorial $n$-strand braid groups} of the manifold $M$ are the fundamental groups
$$B_c(M,n)_{j} = \pi_{1}((M^{dn}_{comb})_{j}), \quad j=1,\dots, q.$$
Here $(M^{dn}_{comb})_{j},\ j=1,\dots,q$ are the connected of the complex $M^{dn}_{comb}$.
\label{def:combinatorial_n_braid}
\end{definition}

\vspace{5mm}

Our main result is the following~\cite{FMN}.
\begin{theorem} \label{thm:gamma_geometrical_invariant}
Let $M^{d}$ be a smooth manifold of dimension $d$ and $B(M_{g},n)_{j}$, $j=1,\dots,p$ be its geometrical, or topological, or combinatorial braid groups. Then for any $j$ there is a well defined mapping $B(M_g,n)_j\to \Gamma_{n}^{d+2}$.
\end{theorem}

The mapping is constructed as follows: given a loop $\alpha$ in general position in $M^{dn}$, we indicate the non-generic configurations on it. Each of these non-generic configurations corresponds to a Pachner move which defines a generator in the group $\Gamma_n^{d+2}$. The product of these generators is the image of the loop $\alpha$ in $\Gamma_n^{d+2}$.

Hence, we get invariants of smooth metrical manifolds and of topological manifolds: these invariants are the corresponding images of the fundamental group of the manifold of triangulations in the corresponding group $\Gamma^k_n$.

One can proceed playing the same game with manifolds having other structures (complex, contact, symplectic, K\"ahler, or related to an action of some group). The definition of the corresponding manifolds of triangulations is left to the reader.

\section{The 3-free links} \label{sec:3-free}

In the present section we give an outline of the theory of {\em 3-free links} and construct an invariant of conjugacy classes of closed braids valued in 3-free links. 

%
%
%

\subsection{The notion of 3-free links}
\label{sec:3-free_links}

As free links may be thought as closures of free braids on some number of strands, the 3-free links may be considered as closures of some sort of free 3-braids, the geometric object behind the groups $G_n^3$. Let us give some explicit definitions.

\begin{definition}
A {\em regular 6-graph} is a disjoint union of regular 6-valent graphs (possibly with loops and multiple edges) and circles. Here we call the circles {\em cyclic edges of the 6-graph}.
\end{definition}

\begin{definition}
A {\em framed 6-graph} is a regular 6-graph such that for each 6-valent vertex the 6 half-edges incident to this vertex are divided into 3 pairs of {\em formally opposite}.
\end{definition}

Let us call two edges $e_0, e_1$ of a framed 6-graph {\em equivalent} if there exists a sequence of edges $e_0=b_1, b_2, \dots, b_n=e_1$ such that for each $i$ the edges $b_i$ and $b_{i+1}$ are opposite. The equivalence class of edges is called a {\em unicursal component} of the graph. A cyclic edge also is a unicursal component. 

\begin{definition}
An {\em oriented framed 6-graph} is a framed 6-graph such that each of unicursal component is oriented.
\end{definition}

The last definition yields that at each 6-valent vertex there are three incoming half-edges and three outgoing half-edges.

\begin{definition}
A {\em 3-free diagram} is an oriented framed regular 6-graph	such that at each vertex three incoming half-edges are ordered.
\end{definition}

In the same way {\em regular 6-graphs with ends} and {\em 3-diagrams with ends} may be defined. In that case we allow the graphs to have 1-valent vertices.

\begin{remark}
When drawing a diagram on the plane we always assume the ordering to be inherited from the plane: the leftmost component is the first, the middle one goes after it, and the rightmost one is the last).
\end{remark}

\begin{figure}
\centering\includegraphics[width=400pt]{ 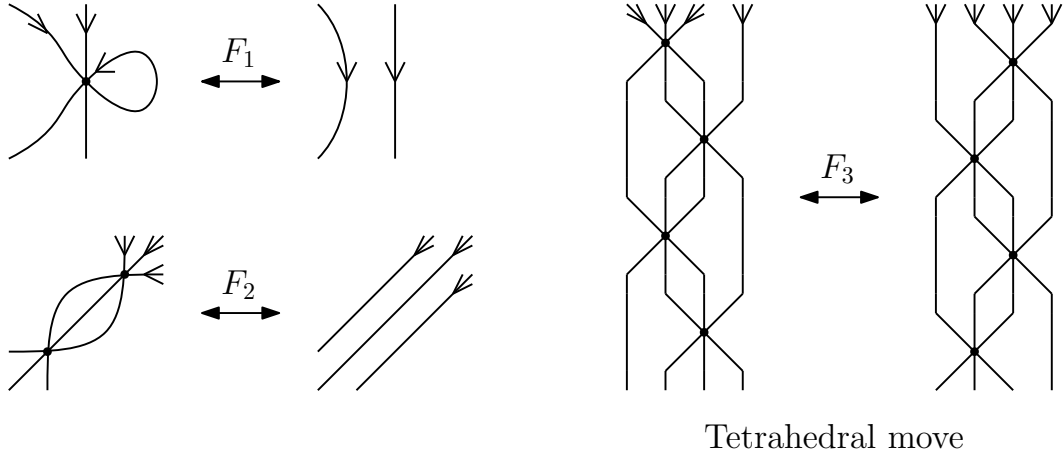}
\caption{3-free moves}
\label{fig:3_moves}
\end{figure}

We consider the set of moves on regular 6-graphs depicted in Fig~\ref{fig:3_moves}.

\begin{definition}
	A {\em 3-free link} (resp., {\em 3-free link with ends}) is an equivalence class of 3-free diagrams (resp., 3-free diagrams with ends) modulo the three 3-free moves depicted in Fig.~\ref{fig:3_moves}.
\end{definition}

Note that these moves do not change the number of unicursal components of a diagram.

\begin{definition}
	A {\em 3-free knot} (resp., {\em 3-free knot with ends}) is 3-free link (resp., with ends) with one unicursal component.
\end{definition}

\subsection{Mapping conjugacy classes of closed braids to 3-free links}
Our first goal is to construct a mapping from the set of conjugacy classes of closed braids to 3-free links (with ends). Note that in case we only have to deal with the second and third Reidemeister moves, and conjugations of the braid.

Consider a closed braid $K$. We may assume that the diagram of $K$ lies in some annulus $A$ on a plane $\Pi$. Fix a line $l$ orthogonal to that plane and intersecting it inside the inner circle of $A$. 

Consider the family of halfplanes whose boundary is the line $l$. Let us naturally parametrise this family by an angle $\varphi\in[0,2\pi]$ and denote the family by $\hat\Pi=\{\Pi_{\varphi}\}$. By a small deformation of $K$ we may assume that the intersection of $K$ and a plane $\Pi_{\varphi}$ is a finite set of points. Let us say that for two angles $\varphi_1, \varphi_2, |\varphi_1-\varphi_2|<\varepsilon,$ a set of points $A_1=K\cap\Pi_{\varphi_1}$ is {\em after} the set of points $A_2=K\cap\Pi_{\varphi_2}$ if $\varphi_2>\varphi_1$ (angles 0 and $2\pi$ are considered equal). 

Now consider the moduli space $M$ of pairs of points on $K$ lying on the same straight ray $m\subset\Pi_{\varphi}\in\hat{\Pi}$ (the origin of the ray may lie on the boundary of $\Pi_{\varphi}$). $M$ is a 1-dimensional manifold with boundary. \newline

 We orient components of this manifold in the following way. For each $x\in M$ consider the corresponding $x_1,x_2\in K$. They lie on some plane $\Pi_\varphi$. It naturally defines two halfspaces of the ambient space $\mathbb{R}^3$. Consider two tangent vectors to $K$ in those points. If their endpoints lie in the same halfspace defined by the plane $\Pi_\varphi$, we orient the tangent vector at $x$ downwards, otherwise --- upwards. If one of the vectors is horizontal, we define orientation to preserve continuity.

Now we need to identify some triples of points of the space $M$. To be precise, we consider such triples of points $(i,j),(i,k),(j,k)$ that the points $i,j,k$ lie on the same straight line $m\in\hat{\Pi}$ (the point $j$ lies between the points $i$ and $k$). We identify them, and define partial ordering so that the component containing the point $(i,k)$ lies {\em between} the other two components. The ordering of the triple (that is, $(i,j), (i,k), (j,k)$ or $(j,k), (i,k), (i,j)$) is defined as follows. 

\begin{figure}
\centering\includegraphics[width=12cm]{ 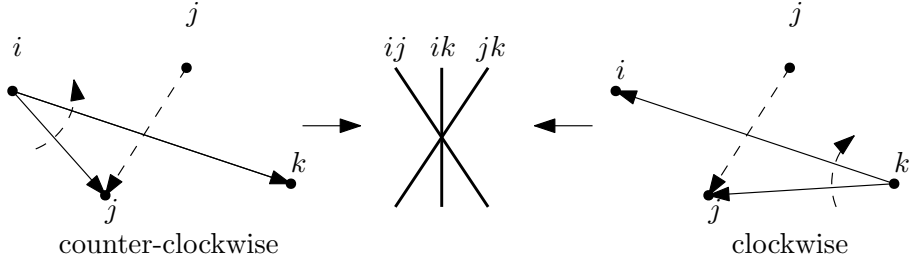}
\caption{If even number of the strands $i,j,k$ go up}
\label{knot_to_3-free1}
\end{figure}

\begin{figure}
\centering\includegraphics[width=12cm]{ 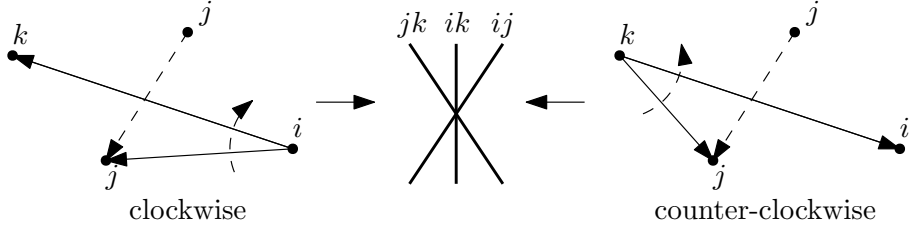}
\caption{If odd number of the strands $i,j,k$ go up}
\label{knot_to_3-free2}
\end{figure}

Consider the points $i',j',k'$ of the closed braid $K$ which lie on the same plane a little {\em after} than the line we are considering. Let $u=\pm 1$ be the sign of the frame $(i',j'),(i'k')$. Now we define $s(i,j,k)$ as $u(-1)^{\uparrow}$, where $\uparrow$ is the number of the points $i,j,k$ where the braid closure is oriented to the ``earlier'' halfspace (``go up''). Finally, we set the order $(i,j), (i,k), (j,k)$ if $s(i,j,k)>0$, and $(j,k), (i,k), (i,j)$ otherwise, see Fig~\ref{knot_to_3-free1}, \ref{knot_to_3-free2}. \newline

This way we have obtained a diagram of a 3-free link, which we denote by $f(K)$.

\begin{theorem}
\label{thm:closed_braids}
	The mapping $f$ is a correct mapping from the set of conjugacy classes of closed braids into the set of 3-free links with ends.
\end{theorem}

This theorem holds for the following reason. We need to check what happens when the knot goes through a codimension 1 singular position.

There are two types of such singularities:

\begin{enumerate}
	\item small deformation of the closed braid $K$ destroys the ray on which three points lie. This situation gives the second move $F_2$ on 3-diagrams;
	\item closed braid goes through a configuration where {\em four} points lie on the same ray. This situation gives the third (tetrahedral) move $F_3$.
\end{enumerate}

Essentially, that means that the closures of two braids with (Reidemeister) equivalent diagrams yield equivalent 3-free link diagrams. \newline

Finally, we need to understand what happens when we deal with closures of two braids which differ by conjugation (a non-Reidemeister equivalence of braids). But note that the closure of a conjugated braid may be transformed into the closure of the original braid via transformations, appearing as Reidemeister moves on the diagram of the closed braid, and such transformations produce only the singularities described above.

\section{Domino tilings and groups $\tau_{N,M}^4$ and $T_{N,M}^4$} \label{sec:domino_tiling}

We have seen that the groups $\Gamma$ naturally act on the set of {\em triangulations} of a given $2$-manifolds. But what about other tilings?

In the present addendum, we review a little the paper by D. Thurston \cite{thurston} which originates from the brilliant idea of Jones \cite{jones} which connects links (in fact, Legendrian links) and domino tilings. The aim of this section is to show that this area of research is related to groups, actually, some groups similar to $\Gamma$'s. 

In fact, this turns out to be a tip of a new iceberg where knots appear threefold:
\begin{enumerate}
\item as Legendrian links which correspond to tilings;
\item as braids which correspond to moving points;
\item through the interrelation between domino tilings, Thompson group and knots.
\end{enumerate}

By comparing 1 and 2 we see that {\em braids act on knots} (in fact, they act on a very specific Legendrian link in the projective space, but this can be easily expanded in many directions). The last item, the Thompson group \cite{jones}, is not even defined in the present paper. Alas! The paper can not be endless! We just wanted to emphasise again and again that groups corresponding to some local transformations ($G, \Gamma, T,$ and $\tau$) play crucial roles everywhere in low-dimensional topology and their relation to new and well known objects is still to be investigated.

In this section we give a brief overview of domino tiling theory, originally studied by D. Thurston \cite{thurston}, which leads to the definition of groups $\tau_{N,M}^4$ and $T_{N,M}^4$, in some regards similar to the groups~$\Gamma_n^4$.

\subsection{Domino tilings}

Assume we have a $2$-surface $\Sigma$ (part of the plane or not, with or without boundary) split into dominoes, that is, quadrilaterals with two ``short'' sides and two ``long'' sides.\index{Domino tiling}

\begin{figure}
\centering\includegraphics[width=250pt]{ 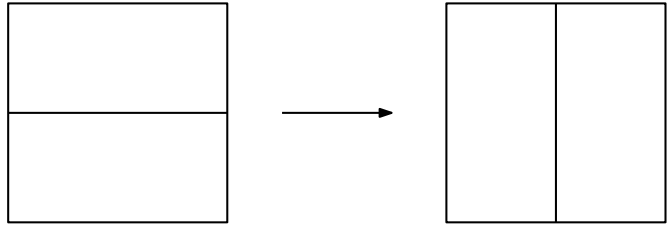}
\caption{The $(2,2)$-move}
\label{2-2-move}
\end{figure}

There are:

\begin{enumerate}
\item Natural moves which transform one tiling into another; the most interesting is the $(2,2)$-move, which is shown in Fig.~\ref{2-2-move};
\item Natural way to associate a projection of (Legendrian) knot in the $1$-tangent surface $T_{S^1} \Sigma$ to each domino tiling.
\end{enumerate}

\begin{figure}
\centering\includegraphics[width=250pt]{ 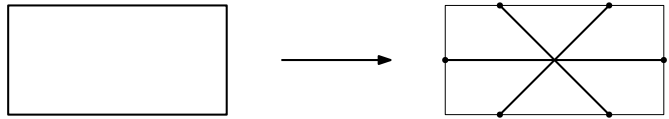}
\caption{Domino and a 3-vertex}
\label{domino_polygon}
\end{figure}

In this construction each domino gives rise to a crossing of three lines. In a sense that means that a domino represents a {\em hexagon}, and the curves connect the opposite vertices of said hexagon. \newline

\begin{figure}
\centering\includegraphics[width=170pt]{ 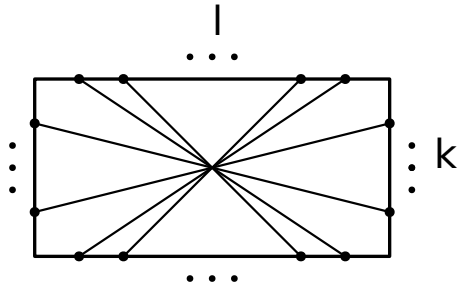}
\caption{A domino giving a $2(k+l)$-gon}
\label{n-gon_domino}
\end{figure}

Naturally, this construction may be extended and generalised. Consider a domino and put $k$ vertices on each of its short sides, and $l$ vertices on each of its long sides. That gives us $2(k+l)-$gon. Now connect the opposite vertices of this polygon, obtaining a $k+l$ crossing, see Fig.~\ref{n-gon_domino}.

This approach leads to a richer set of diagrams, with the usual domino tiling theory being a partial case with $k=1, l=2$.

\begin{figure}
\centering\includegraphics[width=250pt]{ 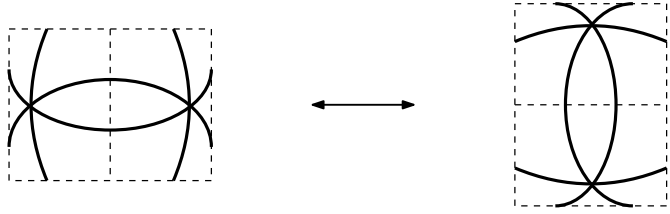}
\caption{The $(2,2)$-move for triple diagrams}
\label{domino_flip}
\end{figure}

The nature of the $(2,2)$-move becomes clearer when we consider the dual picture as shown in Fig.~\ref{domino_flip}: when we associate a triple crossing to each domino placing a vertex on each short side and two vertices on each long side. Then, the $(2,2)$-move (which is essentially just the change of breaking of a square into two dominoes) becomes the following move on diagrams with triple crossings (the square is a bit distorted for better visibility), see Fig~\ref{domino_flip}.

\begin{figure}
\centering\includegraphics[width=250pt]{ 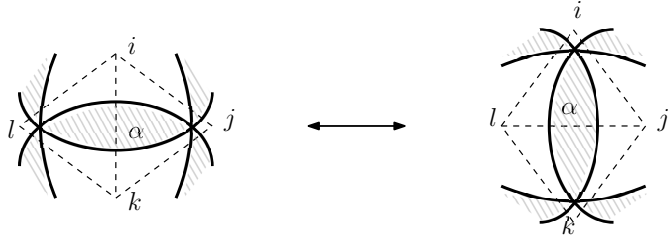}
\caption{The $(2,2)$-move and Ptolemy flip}
\label{domino_flip_square}
\end{figure}

Furthermore, one can consider the ``dual diagram'' to the $(2,2)$-move: consider its checkerboard colouring and attribute vertices to regions of one colour. The result is shown in Fig.~\ref{domino_flip_square}: if we connect the vertices $i,j,k,l$ separated by a single region of the other colour (labelled $\alpha$), the $(2,2)$-move becomes the usual Ptolemy flip for the appearing quadrilateral.

Note that there are five regions of the other colour in this figure, not four. In that regard the move is not ``symmetric'' (it works only with vertices of some one colour in each case).

\begin{figure}
\centering\includegraphics[width=250pt]{ 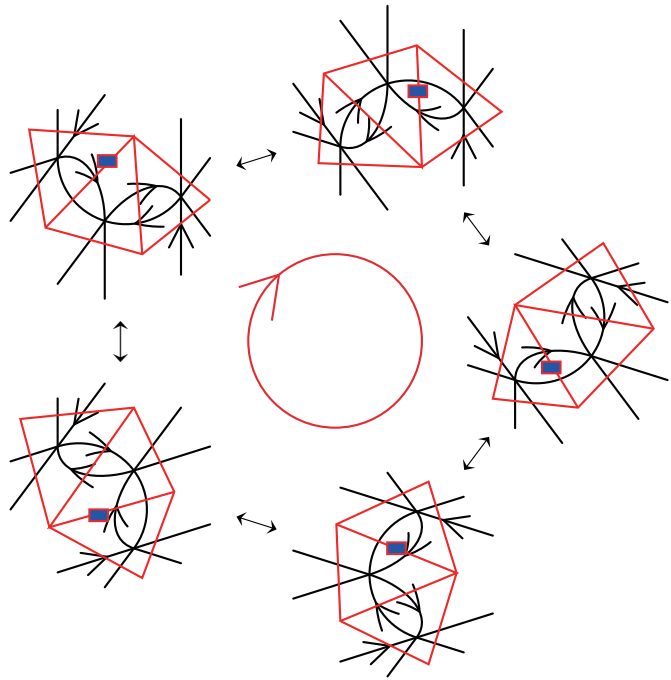}
\caption{Pentagon relation}
\label{fig:pentagon}
\end{figure}

\begin{figure}
\centering\includegraphics[width=300pt]{ 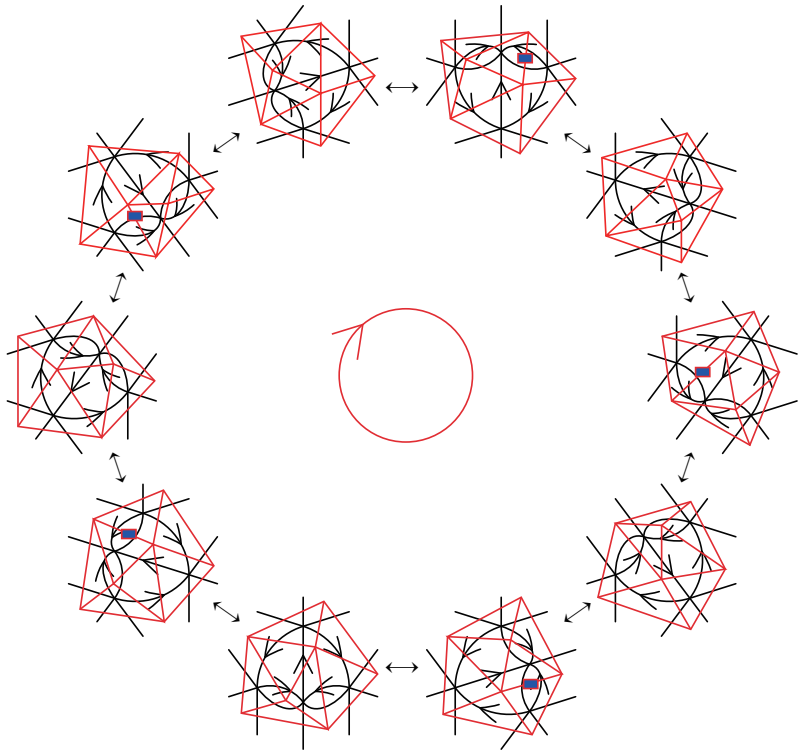}
\caption{Decagon relation}
\label{fig:decagon}
\end{figure}

Consider the following two types of ``relations between relations'', that is, two sequences of $(2,2)$-moves leading back to the beginning (apart from the moves being involutions), see Fig.~\ref{fig:pentagon}, \ref{fig:decagon}.

D. Thurston \cite{thurston} proposed the following conjecture: 

\begin{conjecture}
Consider a cell complex, such that its vertices (0-cells) are triple diagrams, edges (1-cells) are the $(2,2)$-moves, and faces (2-cells) are 
\begin{itemize}
\item squares coming from far commutativity,
\item pentagons as shown in Fig.~\ref{fig:pentagon},
\item decagons as shown in Fig.~\ref{fig:decagon}.	
\end{itemize}
	Then this complex is simply connected.
\end{conjecture}

Now let us study the pentagon and decagon relations in detail. 

We begin with  Fig.~\ref{fig:pentagon}, the pentagon. 
It shows the sequence of $(2,2)-$moves and the dual graph, as described before. Therefore, this relation is just the usual pentagon relation which is known from Ptolemy theory (and from $\Gamma_n^4$ theory). The diagonal being flipped is always marked with a small filled rectangle. 

Note, that this figure cannot appear in a single domino tiling. In other words, there is no domino tiling such that its triple diagram is as shown in any of the five steps of the transformation.

\begin{figure}
\centering\includegraphics[width=250pt]{ 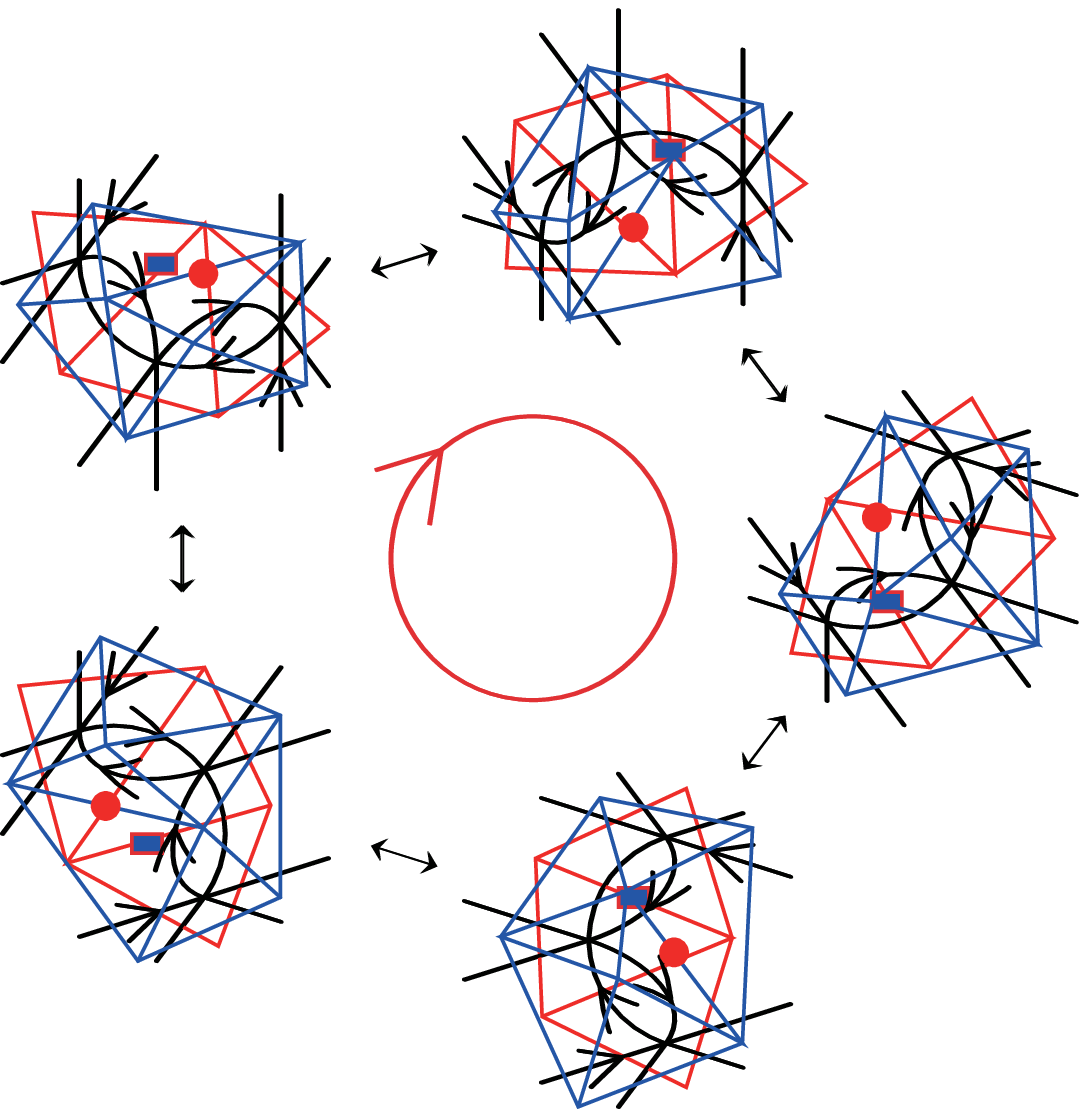}
\caption{Pentagon relation: both dual graphs}
\label{fig:pentagon_both}
\end{figure}

\begin{figure}
\centering\includegraphics[width=300pt]{ 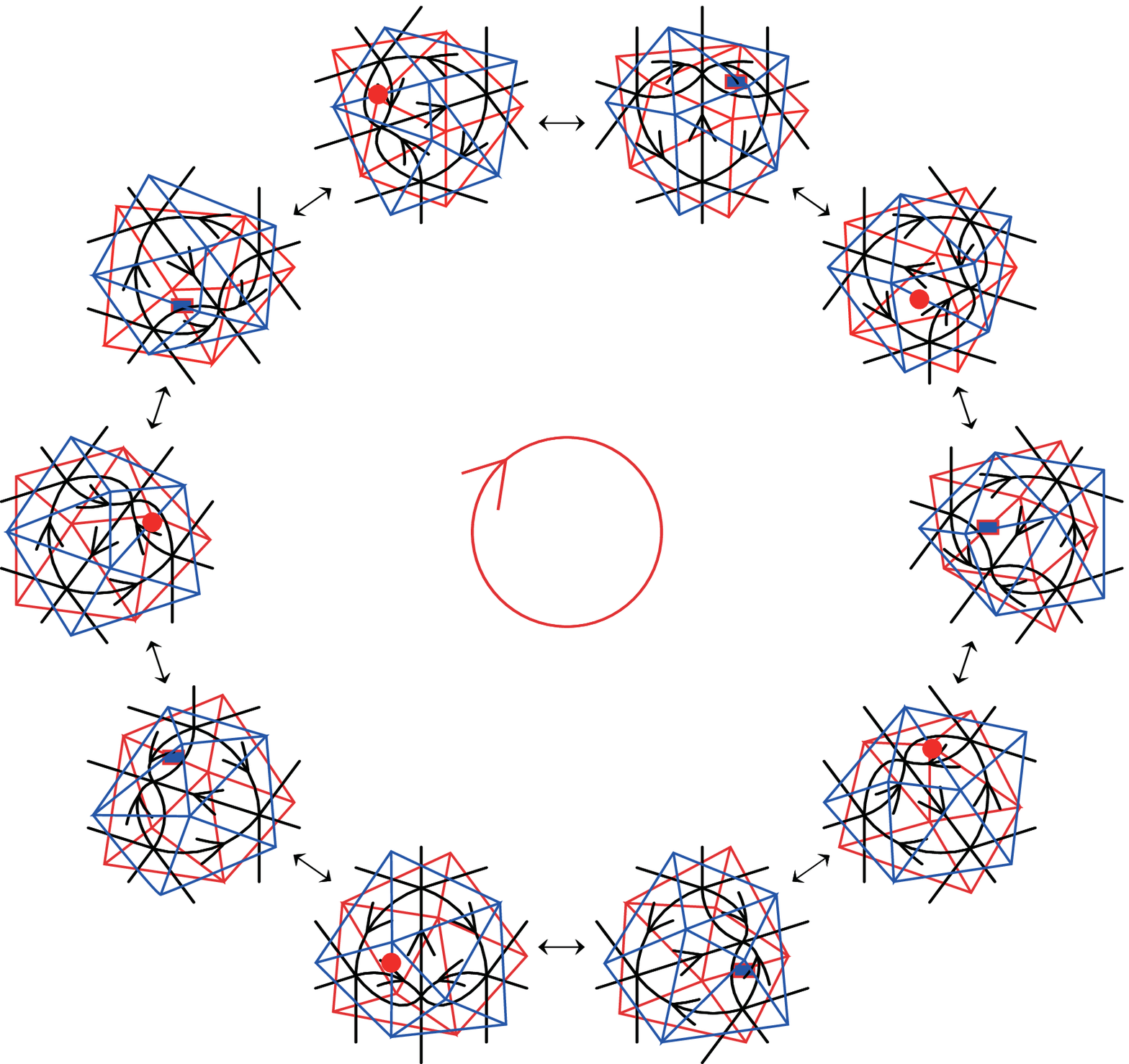}
\caption{Decagon relation: both dual graphs}
\label{fig:decagon_both}
\end{figure}

Fig.~\ref{fig:decagon} is trickier. Again, the sequence of $(2,2)$-moves is clear, but the dual graph is changed by a flip {\em only every second step}. To understand what happens, let us consider the second set of regions (of the ``other colour'' if we consider the checkerboard colouring of the regions) as well, see Fig~\ref{fig:pentagon_both}, \ref{fig:decagon_both}.

Now we can see that each ``missed'' step of the transformation is also a Ptolemy flip but for the {\em second} dual graph (the diagonal being flipped is marked with a circle).

In some sense we can say that the decagon relation is two pentagons interlacing one another. As in case of the pentagon relation, this transformation does not arise from a single domino tiling.

\begin{figure}
\centering\includegraphics[width=250pt]{ 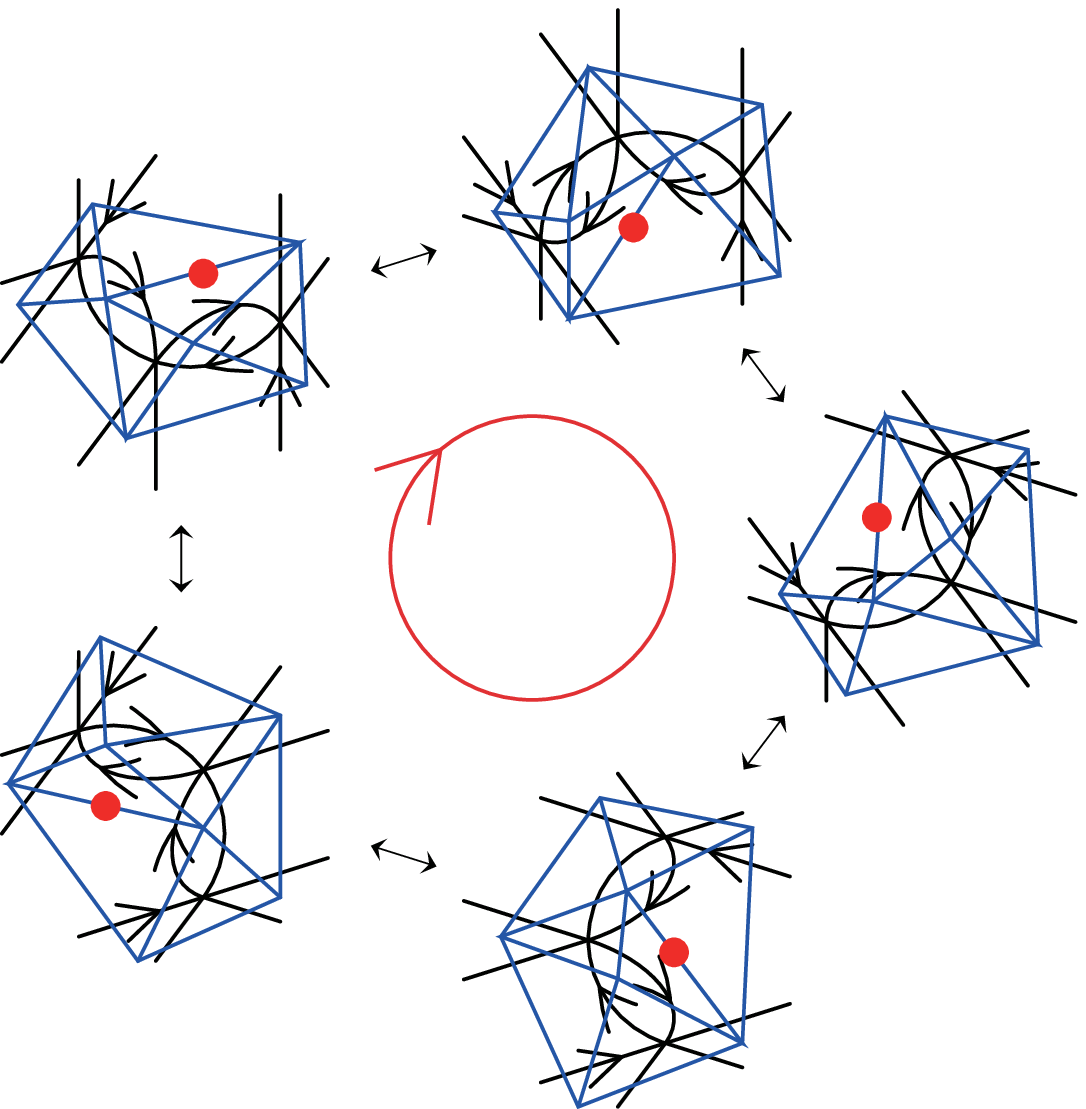}
\caption{Second dual graphs for pentagon}
\label{fig:pentagon_blue}
\end{figure}

\begin{figure}
\centering\includegraphics[width=300pt]{ 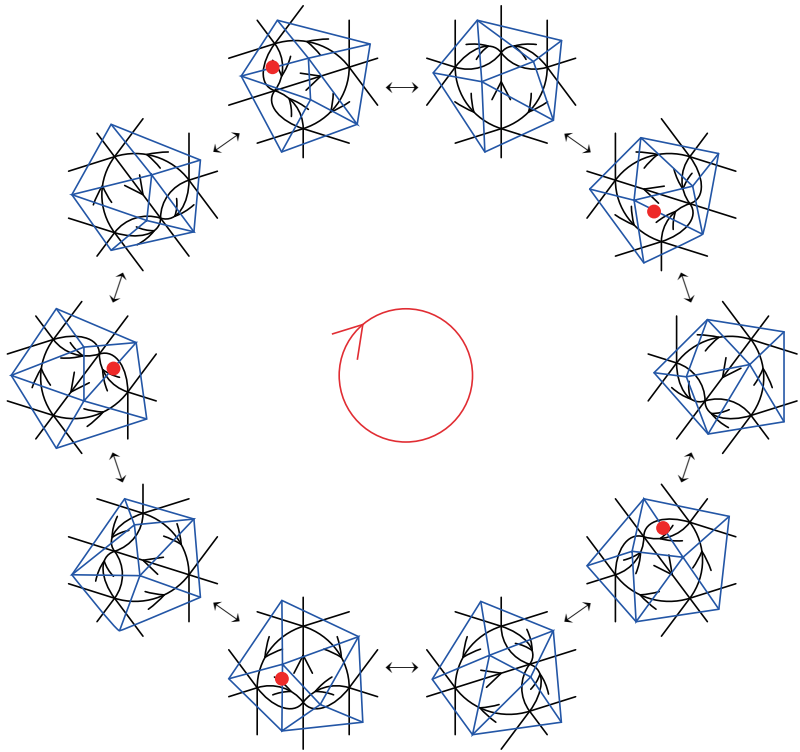}
\caption{Second dual graphs for decagon}
\label{fig:decagon_blue}
\end{figure}

Let us look once again at the second graph, this time without its sibling, see Fig.~\ref{fig:pentagon_blue}, \ref{fig:decagon_blue}.	

Fig.~\ref{fig:decagon_blue} does not give us anything new (again, every second step is missing), but Fig.~\ref{fig:pentagon_blue}, the pentagon relation, stays interesting.

As we can see, this time there is no evident pentagon: there are too many vertices. The graph has 7 vertices: five ``exterior'' ones forming a pentagon, but two ``interior'' ones make the partitioning of the pentagon finer, adding new edges.

Nevertheless, the transformation remains a sequence of Ptolemy flips in the language of the second dual graph. This time, though, it is not the (Ptolemy) pentagon relation, but something new and specific to this theory.

\subsection{The $\tau$-group} \label{sec:tau}

\input{T-group.tex}

\subsection{Braid invariants appearing from tilings}
Let us consider a 2-manifold $M^2$ and a braid $\beta$ in $M^2\times\mathbb{R}$. We may always assume that the braid connects a set of $n$ points on the submanifold $M\times\{1\}$ with a set of $n$ points on the submanifold $M\times\{0\}$ and we shall think of the braid as going down along $\mathbb{R}$. In case of {\em pure braids} these sets of points coincide under the natural identifications of the aforementioned submanifolds.

As the time goes --- as we descend along $\mathbb{R}$ --- the $n$ points in the section of the braid undergo some dynamics. Let us consider the submanifold $\Pi_c=M\times\{c\}$ and the corresponding set of points $\beta_c$. With this set we associate the {\em Vorono\"i diagram} $V_c$ on the submanifold $\Pi_c$. It is a trivalent graph, in other words a 1-foam.

So with every section of the braid we associated a 1-foam. 

Now let us do this operation for {\em each} section of the braid. We obtain a collection of 1-foams, essentially getting a movie of a 2-foam. Let us denote this foam by $F(\beta)$.

\begin{theorem}
\label{thm:braid_foam}
Let $\beta_1, \beta_2$ be two isotopic braids in $M\times\mathbb{R}$, and $F(\beta_1), F(\beta_2)$ be the corresponding 2-foams. Then the foams $F(\beta_1), F(\beta_2)$ may be connected by a sequence of the following moves:
\begin{enumerate}
\item Matveev-Piergallini move,
\item insertion of a 2-foam, spanned by the Vorono\"i diagrams of the steps of the pentagon relation.	
\end{enumerate}

\end{theorem}

\begin{figure}
\centering\includegraphics[width=300pt]{ 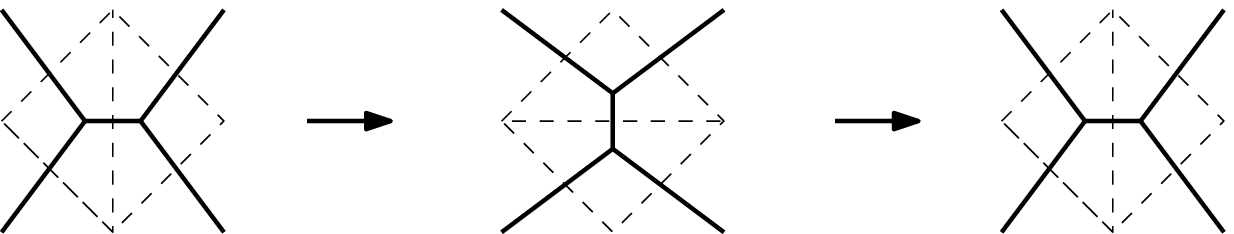}
\caption{The foams giving the Matveev-Piergallini foam skeleton}
\label{matveev-piergallini}
\end{figure}

The Matveev-Piergallini in the language we are now using may be described as an insertion of a foam, spanned by three 1-foams, which are the steps (the $(2,2)$-moves) of a loop in the space of triple diagrams, see Fig.~\ref{matveev-piergallini}.

Note that this move consists of two $(2,2)$-moves. Each move produces a singular point in the constructed 2-foam.

\begin{figure}
\centering\includegraphics[width=220pt]{ 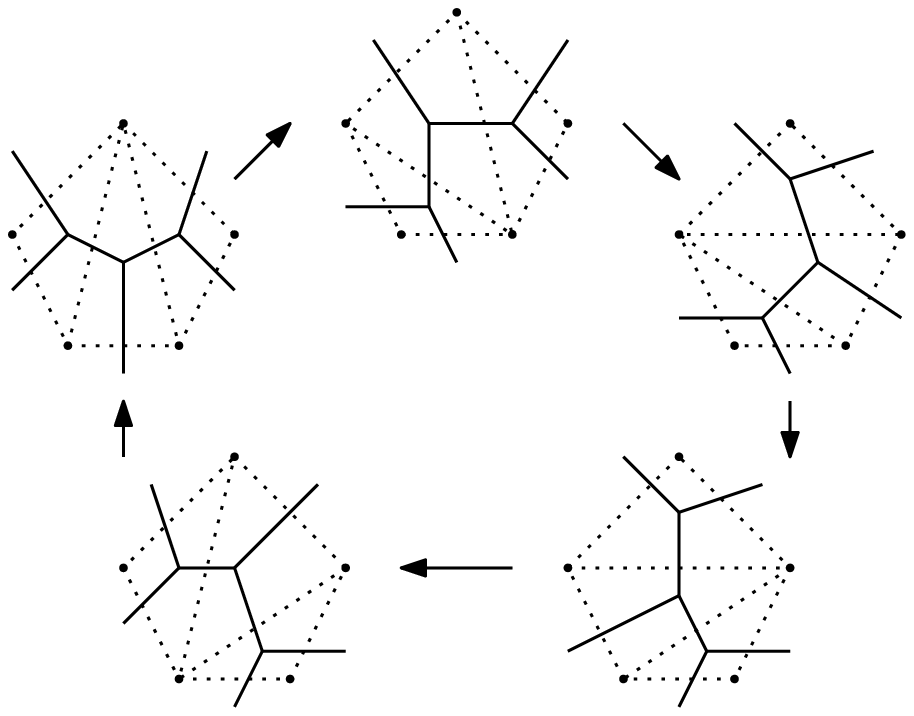}
\caption{Transformation of 1-foams during the pentagon relation steps}
\label{pentagon_voronoi}
\end{figure}

Figure~\ref{pentagon_voronoi} shows the pentagon transformation and the corresponding sketches of Vorono\"i diagrams. The foam mentioned in Theorem~\ref{thm:braid_foam} is constructed by placing these diagrams ``one under the other'' (with the sixth diagram being the same as the first one), and connecting them by a foam (each transition yielding a single singular point of the 2-foam).

Next step is to transform foams into something knot-like, to be precise, into tangles.

Let us go back to the foam $V_c$. It is a trivalent graph, so we perform the following transformation: we double each edge away from a small neighbourhood of each vertex, drawing a parallel line on the right side of the edge, and one --- on the left side of the edge. Finally, at each vertex we connect ``opposite'' lines, producing a triple crossing.

\begin{figure}
\centering\includegraphics[width=300pt]{ 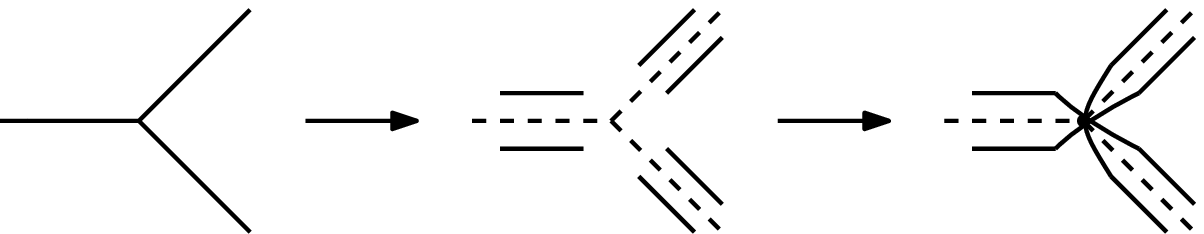}
\caption{From a trivalent vertex to a triple crossing}
\label{trivalent_vertex}
\end{figure}

This operation gives us a 3-free tangle $\hat{T}_c$ for each section of the braid (3-free tangles are defined analogously to 3-free links, see Section~\ref{sec:3-free_links}). Resolving each triple crossing into three usual double crossings, we obtain a (free) tangle $T_c$.

Naturally, doing that on each level we transform the 2-foam $F(\beta)$ first into a ``2-dimensional 3-free tangle'' and then into a usual free 2-tangle $T_2(\beta)$. \newline

During the dynamics of the $n$ points on the plane, the Vorono\"i diagram may undergo the standard ``flips'': their dual are the Ptolemy flips of the Delaunay triangulation. And the corresponding transformation of the 1-foams $V_c$ is the $(2,2)$-move we have discussed. Therefore, {\em all foams $V_c$ for a given braid $\beta$ may be connected by a series of $(2,2)$-moves.} Moreover, the corresponding moves on tangles $T_c$ are just combinations of the second and the third Reidemeister moves.

Therefore, we have

\begin{theorem}
For a given braid $\beta$ and two sections $M\times \{c_1\}, M\times \{c_2\}$, the (free) tangle diagrams $T_{c_1}, T_{c_2}$ are Reidemeister equivalent. 	
\end{theorem}

That means that we have a 1-tangle $T_1(\beta)$ associated to the braid $\beta$. On the other hand, we may say that the braid $\beta$ defines a path in the space of the diagrams of the tangle $T_1(\beta)$. In particular, a {\em pure} braid defines a {\em loop} in that space.

Summing everything up, given a braid $\beta$, we got three objects: 1-dimensional tangle $T_1(\beta)$, 2-dimensional tangle $T_2(\beta)$ and 2-dimensional foam $F(\beta)$. A natural question is: given two isotopic braids $\beta_1, \beta_2$, what may be said about the corresponding objects? The equivalence of 2-foams was established by Theorem~\ref{thm:braid_foam}. The tangles are governed by the following theorem:

\begin{theorem}
Let $\beta_1, \beta_2$ be two isotopic braids. Then
\begin{enumerate}
\item $T_1(\beta_1)=T_1(\beta_2)$ as tangles;
\item $T_2(\beta_1)$ and $T_2(\beta_2)$ are connected by a sequence of moves obtained from the foam moves in Theorem~\ref{thm:braid_foam} by resolution of the 1-cells.	
\end{enumerate}

\end{theorem}

An important partial case is $M=S^2$, in other words, the spherical braids. In that case each foam $V_c$ is closed, and therefore the tangles $T_c$ are in fact {\em links}. Consequently, the tangles $T_1(\beta)$ and $T_2(\beta)$ are not just tangles, but links (1- and 2-dimensional), and the braid $\beta$ defines a path in the space of diagrams of the link $T_1(\beta)$.

\subsection{Further directions: Legendrian knots and cluster algebras}
As was established above, each domino tiling defines a triple diagram --- a collection of lines on the tiled surface. An important partial case is when the tiled surface is a 2-sphere $S^2$. A well-known operation on the curves is {\em lifting them to the unit tangent bundle of the manifold}. That gives us a {\em Legendrian knot (or link)}. To be precise, the knot is obtained by getting the third coordinate from the tangent vector to the curve. In other words, given a curve $\gamma$ on a surface $S$ (say, a sphere), we lift it to a curve $\hat{\gamma}$ using the tangent vector to the curve $\gamma$ as the extra coordinate.

There are natural questions arising from this construction. First of all, which Legendrian knots allow the representation by an {\em oriented} triple diagram with alternating orientations at each triple point? It is know that the answer is not ``all'', but the exact family of ``good'' Legendrian knots is unknown. 

The second question is somewhat inverse: is it true that any two triple diagrams of a Legendrian knot may be connected by a sequence of $(2,2)$-moves? 	\newline

\begin{figure}
\centering\includegraphics[width=300pt]{ 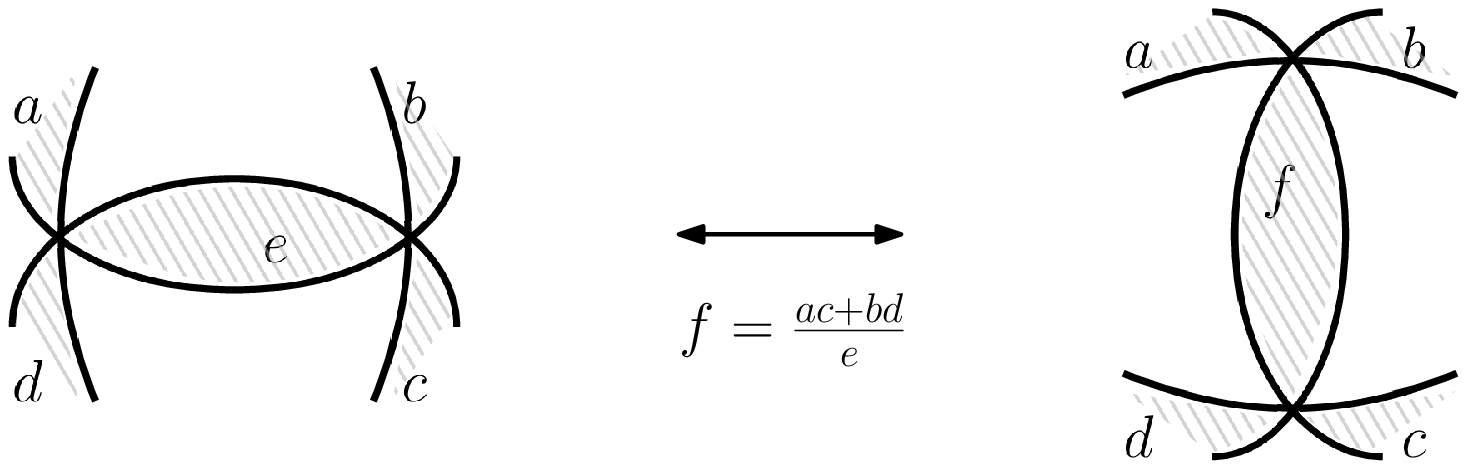}
\caption{Cluster variables change under the $(2,2)$-move}
\label{fig:domino_flip_cluster}
\end{figure}

Domino tilings have an interesting connection to cluster algebras. Consider a triple diagram, given by a domino tiling, and its checkerboard colouring. Associate a variable to each {\em red} region of the diagram. Now, let the diagram be changed by a $(2,2)$-move. If the center region is blue, the variables remain unchanged. If it is red (and we deal with a dual Ptolemy flip), the variables change as shown in Fig.~\ref{fig:domino_flip_cluster}.

Hence, a triple diagram gives a cluster, and the connectivity class of triple diagrams gives a cluster algebra.

Naturally, another algebra may be obtained by attributing variables to the {\em blue} regions of the diagram. \newline

This connection is to be expected, because the $(2,2)$-moves are just Ptolemy flips, and the connection of those with cluster algebra is known. Note that in a similar way $\Gamma_n^4$-groups are connected to cluster algebras: each generator $d_{ijkl}$ corresponds to a Ptolemy flip, and hence gives the cluster variables transformation.

Note further, that like in the $\Gamma-$case the variables were associated with the {\em edges} of the Delaunay triangulation, here the variables are attributed to {\em the regions of the other colour}: for example, if we put variables into red regions, the transformation is given by the {\em blue} dual flips.

%% file: T-group.tex
As we know, the Ptolemy flips give rise to the $\Gamma_n^4$ group, see Section~\ref{sec:gamma_n4}. Since the $(2,2)$-moves of the domino tiling theory translate into flips for dual graphs, it seems natural to consider the group-theoretic interpretation of the underlying geometry. \newline

The domino tiling contains more information than just triangulation and its flips. The reason is simply that the dual triangulation appears from ``half'' the regions (all the regions of one colour) but the second set of regions produces {\em another} triangulation (and it is crucial, as was seen in the decagon relation). So the group which we strive to construct must be different from $\Gamma_n^4$ (though with certain similarities).	

Consider two natural numbers $N, M$ and two sets of indices: $\mathcal{A}, |\mathcal{A}|=N; \, \mathcal{B}, |\mathcal{B}|=M$. \newline

The $\tau_{N,M}^4$ group is defined by the following presentation. The set of generators of the group is $\mathcal{G}=\{a^{\alpha}_{ijkl}, a^i_{\alpha\beta\gamma\delta}\}$, where Latin letters denote the elements of the set $\mathcal{A}$, and Greek letters denote the elements of the set $\mathcal{B}$. The quadruple of indices in both cases is ordered up to cyclic permutation and order reversal, and all indices are distinct; denote the sets of such quadruples of indices by $\bar{\mathcal{A}^4}$ and $\bar{\mathcal{B}^4}$, respectively. \newline

Relations in the group are of the following types:
\begin{enumerate}
	\item[1.] $a^2=1$ for any generator $a$ of the group;
	\item[2.] $a^{\epsilon_1}_{m_1}a^{\epsilon_2}_{m_2}=a^{\epsilon_2}_{m_2}a^{\epsilon_1}_{m_1}$ if $|m_1 \cap m_2|<3$, and $\epsilon_1\notin m_2, \epsilon_2\notin m_1$ (here $m_i$ is a multiindex: $m_i\in \bar{\mathcal{A}^4}$ or $m_i\in\bar{\mathcal{B}^4}$); 
\end{enumerate}
	
\begin{enumerate}
	\item[3.] for any five-element set $\{i,j,k,l,m\}$ (all indices must be either from $\mathcal{A}$ or from $\mathcal{B}$) and for any $\alpha, \beta$ in the second index set, $$a^{\alpha}_{ijkm}a^{\beta}_{jklm}a^{\alpha}_{ijlm}a^{\beta}_{ijkl}a^{\alpha}_{iklm}a^{\beta}_{ijkm}a^{\alpha}_{jklm}a^{\beta}_{ijlm}a^{\alpha}_{ijkl}a^{\beta}_{iklm}=1;$$
	\item[4.] for any two sets $\{i,j,k,l,m,p,q\}$ and $\{\alpha,\beta,\gamma,\delta,\varepsilon,\zeta,\eta\}$ set
	$$a^p_{\alpha\beta\gamma\eta}a^{\gamma}_{jkqp}a^q_{\eta\gamma\varepsilon\zeta}a^{\eta}_{ipqm}a^p_{\beta\delta\gamma\eta}a^{\gamma}_{klqp}a^q_{\alpha\eta\gamma\zeta}a^{\eta}_{ijpq}a^p_{\varepsilon\gamma\eta\delta}a^{\gamma}_{qplm}\cdot$$
	$$\cdot a^q_{\alpha\beta\gamma\eta}a^{\eta}_{jkqp}a^p_{\eta\gamma\varepsilon\zeta}a^{\gamma}_{ipqm}a^q_{\beta\delta\gamma\eta}a^{\eta}_{klqp}a^p_{\alpha\eta\gamma\zeta}a^{\gamma}_{ijpq}a^q_{\varepsilon\gamma\eta\delta}a^{\eta}_{qplm}=1;$$ and vice-versa, replacing Latin and Greek letters (that is, indices from the sets $\mathcal{A}$ and $\mathcal{B}$).
\end{enumerate}

\begin{definition} \label{def:tau}
	The group $\tau_{N,M}^4$ is defined as $$\tau_{N,M}^4=\langle\mathcal{G}\,|\,1,2,3,4\rangle.$$
\end{definition}

Each generator of the group $\tau_{N,M}^4$ may be regarded as a flip in the following sense. Let us return to the dual graphs. If we consider a checkerboard colouring of a plane partitioning and denote ``red'' regions by Latin letters and ``blue'' ones by Greek letters, each generator of the form $a_{ijkl}^\alpha$ corresponds to a flip in the quadrilateral $(i,j,k,l)$ such that the ``interior'' region is labelled with $\alpha$.

\begin{figure}
\centering\includegraphics[width=300pt]{ 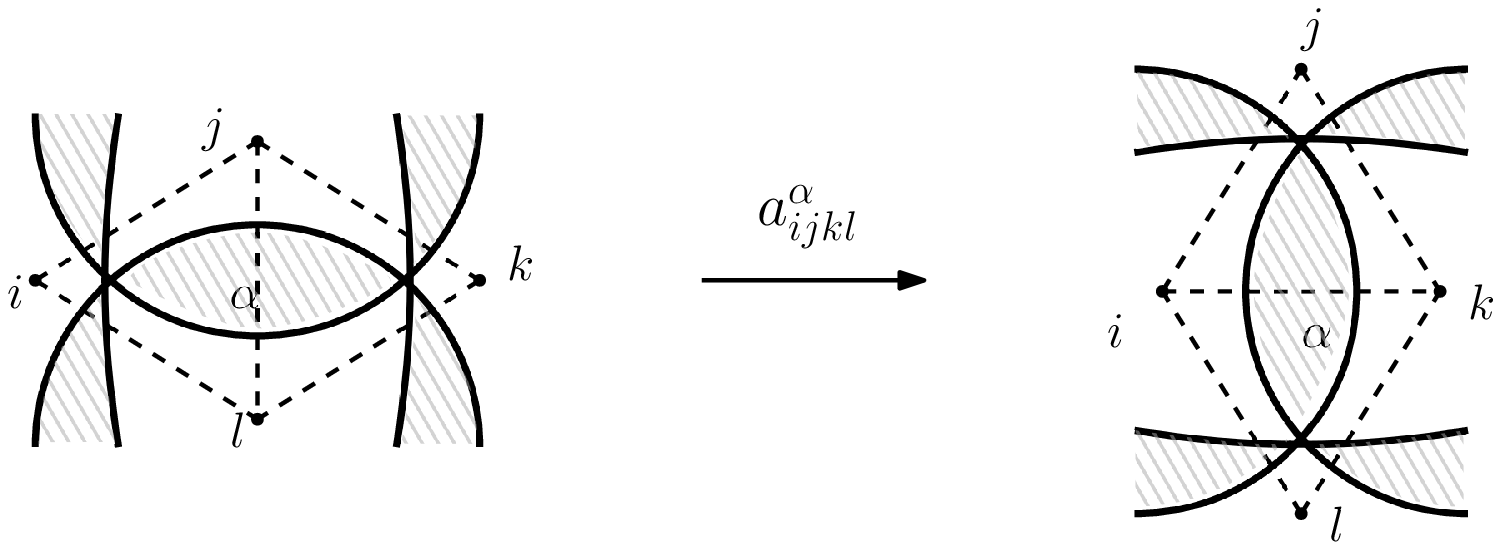}
\caption{The $a_{ijkl}^{\alpha}$ generator interpreted as a dual flip}
\label{dual_flip}
\end{figure}

The relation number 3 of the group $\tau_{N,M}^4$ may be interpreted as a pentagon relation from $\Gamma_N^4$ (and from $\Gamma_M^4$ as well) but such that the five-step process is performed twice. The reason for that is the following: going around the pentagon once interchanges the interior regions $\alpha$ and $\beta$. Thus, to return them to their places, we need to walk around the pentagon once again. \newline

Likewise, the 20-term relation (number 4, page~\pageref{sec:tau}) is interpreted as walking twice around the decagon we have seen earlier (again, walking once interchanges the regions $p$ and $q$, and $\gamma$ and $\eta$, so we need to go around once more, to restore their rightful places).

If one sets $M=1$ (effectively saying ``all Greek letters are the same''), the set of generators becomes equal to the set $\{d_{ijkl}\}$ with distinct $i,j,k,l$. Moreover, the relations $1,2$ become exactly the relations from the group $\Gamma_N^4$, and the relation 3 is just a ``doubled'' pentagonal relation of the group $\Gamma_N^4$. The relation 4 is new, though. 

On the other hand, if we write relations dealing only with exterior vertices the relation 4 disappears (since each term of this relation uses ``interior'' vertices). \newline

Note that if we forget the {\em ordering} of the indices of the generators of the group $\tau_{N,1}^4$ (the resulting group will be denoted by $T_{N,1}^4$), the first three relations give us exactly the relations of the group $G_N^4$. So the ``only exterior vertices matter'' mapping makes $G_N^4$ out of $T_{N,1}^4$.

Let us fix a Greek index $\alpha$ and consider the following mapping: $a_{ijkl}^{\alpha}\to d_{ijkl}, \,$ all other generators are mapped to $1$. In that case the first three relations again give the usual involution, far commutativity and pentagon conditions, while the 20-term relation becomes a 5-term one (but it is not a pentagon relation, for it deals with 7 Latin indices $i,j,k,l,m,p,q$; the indices $p,q$ are present at each term of the relation).

Further connections of the groups $\tau_{N,M}^4$ and $T^4_{N,M}$ with the groups $\Gamma_n^4$ and $G_n^4$ are yet to be discovered.